\documentclass[titlepage, 12pt, a4paper, twoside, openright]{report}
\usepackage{amsmath, amsthm, amsfonts, amscd,epsfig, graphicx}
 
\input xy
\xyoption {all}
\xyoption{2cell} 
\xyoption{v2} 

\renewcommand\a{\alpha}
\renewcommand\b{\beta}
\renewcommand\c{\gamma}
\renewcommand\d{\delta}

\newcommand\s{\sigma}
\newcommand\equivalencerelation{\sim \ }

\newcommand\Reals{\mathbb R}
\renewcommand{\P}{\mathcal{P}}
\newcommand\A{\mathcal{A}} 

\newcommand\Z{\mathbb Z}
\newcommand\Real{\mathbb R}
\newcommand\C{\mathbb C}

\newcommand\F{\mathcal F}

\newcommand\U{\mathcal U}
\newcommand\V{\mathcal V}
\newcommand\Q{\mathcal{Q}} 
 
\newcommand\T{\mathcal{T}}
\newcommand\xsimp{x_{0},\ldots,x_{p}}
\newcommand\gsimp{g_{0},\ldots,g_{p}}
\newcommand\tsimp{t_{1},\ldots,t_{p}}
\newcommand\hsimp{h_{0}[h_{1}|\cdots|h_{p}]}
\newcommand\bhsimp{[h_{1}|\cdots|h_{p}]}

\newcommand\dra{\stackrel{\displaystyle 
\rightarrow}{\rightarrow}}

\newcommand\dirlim{\stackrel{\displaystyle 
Lim}{\rightarrow}}

\newcommand\cstar{\C^\times}

\newcommand\fQ{(f_{1},f_{2})^{-1}Q}

\newcommand\gR{(g_{1},g_{2})^{-1}R}

\newcommand\isom{\stackrel{\simeq}{\rightarrow}}

\newcommand\B{\mathcal{B}}

\newcommand\tensorQ{\pi_{1}^{-1}Q \otimes 
\pi_{3}^{-1}Q}
\newcommand\abYtwo{\pi_{1}^{-1}Y^{[2]}
\times_{X^{[3]}}\pi_{3}^{-1}Y^{[2]}}
\newcommand\abY{\pi_{1}^{-1}Y\times_
{X^{[3]}}\pi_{3}^{-1}Y}
\newcommand\abfQ{(f_{1},f_{2})^{-1}Q}

\newcommand\abfdQ{(f_{1},f_{2})^{-1}\nabla_{Q}}
\newcommand\abFQ{<f_{1},f_{2}>^{-1}Q}
\newcommand\abFdQ{<f_{1},f_{2}>^{-1}\nabla_{Q}}

\newcommand\abfg{(g_{1}\circ f_{1}
,g_{2}\circ f_{2})}

\newcommand\deltatheta{\pi_{1}^{-1}
\theta \otimes \pi_{2}^{-1}
\theta^{*}\otimes \pi_{3}^{-1}\theta 
\otimes \pi_{4}^{-1}
\theta^{*}} 
 
\newcommand\zed{Z_{3}\times_{X^{[3]}}
Y_{23}\times_{X^{[3]}}Y_{12}}

\theoremstyle{plain}
\newtheorem{theorem}{Theorem}[chapter]
\newtheorem{lemma}{Lemma}[chapter]
\newtheorem{proposition}{Proposition}[chapter]
\newtheorem{conjecture}{Conjecture}[chapter] 
\newtheorem{corollary}{Corollary}[chapter] 

\theoremstyle{definition}
\newtheorem{definition}{Definition}[chapter]

\theoremstyle{remark}
\newtheorem{example}{Example}[chapter]
\newtheorem{note}{Note}[chapter]

\title{The Geometry of Bundle Gerbes} 
\author{Daniel Stevenson \\ 
Department of Pure Mathematics \\ 
University of Adelaide} 
\date{\today}

\begin{document}
\UseTwocells 
\UseHalfTwocells 
\null
\pagestyle{empty}
\vspace{3 cm}
\begin{center}
    \Huge\bf The Geometry of Bundle Gerbes
    \end{center}
    
\vspace{4 cm}
\begin{center}
  \bf \Large  Daniel Stevenson
    \end{center}
    
\vspace{4 cm}

\begin{center}
\bf  Thesis submitted for the degree of \\
 Doctor of Philosophy \\
 in the Department of Pure Mathematics \\
 University of Adelaide
    \end{center}
\vspace{1 cm} 
    \begin{center}
 \includegraphics[scale=0.2]{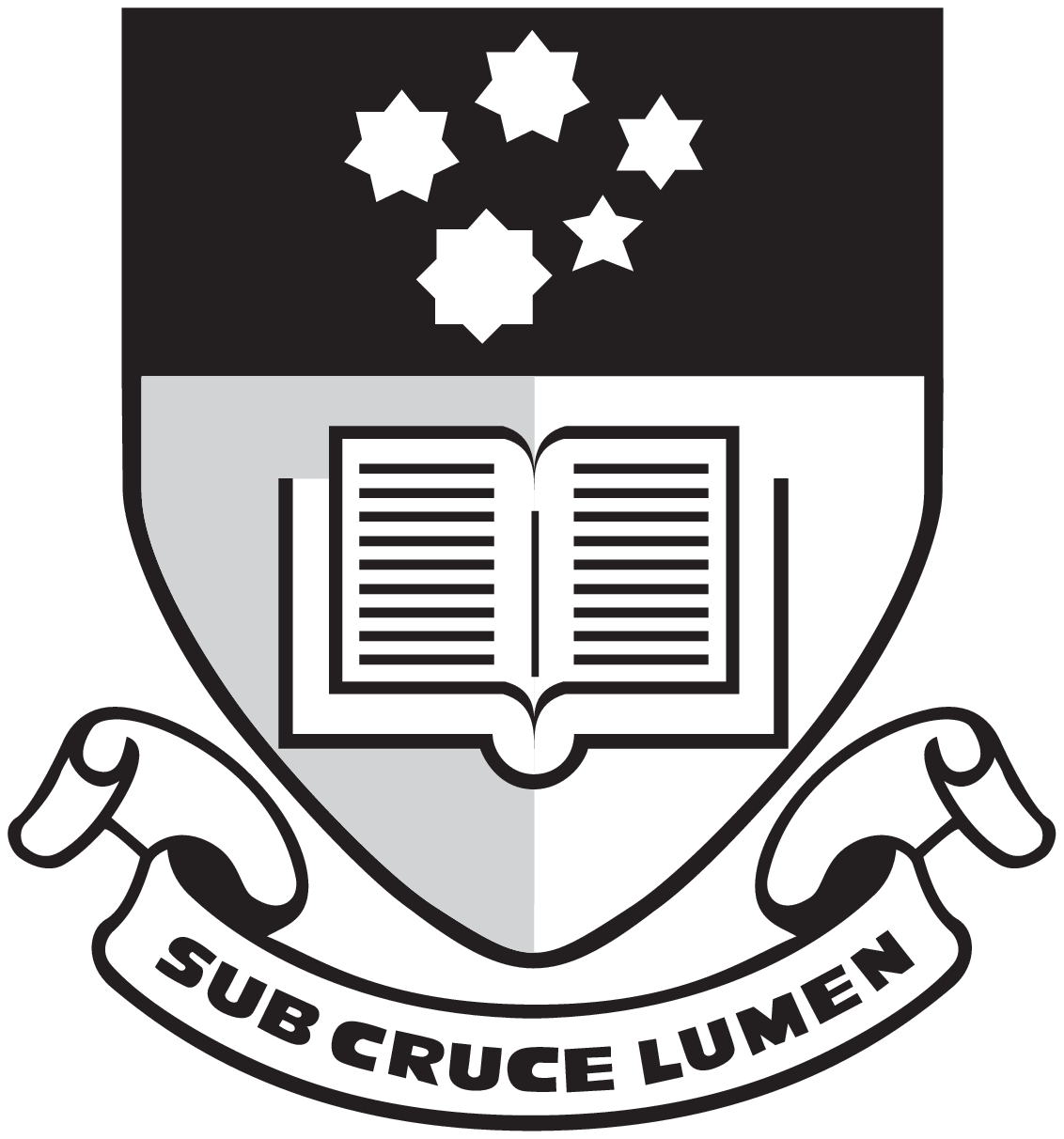}
 \end{center}

 \vspace{1 cm}
 
 \begin{center}
 \bf 9 February 2000
 \end{center}
 
\newpage
\null
\newpage

\pagestyle{plain}
\setcounter{page}{1}
\pagenumbering{roman}

\begin{center} 
   {\bf Abstract} 
\end{center}

This thesis reviews the theory of bundle gerbes 
and then examines the higher dimensional notion 
of a bundle 2-gerbe.  The notion of a bundle 
2-gerbe connection and 2-curving are introduced 
and it is shown that there is a class in $H^{4}(M;\Z)$ 
associated to any bundle 2-gerbe.  

\newpage
\null
\newpage

\begin{center} 
   {\bf Statement of Originality} 
\end{center} 

This thesis contains no material which has been accepted for the 
award of any other degree or diploma at any other university 
or other tertiary institution and, to the best of my 
knowledge and belief, contains no material previously published 
or written by another person, except where due reference 
has been made in the text. 

I give consent to this copy of my thesis, when deposited in 
the University Library, being made available for loan and 
photocopying.  

\vspace{ 3 cm}
\begin{flushright}
    Daniel Stevenson\\
    Adelaide, 9 February, 2000
    \end{flushright}
    
\newpage 
\null
\newpage

\begin{center} 
   {\bf Acknowledgement } 
\end{center} 

I would like to thank my supervisor Michael Murray for his 
patience, generosity and expert guidance in the writing of this thesis.  
He has greatly enhanced my understanding and enjoyment of 
mathematics over the past few years and I would like to say 
that it has been an honour to work with him during this time.  
It is my pleasure also to thank Alan Carey for some 
useful conversations and also Ann Ross for help in administrative matters.  
I would also like to gratefully acknowledge the support 
provided by an Australian Post Graduate Award. 

\newpage 

\tableofcontents

\pagenumbering{arabic}  
\setcounter{page}{1}

\setcounter{chapter}{0} 
\chapter{Introduction.}
\label{chapter:1}
 
There has been much recent interest in
applying the theory of gerbes to differential
geometry \cite{Hit} and  mathematical
physics (\cite{Kal} and \cite{CaMur}).  Gerbes
were originally introduced by Giraud (see \cite{Gir})
and were extensively studied in \cite{Bry}.
Loosely speaking a gerbe is a `sheaf of
categories' except that there is an extra
level of complication when one comes to
discuss sheaf-like glueing laws for objects.
Much of the  importance of gerbes in the previously
cited applications stems from the fact
that equivalence classes of gerbes over
a manifold $M$ are in  a bijective correspondence with $H^{3}(M;\Z)$.
This is entirely similar to the fact that
$H^{2}(M;\Z)$ is in a bijective correspondence
with isomorphism classes of principal $\cstar$
bundles on $M$.

One can think of  $\cstar$ bundles and gerbes as
a geometric `realisation' or `representation' of the
cohomology of the manifold in dimensions two and three
respectively.  From this perspective it is natural
to look for such objects in other dimensions.
   In dimension zero we can consider
the set of all maps $f:M\to \Z$ which are, in fact, in bijective
correspondence with $H^0(M; \Z)$. In dimension one there are various
kinds of objects we can consider which are in bijective correspondence
with $H^1(M; \Z)$. For example we could consider (isomorphism classes
of) principal $\Z$ bundles
over $M$. Or, because the classifying space of $\Z$ is $\cstar$ we could
equally well consider the space of homotopy classes of maps from $M$
into $\C^\times$.   A similar
story is true for $H^{2}(M;\Z)$: one could consider
$H^{2}(M;\Z)$ to be the set of all homotopy classes
of maps $f:M\to B\cstar$, where $B\cstar$ is
the classifying space for the group $\cstar$,
or, we could think of an element of $H^{2}(M;\Z)$
corresponding to a homotopy class of maps
$M\to B\cstar$ as giving rise to a principal
$\cstar$ bundle on $M$.  In the case of
$H^{2}(M;\Z)$ we can go one step further and think
of a principal $\cstar$ bundle on $M$ as a principal $B\Z$
bundle on $M$.  The crucial fact here, as
pointed out in \cite{Seg3}, is that
$B\Z$, can be chosen to be  an abelian topological group
since $\Z$ is.

This last fact is exploited in the paper \cite{Gaj}
by Gajer where the author introduces a notion of a
principal $B^{p}\cstar$ bundle on $M$ and proves
that there is a bijection between $H^{p+2}(M;\Z)$
and the set of all principal $B^{p}\cstar$ bundles
on $M$ (here $B^{p}\cstar$ is the $p$-fold iterated
classifying space $B\cdots B\cstar$ of the group $\cstar$).
Therefore, when one comes to consider $H^{3}(M;\Z)$
there are several possibilities: one option is to
think of $H^{3}(M;\Z)$ as the set of all homotopy
classes of maps $f:M\to BB\cstar$, another is to
think of $H^{3}(M;\Z)$ as the set of equivalence
classes of gerbes on $M$.  Yet another is to follow
Gajer and think of $H^{3}(M;\Z)$
as the set of all principal $B\cstar$ bundles
on $M$ or, alternatively, the set of all principal
$BB\Z$ bundles on $M$.

In \cite{Mur} the notion of a \emph{bundle
gerbe} was introduced.  Bundle gerbes also
give rise to classes in $H^{3}(M;\Z)$.  A
bundle gerbe on $M$ is a triple $(P,X,M)$
where $\pi_{X}:X\to M$ is a surjection which
locally admits sections and where $P$ is a
principal $\cstar$ bundle over the fibre
product $X^{[2]} = X\times_{M}X$.  $P$ is
required to have a product, that is a $\cstar$
bundle isomorphism which on the fibres of
$P$ takes the form
$$
P_{(x_{2},x_{3})}\otimes P_{(x_{1},x_{2})}
\to P_{(x_{1},x_{3})}
$$
for $x_{1}$, $x_{2}$ and $x_{3}$ in the same
fibre of $X$.  It can be shown that every
bundle gerbe on $M$ gives rise to a gerbe
on $M$.  There is a characteristic class
associated to bundle gerbes, the Dixmier-Douady
class.  This is a class in $H^{3}(M;\Z)$
which measures the obstruction to the bundle
gerbe being `trivial' in a certain sense.

In \cite{Bry} the notions of a \emph{connective
structure} on a gerbe and a \emph{curving} for
the connective structure were introduced.  It was
shown that one could associate to a gerbe
equipped with a connective
structure and curving a closed integral three
form $\omega$ on $M$ --- the three curvature
of the curving.  In \cite{Mur}, it
is shown that by considering a connection
on $P$ which is compatible with the product
above, it is possible to introduce the notion
of a bundle gerbe connection and a curving for
a bundle gerbe connection and that, furthermore,
a bundle gerbe on $M$ equipped with a bundle gerbe
connection and curving gives rise to a closed
integral three form on $M$.  This three form
is a representative of the image of the
Dixmier-Douady class of the bundle gerbe in
$H^{3}(M;\Reals)$.  It can be shown that a
bundle gerbe connection and
curving on a bundle gerbe induces a connective
structure and curving on the gerbe associated
to the bundle gerbe.

In the series of papers \cite{BryMcL1} and
\cite{BryMcL2} the authors considered a
certain \emph{2-gerbe} associated to a principal
$G$ bundle $P\to M$ where $G$ was a compact,
simple, simply connected Lie group.  2-gerbes,
introduced in \cite{Bre}, are higher dimensional
analogues of gerbes.  Broadly speaking, a 2-gerbe
on $M$ is a `sheaf of bicategories' on $M$.
Bicategories, defined in \cite{Ben}, are
higher dimensional analogues of categories.
A bicategory consists of a set of objects and for
each pair of objects, a category.  The objects of
this category are then 1-arrows of the bicategory and
the arrows of the category are 2-arrows
of the bicategory.  As a result, the glueing laws
involved in the definition of a sheaf of
bicategories become quite intricate.  One can
show that if one considers a certain class of
2-gerbes (2-gerbes \emph{bound} by $\cstar$) then
one can associate an element of $H^{4}(M;\Z)$ to
each 2-gerbe and moreover it can be shown that
$H^{4}(M;\Z)$ is in a bijective correspondence
with `biequivalence' classes of 2-gerbes.

In \cite{CaMuWa} a related object, a \emph{bundle
2-gerbe}, was introduced.  In the latter part
of this thesis we will study bundle 2-gerbes
using a variation of the definition given in
\cite{CaMuWa}.  A bundle 2-gerbe is a
quadruple $(Q,Y,X,M)$ where $\pi:X\to M$ is a
surjection admitting local sections and where
$(Q,Y,X^{[2]})$ is a bundle gerbe on the
fibre product $X^{[2]}$.  One also requires
that there is a bundle 2-gerbe product.  This
consists of a `product' on $Y$, ie a map
which on the fibres of $Y$ takes the form
$$
m:Y_{(x_{2},x_{3})}\times Y_{(x_{1},x_{2})}
\to Y_{(x_{1},x_{3})},
$$
where $x_{1}$, $x_{2}$ and $x_{3}$ are
points of $X$ lying in the same fibre.  There is
also a product on $Q$ covering the induced
product on $Y^{[2]}$.  The product on $Q$ is a 
$\cstar$ bundle morphism 
$\hat{m}$ given fibrewise by
$$
\hat{m}:Q_{(y_{23},y_{23}^{'})}\otimes
Q_{(y_{12},y_{12}^{'})}\to
Q_{(m(y_{23},y_{12}),m(y_{23}^{'},y_{12}^{'}))}.
$$
This is required to commute with the
bundle gerbe product on $Q$ and satisfy a certain
associativity condition.  One can show that
there is a $\cstar$ valued \v{C}ech 3-cocycle $g_{ijkl}$
representing a class in $H^{4}(M;\Z)$
associated to a bundle 2-gerbe.
In \cite{Mur} it was shown that the
complex
$$
\Omega^{p}(M)\stackrel{\pi^{*}}{\to}
\Omega^{p}(X)\stackrel{\d}{\to}
\Omega^{p}(X^{[2]})\stackrel{\d}{\to} \cdots
$$
was exact, where $\d:\Omega^{p}(X^{[q]})\to
\Omega^{p}(X^{[q+1]})$ is the differential
formed by adding the projection maps
$\pi_{i}^{*}$ with an alternating
sign (here $\pi_{i}:X^{[q]}\to X^{[q-1]}$
deletes the $i^{th}$ factor).  By exploiting
this fact one can define the notion
of a \emph{bundle 2-gerbe connection}
and a \emph{2-curving} for a bundle
2-gerbe.  This
leads to a closed integral four form
$\Theta$ on $M$ which is a representative
in $H^{4}(M;\Reals)$ for the image, in real cohomology,
of the class defined by the
cocycle $g_{ijkl}$.

Associated to a  principal $G$ bundle $P\to M$
for $G$ a compact, simple, simply connected
Lie group, is a   bundle 2-gerbe
$(\tilde{Q},\tilde{P},P,M)$.  If one calculates the cocycle $g_{ijkl}$
associated to the bundle 2-gerbe $\tilde{Q}$
then one recovers the result of \cite{BryMcL1}
giving an explicit formula for the first
Pontryagin class of $P$.  Finally one can also
show that every bundle 2-gerbe on $M$ gives
rise to a 2-gerbe on $M$.

In outline then this thesis is as follows:
In Chapter~\ref{chapter:2} we review the
definitions of simplicial spaces and the
construction of the universal $G$ bundle
$EG\to BG$, we give a description of
the group law on $EG$ and we discuss
some of the results of Gajers paper \cite{Gaj}.
In Chapter~\ref{chapter:3} we begin with some remarks
on principal $\cstar$ bundles before going
on to review the definition of a bundle gerbe
and some of the main results from \cite{Mur}.
We then discuss the notion of
stable isomorphism of bundle gerbes and
the result from \cite{MurSte} that the
set of all stable isomorphism classes
of bundle gerbes on $M$ is in a bijective
correspondence with $H^{3}(M;\Z)$.  Next we discuss
bundle gerbe morphisms and transformations of
bundle gerbe morphisms as well as the related
notion of stable morphisms and stable transformations.
Finally we show that the collection of all
bundle gerbes on $M$ forms a 2-category with bundle gerbe
morphisms as 1-arrows and transformations as
2-arrows.  In Chapter~\ref{chapter:4} we
discuss the singular theory of principal
$\cstar$ bundles in a form which readily
generalises to give a singular theory of
bundle gerbes.  Chapter~\ref{chapter:5} discusses
bundle gerbes from a local viewpoint along
similar lines to the presentation in \cite{Hit} for
gerbes and  Deligne cohomology for bundle gerbes.
Chapter~\ref{chapter:6} is mostly devoted
to a review of the material presented in
\cite{Bry}.  We also discuss the construction
from \cite{MurSte} of a gerbe associated to a
bundle gerbe.  The last section of Chapter~\ref{chapter:6}
contains a discussion of higher glueing
laws.

The remainder of the thesis is concerned with 2-gerbes and
bundle 2-gerbes. In Chapter~\ref{chapter:7} we define the
notions of simplicial bundle gerbe and bundle
2-gerbe.  Chapter~\ref{chapter:10} starts with
the definitions of a bicategory and bigroupoid and
goes on to discuss in some detail the
homotopy bigroupoid $\Pi_{2}(X)$ associated
to a space $X$.  We also show that a bundle
2-gerbe on a point is the same thing as a
bicategory in which the automorphism groups
of each 1-arrow are isomorphic to $\cstar$.
In Chapter~\ref{chapter:11} we give some examples of bundle
2-gerbes, namely the tautological bundle 2-gerbe
from \cite{CaMuWa} and the bundle 2-gerbe
associated to a principal $G$ bundle $P\to M$.
In Chapter~\ref{chapter:8} we define the
notions of bundle 2-gerbe connections and
2-curvings for bundle 2-gerbe connections
and show how a bundle 2-gerbe connection
with a 2-curving on a bundle 2-gerbe gives
rise to a closed, integral four form $\Theta$ on the base
manifold $M$.  We also show that bundle 2-gerbe
connections exist, something that is not
immediately obvious.  Chapter~\ref{chapter:9}
begins with the construction of a $\cstar$
valued \v{C}ech 3-cocycle $g_{ijkl}$ associated
to a bundle 2-gerbe and then goes on to show
how a bundle 2-gerbe equipped with a 2-curving and
bundle 2-gerbe connection gives rise to a
class in the Deligne cohomology group
$H^{3}(M;\underline{\C}^{\times}_{M}\to
\underline{\Omega}^{1}_{M}\to \underline{\Omega}^{2}_{M}
\to \underline{\Omega}^{3}_{M})$.  As a
consequence, we see that the four form
$\Theta$ represents the image in
$H^{4}(M;\Reals)$ of the class in $H^{4}(M;\Z)$
associated to the cocycle $g_{ijkl}$.  Finally
we show the four class of the bundle 2-gerbe
$\tilde{Q}$ arising from a principal $G$ bundle
$P\to M$ is the first Pontryagin class of $P$.

In Chapter~\ref{chapter:12} we discuss the
notion of a trivial bundle 2-gerbe and prove
that if the four class in $H^{4}(M;\Z)$ of
a bundle 2-gerbe vanishes then the bundle 2-gerbe
must be trivial. The final Chapter~\ref{chapter:13}
discusses the relationship of bundle 2-gerbes
to 2-gerbes.

\setcounter{chapter}{1}
\chapter{Review of some simplicial techniques} 
\label{chapter:2} 
\section{Simplicial spaces}
\label{sec:2.1} 
Recall, see \cite{May1}, 
the standard $n$-simplex $\Delta^{n}$ 
in $\Real^{n+1}$ is defined by 
$\Delta^{n} = \{(x_{0},\ldots,x_{n}):x_{0}+\cdots + x_{n} = 1\}$.
If $n \geq 1$ there are face maps 
$\delta^{i}:\Delta^{n-1} \to \Delta^{n}$ for 
$i = 0,1,\ldots,n$ given by 
$$
\delta^{i}(x_{0},\ldots,x_{n-1}) = 
(x_{0},\ldots,x_{i-1},0,x_{i},\ldots,x_{n-1}).
$$ 
For the same $n$ we have degeneracy maps
$\sigma^{i}:\Delta^{n} \to \Delta^{n-1}$, 
$i=0,1,\ldots,n-1$ with 
$$
\sigma^{i}(x_{0},\ldots,x_{n}) = 
(x_{0},\ldots,x_{i-1},x_{i}+x_{i+1},x_{i+2},\ldots,x_{n}).
$$

\begin{definition} 
\label{def:2.1.1}
A \emph{simplicial space} $X$ is a sequence 
$\{X_{n}\}$ of spaces $X_{n}$ for $n=0,1,2,\ldots$ 
together with face operators 
$d_{i}:X_{n} \to X_{n-1}$, $i=0,1,\ldots,n$ for 
$n=1,2,3,\ldots$ and degeneracy operators 
$s_{i}:X_{n} \to X_{n+1}$, $i=0,1,\ldots,n$
for $n=0,1,2,\ldots$ which together satisfy 
the \emph{simplicial identities}

\begin{eqnarray*}
d_{i}d_{j} & = & d_{j-1}d_{i}, \ i < j \\
s_{i}s_{j} & = & s_{j+1}s_{i}, \ i \leq j  \\  
d_{i}s_{j} & = & \begin{cases}
              s_{j-1}d_{i}, &   i < j    \\
              \text{id},  & i = j ,  i = j+1 \\
              s_{j}d_{i-1},  & i > j+1   
              \end{cases}                   
\end{eqnarray*}  
\end{definition}  

Notice that there is a definition of 
simplicial manifold, simplicial group, 
simplicial ring and so on; we simply replace 
the topological category by the smooth 
category or the category of groups and so forth.  
One then requires that 
the face and degeneracy operators are 
morphisms in the appropriate 
category.  Thus a simplicial manifold 
is a sequence of manifolds $\{X_{n}\}$ 
together with smooth maps 
$d_{i}:X_{n} \to X_{n-1}$ for 
$i = 0,\ldots n$ and 
$s_{i}:X_{n}\to X_{n+1}$ for 
$i = 1,\ldots n$  
satisfying the same conditions 
as in Definition~\ref{def:2.1.1} above. 
Similarly, a simplicial group is a 
sequence of groups $\{G_{n}\}$ with face and 
degeneracy operators $d_{i}$ and $s_{i}$ as 
above which are required to be homomorphisms 
between the various groups $G_{n}$.  

A more concise definition is given as 
follows (see \cite{Mac} or \cite{May1}).   
Let $\Delta$ be 
the category whose objects are all 
finite ordered sets $[n] = \{0,1,\ldots,n\}$ 
and whose morphisms $[m]\to [n]$ are all 
nondecreasing monotonic functions 
$[m]\to [n]$.  One then restricts 
attention to the morphisms 
$\d^{i}:[n-1]\to [n]$ and 
$\sigma^{i}:[n]\to [n-1]$ for 
$i = 0,1,\ldots,n$ given by 
\begin{eqnarray*} 
\d^{i}(j) = \begin{cases} 
            j & \text{if}\ j < i,       \\
            j + 1 &\text{if}\ j \geq i, \\ 
            \end{cases} 
\end{eqnarray*} 
and 
\begin{eqnarray*} 
\sigma^{i}(j) = \begin{cases} 
                j &\text{if}\ j \leq i ,  \\ 
                j - 1 & \text{if}\ j > i.  \\ 
                \end{cases} 
\end{eqnarray*} 
One can then prove a Lemma which states that 
every morphism $\mu :[m]\to [n]$ in 
$\text{Hom}(\Delta)$ can be written 
uniquely as 
$$
\mu = \d^{i_{a}}\d^{i_{b}}\cdots \d^{i_{c}}\sigma^{j_{x}}
\sigma^{j_{y}}\cdots \sigma^{j_{z}},
$$
where $n\geq i_{a}>i_{b}>\cdots >i_{c}\geq 0$ 
and $0\leq j_{x}<j_{y}<\cdots <j_{z}< m$ 
and $n = m - z + c$.  (See for instance 
\cite{Mac}). 
Now we can define a \emph{simplicial object} 
in the category of topological spaces,   
a simplicial space in other words, as being a functor 
$F:\Delta^{\circ}\to \mathcal{C}$ where 
$\mathcal{C}$ is the category of topological 
spaces.  Here $\Delta^{\circ}$ 
denotes the category opposite to 
$\Delta$ --- recall the category 
$\mathcal{C}^{\circ}$ opposite 
to $\mathcal{C}$ is the category with 
the same objects as $\mathcal{C}$ but 
with the source and target of each 
arrow in $\mathcal{C}$ interchanged.  
The simplicial identities arise because 
the morphisms $\d^{i}$ and $\sigma^{i}$ satisfy the 
analogous identities in $\Delta$.  
Changing $\mathcal{C}$ to the category of 
smooth manifolds or the category of groups 
recovers the definitions of simplicial 
manifolds and simplicial groups.   

A good source of examples of simplicial 
spaces is the \emph{nerve} of a category.  
(see \cite{Seg2}). 
Let $\mathcal{C}$ be a category.  
We define a simplicial set $N\mathcal{C}$ 
called the nerve of the category $\mathcal{C}$ 
as follows.  $N\mathcal{C}_{0}$ is defined 
to be the set of objects, $Ob(\mathcal{C})$, 
of $\mathcal{C}$ (so we had better only consider 
small categories), $N\mathcal{C}_{1}$ 
is defined to be the set of arrows or 
morphisms $Mor(\mathcal{C})$ of $\mathcal{C}$ and 
in general $N\mathcal{C}_{n}$ is the subset 
$Mor(\mathcal{C})\circ \cdots \circ Mor(\mathcal{C})$ 
of $Mor(\mathcal{C})\times \cdots \times Mor(\mathcal{C})$ 
($n$ factors) consisting of all $n$-tuples 
$(f_{1},\ldots f_{n})$ of morphisms of 
$\mathcal{C}$ with $source(f_{i}) = target(f_{i+1})$.  
The face operators are defined by the source 
and target maps and composition of morphisms 
while the degeneracy operators are defined 
by including an identity morphism.  More 
specifically, $d_{i}:N\mathcal{C}_{n}\to N\mathcal{C}_{n-1}$ 
for $n \geq 1$ is defined by 
$$
d_{i}(f_{1},\ldots,f_{n}) = \begin{cases} 
                            (f_{2},\ldots,f_{n}) &  \text{for 
                             $i = 0$}                         \\
                            (f_{1},\ldots ,f_{i}\circ f_{i+1},
                            \ldots, f_{n})  & \text{for $1 \leq i< n$} \\
                            (f_{1},\ldots,f_{n-1})  & \text{for $i = n$} \\ 
                            \end{cases} 
$$
and the source and target maps for 
$n = 1$.  $s_{i}:N\mathcal{C}_{n} \to N\mathcal{C}_{n+1}$ 
is defined by 
$$
s_{i}(f_{1},\ldots ,f_{n}) = (f_{1},\ldots ,f_{i},\text{id},f_{i+1},
\ldots ,f_{n}),
$$
where $\text{id}$ is the identity arrow at 
$target(f_{i})$ in $\mathcal{C}$.    

We will be especially interested in 
categories $\mathcal{C}$ for which the 
sets of objects and arrows have 
topologies such that the source and 
target maps and the operation of 
composition of arrows are continuous 
in these topologies.   
Segal in \cite{Seg2}  
calls these 
\emph{topological categories}.  Segal 
shows that there is a topological category 
$X_{\mathcal{U}}$ associated to a covering 
$\mathcal{U} = \{U_{i}\}_{i \in I}$ 
of a topological space $X$ (see \cite{Dup}).
The objects of $X_{\mathcal{U}}$ consist 
of pairs $(x,U_{i})$ where $x \in U_{i}$.  
There is a unique morphism $(x,U_{i})\to 
(y,U_{j})$ if and only if $x = y \in U_{ij} = U_{i}\cap U_{j}$.   
The objects of $X_{\mathcal{U}}$ are topologised 
as $Ob(X_{\mathcal{U}}) = \coprod_{i\in I}U_{i}$ 
and the arrows of $X_{\mathcal{U}}$ are 
topologised as $Arr(X_{\mathcal{U}}) = 
\coprod_{i,j\in I, U_{ij} \neq \emptyset}U_{ij}$.   
Clearly 
$X_{\mathcal{U}}$ is an example of a 
topological category.  

Of more interest are the two 
topological categories $G$ and 
$\bar{G}$ defined by Segal which 
are associated to a topological 
group $G$.  Here $G$ 
is the standard category associated 
to a group, so $G$ has one object and 
the morphisms of $G$ are just the set 
of points of $G$.  Clearly $G$ is a 
topological category.  Thus the nerve of the 
topological category $G$ is the simplicial 
space $NG$ with 
$NG_{n} = G\times \cdots \times G$ ($n$ factors) and 
face and degeneracy operators given 
by 
$$
d_{i}(g_{1},\ldots,g_{n}) = \begin{cases} 
                            (g_{2},\ldots,g_{n}), & i = 0     \\
                            (g_{1},\ldots,g_{i}g_{i+1},\ldots g_{n}), 
                             & 1 \leq i < n                             \\
                            (g_{1},\ldots,g_{n-1}), & i = n,           \\
                            \end{cases} 
$$
and 
$$
s_{i}(g_{1},\ldots,g_{n}) = (g_{1},\ldots,g_{i},1,g_{i+1},\ldots,g_{n}). 
$$
$\bar{G}$ is the topological category 
whose objects are the points of $G$ and 
where the unique arrow joining two objects $g_{1}$ and $g_{2}$ is 
$(g_{1},g_{2})$.  Thus the objects have 
a topology induced by $G$ and the arrows have 
a topology induced by $G\times G$.  Therefore 
the nerve of the topological category 
$\bar{G}$ is the simplicial space 
$N\bar{G}$ with $N\bar{G}_{n} = G\times \cdots \times G$ 
($n+1$ factors).  The face and degeneracy 
operators are given by  
$$
d_{i}(g_{0},\ldots,g_{n}) = (g_{0},\ldots,g_{i-1},g_{i+1},\ldots,g_{n})
$$
and 
$$
s_{i}(g_{0},\ldots,g_{n}) = (g_{0},\ldots,g_{i},g_{i},g_{i+1},\ldots,g_{n}).
$$
We shall return to these examples in the 
next section.  Finally, we would like 
to define the notion of a simplicial map 
between simplicial spaces.  
\begin{definition} 
\label{def:2.1.2}  
Let $X$ and $Y$ be simplicial spaces.  
A \emph{simplicial map} $f:X\to Y$ is a sequence 
of maps $f_{n}:X_{n}\to Y_{n}$ which 
commute with the face and degeneracy 
operators.  Thus if $d_{i},s_{i}$ and 
$d_{i}^{\prime},s_{i}^{\prime}$ denote 
the face and degeneracy operators 
of $X_{n}$ and $Y_{n}$ respectively, 
then we have 
$$
f_{n-1}\circ d_{i} = d_{i}^{\prime}\circ f_{n}
$$
and 
$$
f_{n+1}\circ s_{i} = s_{i}^{\prime}\circ f_{n}.
$$
\end{definition} 
If we view simplicial spaces, simplicial groups and so on  
as functors $F:\Delta^{\circ}\to \mathcal{C}$ 
for an appropriate choice of the category 
$\mathcal{C}$, then a simplicial map corresponds 
to a natural transformation between 
functors.  More concretely, if 
$F_{1},F_{2}:\Delta^{\circ}\to \mathcal{C}$ 
represent simplicial objects in some 
category $\mathcal{C}$, then a simplicial 
map from $F_{1}$ to $F_{2}$ 
is a natural transformation 
$\tau :F_{1}\Rightarrow F_{2}$.  The 
definition above can be recovered by 
recalling that a natural transformation 
$\phi:F\Rightarrow G$ between 
functors $F,G:\mathcal{C}\to \mathcal{D}$ 
assigns to each object $c$ of 
$\mathcal{C}$ an arrow $\phi(c):F(c)\to G(c)$ 
in $\mathcal{D}$ such that the diagram 
$$
\xymatrix{ 
F(c_{1}) \ar[r]^{F(\a)} \ar[d]_{\phi(c_{1})} & 
F(c_{2}) \ar[d]^{\phi(c_{2})}                 \\ 
G(c_{1}) \ar[r]^{G(\a)} & G(c_{2}) } 
$$
commutes for all arrows $\a:c_{1}\to c_{2}$ 
in $\mathcal{C}$.  Thus to each object 
$[n]$ of $\Delta^{\circ}$, $\tau$ 
associates an arrow $\tau_{n}:F_{1}([n])\to F_{2}([n])$, 
which, by the commutativity of the diagram 
above, must commute with the face and degeneracy maps 
$F_{1}(\d^{i})$, $F_{1}(\sigma^{i})$, 
$F_{2}(\d^{i})$ and $F_{2}(\sigma^{i})$.  
  
Notice that a functor 
$F:\mathcal{C}_{1}\to \mathcal{C}_{2}$ 
induces a simplicial map 
$NF:N\mathcal{C}_{1}\to N\mathcal{C}_{2}$ 
between simplicial sets by definition 
of a functor.  If $\mathcal{C}_{1}$ 
and $\mathcal{C}_{2}$ are topological 
categories then we of course require that the 
maps on objects and arrows induced by 
$F$ are continuous.  As an example, take 
the topological categories $\bar{G}$ and 
$G$ associated to a topological group $G$.  
Then we can define a functor $p:\bar{G}\to G$ 
by defining $p(g_{1},g_{2}) = g_{1}^{-1}g_{2}$.  
It is easy to check that this defines a 
continuous functor and hence a map of 
simplicial spaces 
$Np:N\bar{G}\to NG$.   

\section{Geometric realisation} 
\label{sec:2.2} 
There is a procedure by which one can 
obtain a single topological space 
from a simplicial space known as                            
\emph{geometric realisation}. 

\begin{definition}
\label{def:2.2.1} 
Let $X$ be a simplicial space.  The 
\emph{geometric realisation}  
of $X$, denoted $|X|$, is  
the quotient  
$$
\coprod \Delta^{n} \times X_{n}/ \equivalencerelation,
$$
where $\equivalencerelation$ is the equivalence relation 
$(\d^{i}(t),x) \equivalencerelation (t,d_{i}(x))$ for 
$t \in \Delta^{n-1},\ x \in X_{n}$ and 
$i = 0,1,\ldots,n$, $n = 1,2,\ldots$ and 
$(\s^{i}(t),x) \equivalencerelation (t,s_{i}(x))$ for 
$t \in \Delta^{n},\ x \in X_{n-1},\ i=0,1,\ldots,n$
and $n=0,1,2,\ldots$.  $|X|$ is made into 
a topological space in the following way 
as explained in \cite{May2}.  Let $F_{q}|X|$ 
denote the image of $\coprod_{i=0}^{q}\Delta^{i}\times X_{i}$ 
in $|X|$ and give $F_{q}|X|$ the quotient 
topology.  Then $F_{q}|X|$ is a closed 
subset of $F_{q+1}|X|$ and $|X|$ is given the topology 
of the union of the $F_{q}|X|$.  
\end{definition} 

We remark that there is another 
realisation associated to a simplicial 
space $X$ called the \emph{fat} realisation 
of $X$ and denoted $||X||$.  $||X||$ is 
formed from the same disjoint union 
$\coprod \Delta^{n}\times X_{n}$ as $|X|$ 
but the equivalence relation 
$\equivalencerelation$ now only identifies 
$(\d^{i}(t),x)\equivalencerelation (t,d_{i}(x))$.  See 
\cite{Seg1} 
for a comparison of the two methods of 
realisation.  In particular, Segal proves 
that under certain conditions on the 
degeneracy operators of a simplicial 
space $X$, there is a homotopy equivalence 
$||X|| \simeq |X|$.  If $X$ is a simplicial 
manifold, then it is sufficient that the 
degeneracy operators $s_{i}$ embed $X_{n-1}$ as 
a submanifold of $X_{n}$,  see \cite{MosPer}.   

Also note that a simplicial map 
$f:X\to Y$ between simplicial spaces $X$ 
and $Y$ induces a map $|f|:|X|\to |Y|$ 
on the geometric realisations of $X$ and 
$Y$.  

We now examine the effect of realisation 
on our examples so far.  Let 
$EG = |N\overline{G}|$ and $BG = |NG|$.  
$EG$ is then the space consisting of all 
equivalence classes $|(x_{0},\ldots,x_{p})(g_{0},\ldots,g_{p})|$ 
with $(x_{0},\ldots,x_{p}) \in \Delta^{p}$ 
and $g_{i} \in G$ under the appropriate 
equivalence relations.  These are

\begin{multline*}
((0,x_{0},\ldots,x_{p-1}),(g_{0},\ldots,g_{p})) \equivalencerelation     
((x_{0},\ldots,x_{p-1}),((g_{1},\ldots,g_{p}))                          \\
\shoveleft{((x_{0},\ldots,x_{p-1},0), 
(g_{0},\ldots,g_{p})) \equivalencerelation 
((x_{0},\ldots,x_{p-1}),(g_{0},\ldots,g_{p-1})) }                      \\ 
\shoveleft{((x_{0},\ldots,x_{i-1},0,x_{i},\ldots,x_{p-1}),
(g_{0},\ldots,g_{p})) \equivalencerelation   
((x_{0},\ldots,x_{p-1}),
(g_{0},\ldots,\hat g_{i},\ldots,g_{p}))}                                  \\    
\shoveleft{((x_{0},\ldots,x_{i-1},x_{i}+x_{i+1},x_{i+2},\ldots,x_{p}),
(g_{0},\ldots,g_{p})) \equivalencerelation}                              \\ 
((x_{0},\ldots,x_{p}),
(g_{0},\ldots,g_{i-1},g_{i},g_{i},\ldots,g_{p}))          
\end{multline*}

$BG$ is then the space consisting of all 
equivalence classes 
$$
|(x_{0},\ldots,x_{p}),[g_{1},\ldots,g_{p}]|
$$ 
with the equivalence relations given by

\begin{multline*}
((0,x_{0},\ldots,x_{p-1}),[g_{1},\ldots,g_{p}]) \equivalencerelation  
((x_{0},\ldots,x_{p-1}),[g_{2},\ldots,g_{p}])                          \\
\shoveleft{((x_{0},\ldots,x_{p-1},0),[g_{1},\ldots, 
g_{p}]) \equivalencerelation 
((x_{0},\ldots,x_{p}),[g_{1},\ldots,g_{p-1}])}                         \\ 
\shoveleft{((x_{0},\ldots,x_{i-1},0,x_{i},\ldots,x_{p-1}),
[g_{1},\ldots,g_{p}]) \equivalencerelation   
((x_{0},\ldots,x_{p}),[g_{1},\ldots
,g_{i-1}g_{i},\ldots,g_{p}])}                                          \\       
\shoveleft{((x_{0},\ldots,x_{i-1},x_{i}+x_{i+1},x_{i+2},\ldots,x_{p}),
[g_{1},\ldots,g_{p-1}]) \equivalencerelation }                          \\ 
((x_{0},\ldots,x_{p}),
[g_{1},\ldots,g_{i-1},1,g_{i},\ldots,g_{p-1}])                          
\end{multline*}
These coordinates are called the \emph{homogenous} 
coordinates on $EG$ and $BG$.   
Recall the simplicial map $Np:N\bar{G}\to NG$ 
defined above.  $Np$ induces a map, also denoted 
$p$, on geometric realisations 
$p:EG\to BG$.  
This map is given in the 
homogenous coordinates by 
$$
p(|(x_{0},\ldots,x_{p}),(g_{0},\ldots,g_{p})|) = 
|(x_{0},\ldots,x_{p}),[g_{0}^{-1}g_{1},\ldots,g_{p-1}^{-1}g_{p}]|.
$$
We will come back to this later.

Of particular interest to us is the following 
fact: the geometric realisation of a 
simplicial group is a group.  For a proof 
of this see \cite{May2}. 
It is not hard to see that 
the simplicial space $N\bar{G}$ is a 
simplicial group while $NG$ is a simplicial 
group if and only if $G$ is abelian.  Thus 
we have the following result (see \cite{Seg3}).  

\begin{proposition}[\cite{Seg3}]  
\label{prop:2.2.2}
$EG$ is a group and contains $G$ as 
a closed subgroup.
\end{proposition} 

To make sense of the group structure on $EG$ it 
is convenient to introduce new coordinates --- the 
so-called non-homogenous coordinates --- (see \cite{Gaj}).  
If 
$(x_{0},\ldots,x_{p}) \in \Delta^{p}$ then let 
$t_{i} = x_{0} + \cdots + x_{i-1}$ for 
$i = 1,\ldots,p$.  Then $0\leq t_{1}\leq \cdots \leq t_{p} \leq 1$.  
An equivalence class $|(x_{0},\ldots,x_{p}),(g_{0},\ldots,g_{p})|$ in 
$EG$ is written in the non-homogenous 
coordinates as $|\tsimp,\hsimp|$ 
where we have set 
$h_{0} = g_{0},\ h_{i} = g_{i-1}^{-1}g_{i},\ i \geq 1$.

Given an element $\xi$ of $EG$ written 
in the non-homogenous coordinates as 
$\xi = |t_{1},\ldots,t_{p},h_{0}[h_{1}|\cdots|h_{p}]|$ we may 
recover the description of $\xi$ in the 
homogenous coordinates by setting 
$x_{i} = t_{i+1} - t_{i}$ (define $t_{0} = 0$ 
and $t_{p+1} = 1$), $g_{0} = h_{0}$, 
$g_{1} = h_{0}h_{1}$, $g_{2} = h_{0}h_{1}h_{2}$, 
etc.  
Then $\xi = |(\xsimp ),(\gsimp )|$.  When we 
use non-homogenous coordinates on $EG$ the 
equivalence relations above take the following 
form (see \cite{Gaj}). 
\begin{multline*}
(0,t_{1},\ldots,t_{p-1},\hsimp ) \equivalencerelation               
(t_{1},\ldots,t_{p-1},h_{0}h_{1}[h_{2}|\cdots|h_{p}])          \\
\shoveleft{(t_{1},\ldots,\hat{t}_{i},\ldots,t_{p},\hsimp ) 
\equivalencerelation (\tsimp ,h_{0}[h_{1}|\cdots|h_{i}|1|h_{i+1}|
\cdots|h_{p}])        }                                          \\  
\shoveleft{(t_{1},\ldots,t_{i-1},t_{i},t_{i},t_{i+1},\ldots,
t_{p-1},\hsimp ) \equivalencerelation }                        \\ 
(t_{1},\ldots,t_{p-1},h_{0}[
h_{1}|\cdots|h_{i}h_{i+1}|\cdots|h_{p}])                      
\end{multline*} 
Similarly an element $\a \in BG$ is expressed 
in the non-homogenous coordinates as 
$\a = |\tsimp, \bhsimp|$ where $t_{i}$ and $h_{i}$ 
are as above.  The above equivalence relations 
then take the form (see \cite{Gaj})   
\begin{multline*} 
(\tsimp , \bhsimp )          \\  
=    \begin{cases} 
     (t_{2},\ldots,t_{p},[h_{2}|\cdots
     |h_{p}]) & \text{for $t_{1}=0$ or 
     $h_{1}=1$}                       \\
     (t_{1},\ldots,\hat t_{i},\ldots,t_{p},
     [h_{1}|\cdots|h_{i}h_{i+1}|\cdots|
     h_{p}]) & \text{for $t_{i}=t_{i+1}$ or 
     $h_{i}=1$}                           \\
     (t_{1},\ldots,t_{p-1},[h_{1}|\cdots|
     h_{p-1}]) & \text{for $t_{p}=1$ or 
     $h_{p}=1$}.                          
     \end{cases}                  
\end{multline*}
The map $p:EG \to BG$ is expressed 
in the non-homogenous coordinates as 
$$ 
|\tsimp,\hsimp| \mapsto    
|\tsimp,\bhsimp|           
$$ 

We can now write down the group law on 
$EG$.  We have 
\begin{multline}\label{eq:EGgrouplaw}   
|t_{1},\ldots,t_{p},h_{0}[h_{1}|\cdots |h_{p}]|\cdot 
|t_{p+1},\ldots,t_{p+q},h_{0}^{\prime}[h_{p+1}|\cdots| 
h_{p+q}]|                                                       \\
= |t_{\sigma(1)},\ldots,t_{\sigma(p+q)},h_{0}h_{0}^{\prime}
[k_{\sigma(1)}|\cdots|k_{\sigma(p+q)}]|,   
\end{multline} 
where 
$$
k_{\sigma(r)} = \begin{cases} 
                h_{\sigma(r)} \text{if $\sigma(r) > p$}           \\ 
                (\prod_{i \in S_{r}}h_{\sigma(i)})^{-1}h_{0}^{\prime -1}     
                 h_{\sigma(r)}h_{0}^{\prime} 
                 (\prod_{i \in S_{r}}h_{\sigma(i)}) 
                \ \text{if $\sigma(r) \leq p$}.                    \\
                \end{cases} 
$$
Here $S_{r} = \{i \leq r : \sigma(i) > p\}$ 
and where $\sigma \in S_{p+q}$ is 
the permutation on $p+q$ letters such 
that 
$$
0 \leq t_{\sigma(1)} \leq \cdots \leq t_{\sigma(p+q)} \leq 1.
$$
One can check that this is well defined.  
It is worth mentioning that Segal gives 
a much more elegant description of this 
group law in \cite{Seg3}.  
Segal takes the viewpoint 
that one may regard $EG$ as the space 
of step functions on $[0,1]$ with values 
in $G$, ie a $G$ valued function constant 
on each half open interval $(t_{i},t_{i+1}]$ 
for some partition 
$0 =t_{0}<t_{1}<\cdots <t_{p}<t_{p+1}=1$ 
of $[0,1]$.  The group law on $EG$ is 
then the pointwise multiplication of these 
step functions.    

The unit of $EG$ is the obvious one and 
an element 
$\xi = |t_{1},\ldots,t_{p},h_{0}[h_{1}|\cdots|h_{p}]|$ 
has inverse given by 
$$
\xi^{-1} = |t_{1},\ldots,t_{p},h_{0}^{-1}[k_{1}|\cdots|
k_{p}]|,
$$
where 
$$
k_{i} = h_{0}h_{1}\cdots h_{p-1}h_{p}^{-1}h_{p-1}^{-1}\cdots 
h_{1}^{-1}h_{0}^{-1}.
$$

It is worth noting that there is a 
method of calculating the cohomology 
of the geometric realisation $|X|$ 
of a simplicial space $X$, at least 
in the case where the homotopy 
equivalence $||X|| \simeq |X|$ is 
available, from the following 
theorem. 

\begin{theorem} 
\label{thm:2.2.3}  
Let $X$ be a simplicial space.  Then 
there is an isomorphism 
$$
H^{\bullet}(||X||;\Real) \cong H^{\bullet}
(S^{\bullet,\bullet}(X);\Real), 
$$
where $H^{\bullet}(S^{\bullet,\bullet}(X);\Real)$ denotes 
the total cohomology of the 
double complex $S^{p,q}(X) = S^{q}(X_{p};\Real)$ 
with boundary operators $d$ and 
$\d$.  Here $d$ is the singular 
coboundary operator and 
$\d:S^{q}(X_{p};\Real)\to S^{q}(X_{p+1};\Real)$ 
is the simplicial coboundary formed by adding 
the pullback face operators 
$d_{i}^{*}$ with an alternating 
sign.  
\end{theorem} 

For a proof see either \cite{Dup}, 
or \cite{BotShuSta}.  
As explained in both of these 
papers, the singular functor $S$ 
can be replaced by the de Rham functor 
$\Omega$ in the above theorem if $X$ 
is a simplicial manifold.   

\section{Classifying spaces} 
\label{sec:2.3} 
Recall \cite{Hus}   
that a topological principal 
$G$ bundle for $G$ a topological 
group is defined as follows.
\begin{definition}[\cite{Hus}]  
\label{def:2.3.1}  
A \emph{topological principal $G$-bundle} 
is a triple $(\pi,P,M)$ where $P$ and 
$M$ are topological spaces such that 
\begin{enumerate} 
\item There is an effective action of 
 $G$ on the right of $P$ 
(effective means that $p\cdot g = p$ 
implies $g = 1$). 
\item $\pi:P\to M$ is a surjection whose 
fibres are the orbits of the $G$-action 
on $P$. 
\item Let $P^{[2]} = P\times_{M}P$ denote 
the fibre product of $P$ with itself.  Then 
the natural map $\tau:P^{[2]}\to G$ is continuous.  
We say that the action of $G$ on $P$ is 
\emph{strongly free}. 
\item There exist local sections of $\pi$, 
ie there is an open cover 
$\mathcal{U} = \{U_{i}\}_{i\in I}$ of $M$ 
together with sections $s_{i}:U_{i}\to P$ 
of $\pi$. 
\end{enumerate} 
\end{definition} 

Recall, see for example \cite{Hus}, that one way 
of defining a \emph{universal} $G$ bundle is to say that 
it is a principal $G$ bundle whose total space is 
contractible.  The base space of such a bundle is 
called a \emph{classifying space} and is unique up to 
homotopy equivalence.  The classifying space has the 
property that isomorphism classes of topological 
principal $G$ bundles are in a bijective correspondence 
with homotopy classes of maps from $M$ into the 
classifying space.  A choice of representative $f$ 
of the homotopy class corresponding to the bundle 
$P$ has the property that the pullback of the 
universal $G$ bundle by $f$ is isomorphic to $P$.  

It is not to difficult to see 
that the space $EG$ constructed 
above is contractible with
contracting homotopy 
$$  
h_{t}(|(x_{0},\ldots,x_{p}),(g_{1},\ldots,g_{p})|)  
 = |(1-t,tx_{0},\ldots,tx_{p}),(1,g_{1},\ldots,g_{p})|, 
$$  
as given in \cite{Dup}.    
Since 
$G$ is a subgroup of $EG$, there is 
a right action of $G$ on $EG$.  
By the formula~\ref{eq:EGgrouplaw} above, we 
have that the right action of $G$ 
on $EG$ is given by 
$$  
|t_{1},\ldots,t_{p},h_{0}[h_{1}|\cdots|h_{p}]|\cdot g   
= |t_{1},\ldots,t_{p},h_{0}g[h_{1}|\cdots|h_{p}]|.   
$$ 
Clearly $p$ is invariant under this right action of $G$. 

In \cite{Seg3},  
it is shown that 
$p:EG\to BG$ has local sections, 
provided that $G$ is a locally 
contractible group, and 
so $p:EG\to BG$ is a principal $G$ 
bundle and therefore a universal 
$G$ bundle.  We shall list some 
properties of $EG$ and $BG$ 
for $G$ abelian in the 
next section.                      
                            
\section{$B^{p}\cstar$ bundles} 
\label{sec:2.4} 
In this section we are interested 
in some properties of $EG$ and 
$BG$ for $G = \cstar$.  Because 
$\cstar$ is abelian, both $E\cstar$ 
and $B\cstar$ are groups --- in 
fact they are abelian groups.  The 
formula of equation~\ref{eq:EGgrouplaw} 
now takes a particularly simple form.  
We have for 
$|t_{1},\ldots,t_{p},z_{0}[z_{1}|\cdots|z_{p}]|$ 
and $|t_{p+1},\ldots,t_{p+q},z_{0}^{'}[z_{p+1}|\cdots|z_{p+q}]|$ 
elements of $E\cstar$ the following 
formula for their product (see \cite{Gaj}).  
\begin{multline*} 
|t_{1},\ldots,t_{p},z_{0}[z_{1}|\cdots|z_{p}]|\cdot 
|t_{p+1},\ldots,t_{p+q},z_{0}^{'}[z_{p+1}|\cdots|z_{p+q}]| \\ 
= |t_{\s(1)},\ldots,t_{\s(p+q)},z_{0}z_{0}^{'}[
z_{\s(1)}|\cdots|z_{\s(p+q)}]|,                    
\end{multline*} 
where $\s \in S_{p+q}$ is a permutation 
on $p+q$ letters such that 
$$
0\leq t_{\s(1)} \leq \cdots \leq t_{\s(p+q)} \leq 1.
$$
The group law on $B\cstar$ is given as follows 
(see \cite{Gaj}).  
If 
$|t_{1},\ldots,t_{p}[z_{1}|\cdots|z_{p}]|$ and 
$|t_{p+1},\ldots,t_{p+q}[z_{p+1}|\cdots|z_{p+q}]|$ 
are elements of $B\cstar$ then their product 
is given by 
\begin{multline*} 
|t_{1},\ldots,t_{p}[z_{1}|\cdots|z_{p}]|\cdot 
|t_{p+1},\ldots,t_{p+q}[z_{p+1}|\cdots|z_{p+q}]|    \\ 
= |t_{\s(1)},\ldots,t_{\s(p+q)}[z_{\s(1)}|\cdots| 
z_{\s(p+q)}]|,                                     
\end{multline*} 
where $\s$ is as above.   
Clearly the projection 
$p:E\cstar\to B\cstar$ is now a group 
homomorphism.  In fact there is a short 
exact sequence of groups 
$$
1\to \cstar \to E\cstar \to B\cstar \to 1.
$$
This is a generalisation of the 
exponential short exact sequence
$$
1\to \Z\to \C\to \cstar \to 1.
$$
Since $B\cstar$ is now an abelian 
topological group, we can iterate the 
classifying space construction to 
form the abelian topological group 
$BB\cstar = B^{2}\cstar$ and in general 
$B\cdots B\cstar = B^{p}\cstar$.  
Again we have the iterated exponential 
short exact sequence 
$$
1\to B^{p}\cstar \to EB^{p}\cstar\to B^{p+1}\cstar \to 1.
$$
In \cite{Gaj},  
Gajer develops a theory of smooth principal 
$B^{p}\cstar$ bundles.  
To do this he introduces a differentiable 
space structure on $EB^{p}\cstar$ and 
$B^{p}\cstar$ and defines the notion 
of a smooth (in this differentiable space 
sense) $B^{p}\cstar$ bundle.   
It is important to note that this 
differentiable space structure is not 
the same as a differentiable structure on 
a smooth manifold.  Indeed, this 
differentiable space structure on a 
Hausdorff space $M$ is defined 
in \cite{Che}  
(see also \cite{Mos}) by requiring that there is a family 
of continuous maps $\c :U\to M$
called \emph{plots}, each of 
whose domain is a convex set in some 
Euclidean space (the dimension of which 
need not be fixed), which satisfy 
\begin{itemize} 
\item if $\c:U\to M$ is a plot,  
$V$ a convex set and $\theta:V\to U$ a 
smooth map, then $\c \circ \theta:V\to M$ 
is also a plot,  
\item every constant map from a convex set 
into $M$ is also a plot, 
\item if $U$ is a convex set and 
$\{U_{i}\}_{i\in I}$ is a covering of 
$U$ by convex sets $U_{i}$ each of which is open 
in $U$, and $\c:U\to M$ is a map 
such that each restriction 
$\c_{i} = \c|_{U_{i}}$ is a plot, then 
$\c$ is also a plot.
\end{itemize} 
Clearly every smooth manifold has a 
differentiable space structure, but not 
every Hausdorff space with a differentiable 
space structure is a smooth manifold.       

Gajer then goes on to show how every $B^{p}\cstar$ 
bundle is induced by pullback from the 
universal $B^{p}\cstar$ bundle 
$EB^{p}\cstar \to B^{p+1}\cstar$ and 
shows that there is an bijection  
between isomorphism classes of $B^{p}\cstar$ 
bundles and $H^{p+2}(M;\Z)$.  
For us the most important thing about all 
of this is that $B^{2}\cstar$ provides 
us with an explicit realisation of an  
Eilenberg-Maclane space $K(\Z,3)$.  
We see this as follows.  From the exact 
homotopy sequence of the fibering 
$B^{p}\cstar\to EB^{p}\cstar\to B^{p+1}\cstar$ 
and the contractibility of $EB^{p}\cstar$ we 
get $\pi_{q}(B^{p+1}\cstar)=\pi_{q-1}(B^{p}\cstar)$.  
Hence by induction we get $\pi_{q}(B^{p}\cstar) = 
\pi_{q-p}(\cstar)$ if $q\geq p$ and $\pi_{q}(
B^{p}\cstar) = 0$ if $q<p$.  Hence 
$$ 
\pi_{q}(B^{p}\cstar) = \begin{cases} 
                       \Z \ \text{if}\ q = p + 1      \\ 
                       0 \ \text{otherwise}.            \\ 
                       \end{cases} 
$$ 
So $B^{p}\cstar$ is a $K(\Z,p+1)$.    
A theorem 
of algebraic topology  
asserts that 
for a given space $M$ and a cohomology 
class $\omega \in H^{3}(M;\Z)$ there 
is a map $f:M\to K(\Z,3)$, unique up 
to homotopy, such that $\omega$ is the 
pullback $f^{*}\iota$ of the fundamental 
class $\iota \in H^{3}(K(\Z;3);\Z) = \Z$ 
(see for example \cite{Spa}).  
Thus every class $\omega \in H^{3}(M;\Z)$ is the 
pullback of the generator of $H^{3}(B^{2}\cstar;\Z)$ 
by some map $M\to B^{2}\cstar$ which is 
unique up to homotopy.

\setcounter{chapter}{2}
\chapter{Review of bundle gerbes} 
\label{chapter:3} 
\section{Properties of $\cstar$ bundles} 
\label{sec:3.1} 

Let $\underline{G}_{M}$ denotes 
the sheaf of smooth $G$ valued functions 
on $M$ for $G$ an abelian Lie group.  
For the definitions and basic 
properties of sheaves  
see Chapter 1 of Brylinski's book \cite{Bry}.  For us, 
it is sufficient to know that a sheaf 
$\underline{\mathcal{S}}$  
of abelian groups on a manifold 
$M$ gives rise to sheaf cohomology 
groups $H^{i}(M;\underline{\mathcal{S}})$ for 
$i = 0,1, \ldots$ and that given an exact 
sequence of sheaves of abelian groups 
$$
0 \to \underline{\mathcal{A}} \to 
\underline{\mathcal{B}} \to 
\underline{\mathcal{C}} \to 0,
$$
on $M$, there is a long exact sequence of 
sheaf cohomology groups 
$$
\cdots \stackrel{\d}{\to} H^{n}(M;\underline{\mathcal{A}}) 
\to H^{n}(M;\underline{\mathcal{B}}) \to 
H^{n}(M;\underline{\mathcal{C}}) 
\stackrel{\d}{\to} H^{n+1}(M;\underline{
\mathcal{A}}) \to \cdots
$$
The homomorphism 
$\d:H^{n}(M;\underline{\mathcal{C}})\to H^{n+1}(M;\underline{\A})$ 
is called the connecting homomorphism.  
In \cite{Bry}[Corollary 1.1.9], Brylinski 
proves that 
the sequence 
of sheaves of abelian groups on $M$, 
$$
0\to \underline{\Z}_{M} \to \underline{\C}_{M} 
\to \underline{\C}^{\times}_{M}\to 0   
$$
is exact.  Hence there is a long exact sequence 
in sheaf cohomology
\begin{equation} 
\label{eq:longexactsheafcohom} 
\cdots \stackrel{\d}{\to} H^{n}(M;\underline{\Z}_{M})\to 
H^{n}(M;\underline{\C}_{M})\to 
H^{n}(M;\underline{\C}^{\times}_{M}) 
\stackrel{\d}{\to} H^{n+1}(M;\underline{\Z}_{M})\to 
\cdots 
\end{equation} 
In Chapter 1 of his book \cite{Bry}, 
Brylinski defines the notion of a soft 
sheaf of abelian groups and proves 
that the sheaf cohomology groups 
associated to a soft sheaf are all 
zero.  It is a fact that $\underline{\C}_{M}$ 
is a soft sheaf.  Since the sequence of 
sheaf cohomology groups~\ref{eq:longexactsheafcohom} 
is exact, we have in particular that 
$H^{1}(M;\underline{\C}^{\times}_{M})\cong H^{2}(M;\underline{\Z}_{M})$, 
the isomorphism being provided by the 
connecting homomorphism $\d$.  
Brylinski in \cite{Bry} proves that if $M$ 
is a paracompact manifold, as we will 
always assume, then \v{C}ech cohomology 
is canonically isomorphic to sheaf 
cohomology.  Thus we have isomorphisms 
$\check{H}^{1}(M;\underline{\C}^{\times}_{M})\cong 
H^{1}(M;\underline{\C}^{\times}_{M})$  
and $\check{H}^{2}(M;\underline{\Z}_{M})\cong 
H^{2}(M;\underline{\Z}_{M})$.  
Finally, it is shown in \cite{BotTu} that 
the \v{C}ech cohomology groups 
$\check{H}^{n}(M;\underline{\Z}_{M})$ and the 
singular cohomology groups 
$H^{n}(M;\Z)$ coincide.  Thus we 
have an isomorphism 
$\check{H}^{1}(M;\underline{\C}^{\times}_{M})\cong H^{2}(M;\Z)$.  
 
Recall that there are several ways of 
characterising principal $\cstar$ bundles.  
One such way is to describe a $\cstar$ bundle 
$P\to M$ by its transition functions 
$g_{ij}:U_{ij}\to \cstar$ relative to some 
open cover $\U =\{U_{i}\}_{i\in I}$ of $M$.  
These transition functions are then 
representatives of a class in the \v{C}ech 
cohomology group $\check{H}^{1}(M;\underline{\C}^{\times}_{M})$.  
By the above discussion, the transition 
functions $g_{ij}$ define a class in 
$H^{2}(M;\Z)$.   
This is the 
Chern class of the $\cstar$ 
bundle $P$. It is a standard 
fact, see \cite{KoNo}, that given smooth maps 
$g_{ij}:U_{ij}\to \cstar$ satisfying 
the cocycle condition
$g_{jk}g_{ik}^{-1}g_{ij} = 1$, 
one can construct a $\cstar$ 
bundle $\tilde{P}\to M$ 
with Chern class represented 
by $g_{ij}$ in 
$\check{H}^{1}(M;\underline{\C}^{\times}_{M})$.    

If $P_{1}$ and $P_{2}$ are two 
$\cstar$ bundles on $M$ described 
by transition functions 
$g_{ij}:U_{ij}\to \cstar$ and 
$h_{ij}:U_{ij}\to \cstar$ respectively, 
where the transition functions are 
both defined relative to the same
cover $\mathcal{U} = \{U_{i}\}_{i \in I}$ 
of $M$, then we can form a 
$\cstar$ bundle with 
transition functions equal to 
$g_{ij}h_{ij}$.  The reason this is 
possible is because $\cstar$ is 
abelian.  The resulting bundle 
we call the contracted product 
of $P_{1}$ and $P_{2}$ and denote 
by $P_{1}\otimes P_{2}$.   
Another way of constructing 
this bundle, see \cite{Bry}[Chapter 2], is to first form 
the fibre product $P_{1}\times_{M}P_{2}$ 
and then factor out by the subgroup 
of $\cstar \times \cstar$ 
consisting of all pairs $(z,z^{-1})$ 
for $z \in \cstar$.  The resulting 
$\cstar$ bundle is $P_{1}\otimes P_{2}$ 
and one can check that 
$P_{1}\otimes P_{2}$ has transition 
functions $g_{ij}h_{ij}$.   

Given a $\cstar$ bundle $P\to M$ 
with transition functions $g_{ij}$ 
then we can construct a $\cstar$ bundle 
$P^{*}\to M$ with transition functions 
$g_{ij}^{-1}$.  Again the reason this 
is possible is because $\cstar$ is 
a commutative group.  This bundle is easily 
seen to be the one obtained from $P$ 
by changing the action of $\cstar$ on 
$P$ to its inverse, that is the $\cstar$ 
bundle associated to $P$ via the 
isomorphism $\cstar \to \cstar$, $z\mapsto z^{-1}$.  
Notice that if we have a $\cstar$ bundle 
$P$ on $M$, then the $\cstar$ bundle 
$P\otimes P^{*}$ is canonically trivial (look at the 
transition functions). 

It is often more convenient to work 
with complex line bundles than with 
$\cstar$ principal bundles.  Given 
a principal $\cstar$ bundle $P\to M$ 
we may form a complex line bundle 
$L_{P}$ on $M$ by the associated 
bundle construction.  Recall, see \cite{KoNo}, 
that $L_{P}$ is formed from the 
product $P\times \C$ by factoring 
out by the equivalence relation 
$\equivalencerelation$ defined by 
$(p,u)\equivalencerelation (q,v)$ 
if and only if there is a $z\in \cstar$ 
such that $q = p\cdot z$ and $v = uz^{-1}$.  
Conversely, given a complex line bundle 
$L$ on $M$ we can form a principal 
$\cstar$ bundle on $M$ by forming the 
frame bundle associated to $L$ (see 
\cite{KoNo}).   
This correspondence gives an 
equivalence of categories between the 
category of principal $\cstar$ bundles 
on $M$ and the category of complex 
line bundles on $M$, see \cite{Bry}[Chapter 2].
Under this correspondence, 
the operation of contracted product 
on principal $\cstar$ bundles 
becomes the operation of tensor 
product on complex line bundles while 
the `inverse' of a principal $\cstar$ 
bundle corresponds to taking the dual of 
the associated complex line bundle. 

Recall, see for example \cite{KoNo}, 
that the pullback of a $\cstar$ 
bundle $\pi:P\to M$ by a smooth map 
$f:N\to M$ is the $\cstar$ bundle 
$f^{-1}\pi:f^{-1}P\to N$ where $f^{-1}P$ 
is the submanifold of $N\times P$ 
consisting of all pairs $(n,p)$ 
such that $f(n) = \pi(p)$.  
Every $\cstar$ bundle $P$ on $M$ is 
induced by pullback from 
the universal bundle $E\cstar \to B\cstar$ 
with a map $g:M\to B\cstar$, the 
classifying map of $P$.  This 
classifying map $g$ is 
unique up to homotopy.  If $P$ has 
transition functions $g_{ij}:U_{ij}\to \cstar$ 
relative to some open cover $\U = \{U_{i}\}_{i\in I}$ 
of $M$, then a choice of the 
classifying map $g:M\to B\cstar$ is given by 
\begin{equation} 
\label{eq:classifyingmap} 
g(m) = |\phi_{1}(m),\ldots,\phi_{p}(m), 
[g_{i_{0}i_{1}}(m)|\cdots|g_{i_{p-1}i_{p}}]|,
\end{equation} 
see for example \cite{Gaj}. 
Here $\phi_{n}(m) = \psi_{i_{0}}(m)+\cdots +\psi_{i_{n}}(m)$, 
where $\{\psi_{i}\}_{i\in I}$ is a 
partition of unity subordinated to the 
cover $\mathcal{U}$.  
We have the following proposition. 
\begin{proposition} 
\label{prop:3.1.1} 
Let $P$ and $Q$ be principal 
$\cstar$ bundles on a smooth manifold 
$M$.  Suppose that $P$ has a classifying 
map $f:M\to B\cstar$ and that $Q$ has 
a classifying map $g:M\to B\cstar$.  
Then the contracted product $P\otimes Q$ 
has a choice of classifying map given by 
$f\cdot g:M\to B\cstar$ and $P^{*}$ has 
a choice of classifying map given by 
$f^{-1}:M\to B\cstar$.  Here $f^{-1}$ 
means the map $\text{inv}\circ f$, 
where $\text{inv}:B\cstar\to B\cstar$ 
is the map $\xi\mapsto \xi^{-1}$.  
\end{proposition} 
\begin{proof}  
First of all, choose an open covering 
$\mathcal{U} = \{U_{i}\}_{i\in I}$ 
of $M$ such that $P$ and $Q$ have 
transition functions $g_{ij}:U_{ij}\to \cstar$ 
and $h_{ij}:U_{ij}\to \cstar$ respectively, 
both relative to the same open cover 
$\mathcal{U}$.  Then, as we have already 
seen, choices for $f$ and $g$ are 
given by 
$$
f(m) = |\phi_{1}(m),\ldots,\phi_{p}(m),[g_{i_{0}i_{1}}(m)|
\cdots |g_{i_{p-1}i_{p}}(m)]|,
$$
and 
$$
g(m) = |\phi_{1}(m),\ldots,\phi_{p}(m),[h_{i_{0}i_{1}}(m)|
\cdots |h_{i_{p-1}i_{p}}(m)]|, 
$$
where the $\phi_{n}$ are defined as above.  
If one performs the calculation, 
using the group law in the abelian 
group $B\cstar$ given in Section~\ref{sec:2.4}, 
then one finds that 
$$
f(m)g(m) = |\phi_{1}(m),\ldots,\phi_{p}(m),[g_{i_{0}i_{1}}(m)
h_{i_{0}i_{1}}(m)|\cdots |g_{i_{p-1}i_{p}}(m)h_{i_{p-1}i_{p}}(m)]|. 
$$
But this is precisely the classifying map we would 
have obtained for $P\otimes Q$ 
using the recipe of equation~\ref{eq:classifyingmap}.  
The other statement of the proposition 
is proved similarly.      
\end{proof} 
 
Recall, see \cite{KoNo}, that one 
way of defining a \emph{connection} on 
a principal $\cstar$ bundle 
$P\to M$ is to specify a family 
of one forms $\{\omega_{i}\}_{i\in I}$  
with $\omega_{i} \in \Omega^{1}(U_{i})$, 
for some open cover $\mathcal{U} = \{U_{i}\}_{i\in I}$ 
of $M$ relative to which $P$ has 
transition functions $g_{ij}:U_{ij}\to \cstar$, 
and to require the $\omega_{i}$ to 
satisfy 
$$
\omega_{j} = \omega_{i} + \frac{dg_{ij}}{g_{ij}} 
$$
on each overlap $U_{ij} = U_{i}\cap U_{j}$. 
A connection on a complex line 
bundle $L \to M$ is a linear map 
$\nabla:\Gamma(M,L)\to \Gamma(M,T^{*}M\otimes L)$ 
which satisfies the Liebniz rule 
$\nabla(fs) = df\otimes s + f\nabla(s)$ 
for all $s\in \Gamma(M,L)$ and all 
$f\in C^{\infty}(M,\C)$.  If $L_{1}$ 
and $L_{2}$ are complex line bundles 
on $M$ with connections $\nabla_{1}$ 
and $\nabla_{2}$ respectively, then 
the tensor product bundle $L_{1}\otimes L_{2}$ 
has a connection $\nabla_{1} + \nabla_{2}$ 
defined by (see \cite{Bry}[Chapter 2]) 
$$
(\nabla_{1} + \nabla_{2})(s_{1}\otimes s_{2}) = 
\nabla_{1}(s_{1})\otimes s_{2} + s_{1}\otimes 
\nabla_{2}(s_{2}).  
$$
Recall also that if $\nabla$ is a 
connection on a complex line bundle 
$L$ on $M$, then given any complex valued one form 
$A$ on $M$, we can define a new 
connection $\nabla + A$ on $L$ by 
$$
(\nabla + A)(s) = \nabla(s) + A\otimes s 
$$
for a section $s$ of $L$.  
Connections on complex line bundles are 
easier to get a handle on than connections 
on principal $\cstar$ bundles and so in 
most calculations we will work with 
connections on complex line bundles.  
This is justified by Proposition 2.2.5 
of \cite{Bry}.  

Let $P\to M$ be a $\cstar$ bundle,  
and let $f:N\to M$ and 
$g:O\to N$ be smooth 
maps.  
There is a canonical natural
isomorphism 
$\varphi_{f,g}:(f\circ g)^{-1}P\rightarrow g^{-1}f^{-1}P$, 
where by natural we mean natural with respect to 
isomorphisms $\psi:P_{1}\to P_{2}$.  
It has the property that if 
$\nabla_{L}$ is a connection on the 
associated line bundle $L = L_{P}$ 
over $M$, then 
$$
(f\circ g)^{-1} 
\nabla_{L} 
= \tilde{\varphi}_{f,g}\circ 
g^{-1}f^{-1}\nabla_{L}\circ 
\tilde{\varphi}_{f,g}, 
$$
where $\tilde{\varphi}_{f,g}$ 
denotes the isomorphism of line 
bundles induced by $\varphi_{f,g}$.  
This means that we can freely 
identify 
$(f\circ g)^{-1}P = g^{-1}f^{-1}P$ 
for any calculations. 

\section{Bundle gerbes and the Dixmier-Douady class}
\label{sec:3.2}                                                       
Suppose $Y \stackrel{\pi}{\to} M$ 
is a surjection admitting local 
sections.  We form the fibre 
product $Y \times_{M} Y =Y^{[2]}$ 
in the usual way.  That is, $Y^{[2]}$ 
is the set consisting of 
all pairs $(y_{1},y_{2}) \in Y \times Y$  
such that $\pi(y_{1}) = \pi(y_{2}).$  
We denote by $Y^{[3]}$ the triple 
fibre product $Y \times_{M} Y \times_{M} Y$ 
and in general we put 
$Y^{[p]}$ equal to the $p$-fold fibre product 
$Y \times_{M} \cdots \times_{M} Y.$  
Notice that we have 
projection maps $\pi_{i}:Y^{[p]} \to Y^{[p-1]}$ for 
$i = 1,\ldots,p$, obtained by omitting 
the $i^{th}$ factor.   
It is easy 
to see that in fact we have a 
simplicial space $Y = \{Y_{p}\}$ 
with $Y_{p} = Y^{[p+1]}$ and the obvious face and 
degeneracy operators.

\begin{definition}[\cite{Mur}]
\label{def:3.2.1} 
A \emph{bundle gerbe} consists of a triple 
of manifolds $(P,X,M)$ with a surjection   
$\pi : X \to M$ admitting
local sections and a $\cstar$-bundle 
$P \to X^{[2]}$,  
such that $P$ has a  
product; that is an isomorphism of  $\cstar$-bundles 
$m: \pi_{1}^{-1}P\otimes \pi_{3}^{-1}P \to \pi_{2}^{-1}P$
over $X^{[3]}$, covering the identity and such 
that the following diagram of $\cstar$ 
bundle isomorphisms over $X^{[4]}$ commutes: 
$$ 
\xymatrix{ 
\pi_{2}^{-1}\pi_{1}^{-1}P\otimes \pi_{4}^{-1}(
\pi_{1}^{-1}P\otimes \pi_{3}^{-1}P) 
\ar @2{-}[r] \ar[d]_{1\otimes \pi_{4}^{-1}m} & 
\pi_{1}^{-1}(\pi_{1}^{-1}P\otimes \pi_{3}^{-1}
P)\otimes \pi_{3}^{-1}\pi_{3}^{-1}P 
\ar[d]^{\pi_{1}^{-1}m\otimes 1}                  \\ 
\pi_{2}^{-1}\pi_{1}^{-1}P\otimes \pi_{4}^{-1}
\pi_{2}^{-1}P \ar @2{-}[d] & \pi_{1}^{-1}
\pi_{2}^{-1}P\otimes \pi_{3}^{-1}\pi_{3}^{-1}P 
\ar @2{-}[d]                                       \\ 
\pi_{2}^{-1}(\pi_{1}^{-1}P\otimes \pi_{3}^{-1}P) 
\ar[d]_{\pi_{2}^{-1}m} & \pi_{3}^{-1}(\pi_{1}
^{-1}P\otimes \pi_{3}^{-1}P) \ar[d]^{\pi_{3}
^{-1}m}                                            \\ 
\pi_{2}^{-1}\pi_{2}^{-1}P \ar @2{-}[r] & 
\pi_{3}^{-1}\pi_{2}^{-1}P.                           } 
$$ 
Here we have used the remarks made at the 
end of Section~\ref{sec:3.1} and the 
simplicial identities satisfied by the $\pi_{i}$.  
This last condition gives an `associativity' 
constraint on the product $m$, in the sense 
that we have 
$m(u\otimes m(v\otimes w)) = m(m(u\otimes v)\otimes w)$ 
for $u \in P_{(x_{3},x_{4})}$, 
$v \in P_{(x_{2},x_{3})}$ and 
$w \in P_{(x_{1},x_{2})}$, where 
$(x_{1},x_{2},x_{3},x_{4}) \in X^{[4]}$. 
\end{definition}

We will sometimes denote a bundle gerbe
$(P,X,M)$ by a single letter $P$.  A bundle 
gerbe is typically depicted as in the diagram
below.
$$
\begin{array}{ccc}
P         &          &        \\
\downarrow &         &         \\
Y^{[2]}   & \stackrel{\pi_{i}}{\dra} & Y \\
          &               &      \downarrow \pi \\
          &               &        M         \\
\end{array}
$$
We will need the definition of a \emph{simplicial 
line bundle} --- see \cite{BryMcL1}.  
\begin{definition}[\cite{BryMcL1}]
\label{def:3.2.2} 
Let $X = \{X_{p}\}$ be a 
simplicial manifold.  Let $L$ 
be a line bundle on $X_{p}$.  
Denote by $\d(L)$ the line 
bundle 
$$
d_{0}^{-1}L\otimes d_{1}^{-1}L^{*}\otimes 
d_{2}^{-1}L\otimes \cdots   
$$
on $X_{p+1}$, where the last factor in 
the above tensor product line bundle 
is $d_{p+1}^{-1}L$ if $p$ is even and 
$d_{p+1}^{-1}L^{*}$ if $p$ is odd or zero, 
where $d_{i}:X_{p+1}\to X_{p}$ 
are the face operators.  Notice that the 
line bundle $\d\d(L)$ on $X_{p+2}$ is 
canonically trivial.  This is a result of the 
commutation relations satisfied by the 
$d_{i}$ which allow us to contract off 
pairs of line bundles in the tensor 
product $\d\d(L)$.  Clearly $\d$ gives 
rise to a functor from the category of 
line bundles on $X_{p}$ to the category 
of line bundles on $X_{p+1}$.  If $s$ is a section of $L$, 
then it induces a section $\d(s)$ of 
$\d(L)$ defined by 
$$
d_{0}^{-1}s\otimes d_{1}^{-1}s^{*}\otimes 
d_{2}^{-1}s\cdots    
$$
where, as above, the last factor in this 
induced section is $d_{p+1}^{-1}s$ if $p$ is even 
and $d_{p+1}^{-1}s^{*}$ if $p$ is odd or zero.   
A \emph{simplicial line 
bundle} on $X$ consists of a line 
bundle $L$ on $X_{1}$, together with 
a non-vanishing section $s$ of $\d(L)$ 
on $X_{2}$ such that $\d(s)$ equals the 
canonical non-vanishing section 
$\underline{1}$ of $\d\d(L)$ on $X_{3}$.  
\end{definition} 

It is easy to see that the bundle 
gerbe product $m_{P}$ on a bundle 
gerbe $(P,X,M)$ forces the $\cstar$ 
bundle 
$$
\pi_{1}^{-1}P\otimes \pi_{2}^{-1}P^{*}
\otimes \pi_{3}^{-1}P\to X^{[3]} 
$$
to be trivial, and indeed one can define 
a section $s$ of this bundle as follows.  
Choose $u\in P_{(x_{2},x_{3})}$ and 
$v\in P_{(x_{1},x_{2})}$, then $m_{P}(u\otimes 
v)\in P_{(x_{1},x_{3})}$.  We put 
$$
s(x_{1},x_{2},x_{3}) = u\otimes m_{P}(u\otimes v)^{*}
\otimes v. 
$$
One can check that this is well defined and 
one can also show that the associativity of 
$m_{P}$ is equivalent to the following 
condition on $s$: 
$$
\label{eq:coherency} 
\pi_{1}^{-1}s\otimes \pi_{2}^{-1}s^{*}\otimes 
\pi_{3}^{-1}s\otimes \pi_{4}^{-1}s^{*} = 
\underline{1}. 
$$
It follows that the bundle gerbe $P$ can be 
viewed as a simplicial line bundle on the 
simplicial manifold $X =\{X_{p}\}$ with 
$X_{p} = X^{[p+1]}$.  
Simplicial line bundles on the 
simplicial manifold $NG$ associated 
to the classifying space of a Lie 
group $G$ take a particularly 
simple form.  We have the following 
result. 

\begin{proposition}[\cite{BryMcL1}] 
\label{prop:3.2.3} 
A simplicial line bundle on the 
simplicial manifold $NG$ is a central 
extension of $G$ by $\cstar$. 
\end{proposition} 

Recall from \cite{Mur} that a bundle 
gerbe $(P,X,M)$ has two important 
algebraic structures; an identity and an 
inverse.  The identity is a section $e$ of the 
pullback $\cstar$ bundle $\Delta^{-1}P$ 
over $X$, where $\Delta:X\to X^{[2]}$ is the 
inclusion of the diagonal.  The inverse is 
an isomorphism given fibrewise by 
$P_{(x_{1},x_{2})}\to P_{(x_{2},x_{1})}$ and 
denoted $p\mapsto p^{-1}$.  The identity section 
and inverse satisfy the properties one would 
expect; for instance $m_{P}(p\otimes e)= p$ 
and $m_{P}(p\otimes p^{-1}) = e$.  
 
It is useful sometimes, particularly 
when we discuss bundle gerbe 
connections and curvature, to have a 
slightly different formulation of 
the definition of a bundle gerbe (see 
\cite{CaMur}).  
In Definition~\ref{def:3.2.1} above 
we replace the $\cstar$ bundle $P\to X^{[2]}$ 
by a line bundle $L$.  We then require 
that there is a line bundle isomorphism 
$m_{L}:\pi_{1}^{-1}L\otimes \pi_{3}^{-1}L\to \pi_{2}^{-1}L$. 
The product on $L$ defined by $m_{L}$ is 
then required to be associative in the 
sense described above.  The two 
definitions are easily seen to be 
equivalent.  

\begin{example} 
\label{ex:3.2.4} 
The basic example of a bundle gerbe is 
the \emph{tautological bundle gerbe} (see \cite{Mur}  
and also \cite{CaMuWa}).  
Let $M$ be a 
smooth manifold with a base point 
$m_{0}$ and suppose $\omega$ is a 
closed, integral three form on $M$.  Assume 
that $M$ is 2-connected.  Form the path 
fibration $\pi:\mathcal{P}M\to M$, where 
$\mathcal{P}M$ is the path space of $M$ 
so $\mathcal{P}M$ consists of all piecewise smooth 
paths $\c:[0,1]\to M$ with $\c(0) = m_{0}$.  
The projection $\pi$ is given by 
$\pi(\c) = \c(1)$.  The fibre product 
$\mathcal{P}M^{[2]}$ consists of all pairs of paths 
$(\c_{1},\c_{2})$ in $\mathcal{P}M$ 
which have the same endpoint.  
There is an identification of  
$\mathcal{P}M^{[2]}$ with $\Omega M$, the 
based, piecewise smooth loop space of $M$, which is given 
by sending a pair of paths $(\c_{1},\c_{2})$ 
to the loop obtained by going along 
$\c_{1}$ at double speed and then going 
along the path $\c_{2}$ at double speed 
but in the reverse direction.  Define 
a closed integral two form $F$ on 
$\mathcal{P}M^{[2]} = \Omega M$ by 
pulling back $\omega$ to $\Omega M\times [0,1]$ 
via the evaluation map $\text{ev}:\Omega M\times 
[0,1]\to M$ and integrating out $[0,1]$. Here 
the evaluation map $\text{ev}$ is defined by 
$\text{ev}(\c,t) = \c(t)$ for $\c\in \Omega M$ 
and $t\in [0,1]$.  We have  
$\pi_{0}(\mathcal{P}M^{[2]}) = 
\pi_{0}(\Omega M) = \pi_{1}(M) = 0$.  
Note also that $\mathcal{P}M^{[2]}$ is 
simply connected.   
We would like to construct the tautological 
$\cstar$ bundle on $\mathcal{P}M^{[2]}$ 
with Chern class $F$.   
So let $(\c_{1},\c_{2})$ 
be a point in $\mathcal{P}M^{[2]}$.   
A path from  
$(m_{0},m_{0})$ to $(\c_{1},\c_{2})$ 
is the same thing as a piecewise smooth homotopy 
$\mu:[0,1]\times [0,1]\to M$ with 
endpoints fixed from $\c_{1}$ to 
$\c_{2}$.  This means $\mu$ satisfies 
$\mu(0,t) = \c_{1}(t)$, $\mu(1,t) = \c_{2}(t)$, 
$\mu(s,0) = m_{0}$ and $\mu(s,1) = \c_{1}(1) = \c_{2}(1)$.  
Introduce an equivalence relation 
$\equivalencerelation$ on the space 
$\mathcal{P}(\mathcal{P}M^{[2]})\times \cstar$, 
where we identify paths in $\mathcal{P}M^{[2]}$ 
with homotopies $\mu$ as above, by 
saying $(\mu,z)$ is equivalent to 
$(\mu^{\prime},w)$ if and only if, for all 
piecewise smooth homotopies $H:[0,1]\times [0,1]\times [0,1]\to M$ 
with endpoints fixed between 
$\mu$ and $\mu^{\prime}$ we have 
$$
w = z\exp(\int_{I^{3}}H^{*}\omega).
$$
The homotopy $H$ is required to satisfy 
$H(0,s,t) = \mu(s,t)$, $H(1,s,t) = \mu^{\prime}(s,t)$, 
and for each fixed $r$, the map 
$(s,t)\mapsto H(r,s,t)$ is a piecewise smooth homotopy 
with endpoints fixed between $\c_{1}$ and 
$\c_{2}$.  We let 
$$
\mathcal{Q} = (\mathcal{P}(\mathcal{P}M^{[2]})
\times \cstar)/\equivalencerelation.
$$
$\mathcal{Q}$ is a $\cstar$ bundle over 
$\mathcal{P}M^{[2]}$.  The triple 
$(\mathcal{Q},\mathcal{P}M,M)$ defines 
a bundle gerbe with bundle gerbe 
product given by 
$$
[\mu,z]\otimes [\nu,w] \mapsto [\mu \circ \nu,zw]
$$
where the square brackets denote 
equivalence classes and where $\mu$ 
is a piecewise smooth homotopy with endpoints fixed 
between paths $\c_{2}$ and $\c_{3}$, 
$\nu$ is a piecewise smooth homotopy with endpoints 
fixed between paths $\c_{1}$ and $\c_{2}$, 
and where $\mu\circ \nu$ denotes the 
piecewise smooth homotopy given by 
$$
(\mu \circ \nu)(s,t) = \begin{cases} 
                      \nu(2s,t), s \in [0,1/2],       \\
                      \mu(2s-1,t), s \in [1/2,1].     \\ 
                      \end{cases} 
$$
The proof that this is well defined and 
associative can be found in \cite{CaMuWa}. 
\end{example} 

A particularly important example of 
this construction is when $M$ is $G$, a 
simple, simply connected compact Lie 
group.  In this case we can identify 
$\mathcal{P}G^{[2]}$ with the group of 
based, piecewise smooth loops in $G$, 
$\Omega G$, by taking a pair of paths 
$(\c_{1},\c_{2})$ in $\mathcal{P}G^{[2]}$ 
and forming the loop $\c_{2}^{o}\circ \c_{1}$, 
by which we mean the loop which follows 
$\c_{1}$ at double speed and then follows 
the path $\c_{2}$ from its endpoint to its 
start point at double speed.  We can also 
identify the $\cstar$ bundle $\mathcal{Q}\to 
\Omega G$ with the Kac-Moody group 
$\widehat{\Omega G}$ which forms part of the 
central extension of the loop group 
$\cstar \to \widehat{\Omega G}\to \Omega G$.  
This is an example of the following 
general method of constructing bundle gerbes 
from central extensions.  

\begin{example} 
\label{ex:3.2.5} 
A frequently occurring example of a bundle 
gerbe is the so-called \emph{lifting} 
bundle gerbe associated to a principal 
$G$ bundle $\pi:P\to M$ where $G$ is 
a Lie group with a central  
extension  
$$
\cstar \to \hat{G}\to G.
$$
The lifting bundle gerbe is then formed 
by pulling back the $\cstar$ bundle 
$\hat{G}$ to $P^{[2]}$ via the 
natural map $\tau:P^{[2]}\to G$ defined by 
$p_{2} = p_{1}\cdot \tau(p_{1},p_{2})$ for 
$p_{1}$ and $p_{2}$ in the same fibre of 
$P$.  If we let this bundle be $\hat{P}$ 
then it is easy to see that the triple 
$(\hat{P},P,M)$ is a bundle gerbe.  
For more details see \cite{Mur}.  
\end{example} 
 
Various operations are possible with bundle gerbes.  
For example, given a map $f:N \to M$ we can 
pullback a bundle gerbe $P = (P,X,M)$  
to obtain a new bundle gerbe 
$f^{-1}P = (f^{-1}P,f^{-1}X,N)$ called the pullback
of $P$ by $f$.  Thus $f^{-1}X$ is the 
manifold obtained by pulling back 
$\pi:X\to M$ with $f$ as shown in the 
diagram below.  
$$ 
\xymatrix{ 
f^{-1}X \ar[d]_{f^{-1}\pi} \ar[r]^{\bar{f}} & 
X \ar[d] \pi                                    \\ 
N \ar[r]^{f} & M              } 
$$ 
Then $f^{-1}P \to (f^{-1}X)^{[2]}$ is the 
pullback of $P\to X^{[2]}$ by the 
map $\bar{f}^{[2]}$.  Clearly $f^{-1}P$ 
inherits a bundle gerbe product from 
$P$.    
Suppose we have a second 
map $g:O\to N$.  Then we have 
two bundle gerbes $g^{-1}f^{-1}P$ 
and $(f\circ g)^{-1}P$ over 
$O$.  We will make the convention 
once and for all that these 
two bundle gerbes are the same.     

Given two bundle 
gerbes $(P,X,M)$ and $(Q,Y,M)$ we can form
their product.  This is a bundle gerbe 
over $M$ which is described as follows.  
We first form the fibre 
product $X \times_{M} Y \to M$.  There are 
two projection maps $p_{1}:X\times_{M}Y\to X$ 
and $p_{2}:X\times_{M}Y\to Y$.  We use these 
to form the pullback $\cstar$ bundles 
$(p_{1}^{[2]})^{-1}P$ and $(p_{2}^{[2]})^{-1}Q$ over 
$(X\times_{M}Y)^{[2]}$  and then 
form the contracted product  
$P \otimes Q = (p_{1}^{[2]})^{-1}P\otimes (p_{2}^{[2]})^{-1}Q$ 
over $(X\times_{M}Y)^{[2]}$.  
This is the bundle whose fibre 
at $((x_{1},y_{1}),(x_{2},y_{2}))$
is
$$
P_{(x_{1},x_{2})} \otimes Q_{(y_{1},y_{2})}.
$$
The triple $(P\otimes Q,X\times_{M}Y,M)$ 
is a bundle gerbe on $M$ which we call 
the product of the bundle gerbes $(P,X,M)$ 
and $(Q,Y,M)$.  If $X = Y$ then we can form 
the product above in a slightly 
different way.  If $P$ and 
$Q$ are as above, then we define 
a $\cstar$ bundle $P\otimes Q\to X^{[2]}$ 
by forming the contracted product. 
The triple $(P\otimes Q,X,M)$ is 
then a bundle gerbe over $M$. 

Given a bundle gerbe $(P,X,M)$ we can 
construct a bundle gerbe $(P^{*},X,M)$ 
called the `inverse' of $P$.  
$P^{*}\to X^{[2]}$ is the $\cstar$ 
bundle `inverse' to $P\to X^{[2]}$, 
see Section~\ref{sec:3.1}.   
$P^{*}$ inherits a bundle gerbe product  
via the bundle gerbe product on $P$ using 
the isomorphism $P\otimes P^{*}\cong X^{[2]}\times \cstar$.  

Recall from \cite{Mur} that there is a 
characteristic class associated to a bundle gerbe 
$(P,X,M)$ called the Dixmier-Douady class 
of $P$.  This is a class in $H^{3}(M;\Z)$.  
A \v{C}ech representative for this class (or 
rather a representative for its image in the \v{C}ech cohomology 
group $\check{H}^{2}(M;\underline{\C}^{\times}_{M})$ under 
the isomorphism 
$\check{H}^{2}(M;\underline{\C}^{\times}_{M})\cong H^{3}(M;\Z)$ ---   
see Section~\ref{sec:3.1}) 
can be constructed as follows --- see \cite{Mur}. 
Choose an open cover 
$\mathcal{U} = \{U_{i}\}_{i\in I}$ of $M$, all 
of whose intersections are either empty 
or contractible and such that there exist 
local sections $s_{i}:U_{i}\to X$ of $\pi$.  
Form maps $(s_{i},s_{j}):U_{ij}\to X^{[2]}$ 
which map $m\in U_{ij}$ to 
$(s_{i}(m),s_{j}(m))\in X^{[2]}$.  Denote by 
$P_{ij}$ the pullback $\cstar$ bundle 
$(s_{i},s_{j})^{-1}P$ on $U_{ij}$.  Since 
$U_{ij}$ is contractible, there exists a 
section $\sigma_{ij}:U_{ij}\to P_{ij}$ 
of $P_{ij}\to U_{ij}$.  Define a map 
$g_{ijk}:U_{ijk}\to \cstar$ by 
$$
m_{P}(\sigma_{jk}\otimes \sigma_{ij}) = 
\sigma_{ik}\cdot g_{ijk}.
$$
The associativity of the bundle gerbe 
product $m_{P}$ implies that $g_{ijk}$ 
satisfies the cocycle condition 
$$
g_{jkl}g_{ikl}^{-1}g_{ijl}g_{ijk}^{-1} = 1
$$
on $U_{ijkl}$. $g_{ijk}$ is then a 
representative for a class in the \v{C}ech 
cohomology group $\check{H}^{2}(M;\underline{\C}^{\times}_{M})$ 
and it therefore defines a class $[g_{ijk}]$ in $H^{3}(M;\Z)$ 
under the isomorphism 
$\check{H}^{2}(M;\underline{\C}^{\times}_{M})\simeq H^{3}(M;\Z)$ 
mentioned in Section~\ref{sec:3.1}.  
This class in $H^{3}(M;\Z)$ is the 
Dixmier-Douady class of the bundle gerbe 
$P$ and is usually denoted by 
$DD(P)$ or $DD(P,X)$.  

Suppose we have our usual map $\pi:Y\to M$ 
admitting local sections and also a $\cstar$ 
bundle $Q\to Y$.  Then we can construct a bundle 
gerbe $(\d(Q),Y,M)$ where $\d(Q)\to Y^{[2]}$ is 
the $\cstar$ bundle $\d(Q) = \pi_{1}^{-1}Q
\otimes \pi_{2}^{-1}Q^{*}$ of Definition~\ref{def:3.2.2}.  
$\d(Q)$ has an associative product because we 
can regard such a product as a section of $\d\d(Q)$ 
which satisfies the coherency condition~\ref{eq:coherency} 
in $\d\d\d(Q)$.  Such a section is provided by the 
canonical trivialisation $\underline{1}$ of $\d\d(Q)$ 
and one can check that $\d(\underline{1})$ matches 
the canonical trivialisation of $\d\d\d(Q)$.   

\begin{definition}[\cite{Mur}] 
\label{def:3.2.6} 
A bundle gerbe $(P,X,M)$ is said to be 
\emph{trivial} if there exists a $\cstar$ 
bundle $Q\to X$ together with an isomorphism 
$\phi:\d(Q)\to P$ such that the following 
diagram commutes: 
$$
\xymatrix{ 
\pi_{1}^{-1}\d(Q)\otimes \pi_{3}^{-1}\d(Q) 
\ar[r] \ar[d]_-{\pi_{1}^{-1}\phi\otimes 
\pi_{3}^{-1}\phi} & \pi_{2}^{-1}\d(Q) 
\ar[d]^-{\pi_{2}^{-1}\phi}                 \\ 
\pi_{1}^{-1}P\otimes \pi_{3}^{-1}P 
\ar[r]^-{m_{P}} & \pi_{2}^{-1}P,             } 
$$
where the map $\pi_{1}^{-1}\d(Q)\otimes 
\pi_{3}^{-1}\d(Q)\to \pi_{2}^{-1}\d(Q)$ 
denotes the bundle gerbe product in 
the bundle gerbe $(\d(Q),X,M)$.  The pair 
$(Q,\phi)$ is called a \emph{trivialisation} 
of the bundle gerbe $P$.  
\end{definition} 

The Dixmier-Douady class behaves as one 
would expect under operations such as 
pullbacks and the like.  It also has 
the property that a bundle gerbe $P$ is 
trivial if and only if its Dixmier-Douady 
class is zero.  More precisely, 
we have the following proposition.  

\begin{proposition}[\cite{Mur}] 
\label{prop:3.2.7} 
Let $(P,X,M)$ and $(Q,Y,M)$ be 
bundle gerbes with Dixmier-Douady 
classes $DD(P)$ and 
$DD(Q)$ respectively.  Then we have 
\begin{enumerate} 
\item The product bundle gerbe 
$(P\otimes_{M} Q,X\times_{M}Y,M)$ has 
Dixmier-Douady class 
$$
DD(P\otimes_{M} Q) = DD(P) + DD(Q). 
$$
\item If $P^{*}$ denotes the bundle gerbe 
inverse to $P$ then it has 
Dixmier-Douady class $DD(P^{*})$ 
equal to $-DD(P)$.  
\item If 
$f:N\to M$ is a smooth map then the 
Dixmier-Douady class $DD(f^{-1}P)$ 
of the pullback bundle gerbe 
$(f^{-1}P,f^{-1}X,N)$ is equal to 
$f^{*}DD(P)$.  
\item The Dixmier-Douady class $DD(P)$ 
is the obstruction to $P$ being trivial.  
\end{enumerate} 
\end{proposition} 

\begin{proposition}[\cite{Mur}]  
\label{prop:3.2.8} 
Suppose $(P,X,M)$ is the lifting bundle 
gerbe for some principal $G$ bundle 
$X\to M$ where $G$ is part of a central 
extension 
$$
\cstar\to \hat{G}\to G, 
$$
see example~\ref{ex:3.2.5}.  Then the 
Dixmier-Douady class of $P$ is the 
obstruction to lifting the structure group 
of $X$ to $\hat{G}$.  
\end{proposition} 

Finally we note that it is possible to 
consider `$A$ bundle gerbes' where $A$ 
is an abelian Lie group.  For this we need 
to know that there is a well defined notion 
of contracted product of principal $A$ bundles.  
To define the contracted product 
$P_{1}\otimes_{A}P_{2}$ of two principal 
$A$ bundles $P_{1}$ and $P_{2}$ we simply mimic 
the procedure for the case $A = \cstar$, ie we 
set $P_{1}\otimes_{A}P_{2}$ to be the quotient 
of the $A\times A$ bundle $P_{1}\times_{M}P_{2}$ 
by the anti-diagonal subgroup of $A\times A$ 
consisting of all pairs $(a,a^{-1})$.  Now we 
can define an $A$-bundle gerbe to be a triple 
$(P,X,M)$ where $P\to X^{[2]}$ is a principal 
$A$ bundle with an associative product $\pi_{1}^{-1}
P\otimes_{A}\pi_{3}^{-1}P\to \pi_{2}^{-1}P$.  
Such an $A$ bundle gerbe will then give rise to 
class in the \v{C}ech cohomology group $\check{H}^{2}(M;\underline{A}_{M})$ 
rather than $H^{3}(M;\Z)$.  
  
\section{Bundle gerbe connections and 
the extended \\ Mayer-Vietoris sequence} 
\label{sec:3.3}       
Let $\pi:Y\to M$ be a surjection 
admitting local sections.  Then we 
have seen that we can define a 
simplicial space $Y = \{Y_{p}\}$ 
with $Y_{p} = Y^{[p+1]}$ and face 
operators $\pi_{i}:Y^{[p+1]}\to Y^{[p]}$ 
given by omitting the $i$th factor 
as mentioned above.  Then for 
$p = 0,1,\ldots$ we can form a 
complex $K^{\bullet}$ with 
$K^{q} = \Omega^{p}(Y^{[q]})$ and 
$\d:K^{q}\to K^{q+1}$ given by the 
alternating sum of pullbacks 
$\pi_{i}^{*}$.  The simplicial 
identities satisfied by the 
$\pi_{i}$ ensure that $\d^{2} = 0$.  

\begin{proposition}[\cite{Mur}] 
\label{prop:3.3.1}  
The complex 
\begin{equation}  
\label{eq:genmayviet1}  
0\to \Omega^{p}(M)\stackrel{\pi^{*}}{\to} 
\Omega^{p}(Y)\stackrel{\d}{\to} \cdots 
\stackrel{\d}{\to} 
\Omega^{p}(Y^{[q]})\stackrel{\d}{\to} \cdots 
\end{equation} 
is exact for $p = 0,1,2,\ldots$.
\end{proposition} 
  
We also have an analogue of this 
proposition with the de Rham 
functor $\Omega$ replaced by the 
singular functor $S$.  For the 
definitions and basic properties 
of singular cohomology see for 
example \cite{Spa}.  Again we make 
$S^{p}(Y^{[q]};\Z)$ a complex 
for fixed $p$ by defining 
$\d:S^{p}(Y^{[q]};\Z)\to S^{p}(Y^{[q+1]};\Z)$ 
by the alternating sum of pullbacks 
$\pi_{i}^{*}$.  As before $\d^{2} = 0$ 
and we have the following 
proposition.  

\begin{proposition} 
\label{prop:3.3.2} 
For fixed $p = 0,1,2,\ldots$ the 
complex 
\begin{equation} 
\label{eq:genmayviet2}        
0\to S^{p}(M;\Z)\stackrel{\pi^{*}}{\to} 
S^{p}(Y;\Z)\stackrel{\d}{\to}\cdots 
\cdots \stackrel{\d}{\to} 
S^{p}(Y^{[q]};\Z)\stackrel{\d}{\to} 
\cdots
\end{equation}   
is acyclic.
\end{proposition} 
   
The following proof is due to 
M. Murray.
\begin{proof} 
Consider first the case of a map 
$\pi:Z\to X$ which is onto.  
Here $Z$ and $X$ only have to be sets 
and $\pi$ is a function.  Denote 
as usual by $Z^{[p]}$ the $p$-fold fibre 
product over $\pi$.  Then define 
$$
\d : \text{Map}(Z^{[p]},\Z)\to 
\text{Map}(Z^{[p+1]},\Z) 
$$
by $\d = \sum_{i=0}^{p}(-1)^{i}\pi_{i}^{*}$ 
where 
$\pi_{i}(z_{1},\ldots,z_{p}) = (z_{1},\ldots,\hat{z}_{i},\ldots,z_{p})$. 
In the case that $X$ is a point there 
is a standard proof that this complex 
has no cohomology, see for example 
\cite{Spa}.  This can be adapted 
completely to the present case by 
choosing a section $s$ of $\pi$.  
Notice that because we are dealing 
with sets the question of existence 
of a section depends only on the 
fact that $\pi$ is onto.  Indeed if 
$f\in \text{Map}(Z^{[p]},\Z)$ and 
$\d(f) = 0$ then define 
$F\in \text{Map}(Z^{[p-1]},\Z)$ by 
$$
F(z_{1},\ldots,z_{p-1}) = 
f(s(\pi(z_{1})),z_{1},\ldots,z_{p-1}) 
$$
and it is straightforward to check that 
$\d(F) = f$.  Notice that if 
$f:Z\to \Z$ then $\d(f) = 0$ means that 
$f$ is constant on the fibres of $\pi$ 
so $f$ descends to a map $F:M\to \Z$.  
Clearly $f = \pi^{*}F$.   
Now return to the case at hand and 
define $X$ to be the space of all 
singular simplices in $M$ and $Z$ 
to be the space of all singular 
simplices in $Y$.  $\pi:Y\to M$ 
induces a surjective map 
$\pi_{*}:Z\to X$.  Notice that the 
space of all singular simplices in 
$Y^{[p]}$ is just $Z^{[p]}$ with respect 
to $\pi_{*}$.  To prove the 
proposition we just have to notice 
that the space of all singular 
cochains in $Y^{[p]}$ is just 
$\text{Map}(Z^{[p]},\Z)$.  
\end{proof} 

\begin{note} 
One could think of the 
complexes~\ref{eq:genmayviet1} 
and~\ref{eq:genmayviet2} above 
as generalisations of the 
Mayer-Vietoris exact sequence 
of an open cover (see \cite{BotTu}).  
One sees this by choosing an 
open cover $\{U_{i}\}_{i\in I}$ 
of $M$ and letting $Y$ be the 
disjoint union 
$$
\coprod_{i\in I}U_{i}.
$$
Then the fibre products $Y^{[2]}$ 
and $Y^{[3]}$ correspond to 
disjoint unions of double intersections 
$\coprod U_{ij}$ and disjoint unions 
of triple intersections 
$\coprod U_{ijk}$ and so on.  Thus we 
can identify the exact 
sequences~\ref{eq:genmayviet1} 
and~\ref{eq:genmayviet2} with 
the de Rham and singular Mayer-Vietoris 
sequences of an open cover respectively.  

In fact we can think of surjections 
$\pi:Y\to M$ admitting local sections 
as being `generalised open sets' 
of $M$.  This is the viewpoint adopted 
in Chapter 5 of \cite{Bry} where open sets 
are replaced by local homeomorphisms.  
In this framework fibre products 
correspond to intersections of open sets 
and the relationship between the extended  
Mayer-Vietoris sequences~\ref{eq:genmayviet1} 
and~\ref{eq:genmayviet2} and the more familiar 
Mayer-Vietoris sequences  
associated to an open covering of $M$ 
discussed in \cite{BotTu} becomes more 
apparent.  
\end{note} 

Recall from \cite{Mur} that 
there is a notion of \emph{bundle gerbe 
connection} and \emph{curving} for a bundle gerbe connection.  
We will briefly recall these definitions.  
Let $(P,X,M)$ be a bundle gerbe and 
let $L\to X^{[2]}$ denote the line bundle 
on $X^{[2]}$ associated to $P$.  We let 
$m_{L}$ denote the line bundle isomorphism 
induced by the bundle gerbe product $m_{P}$ 
on $P$.  As we noted earlier, $m_{L}$ 
defines a product on $L$ which is 
associative.  We make the following 
definition. 

\begin{definition}[\cite{Mur}]  
\label{def:3.3.3} 
A \emph{bundle gerbe connection} on $L$ 
is a connection $\nabla_{L}$ on $L$ 
such that 
\begin{equation} 
\label{eq:bundlegerbeconnection} 
\pi_{1}^{-1}\nabla_{L} + \pi_{3}^{-1}
\nabla_{L} = m_{L}^{-1}\circ \pi_{2}^{-1}
\nabla_{L}\circ m_{L}
\end{equation} 
on the line bundle 
$\pi_{1}^{-1}L\otimes \pi_{3}^{-1}L$ 
over $X^{[3]}$. 
\end{definition} 

$\nabla_{L}$ induces a connection 
$\nabla_{P}$ on $P$ and we say that 
$\nabla_{P}$ is a bundle gerbe 
connection on $P$.  Bundle gerbe 
connections always exist since we can 
choose a connection $\nabla_{L}$ on 
$L$ and then notice that the two 
connections 
$\pi_{1}^{-1}\nabla_{L} + \pi_{3}^{-1}\nabla_{L}$ 
and $m_{L}^{-1}\circ \pi_{2}^{-1}\nabla_{L}\circ m_{L}$ on 
$\pi_{1}^{-1}L\otimes \pi_{3}^{-1}L$ 
differ by the pullback of a one form 
$A$ on the base $X^{[3]}$.  The 
associativity condition on $m_{L}$ 
forces $\d(A) = 0$ and so we conclude 
by Proposition~\ref{prop:3.3.1} 
that there exists a one form $B$ on 
$X^{[2]}$ such that $\d(B) = A$.  The 
connection $\nabla_{L} - B$ on $L$ 
defines a bundle gerbe connection on 
$P$.  It is shown in 
\cite{Mur} that the space 
of all bundle gerbe connections on 
$P$ is an affine space for $\Omega^{1}(X)$.  

If we let $F_{\nabla_{L}}$ denote 
the curvature of $\nabla_{L}$ then 
equation~\ref{eq:bundlegerbeconnection} 
implies that we must have $\d(F_{\nabla_{L}}) = 0$, 
and so by Proposition~\ref{prop:3.3.1} 
again we see that there is a two 
form $f$ on $X$ such that 
$F_{\nabla_{L}} = \d(f)$.   
Note that $f$ is unique only up to pullbacks 
of two forms on $M$.  A choice of 
$f$ is called a \emph{curving} for the 
bundle gerbe connection $\nabla_{L}$.   
Therefore we have $\d(df) =0$ and 
hence there exists a closed three 
form $\omega$ on $M$ such that 
$df = \pi^{*}\omega$.  Furthermore 
it can be shown that $\omega$ is an 
integral three form, ie the integral 
of $\omega$ over any closed three 
manifold in $M$ is $2\pi i$ times 
an integer.  We call $\omega$ the 
\emph{three curvature} of the 
bundle gerbe connection and curving.  
In \cite{Mur} it is shown that 
$\omega$ is a de Rham representative 
for the image of the Dixmier-Douady class in 
$H^{3}(M;\Reals)$.  Note 
that a bundle gerbe can have a non-zero 
\v{C}ech representative for its Dixmier-Douady 
class but, due to the phenomenon of 
torsion, the de Rham representative 
for the image of the Dixmier-Douady 
class inside $H^{3}(M;\Reals)$ can be zero.   

\begin{example} 
\label{ex:tautbgconn} 
Recall from example~\ref{ex:3.2.4} the 
tautological bundle gerbe $(\mathcal{Q},
\mathcal{P}M,M)$ associated to 
a closed integral three form $\omega$ on 
a 2-connected manifold $M$.  In \cite{Mur} 
it is shown that there is a connection $\nabla$ on the 
$\cstar$ bundle $\mathcal{Q}\to \mathcal{P}M^{[2]}$ 
which is compatible with the bundle gerbe product 
on $\mathcal{Q}$ and whose curvature equals 
$\d(f)$ where $f$ is the two form on $\mathcal{P}M$ 
defined by $f = \int_{I}\text{ev}^{*}\omega$ (here 
$\text{ev}:\mathcal{P}M\times I\to M$ is the 
evaluation map sending $(\c,t)$ to $\c(t)$).  
Therefore $\nabla$ and $f$ is an example 
of a bundle gerbe connection and curving on 
$\mathcal{Q}$.  
\end{example} 
     
\section{Bundle gerbe morphisms and  the 
2-category of bundle gerbes} 
\label{sec:3.4} 
Recall the following definition from 
\cite{Mur}.   

\begin{definition}[\cite{Mur}]   
\label{def:3.4.1}
Let $(P,X,M)$ and $(Q,Y,N)$ be bundle 
gerbes on $M$ and $N$ respectively.  A 
\emph{bundle gerbe morphism} $\bar{f}:P\to Q$ 
consists of a triple 
$\bar{f} = (\hat{f},f,f^{'})$ of 
smooth maps $f^{'}:M\to N$, $f:X\to Y$ 
covering $f^{'}$ and $\hat{f}:P\to Q$ 
covering the induced map 
$f^{[2]}:X^{[2]}\to Y^{[2]}$.  The 
map $\hat{f}$ is required to commute 
with the bundle gerbe products on $P$ 
and $Q$.  
\end{definition} 

As an example, consider the situation 
in which we have a smooth map 
$f:N\to M$ and a bundle gerbe $(P,X,M)$.  
We can pullback the bundle gerbe $P$ 
to a bundle gerbe $(f^{-1}P,f^{-1}X,N)$ 
on $N$.  By definition of the pullback 
construction there is a bundle 
gerbe morphism 
$\bar{f}:f^{-1}P\to P$.  

We will almost always be interested in 
the case where $M = N$ and $f^{'}$ is the 
identity.  In this case, a bundle gerbe 
morphism $\bar{f}:P\to Q$  
is a triple 
$(\hat{f},f,\text{id}_{M})$.  

\begin{definition}[\cite{MurSte}] 
\label{def:3.4.2}
Let $(P,X,M)$ and $(Q,Y,M)$ be bundle 
gerbes.  A \emph{bundle gerbe isomorphism} 
$\bar{f}:P\to Q$ is a bundle gerbe morphism 
$\bar{f} = (\hat{f},f,\text{id}_{M})$ 
from $P$ to $Q$ such that $f:X\to Y$ is a 
diffeomorphism and $\hat{f}:P\to Q$ covering 
$f^{[2]}:X^{[2]}\to Y^{[2]}$ is a $\cstar$ 
bundle isomorphism. 
\end{definition} 

Unfortunately, bundle gerbe isomorphisms 
seldom appear in practice.  
Given bundle gerbes $(P,X,M)$ 
and $(Q,Y,M)$, it is not clear when 
a bundle gerbe morphism between them 
exists.  Obviously, a necessary 
condition is that $DD(P) = DD(Q)$.  
One way of dealing with this 
problem is to consider stable 
isomorphism classes of bundle gerbes 
on $M$.  

\begin{definition}[\cite{CaMicMur}, \cite{MurSte}]  
\label{def:3.4.3} 
We say that bundle gerbes $(P,X,M)$ 
and $(Q,Y,M)$ are \emph{stably isomorphic} 
if there exist $\cstar$ bundles 
$T_{1}$ and $T_{2}$ over $X\times_{M}Y$ 
such that there is an isomorphism 
$$
(p_{1}^{[2]})^{-1}P\otimes \d(T_{1}) 
\cong \d(T_{2})\otimes (p_{2}^{[2]})^{-1}Q
$$
covering the identity on 
$(X\times_{M}Y)^{[2]}$ and commuting 
with the bundle gerbe products.  (Here 
$p_{1}$ and $p_{2}$ are the two 
projections of $X\times_{M}Y$ onto 
$X$ and $Y$ respectively).  
\end{definition} 

Suppose bundle gerbes $(P,X,M)$ and 
$(Q,Y,M)$ have $DD(P)= DD(Q)$.  Then 
the product bundle gerbe 
$(P\otimes Q^{*},X\times_{M}Y,M)$ 
has $DD(P\otimes Q^{*}) = 0$.  
Therefore there is a $\cstar$ bundle 
$T\to X\times_{M}Y$ together with 
an isomorphism $P\otimes Q^{*}\to \d(T)$ 
covering the identity on 
$(X\times_{M}Y)^{[2]}$ and commuting 
with the respective bundle gerbe 
products.  Now form the bundle gerbe 
$(P\otimes Q^{*})\otimes (p_{2}^{[2]})^{-1}Q
\to (X\times_{M})^{[2]}$.  
The isomorphism above translates into 
$$
(P\otimes Q^{*})\otimes (p_{2}^{[2]})^{-1}Q 
\to \d(T)\otimes (p_{2}^{[2]})^{-1}Q.  
$$
We have 
$$
(P\otimes Q^{*})\otimes (p_{2}^{[2]})^{-1}Q 
= (p_{1}^{[2]})^{-1}P \otimes (p_{2}^{[2]})^{-1}
(Q^{*}\otimes Q).
$$
The bundle gerbe $Q^{*}\otimes Q\to Y^{[2]}$ 
has zero Dixmier-Douady class, hence it is 
trivial. Therefore its pullback to 
$(X\times_{M}Y)^{[2]}$ is trivial and so 
$P$ and $Q$ are stably isomorphic.  
Conversely, if $P$ and $Q$ are stably 
isomorphic, then it is easy to see that 
$DD(P) = DD(Q)$.  Thus we have the 
proposition (see \cite{MurSte}):  

\begin{proposition}[\cite{MurSte}]  
\label{prop:3.4.4} 
Two bundle gerbes $(P,X,M)$ and 
$(Q,Y,M)$ have the same 
Dixmier-Douady class if and only if 
they are stably isomorphic.
\end{proposition} 

This shows that the notion of 
stable isomorphism is an 
equivalence relation and so we can form 
the set of equivalence classes.  
Taking products of bundle gerbes 
clearly respects stable isomorphism 
and indeed there is a group 
structure on the set of stable 
isomorphism classes of bundle 
gerbes.  In Chapter~\ref{chapter:4}  
we shall see that there is an isomorphism 
between stable isomorphism classes 
of bundle gerbes and $H^{3}(M;\Z)$.  

There is another way of looking at 
the question of whether a bundle gerbe 
morphism exists between two given bundle 
gerbes or not. 
Before we explain this, we need to digress 
for a moment.  Let $P\to X$ be a 
principal $\cstar$ bundle on $X$ and 
let $\pi:X\to M$ be an onto map 
admitting local sections.  We have 
the following lemma  
(see \cite{Bry}[pages 185--186]).

\begin{lemma}[\cite{Bry}]      
\label{lemma:3.4.5} 
A necessary and sufficient condition 
for $P$ to descend to a bundle on $M$ 
is that there is an isomorphism 
$\phi:\pi_{1}^{-1}P\to \pi_{2}^{-1}P$ 
covering the identity on $X^{[2]}$ 
which satisfies the gluing law 
\begin{equation} 
\label{eq:descentisomorphism} 
\pi_{3}^{-1}\phi\circ \pi_{1}^{-1}\phi = \pi_{2}^{-1}\phi
\end{equation}  
over $X^{[3]}$.   
$\phi$ is known as 
a descent isomorphism.  $P$ descends 
to a bundle on $M$ means that there is 
a $\cstar$ bundle $Q\to M$ together 
with an isomorphism $P\simeq \pi^{-1}Q$.    
\end{lemma} 

Note that another way of looking at this 
is that a necessary and sufficient 
condition for $P$ to descend to a 
bundle on $M$ is that there is a 
section $\hat{\phi}$ of $\d(P)\to X^{[2]}$ 
such that $\d(\hat{\phi}) = \underline{1}$, 
where $\underline{1}$ is the 
canonically trivialising section of 
$\d\d(P)$, see Definition~\ref{def:3.2.2}. 
We have the following slight generalisation 
of this lemma.  

Suppose $\pi_{X}:X\to M$ and $\pi_{Y}:Y\to M$ 
are surjections admitting local sections.  
Then $X\times_{M}Y$ is also a surjection 
admitting local sections.  Given a $\cstar$ 
bundle $P$ on $X\times_{M}Y^{[p]}$, let 
$\d_{Y}(P)$ denote the $\cstar$ bundle 
on $X\times_{M}Y^{[p+1]}$ whose fibre at a 
point $(x,y_{1},y_{2},\ldots,y_{p+1})$ of 
$X\times_{M}Y^{[p+1]}$ is 
$$
P_{(x,y_{1})}\otimes P_{(x,y_{2})}^{*}\otimes 
\cdots, 
$$
where the last factor is $P_{(x,y_{p+1})}$ 
if $p$ is even and $P_{(x,y_{p+1})}^{*}$ if 
$p$ is odd --- compare with Definition~\ref{def:3.2.2}.  
As in Definition~\ref{def:3.2.2} it is easy to 
check that $\d_{Y}\d_{Y}(P)$ is canonically 
trivialised.  We will denote this canonical 
trivialisation by $\underline{1}$.  

\begin{lemma} 
\label{lemma:gendescent} 
Suppose $\pi_{X}$ and $\pi_{Y}$ are as above and 
that $P$ is a $\cstar$ bundle on $X\times_{M}Y$ 
such that there is a section $\hat{\phi}$ of 
$\d_{Y}(P)$ on $X\times_{M}Y^{[2]}$ which 
satisfies the coherency condition $\d_{Y}(\hat{\phi}) 
= \underline{1}$ on $X\times_{M}Y^{[3]}$.  
Then $P$ descends to a bundle $Q$ on $X$, ie 
there is an isomorphism $P\simeq p_{1}^{-1}Q$, 
where $p_{1}:X\times_{M}Y\to X$ is the projection 
onto the first factor.  
\end{lemma} 

\begin{proof} 
Choose an open cover $\U = \{U_{\a}\}_{\a\in \Sigma}$ 
of $M$ such that there exist local sections 
$t_{\a}:U_{\a}\to Y$ of $\pi_{Y}$.  Let $X_{\a} = 
\pi_{X}^{-1}(U_{\a})$.  Define a $\cstar$ bundle 
$P_{\a}$ on $X_{\a}$ by $P_{\a} = \hat{t}_{\a}^{-1}P$, 
where $\hat{t}_{\a}:X_{\a}\to X\times_{M}Y$ sends 
$x\in X_{\a}$ to $(x,t_{\a}(m))$, where $m = \pi_{X}(x) 
\in U_{\a}$.  The section $\hat{\phi}$ of 
$\d_{Y}(P)$ induces a section $\hat{\phi}_{\a\b}$ 
of $P_{\b}\otimes P_{\a}^{*}$ and hence an 
isomorphism $\phi_{\a\b}:P_{\a}\to P_{\b}$.  
The coherency condition on $\hat{\phi}$ 
ensures that the $\phi_{\a\b}$ satisfy the 
glueing condition $\phi_{\b\c}\phi_{\a\b} = 
\phi_{\a\c}$.  Hence we can use the standard 
clutching construction to form a $\cstar$ bundle $Q$ 
from the $P_{\a}$ and one can check that 
there is an isomorphism $P\simeq p_{1}^{-1}Q$.  
\end{proof} 

As an example of this, suppose we have 
a closed, integral two form $F$ on a 
non-simply connected manifold $M$.  If 
$\tilde{M}$ denotes the simply connected 
universal covering space of $M$ then 
we can lift $F$ to a closed integral 
two form on $\tilde{M}$ which we denote 
by $\tilde{F}$.  We cannot construct the 
tautological line bundle on $M$ directly, 
but we can on $\tilde{M}$.  Let $\tilde{L}$ 
denote the line bundle on $\tilde{M}$ 
with curvature $\tilde{F}$ constructed by 
the tautological method.  We know that 
$\tilde{L}$ descends to a line bundle 
$L$ on $M$.  In \cite{CaMuWa} 
it is shown that if $\pi_{1}(M)$ has no 
non-trivial extensions by $\cstar$, then 
there is an isomorphism 
$\phi:\pi_{1}^{-1}\tilde{L}\to \pi_{2}^{-1}\tilde{L}$ 
satisfying~\ref{eq:descentisomorphism} 
above.  

Suppose we are given two bundle gerbes 
$(P,X,M)$ and $(Q,Y,M)$ together with 
bundle gerbe morphisms 
$\bar{f},\bar{g}:P\to Q$.   
Let $\bar{f} = (\hat{f},f,\text{id}_{M})$ 
and $\bar{g}= (\hat{g},g,\text{id}_{M})$.  Form 
a $\cstar$ bundle $\hat{D}\to X$ by 
letting $\hat{D} = (f,g)^{-1}Q$, where 
$(f,g):X\to Y^{[2]}$ is the map 
which sends $x\in X$ to 
$(f(x),g(x))\in Y^{[2]}$.  We have 
the following series of lemmas. 

\begin{lemma} 
\label{lemma:3.4.6}
$\hat{D}$ descends to a bundle $D_{\bar{f},\bar{g}}\to M$, 
that is, we have $\hat{D}\cong \pi_{X}^{-1}D_{\bar{f},\bar{g}}$, 
where $\pi_{X}$ denotes the projection 
$X\to M$. 
\end{lemma} 

\begin{proof}  
We first define a map 
$\phi:\pi_{1}^{-1}\hat{D}\to \pi_{2}^{-1}\hat{D}$.  
If $\xi \in \pi_{1}^{-1}\hat{D}_{(x_{1},x_{2})} = Q_{(f(x_{2}),g(x_{2}))}$ 
then choose $p \in P_{(x_{1},x_{2})}$ and 
define 
$$
\phi(\xi) = m_{Q}(m_{Q}(\hat{g}(p)^{*}\otimes \xi)\otimes 
\hat{f}(p)).
$$
Here $m_{Q}$ is the bundle gerbe 
multiplication in $Q$. 
Clearly $\phi$ is independent of the 
choice of $p$.  We need to show that 
$\phi$ is a descent isomorphism, ie 
that the following diagram commutes
$$
\xymatrix{ 
\pi_{1}^{-1}\pi_{1}^{-1}\hat{D} = 
\pi_{2}^{-1}\pi_{1}^{-1}\hat{D} \ar[dr]_{\pi_{1}
^{-1}\phi} \ar[rr]^{\pi_{2}^{-1}\phi} & & 
\pi_{2}^{-1}\pi_{2}^{-1}\hat{D} = \pi_{3}^{-1}
\pi_{2}^{-1}\hat{D}                            \\ 
& \pi_{1}^{-1}\pi_{2}^{-1}\hat{D} = \pi_{3}^{-1}
\pi_{1}^{-1}\hat{D} \ar[ur]_{\pi_{3}^{-1}\phi} } 
$$
So let $\xi \in Q_{(f(x_{3}),g(x_{3}))}$, 
$p \in P_{(x_{1},x_{2})}$, and 
$p^{\prime} \in P_{(x_{2},x_{3})}$.  
Then 
\begin{eqnarray*} 
\pi_{3}^{-1}\phi \circ \pi_{1}^{-1}\phi(\xi) 
& = & \pi_{3}^{-}\phi(m_{Q}(m_{Q}(\hat{g}(p^{\prime})^{*}
\otimes \xi)\otimes \hat{f}(p^{\prime})))         \\
& = & m_{Q}(m_{Q}(\hat{g}(p)^{*}\otimes m_{Q}(
m_{Q}(\hat{g}(p^{\prime})^{*}\otimes \xi)\otimes \hat{f}(
p^{\prime})))\otimes \hat{f}(p))                  \\
& = & \vdots       \\ 
& = & m_{Q}(m_{Q}(\hat{g}(m_{P}(p^{\prime}\otimes p))^{*} 
\otimes \xi)\otimes \hat{f}(m_{P}(p^{\prime}\otimes p))) \\
& = & \pi_{2}^{-1}\phi(\xi).                       \\ 
\end{eqnarray*} 
Therefore $\hat{D}$ descends to a bundle 
$D_{\bar{f},\bar{g}}$.  
\end{proof} 

\begin{lemma} 
\label{lemma:3.4.7}
Suppose we are given bundle gerbes 
$(P,X,M)$ , $(Q,Y,M)$ and three 
bundle gerbe morphisms 
$\bar{f}_{i} = (\hat{f}_{i},f_{i},\text{id}_{M}):P\to Q$ 
for $i = 1,2,3$.  Then we have 
the following isomorphism of 
$\cstar$ bundles on $M$ 
$$
D_{\bar{f}_{2},\bar{f}_{3}}\otimes D_{\bar{f}_{1},\bar{f}_{2}
}\cong D_{\bar{f}_{1},\bar{f}_{3}}.
$$
\end{lemma} 

\begin{proof} 
The bundle gerbe product in $Q$ 
defines an isomorphism 
$\psi:(f_{2},f_{3})^{-1}Q\otimes 
(f_{1},f_{2})^{-1}Q\to (f_{1},f_{3})^{-1}Q$.
We need to check that this isomorphism 
is compatible with the descent 
isomorphism defined in Lemma~\ref{lemma:3.4.6}.
That is, we need to check that the following 
diagram over $X^{[2]}$ commutes:
$$
\xymatrix{ 
\pi_{1}^{-1}(f_{2},f_{3})^{-1}Q\otimes 
\pi_{1}^{-1}(f_{1},f_{2})^{-1}Q       
\ar[r]^-{\pi_{1}^{-1}\psi} \ar[d]_-{\phi_{23}\otimes 
\phi_{12}} & \pi_{1}^{-1}(f_{1},f_{3})^{-1}Q 
\ar[d]^-{\phi_{13}}                           \\ 
\pi_{2}^{-1}(f_{2},f_{3})^{-1}Q\otimes 
\pi_{2}^{-1}(f_{1},f_{2})^{-1}Q     
\ar[r]^-{\pi_{2}^{-1}\psi} & 
\pi_{2}^{-1}(f_{1},f_{3})^{-1}Q, } 
$$
where 
$\phi_{ij}:\pi_{1}^{-1}(f_{i},f_{j})^{-1}Q\to 
\pi_{2}^{-1}(f_{i},f_{j})^{-1}Q$ 
is the descent isomorphism 
constructed in Lemma~\ref{lemma:3.4.6}. 
Let $\xi_{23}\in \pi_{1}^{-1}(f_{2},f_{3})^{-1}Q_{(x_{1},x_{2})}$, 
$\xi_{12}\in \pi_{1}^{-1}(f_{1},f_{2})^{-1}Q_{(x_{1},x_{2})}$, 
and $p \in P_{(x_{1},x_{2})}$.  Then 
\begin{eqnarray*} 
& & \pi_{2}^{-1}\psi\circ (\phi_{23}\otimes \phi_{12}) 
(\xi_{23}\otimes \xi_{12})                        \\                        
& = &  \pi_{2}^{-1}\psi(m_{Q}(m_{Q}(\hat{f}_{3}(p)^{*}    
\otimes \xi_{23})\otimes \hat{f}_{2}(p))\otimes        
 m_{Q}(m_{Q}(\hat{f}_{2}(p)^{*}\otimes \xi_{12})\otimes 
\hat{f}_{1}(p)))                                      \\
& = &  m_{Q}(m_{Q}(m_{Q}(\hat{f}_{3}(p)^{*}\otimes 
\xi_{23})\otimes \hat{f}_{2}(p))\otimes 
m_{Q}(m_{Q}(\hat{f}_{2}(p)^{*}\otimes \xi_{12})
\otimes \hat{f}_{1}(p)))                            \\
& = &  \cdots                                         \\
& = &  m_{Q}(m_{Q}(\hat{f}_{3}(p)^{*}\otimes 
m_{Q}(\xi_{23}\otimes \xi_{12}))\otimes 
\hat{f}_{1}(p))                                   \\
& = & \phi_{13}\circ \pi_{1}^{-1}\psi 
(\xi_{23}\otimes \xi_{12}).                     \\
\end{eqnarray*} 
\end{proof} 

\begin{lemma}
\label{lemma:3.4.8}
Given three bundle gerbes $(P,X,M)$ and 
$(Q,Y,M)$ and $(R,Z,M)$ together with 
pairs of bundle gerbe morphisms 
\begin{eqnarray*} 
\bar{f}_{i} = (\hat{f}_{i},f_{i},\text{id}_{M}):P\to Q \\
\bar{g}_{i} = (\hat{g}_{i},g_{i},\text{id}_{M}):Q\to R \\
\end{eqnarray*} 
for $i =1,2$.  Then we have the following 
isomorphism of $\cstar$ bundles on $M$.
$$
D_{\bar{g}_{1}\circ \bar{f}_{1},\bar{g}_{2}\circ \bar{f}_{2}}\cong 
D_{\bar{f}_{1},\bar{f}_{2}}\otimes D_{\bar{g}_{1},\bar{g}_{2}}.
$$
\end{lemma}

\begin{proof} 
First of all define a map 
$<g_{1},g_{2}>:Y^{[2]}\to Z^{[2]}$ 
by mapping a point $(y_{1},y_{2})$ 
of $Y^{[2]}$ to the point 
$(g_{1}(y_{1}),g_{2}(y_{2}))$ of 
$Z^{[2]}$.  We then have 
$(g_{1}\circ f_{1},g_{2}\circ f_{2}) = <g_{1},g_{2}>\circ (f_{1},f_{2})$.  
Note that 
$<g_{1},g_{2}>^{-1}R\cong (g_{2}^{[2]})^{-1}R\otimes \pi_{2}^{-1}\gR$.
We therefore have the following series 
of isomorphisms of $\cstar$ bundles
\begin{eqnarray*} 
(g_{1}\circ f_{1},g_{2}\circ f_{2})^{-1}R & = & 
(<g_{1},g_{2}>\circ (f_{1},f_{2}))^{-1}R    \\ 
& = & (f_{1},f_{2})^{-1}<g_{1},g_{2}>^{-1}R   \\ 
& =& (f_{1},f_{2})^{-1}((g_{2}^{[2]})^{-1}R\otimes 
\pi_{2}^{-1}(g_{1},g_{2})^{-1}R)         \\ 
& = & (f_{1},f_{2})^{-1}(g_{2}^{[2]})^{-1}R\otimes 
f_{2}^{-1}\pi_{Z}^{-1}D_{\bar{g}_{1},\bar{g}_{2}}       \\ 
& =& (f_{1},f_{2})^{-1}(g_{2}^{[2]})^{-1}R\otimes 
\pi_{X}^{-1}D_{\bar{g}_{1},\bar{g}_{2}},                  \\
\end{eqnarray*} 
where $\pi_{X}$ denotes the projection 
$\pi_{X}:X\to M$ and so on.  The $\cstar$ 
bundle map $\hat{g}_{2}:Q\to R$ pulls back 
to define an isomorphism 
$$ 
\psi = (f_{1},f_{2})^{-1}\hat{g}_{2}:\fQ\to 
(f_{1},f_{2})^{-1}(g_{2}^{[2]})^{-1}R.  
$$ 
We need to show that $\psi$ is 
compatible with the descent 
isomorphisms:
$$
\xymatrix{ 
\pi_{1}^{-1}(f_{1},f_{2})^{-1}Q \ar[r]^-{\pi_{1}
^{-1}\psi} \ar[d]_-{\phi_{Q}} & \pi_{1}^{-1}
(f_{1},f_{2})^{-1}(g_{2}^{[2]})^{-1}R 
\ar[d]^-{\phi_{R}}                                \\ 
\pi_{2}^{-1}(f_{1},f_{2})^{-1}Q \ar[r]^-{\pi_{2}
^{-1}\psi} & \pi_{2}^{-1}(f_{1},f_{2})^{-1}(
g_{2}^{[2]})^{-1}R.                              } 
$$ 
This is clearly true, since $\psi$ is 
induced by the bundle gerbe morphism 
$\hat{g}_{2}$.  Therefore there is an isomorphism 
$\pi_{X}^{-1}D_{\bar{g}_{1}\circ \bar{f}_{1},\bar{g}_{2}\circ \bar{f}_{2}}
\cong \pi_{X}^{-1}(D_{\bar{f}_{1},\bar{f}_{2}}\otimes 
D_{\bar{g}_{1},\bar{g}_{2}})$ 
which is compatible with the descent 
isomorphisms.  Therefore we must 
have 
$D_{\bar{g}_{1}\circ \bar{f}_{1},\bar{g}_{2}\circ \bar{f}_{2}}\cong 
D_{\bar{f}_{1},\bar{f}_{2}}\otimes D_{\bar{g}_{1},\bar{g}_{2}}$.
\end{proof} 

Motivated by the results above, we make the 
following definitions.  

\begin{definition} 
\label{def:3.4.9} 
Let $(P,X,M)$ and $(Q,Y,M)$ be bundle 
gerbes and suppose $\bar{f} = (\hat{f},
f,\text{id}_{M})$ and $\bar{g} = (\hat{
g},g,\text{id}_{M})$ are bundle gerbe 
morphisms $\bar{f},\bar{g}:P\to Q$.  
A \emph{transformation} of bundle gerbe 
morphisms $\theta:\bar{f}\Rightarrow 
\bar{g}$ is a section of the $\cstar$ 
bundle $D_{\bar{f},\bar{g}}$ on $M$. 
\end{definition} 

\begin{definition} 
\label{def:3.4.10} 
Let $(P,X,M)$ and $(Q,Y,M)$ be bundle 
gerbes.  A bundle gerbe \emph{equivalence} 
consists of a pair of bundle gerbe 
morphisms $\bar{f} = (\hat{f},f,\text{id}_{M})$ 
and $\bar{g} = (\hat{g},g,\text{id}_{M})$ 
with $\bar{f}:P\to Q$ and $\bar{g}:Q\to P$ 
such that there exist transformations 
of bundle gerbe morphisms $\theta_{1}:\bar{g}
\circ \bar{f}\Rightarrow \text{id}_{P}$ 
and $\theta_{2}:\bar{f}\circ \bar{g}
\Rightarrow \text{id}_{Q}$.  
\end{definition} 

Now we can look at bundle gerbe 
morphisms in a different light.  We can view 
the collection of all bundle gerbes on a manifold 
$M$ as objects of a \emph{2-category} $\mathcal{B}$ 
(for the definition of a 2-category we refer 
to Section~\ref{sec:10.1}).  
For given two bundle gerbes $P$ and $Q$ 
on $M$ we can define a category 
$\text{Hom}(P,Q)$ whose objects 
consist of all bundle gerbe morphisms 
from $P$ to $Q$.  Given two objects 
$\bar{f}_{1}$ and $\bar{f}_{2}$ 
of $\text{Hom}(P,Q)$, that is 1-arrows of 
$\mathcal{B}$, the arrows from 
$\bar{f}_{1}$ to $\bar{f}_{2}$, 
ie 2-arrows of $\mathcal{B}$, consist 
of all transformations $\theta_{f_{1},f_{2}}$ 
from $\bar{f}_{1}$ to $\bar{f}_{2}$.  
Composition of arrows in $\text{Hom}(P,Q)$ 
is by tensor product, with the identifications 
of Lemma~\ref{lemma:3.4.7}.  
(For details on 2-categories and 
bicategories see Chapter~\ref{chapter:10}).  
We have the following proposition.

\begin{proposition} 
\label{prop:3.4.11}
$\mathcal{B}$ as defined above 
is a 2-category.
\end{proposition} 

All we really need to do is show that there is 
composition functor, in other words if $P$, $Q$ 
and $R$ are three objects of $\mathcal{B}$, 
then there is a functor 
$$
m:\text{Hom}(Q,R)\times \text{Hom}(P,Q) 
\to \text{Hom}(P,R)
$$
which is associative.  If we define the 
action of $m$ on objects to be composition 
of bundle gerbe morphisms and the action of 
$m$ on 2-arrows to be tensor product, with 
the identifications of Lemma~\ref{lemma:3.4.8}, 
then this condition is clearly satisfied.  

As mentioned in \cite{Mur}, an annoying feature 
of bundle gerbes is that two bundle gerbes can 
have the same Dixmier-Douady class and yet fail 
to be isomorphic by a bundle gerbe isomorphism 
nor even equivalent by a bundle gerbe equivalence.  
There is a way to circumvent this problem by 
introducing the notion of a \emph{stable morphism} 
between bundle gerbes.  The following definition is due to 
M. Murray.  

\begin{definition} 
\label{def:3.4.12} 
Let $(P,X,M)$ and $(Q,Y,M)$ be bundle gerbes.  
A \emph{stable morphism} $f:P\to Q$ is a choice 
of trivialisation $(L_{f},\phi_{f})$ of the bundle gerbe 
$P^{*}\otimes Q$, ie a choice of $\cstar$ bundle 
$L_{f}\to X\times_{M}Y$ and a choice of 
isomorphism $\phi_{f}:\d(L_{f})\to P^{*}\otimes Q$ 
such that the diagram below commutes: 
$$
\xymatrix{ 
\pi_{1}^{-1}\d(L_{f})\otimes 
\pi_{3}^{-1}\d(L_{f}) \ar[rr]^-{\pi_{1}^{-1}
\phi_{f}\otimes \pi_{3}^{-1}\phi_{f}} 
\ar[d] & & \pi_{1}^{-1}(P^{*}\otimes Q)\otimes 
\pi_{3}^{-1}(P^{*}\otimes Q) 
\ar[d]                                         \\ 
\pi_{2}^{-1}\d(L_{f}) \ar[rr]^-{\pi_{2}^{-1}\phi_{f}} 
& & \pi_{2}^{-1}(P^{*}\otimes Q),                 } 
$$
see Definition~\ref{def:3.2.6}.   
Suppose we have two stable morphisms 
$f,g:P\to Q$ corresponding to choices of trivialisations 
$(L_{f},\phi_{f})$ and $(L_{g},\phi_{g})$ 
of $P^{*}\otimes Q$ respectively.  
Then we have an isomorphism $\phi_{g}^{-1}\circ 
\phi_{f}:\d(L_{f})\to \d(L_{g})$ which corresponds 
to a section $\hat{\phi}_{f,g}\in \Gamma(\d(L_{f}^{*}
\otimes L_{g}))$.  Because the isomorphisms $\phi_{f}$ 
and $\phi_{g}$ commute with the bundle gerbe 
products, the section $\hat{\phi}_{f,g}$ will 
satisfy the coherency condition $\d(\hat{\phi}_{f,g}) = 
\underline{1}$.  Therefore there is a 
$\cstar$ bundle $D_{f,g}$ on $M$ and 
an isomorphism $L_{f}^{*}\otimes L_{g} 
\simeq D_{f,g}$.  We define a transformation 
$\theta:f\Rightarrow g$ to be a section $\hat{\theta}$ of 
$L_{f}^{*}\otimes L_{g}$ such that $\d(\hat{\theta}) = 
\hat{\phi}_{f,g}$.  Thus $\hat{\theta}$ will 
descend to a section $\theta$ of $D_{f,g}$.   
\end{definition} 

\begin{note} 
\label{note:3.4.13} 
\begin{enumerate} 
\item Note that a bundle gerbe morphism 
$\bar{f}=(\hat{f},f,\text{id}_{M}):P\to Q$ 
gives rise to a stable morphism $f:P\to Q$ 
because we can define a trivialisation $L_{f}$ 
of $P^{*}\otimes Q$ by letting $L_{f}$ have 
fibre $L_{f}(x,y)$ at $(x,y)\in X\times_{M}Y$ 
equal to $Q_{(f(x),y)}$ --- the bundle gerbe 
product on $Q$ and the isomorphism $P\simeq 
(f^{[2]})^{-1}Q$ provides an isomorphism 
$Q_{(f(x_{2}),y_{2})}\otimes Q_{(f(x_{1}),
y_{1})}^{*}\simeq Q_{(y_{1},y_{2})}\otimes 
P_{(x_{1},x_{2})}^{*}$.  Suppose we have stable 
morphisms $f:P\to Q$ and $g:P\to Q$ arising 
from bundle gerbe morphisms $\bar{f},\bar{g}:
P\to Q$.  Then the $\cstar$ bundle $D_{f,g}$ of 
Definition~\ref{def:3.4.12} above is the $\cstar$ 
bundle $D_{\bar{f},\bar{g}}$ of Lemma~\ref{lemma:3.4.6}.  

\item Note that a stable morphism exists 
if and only if $P$ and $Q$ have the same 
Dixmier-Douady class.  

\item Suppose we have bundle gerbes $(P,X,M)$ 
and $(Q,Y,M)$ and stable morphisms $f,g,h:P\to Q$ 
corresponding to trivialisations $(L_{f},\phi_{f})$, 
$(L_{g},\phi_{g})$ and $(L_{h},\phi_{h})$ 
respectively.  Suppose also that we have 
transformations $\theta_{f,g}:f\Rightarrow g$ 
and $\theta_{g,h}:g\Rightarrow h$ corresponding 
to sections $\hat{\theta}_{f,g}\in \Gamma(L_{f}^{*}
\otimes L_{g})$ and $\hat{\theta}_{g,h}\in 
\Gamma(L_{g}^{*}\otimes L_{h})$ respectively.  
We would like to define a composed transformation 
$\theta_{g,h}\theta_{f,g}:f\Rightarrow h$.  We do 
this as follows.  Let $\theta_{g,h}\theta_{f,g}$ 
denote the image of the section $\hat{\theta}_{f,g}
\otimes \hat{\theta}_{g,h}$ of $L_{f}^{*}\otimes 
L_{g}\otimes L_{g}^{*}\otimes L_{h}$ under the 
isomorphism $\text{can}:L_{f}^{*}\otimes L_{g}\otimes 
L_{g}^{*}\otimes L_{h}\to L_{f}^{*}\otimes L_{h}$ 
induced by contraction.  To show that $\theta_{g,h}
\theta_{f,g}$ defines a transformation we need 
to show that $\d(\theta_{g,h}\theta_{f,g}) = 
\hat{\phi}_{f,h}$.  We have $\d(\hat{\theta}_{f,g}
\otimes \hat{\theta}_{g,h}) = \hat{\phi}_{f,g}
\otimes \hat{\phi}_{g,h}$ and $\hat{\phi}_{f,g}
\otimes \hat{\phi}_{g,h}$ is mapped to $\hat{\phi}_{f,h}$ 
under the isomorphism $\d(L_{f})^{*}\otimes \d(L_{g})
\otimes \d(L_{g})^{*}\otimes \d(L_{h})\to \d(L_{f})^{*}
\otimes \d(L_{h}) = \d(L_{f}^{*}\otimes L_{h})$.  
Therefore we can compose transformations.  This 
operation of composition is clearly associative.  
Note that there is an identity transformation 
$1_{f}:f\Rightarrow f$ and so it follows that we have 
defined a category $\text{Hom}_{s}(P,Q)$ whose objects 
consist of the stable morphisms $P\to Q$ and whose 
arrows are the transformations between these 
stable morphisms.  Finally note that all the arrows 
in $\text{Hom}_{s}(P,Q)$ are invertible.     
\end{enumerate} 
\end{note} 

Suppose we have bundle gerbes $(P,X,M)$, 
$(Q,Y,M)$ and $(R,Z,M)$ and stable morphisms 
$f:P\to Q$ and $g:Q\to R$ corresponding to 
trivialisations $(L_{f},\phi_{f})$ and 
$(L_{g},\phi_{g})$ of $P^{*}\otimes Q$ and 
$Q^{*}\otimes R$ respectively.  So $\phi_{f}$ 
is an isomorphism $P^{*}\otimes Q\simeq \d(L_{f}) 
=\pi_{1}^{-1}L_{f}\otimes \pi_{2}^{-1}L_{f}^{*}$ 
and $\phi_{g}$ is an isomorphism $Q^{*}\otimes R\simeq 
\pi_{1}^{-1}L_{g}\otimes \pi_{2}^{-1}L_{g}^{*}$.  
Both isomorphisms respect the bundle gerbe products.  
Fibrewise, $\phi_{f}$ is an isomorphism 
$$
P^{*}(x_{1},x_{2})\otimes Q(y_{1},y_{2}) 
\simeq L_{f}(x_{2},y_{2})\otimes L_{f}^{*}(x_{1},y_{1}) 
$$
which can be rewritten as 
$$
P(x_{1},x_{2})\otimes L_{f}(x_{2},y_{2})\simeq 
Q(y_{1},y_{2})\otimes L_{f}(x_{1},y_{1}).  
$$
Similarly we can write 
$$
Q(y_{1},y_{2})\otimes L_{g}(y_{2},z_{2})\simeq 
R(z_{1},z_{2})\otimes L_{g}(y_{1},z_{1}).  
$$
We want to define a trivialisation $(L_{g\circ f},
\phi_{g\circ f})$ of $P^{*}\otimes R$ corresponding 
to a stable morphism $g\circ f:P\to R$.  Consider the 
$\cstar$ bundle $L$ on $X\times_{M}Z\times_{M}Y^{[2]}$ 
whose fibre at a point $(x,z,y_{1},y_{2})\in X\times_{M}
Z\times_{M}Y^{[2]}$ is 
$$
L(x,z,y_{1},y_{2}) = L_{f}(x,y_{1})\otimes L_{g}(y_{2},z)
\otimes Q(y_{1},y_{2}).  
$$
We will show there is a section $\hat{\phi}$ of 
$\d_{Y^{[2]}}(L)$ on $X\times_{M}Z\times_{M}(Y^{[2]})^{[2]}$ 
satisfying $\d_{Y^{[2]}}(\hat{\phi}) = \underline{1}$.  
We construct $\hat{\phi}$ as an isomorphism 
given fibrewise by $L(x,z,y_{1},y_{2})\to 
L(x,z,y_{3},y_{4})$.  
\begin{eqnarray*} 
&   & L_{f}(x,y_{1})\otimes L_{g}(y_{2},z)\otimes 
Q(y_{1},y_{2})                                      \\ 
& \simeq & L_{f}(x,y_{1})\otimes L_{g}(y_{2},z)\otimes 
Q(y_{4},y_{2})\otimes Q(y_{3},y_{4})\otimes Q(y_{1},y_{3}) \\ 
& \simeq & L_{f}(x,y_{3})\otimes P(x,x)\otimes 
L_{g}(y_{2},z)\otimes Q(y_{4},y_{2})\otimes Q(y_{3},y_{4}) \\ 
& \simeq & L_{f}(x,y_{3})\otimes L_{g}(y_{2},z)
\otimes Q(y_{4},y_{2})\otimes Q(y_{3},y_{4})               \\ 
& \simeq & L_{f}(x,y_{3})\otimes L_{g}(y_{4},z)\otimes 
R(z,z)\otimes Q(y_{3},y_{4})                                \\ 
& \simeq & L_{f}(x,y_{3})\otimes L_{g}(y_{4},z)\otimes 
Q(y_{3},y_{4}).           
\end{eqnarray*} 
One can check that this isomorphism is coherent given 
a third point $(y_{5},y_{6})$ of $Y^{[2]}$.  
Hence by Lemma~\ref{lemma:gendescent} $L$ descends 
to a $\cstar$ bundle $L_{g\circ f}$ on $X\times_{M}Z$.  
We now need to show that $L_{g\circ f}$ trivialises 
$P^{*}\otimes R$.  Consider $\d(L)$ on $(X\times_{M}
Z\times_{M}Y^{[2]})^{[2]}$.  It has fibre at a 
point $((x_{1},z_{1},y_{1},y_{2}),(x_{2},z_{2},y_{1}^{'},
y_{2}^{'}))$ of $(X\times_{M}Z\times_{M}Y^{[2]})^{[2]}$ 
equal to
\begin{eqnarray*} 
&     & \d(L)((x_{1},z_{1},y_{1},y_{2}),(x_{2},z_{2},
y_{1}^{'},y_{2}^{'}))                                   \\ 
& = & L(x_{2},z_{2},y_{1}^{'},y_{2}^{'})\otimes 
L^{*}(x_{1},z_{1},y_{1},y_{2})                          \\ 
& = & L_{f}(x_{2},y_{1}^{'})\otimes L_{g}(y_{2}^{'},
z_{2})\otimes Q(y_{1}^{'},y_{2}^{'})\otimes L_{f}^{*}
(x_{1},y_{1})\otimes L_{g}^{*}(y_{2},z_{1})\otimes 
Q^{*}(y_{1},y_{2})                                               \\  
& \simeq & P^{*}(x_{1},x_{2})\otimes Q(y_{1},y_{1}^{'})
\otimes Q(y_{2},y_{2}^{'})\otimes R(z_{1},z_{2})\otimes 
Q(y_{1}^{'},y_{2}^{'})\otimes Q^{*}(y_{1},y_{2})              \\ 
& \simeq & P^{*}(x_{1},x_{2})\otimes R(z_{1},z_{2}).      
\end{eqnarray*} 
This isomorphism commutes with the bundle gerbe 
products on $\d(L)$ and $P^{*}\otimes R$ and is 
also compatible with the descent isomorphism for 
$\d(L)$ induced from the descent isomorphism for $L$.  
Hence this isomorphism descends to $(X\times_{M}Z)^{[2]}$ 
to provide a trivialisation $\phi_{g\circ f}:\d(L_{g\circ f})\simeq 
P^{*}\otimes R$.  We define the composite stable 
morphism $g\circ f:P\to R$ to be the 
trivialisation $(L_{g\circ f},\phi_{g\circ f})$.  

Given bundle gerbes $(P,X,M)$, $(Q,Y,M)$ and 
$(R,Z,M)$ and stable morphisms $f_{1},f_{2}:P\to Q$ 
and $g_{1},g_{2}:Q\to R$ together with transformations 
$\theta:f_{1}\Rightarrow f_{2}$ and $\rho:g_{1}\Rightarrow 
g_{2}$, one can define a composed transformation 
$\rho\circ \theta:g_{1}\circ f_{1}\Rightarrow g_{2}\circ f_{2}$ 
using a similar technique to that in the preceding 
paragraph (except that we now descend a morphism).  
One can show that this operation of composing 
transformations is compatible with the operations 
of composing transformations in the categories 
$\text{Hom}_{s}(P,Q)$ and $\text{Hom}_{s}(Q,R)$.  It follows that 
we have defined a functor 
$$
\circ :\text{Hom}_{s}(P,Q)\times \text{Hom}_{s}(Q,R)\to 
\text{Hom}_{s}(P,R).  
$$
We would like to know whether or not this 
operation of composition is associative.  
Suppose we have bundle gerbes $(P,X,M)$, $(Q,Y,M)$, 
$(R,Z,M)$ and $(S,W,M)$ together with stable 
morphisms $f:P\to Q$, $g:Q\to R$ and $h:R\to S$.  
Let $p:X\times_{M}Y^{[2]}\times_{M}Z^{[2]}\times_{M}W
\to X\times_{M}W$ denote the natural projection obtained 
by omitting the factors $Y^{[2]}$ and $Z^{[2]}$ in 
the fibre product.  Then there is an isomorphism 
$$
p^{-1}L_{h\circ (g\circ f)}(x,y_{1},y_{2},z_{1},z_{2},w)\simeq 
L_{h}(z_{2},w)\otimes L_{g}(y_{2},z_{1})\otimes L_{f}(x,y_{1})
\otimes Q(y_{1},y_{2})\otimes R(z_{1},z_{2}) 
$$
of fibres at a point $(x,y_{1},y_{2},z_{1},z_{2},w)$.  
At the same point there is an isomorphism 
$$
p^{-1}L_{(h\circ g)\circ f}(x,y_{1},y_{2},z_{1},z_{2},w)\simeq 
L_{h}(z_{2},w)\otimes L_{g}(y_{2},z_{1})\otimes 
L_{f}(x,y_{1})\otimes Q(y_{1},y_{2})\otimes R(z_{1},z_{2}).  
$$
It follows from this that $p^{-1}L_{h\circ (g\circ f)}
\simeq p^{-1}L_{(h\circ g)\circ f}$ and hence 
that $L_{h\circ (g\circ f)}\simeq L_{(h\circ g)\circ f}$.  
Furthermore this isomorphism is compatible with 
the isomorphisms $\phi_{h\circ (g\circ f)}$ and 
$\phi_{(h\circ g)\circ f}$ and hence defines a 
transformation $\theta_{f,g,h}:h\circ (g\circ f)
\Rightarrow (h\circ g)\circ f$.  One can show that this 
transformation $\theta_{f,g,h}$ satisfies a  
coherency condition given a fifth bundle gerbe $(T,V,M)$ 
and a fourth stable morphism $k:S\to T$.  This coherency 
condition is 
$$
(\theta_{g,h,k}\circ 1_{f})\theta_{f,h\circ g,k}
(1_{k}\circ \theta_{f,g,h}) = 
\theta_{f,g,k\circ h}\theta_{g\circ f,h,k}.  
$$
Given a stable morphism $f:P\to Q$ between 
bundle gerbes $P = (P,X,M)$ and $Q = (Q,Y,M)$, 
there is a natural stable morphism $\tilde{f}:Q\to P$.  
If $f$ is given by a trivialisation $(L_{f},\phi_{f})$ 
of $P^{*}\otimes Q$, then $\tilde{f}$ corresponds 
to the trivialisation $(L_{f}^{*},\phi_{f}^{*})$ 
of $Q^{*}\otimes P$.  Let us calculate the 
trivialisation of $P^{*}\otimes P$ corresponding 
to the stable morphism $\tilde{f}\circ f$.  First 
we form the $\cstar$ bundle $L$ on $X\times_{M}X\times_{M}
Y^{[2]}$ whose fibre at a point $(x_{1},x_{2},y_{1},y_{2})$ 
of $X^{[2]}\times_{M}Y^{[2]}$ is 
$$
L(x_{1},x_{2},y_{1},y_{2}) = L_{f}(x_{1},y_{1})\otimes 
L_{f}^{*}(x_{2},y_{2})\otimes Q(y_{1},y_{2}).  
$$
Since $(L_{f},\phi_{f})$ is a trivialisation of 
$P^{*}\otimes Q$, we see there is an isomorphism 
$$
L(x_{1},x_{2},y_{1},y_{2}) \simeq P(x_{1},x_{2}).  
$$
This isomorphism preserves the descent isomorphism for 
$L$ and so we see there is an isomorphism 
$L_{\tilde{f}\circ f}\simeq P$.  This is a transformation 
$\tilde{f}\circ f\Rightarrow 1_{P}$, where $1_{P}$ 
is the stable morphism $P\to P$ induced by the 
identity bundle gerbe morphism $P\to P$.  Similarly 
there is a transformation $f\circ \tilde{f}\Rightarrow 1_{Q}$.  
These transformations are compatible with the 
transformations $\theta_{f,g,h}$ in the sense 
of Definition~\ref{def:10.1.2}.  We have the 
following Proposition.  

\begin{proposition} 
\label{prop:stablebigrpd} 
Given a manifold $M$, we can form a bicategory 
$\mathcal{B}_{s}$ whose objects are the bundle gerbes on $M$ and where 
the 1-arrows are the objects of the categories 
$\text{Hom}_{s}(P,Q)$ and where the 2-arrows are the 
arrows of the categories $\text{Hom}_{s}(P,Q)$.  This 
bicategory is in fact a bigroupoid.  
\end{proposition}  

We refer to Definitions~\ref{def:10.1.1} 
and~\ref{def:10.1.2} for the precise meanings of 
the terms \emph{bicategory} and \emph{bigroupoid}.

\setcounter{chapter}{3}
\chapter{Singular theory of bundle gerbes}  
\label{chapter:4} 
\section{Construction of a classifying map}
\label{sec:4.1} 

\begin{definition} 
The \emph{universal 
bundle gerbe} $(\widetilde{EB\cstar},EB\cstar,BB\cstar)$ 
is the lifting bundle 
gerbe associated to the principal $B\cstar$ 
bundle $EB\cstar\to BB\cstar$ via the short exact 
sequence of abelian groups 
$$
\cstar\to E\cstar\to B\cstar.  
$$
So $\widetilde{EB\cstar}\to (EB\cstar)^{[2]}$ 
is the pullback $\tau^{-1}E\cstar$ of the universal $\cstar$ 
bundle $E\cstar\to B\cstar$ via the 
canonical map $\tau:(EB\cstar)^{[2]}\to B\cstar$.  
\end{definition} 

Let $\widetilde{EB\cstar}$ have Dixmier-Douady class 
$\omega$ in $H^{3}(BB\cstar;\Z) = \Z$.  
Suppose $(P,Y,M)$ is a bundle gerbe on $M$ 
with Dixmier-Douady class $DD(P)$ in 
$H^{3}(M;\Z)$.  
We will review a construction of 
\cite{Gaj} which shows how to construct  
a map $g:M \to BB\cstar$ from a \v{C}ech 
cocycle $g_{ijk}$ representing the Dixmier-Douady 
class of $P$.  It will follow from this 
that the two classes $DD(P)$ and $g^{*}\omega$ 
are equal.     
 
Suppose the class $DD(P)$ is 
represented by a \v{C}ech 2-cocycle  
$g_{ijk}:U_{ijk}\to \cstar$, relative 
to a Leray cover $\U = \{U_{i}\}_{i \in I}$ 
of $M$, in the \v{C}ech cohomology group 
$\check{H}^{2}(M;\underline{\C}^{\times}_{M})$ via 
the isomorphism induced by the short 
exact sequence of sheaves of abelian 
groups on $M$ 
$$
1\to \underline{\Z}_{M} \to 
\underline{\C}_{M} \to 
\underline{\C}^{\times}_{M} \to 1, 
$$
see Section~\ref{sec:3.1}.  Thus $g_{ijk}$ 
satisfies the \v{C}ech 2-cocycle 
condition 
\begin{equation} 
\label{eq:4.1.1}
g_{jkl}g_{ikl}^{-1}g_{ijl}g_{ijk}^{-1} = 1 
\end{equation} 
on $U_{ijkl}$.  We will construct maps 
$g_{ij}:U_{ij}\to B\cstar$ which satisfy 
the \v{C}ech 1-cocycle condition 
$g_{jk}g_{ik}^{-1}g_{ij} = 1$ on $U_{ijk}$.  
These will be transition functions for a 
$B\cstar$ bundle on $M$, whose classifying 
map will be the map $g:M\to BB\cstar$ 
referred to above.  

We will make the assumption that every 
manifold $M$ we deal with has a locally finite, 
countable Leray cover $\U$ such that there is a  
smooth partition of unity $\{\phi_{i}\}_{i \in I}$ 
subordinate to the open cover $\U$.  For 
each $m \in M$ form the subset 
$I_m = \{i_{0},i_{1},\ldots,i_{n(m)}\}$ of $I$ 
consisting of those $i \in I$ such that 
$\phi_{i}(m)$ is not zero.   
Since our Leray covering is locally 
finite, this subset will always be finite.  
For $r = 0,1,\ldots,n(m)$, let 
$$
\psi_{r}(m) = \phi_{i_{0}}(m) + \cdots +   
\phi_{i_{r-1}}(m).  
$$

Then we will have $\psi_{0}(m)\leq \cdots \leq \psi_{n(m)}(m)$ 
and using the non-homogenous 
coordinates on $B\cstar$, see Section~\ref{sec:2.2} 
and \cite{Gaj}, we can define 
$g_{ij}:U_{ij} \to B\cstar$ by 
$$
g_{ij}(m) = |\psi_{1}(m),\ldots,\psi_{n}(m),  
[g_{iji_{0}}(m)^{-1}g_{iji_{1}}(m)|\cdots | 
g_{iji_{n-1}}(m)^{-1}g_{iji_{n}}(m)]|
$$
where $m\in U_{ij}$, $n = n(m)$ and $\{i_{0},\ldots,i_{n(m)}\}$ 
is the set $I_{m}$ above.  
Notice that this makes sense, that is 
respects the face and degeneracy operators, precisely 
because each $g_{iji_{r}}^{-1}g_{iji_{s}}$ 
defines a trivial cocycle.  By Lemma 1.2 of 
\cite{Gaj} we have that $g_{ij}:U_{ij}\to B\cstar$ 
is smooth --- that is smooth in the 
differentiable space sense of \cite{Che} and 
\cite{Mos} --- see also Section~\ref{sec:2.4} ---  
and hence is continuous.  
      
\begin{lemma} 
\label{lemma:4.1.1} 
The maps $g_{ij}:U_{ij}\to B\cstar$ satisfy 
the 1-cocycle condition 
$$
g_{jk}(m)g_{ik}(m)^{-1}g_{ij}(m) = 1 
$$
on $U_{ijk}$. 
\end{lemma} 

\begin{proof} To show this we need to use the 
group multiplication in $B\cstar$ using the 
non-homogenous coordinates on $B\cstar$ as 
reviewed in \cite{Gaj} and see also Section~\ref{sec:2.2}.  
We get 
\begin{eqnarray*} 
&   & g_{ij}g_{jk}                            \\ 
& = & |\psi_{0},\ldots,\psi_{n},[g_{iji_{0}}^{-1}
      g_{iji_{1}}g_{jki_{0}}^{-1}g_{jki_{1}}|\cdots | 
      g_{iji_{n-1}}^{-1}g_{iji_{n}}g_{jki_{n-1}}^{-1} 
      g_{jki_{n}}]|                                  \\ 
& = & |\psi_{0},\ldots,\psi_{n},[g_{iki_{0}}^{-1}
      g_{iki_{1}}|\cdots|g_{iki_{n-1}}g_{iki_{n}}]|  \\ 
& = & g_{ik},                                         
\end{eqnarray*}    
where we have used the 2-cocycle condition~\ref{eq:4.1.1} 
to write $g_{iji_{0}}^{-1}g_{iji_{1}}g_{jki_{0}}^{-1}g_{jki_{1}}$ 
equal to $g_{iki_{0}}^{-1}g_{iki_{1}}$ and so on.  
\end{proof} 

Hence the $g_{ij}$ represent transition functions 
for a $B\cstar$ bundle $X_{P}$ on $M$. 
We want to show that the class in $H^{3}(M;\Z)$ 
determined by $X_{P}$ and the Dixmier-Douady 
class $DD(P)$ of 
the bundle gerbe $P$ are the same.  

\begin{lemma} 
\label{lemma:4.1.2} 
Let $(\tilde{X}_{P},X_{P},M)$ denote the 
lifting bundle gerbe associated to $X_{P}$ 
via the short exact sequence of groups 
$$
1\to \cstar\to E\cstar\to B\cstar\to 1. 
$$
Then $DD(\tilde{X}_{P}) = DD(P)$. 
\end{lemma} 

\begin{proof} 
We will calculate the Dixmier-Douady class of 
$\tilde{X}_{P}$ relative to the same open cover $\U$ 
of $M$ used above.  We choose sections $s_{i}$ 
of $X_{P}\to M$ above each $U_{i}$ and form the 
pullback $\cstar$ bundles $(s_{i},s_{j})^{-1}\tilde{X}_{P}$ 
over $U_{ij}$.  Since the sections $s_{i}$ and 
$s_{j}$ are related by the transition cocycles 
$g_{ij}$, we see that a classifying map for 
$(s_{i},s_{j})^{-1}\tilde{X}_{P}$ is given 
by $g_{ij}:U_{ij}\to B\cstar$.  

Thus choosing sections $\sigma_{ij}$ of 
$(s_{i},s_{j})^{-1}\tilde{X}_{P}\to U_{ij}$ 
corresponds to choosing lifts $\hat{g}_{ij}:U_{ij}\to E\cstar$ 
of the $g_{ij}$.  We have a canonical choice 
for the $\hat{g}_{ij}$ given in the 
non-homogenous coordinates by 
\begin{eqnarray*} 
\hat{g}_{ij}(m) & = & |\psi_{0}(m),\ldots,\psi_{n}(m), 
g_{iji_{0}}(m)[g_{iji_{0}}(m)^{-1}g_{iji_{1}}(m)|\cdots \\ 
&  &   |g_{iji_{n-1}}(m)^{-1}g_{iji_{n}}(m)]|,       
\end{eqnarray*}    
where again $n = n(m)$ and the $\psi_{r}$ are defined as above.  
The $\hat{g}_{ij}$ need not be cocycles and 
indeed their failure to be is measured exactly by the 
2-cocycle $g_{ijk}$ for a calculation similar to 
that in the proof of Lemma~\ref{lemma:4.1.1} 
above shows that 
$$
\hat{g}_{jk}\hat{g}_{ik}^{-1}\hat{g}_{ij} = 
i(g_{ijk}), 
$$
where $i:\cstar\to E\cstar$ is the homomorphism 
induced by the inclusion of the subgroup $\cstar$.  
In other words, the lifting bundle gerbe $\tilde{X}_{P}$ 
has Dixmier-Douady class equal to $[g_{ijk}]$.  
\end{proof}

\section{The singular theory of $\cstar$ bundles}  
\label{sec:4.2} 
In this section we will show 
that there is an analogue of the curvature two form $F$ 
associated to a connection one form $A$ on 
a principal $\cstar$ bundle $\pi:P\to X$ in 
singular cohomology.  To be more precise, 
we will construct a singular 1-cochain $A\in S^{1}(P;\Z)$ 
and a singular 2-cocycle $F \in S^{2}(X;\Z)$ 
such that $\pi^{*}F = dA$.  
 
We first digress for 
a moment.  Recall from Chapter~\ref{chapter:3} 
that we can replace the structure group $\cstar$ 
of a bundle gerbe $P\to X^{[2]}$ with any abelian 
Lie group $A$ and still get a well defined notion 
of bundle gerbe.  In particular we can replace 
$\cstar$ by $\Z$ and get the notion of a `$\Z$ 
bundle gerbe'.  $\Z$ bundle gerbes arise whenever 
we have a $\cstar$ bundle $P\to X$, since we have 
the canonical map $\tau:P^{[2]}\to \cstar$ defined 
by $p_{2} = p_{1}\tau(p_{1},p_{2})$ for $p_{1}$ 
and $p_{2}$ points in the same fibre of $P$.  Then we 
can pullback the universal $\Z$ bundle $\C\to \cstar$ 
with the map $\tau$ to get a $\Z$ bundle $\tilde{P}\to 
P^{[2]}$.  Because the map $\tau$ satisfies $\tau(p_{2},
p_{3})\tau(p_{1},p_{2}) = \tau(p_{1},p_{3})$ we clearly 
get an associative product on $\tilde{P}$ and hence the 
triple $(\tilde{P},P,X)$ defines a $\Z$ bundle gerbe.  
Of course this is a rather artificial construction, 
but it does have the useful property that if one 
calculates the class in $H^{2}(X;\Z)$ defined by the  
$\Z$ bundle gerbe $\tilde{P}$ then one finds that it is equal to the Chern 
class of the $\cstar$ bundle $P$ (as one would hope).    
To see this we first choose an open cover $\{U_{i}\}$ 
of $X$ such that there exist local sections $s_{i}:U_{i}\to P$ 
of $P\to X$.  We then form the pullback $\Z$ bundles 
$(s_{i},s_{j})^{-1}\tilde{P} = \tilde{P}_{ij}$ and 
choose sections $\sigma_{ij}:U_{ij}\to \tilde{P}_{ij}$ --- 
this is like choosing lifts $z_{ij}:U_{ij}\to \C$ of the 
transition functions $g_{ij}:U_{ij}\to \cstar$ of $P$.  
Finally we can define an integer valued \v{C}ech 2-cocycle 
$n_{ijk}:U_{ijk}\to \Z$ by $m(\sigma_{jk}\otimes_{\Z} 
\sigma_{ij}) = \sigma_{ik} + n_{ijk}$ --- this translates into 
the condition $z_{jk} + z_{ij} = z_{ik} + n_{ijk}$.  
From this it is clear that $[n_{ijk}]$ is the 
image of $[g_{ij}]$ under the isomorphism 
$\check{H}^{1}(M;\underline{\C}^{\times}_{M})
\simeq \check{H}^{2}(M;\underline{\Z}_{M})$ 
induced from the long exact sequence in sheaf 
cohomology arising from the short exact sequence 
of sheaves of abelian groups 
$$
\underline{\Z}_{M}\to \underline{\C}_{M}\to 
\underline{\C}^{\times}_{M}.  
$$
Therefore the class $[n_{ijk}]$ in $H^{2}(M;\Z)$ is the Chern 
class of the $\cstar$ bundle $P$.     
       
We have the following lemma.  
\begin{lemma} 
\label{lemma:4.2.1} 
There exists a singular cocycle 
$c\in Z^{1}(\cstar;\Z)$ and a 
singular cochain $b\in S^{0}(\cstar\times \cstar;\Z)$ 
such that 
\begin{equation} 
\label{eq:4.2.1} 
d_{0}^{*}c - d_{1}^{*}c + d_{2}^{*}c = db 
\end{equation} 
and 
\begin{equation} 
\label{eq:4.2.2}  
d_{0}^{*}b - d_{1}^{*}b + d_{2}^{*}b - d_{3}^{*}b = 0.
\end{equation} 
\end{lemma} 
Here the $d_{i}$ denote the face 
operators in the simplicial 
manifold $N\cstar$.  
\begin{proof} 
First of all note that a singular 0-chain 
in $S_{0}(\Z;\Z)$ is a linear combination 
of 0-simplices $\sigma:\Delta^{0}\to \Z$ 
with integer coefficients.  Since $\Delta^{0}$ 
is just a point, it follows that we may 
identify a singular 0-chain of $S_{0}(\Z;\Z)$ 
with a formal linear combination $n_{1}m_{1} + 
\cdots + n_{k}m_{k}$.  A singular 0-cochain is 
a linear map $S_{0}(\Z;\Z)\to \Z$.  There is 
a canonical choice of a singular 0-cochain 
in $S^{0}(\Z;\Z)$, namely the linear map 
which sends a formal linear combination 
$n_{1}m_{1} + \cdots + n_{k}m_{k}$ to the 
integer $n_{1}m_{1} + \cdots + n_{k}m_{k}$.  
Let us call this singular 0-cochain $a$.  It 
has the property that $\text{pr}_{1}^{*}a + \text{pr}_{2}^{*}a 
= m^{*}a$, where $\text{pr}_{1}$ and $\text{pr}_{2}$ are the 
projections on the first and second factors of 
$\Z\times \Z$ respectively, and where $m:\Z\times 
\Z\to \Z$ is addition.  This is because a singular 
0-chain of $\Z\times \Z$ may be identified in the 
same way as above with a formal linear combination 
$n_{1}(m_{1},m_{1}^{'}) + \cdots + n_{k}(m_{k},m_{k}^{'})$ 
and it is easy to see that the value of $m^{*}a$ 
on such a 0-chain is equal to the value of 
$\text{pr}_{1}^{*}a + \text{pr}_{2}^{*}a$ on such a 0-chain.  
We also have $da = 0$, where $d$ is the singular 
coboundary, since a singular 1-chain of $S_{1}(\Z;\Z)$ 
must be a linear combination of constant maps 
$\Delta^{1}\to \Z$ with integer coefficients 
from which it follows that $da$ evaluated at 
such a 1-chain must be zero.    
         
Next we have the short exact sequence of groups 
$\Z\to \C\to \cstar$, where $\Z\to \C$ is the 
inclusion and $\C\to \cstar$ is the map 
sending $z\in \C$ to $\exp (2\pi\sqrt{-1}z)$.  
We will think of $\C$ as a `principal $\Z$ 
bundle' on $\cstar$.  Form the fibre product 
$\C^{[2]}$.  We then have a canonical map 
$\tau:\C^{[2]}\to \Z$.  Form a 0-cochain 
$\tau^{*}a$ of $S^{0}(\C^{[2]};\Z)$ by pullback with 
$\tau$.  The property $\text{pr}_{1}^{*}a + \text{pr}_{2}^{*}a 
= m^{*}a$ of $a$ shows that $\d(\tau^{*}a) = 0$ 
in $S^{0}(\C^{[3]};\Z)$, where $\d:S^{0}(\C^{[2]};\Z) 
\to S^{0}(\C^{[3]};\Z)$ is the coboundary 
defined in Proposition~\ref{prop:3.3.2} 
by adding the pullback maps $\pi_{i}^{*}$ 
with alternating signs.  It follows 
from Proposition~\ref{prop:3.3.2} that there exists 
$\hat{c} \in S^{0}(\C;\Z)$ such that 
$\tau^{*}a = \d(\hat{c})$.  Since $da = 0$, 
we have $\d(d\hat{c}) = 0$ so by Proposition~\ref{prop:3.3.2} 
again we get $d\hat{c} = p^{*}c$ for some 
$c\in S^{1}(\cstar;\Z)$, where $p:\C\to \cstar$ 
is the projection.     

Denote multiplication in $\cstar$ by $m$ 
and denote multiplication (addition) in 
$\C$ by $\hat{m}$.  We have a commutative 
diagram 
$$
\xymatrix{ 
(\C\times \C)^{[2]} \ar[d]_-{\hat{m}^{[2]}} 
\ar[rr]^-{(\tau\circ \text{pr}_{1}^{[2]},\tau\circ 
\text{pr}_{2}^{[2]})} & & \Z\times \Z \ar[d]          \\ 
\C^{[2]} \ar[rr]^-{\tau} & & \Z,                 } 
$$
where now $\text{pr}_{1}:\C\times \C\to \C$ and 
$\text{pr}_{2}:\C\times \C\to \C$ denote the 
projections on the first and second factors 
in $\C\times \C$ respectively and 
$\Z\times \Z\to \Z$ is addition in $\Z$.  From 
this diagram we get $(\hat{m}^{[2]})^{*}\tau^{*}
a = (\tau\circ \text{pr}_{1}^{[2]})^{*}a + (\tau\circ 
\text{pr}_{2}^{[2]})^{*}a$ and since $\tau^{*}a = \d(\hat{c})$, 
this becomes $\d(\hat{m}^{*}\hat{c}) = \d(\text{pr}_{1}^{*}
\hat{c} + \text{pr}_{2}^{*}\hat{c})$.  Using Proposition~
\ref{prop:3.3.2} again, we see that there exists 
$b\in S^{0}(\cstar\times \cstar;\Z)$ such that 
$(p\times p)^{*}b = \text{pr}_{1}^{*}\hat{c} - \hat{m}^{*}
\hat{c} + \text{pr}_{2}^{*}\hat{c}$, where $p\times p$ 
is the projection $\C\times \C\to \cstar\times 
\cstar$.  From this last equation, it is easy 
to see that $b$ satisfies equation~\ref{eq:4.2.2} 
above.  Also from the last equation we see that 
$(p\times p)^{*}db = (p\times p)^{*}(\text{pr}_{1}^{*}c 
- m^{*}c + \text{pr}_{2}^{*}c)$, in other words 
equation~\ref{eq:4.2.1} is satisfied.      
\end{proof} 
\begin{note} 
We will briefly check here that the singular 
cochains $\hat{c}$ and $c$ defined above are 
not identically zero and that $c$ is not 
the coboundary of a 0-cochain.  To solve 
$\tau^{*}a = \d(\hat{c})$ for $\hat{c}$, 
recall from the proof of Proposition~\ref{prop:3.3.2} 
that we put $\hat{c}$ equal to the 0-cochain 
that sends a 0-simplex $\sigma:\Delta^{0}\to \C$ 
to $\tau^{*}a(\sigma, s(p(\sigma)))$, where 
$s:S_{0}(\cstar;\Z)\to S_{0}(\C;\Z)$ is a 
section of the set map $p_{*}:S_{0}(\C;\Z)
\to S_{0}(\cstar;\Z)$.  Since we may identify 
$S_{0}(\C;\Z)$ and $S_{0}(\cstar;\Z)$ with 
the free abelian groups on the points of 
$\C$ and $\cstar$ respectively, choosing such 
a section $s$ amounts to choosing a section of 
$p:\C\to \cstar$ where we regard $\C$ and $\cstar$ 
as sets so $s$ need not be continuous.  The 
obvious candidate in this case is the map 
$s:\cstar\to \C$ sending $r\exp (\sqrt{-1}\theta)$ 
(where $0\leq \theta < 2\pi$) to 
$\frac{1}{2\pi \sqrt{-1}}(\log r + 
\sqrt{-1}\theta)$.  Now consider the 0-simplex 
$1$ sending the point $1$ of $\Delta^{0}$ to $1\in \C$.  
We have $\hat{c}(1) = a(\tau(\frac{1}{2\pi \sqrt{-1}}
\log(1),1))) = a(\tau(0,1)) = a(1) = 1$.  
In a similar way one can show that $d\hat{c}$ 
is not identically zero, and so $c$ is not 
identically zero.  We will now check that 
we do not have $c = de$, for some $e\in S^{0}(
\cstar;\Z)$.  Suppose we did, then we would have 
$d(\hat{c} - p^{*}e) = 0$ in $S^{1}(\C;\Z)$.  
This means that $\hat{c} - p^{*}e$ must be 
constant on the 0-simplexes $\Delta^{0}\to \C$.  
Let this constant value be $n$ say.  Then 
$e+n$ defines a new 0-cochain of $S^{0}(\cstar;\Z)$ 
and we have $d(e+n) = c$.  Also we have 
$\hat{c} = p^{*}(e+n) = p^{*}e + n$ and so $\d(\hat{c}) = 0$, 
which contradicts the fact that $\d(\hat{c}) = 
\tau^{*}a \neq 0$.  In fact, one can give an explicit 
formula for $c$.  First of all suppose we have 
chosen a set map $s:\cstar\to \C$ which is also 
a section of $p:\C\to \cstar$.  We need to find a section 
of $p_{*}:S_{1}(\C;\Z)\to S_{1}(\cstar;\Z)$.  There is a 
canonical way of doing this.  Let $\sigma:\Delta^{1}\to \cstar$ 
be a 1-simplex.  Then there is a unique way of lifting 
$\sigma$ to a 1-simplex $\tilde{\sigma}:\Delta^{1}\to \C$ 
with $p\circ \tilde{\sigma} = \sigma$ and $\tilde{\sigma}(0) = 
s(\sigma(0))$.  Now we can write down a formula for $c$: 
$$
c_{\sigma} = \tilde{\sigma}(1) - s(\sigma(1)), 
$$ 
where $\sigma:\Delta^{1}\to \cstar$ is a 1-simplex.  It 
is not hard to check that $c$ defined in this manner 
is closed and it value on the canonical generator of 
$H_{1}(\cstar;\Z)$, a triangle in $\cstar$ around $0$, 
is the winding number of the triangle ie $1$.  Hence 
$c$ represents the generator or fundamental class of 
$H^{1}(\cstar;\Z)$.  Similarly 
one can write down a formula for $b$, 
$$
b_{\sigma_{1},\sigma_{2}} = s(\sigma_{1}\sigma_{2}) - s(\sigma_{1}) 
- s(\sigma_{2}), 
$$
where $\sigma_{1}:\Delta^{0}\to \cstar$ and $\sigma_{2}:
\Delta^{0}\to \cstar$ are 0-simplexes.            
\end{note} 
        
Now suppose we are given a principal $\cstar$  
bundle $P\to X$.  Let $\pi_{P}$ 
denote the projection $P\to X$.  
Let, as usual, $P^{[2]}$ denote 
the fibre product of $P$ with 
itself over $X$ and let 
$P^{[p]}$ denote the $p$-fold 
such fibre product.  We have the 
natural map 
$\tau:P^{[2]}\to \cstar$ defined by 
$p_{2} = p_{1}\cdot \tau(p_{1},p_{2})$ 
for $p_{1}$ and $p_{2}$ in the same 
fibre of $P$.  $\tau$ satisfies 
\begin{equation} 
\label{eq:4.2.3} 
\tau(p_{2},p_{3})\tau(p_{1},p_{2}) = 
\tau(p_{1},p_{3}) 
\end{equation} 
for $(p_{1},p_{2},p_{3})\in P^{[3]}$.   
We have the 
following commutative diagram: 
$$
\xymatrix{ 
\tilde{P} \ar[d]_-{\pi_{\tilde{P}}} 
\ar[r]^-{\tilde{\tau}} & \C \ar[d]^-{p}      \\  
P^{[2]} \ar[r]^-{\tau} &           \cstar.    } 
$$
Let $\a = \tau^{*}c$ and let $\hat{\a} = 
\tilde{\tau}^{*}\hat{c}$ so $d\hat{\a} = 
\pi_{\tilde{P}}^{*}\a$.      
Notice that even if 
$P$ has zero Chern class in $H^{2}(X;\Z)$, 
in which case the structure group of 
$P$ will reduce to $\Z_{n}$ for some $n$, 
$\a$ will be non zero since the map 
$\tau:P^{[2]}\to \cstar$ is onto $\cstar$.   
Then equations~\ref{eq:4.2.1} 
and~\ref{eq:4.2.3} imply that 
\begin{equation} 
\label{eq:4.2.4} 
\pi_{1}^{*}\a - \pi_{2}^{*}\a + \pi_{3}^{*}\a = d\b, 
\end{equation} 
where $\b\in S^{0}(P^{[3]};\Z)$ is defined 
by $\b = (\tau\circ \pi_{1},\tau\circ \pi_{3})^{*}b$.  
Another way of seeing this is to note that 
we have the commutative diagram 
$$
\xymatrix{ 
\pi_{1}^{-1}\tilde{P} \times_{P^{[3]}} 
\pi_{3}^{-1}\tilde{P} \ar[d]_-{\tilde{m}_{\tilde{P}}} 
\ar[rrr]^-{(\pi_{1}^{-1}\tilde{\tau}\circ \text{pr}_{1},
\pi_{3}^{-1}\tilde{\tau}\circ \text{pr}_{2})} & & & \C
\times \C \ar[d]^-{\hat{m}}                   \\ 
\pi_{2}^{-1}\tilde{P} \ar[rrr]^-{\pi_{2}^{-1}
\tilde{\tau}} & & & \C ,          } 
$$
where $\tilde{m}_{\tilde{P}}:\pi_{1}^{-1}\tilde{P}
\times_{P^{[3]}}\pi_{3}^{-1}\tilde{P}\to \pi_{2}^{-1}
\tilde{P}$ is induced by $m_{\tilde{P}}$.   
From the commutativity of this diagram and 
the equation $\text{pr}_{1}^{*}\hat{c} + \text{pr}_{2}^{*}
\hat{c} = \hat{m}^{*}\hat{c} + (p\times p)
^{*}b$ in $S^{0}(\C\times \C;\Z)$, we get the 
following equation in $S^{0}(\pi_{1}^{-1}
\tilde{P}\times_{P^{[3]}}\pi_{3}^{-1}\tilde{P};\Z)$: 
\begin{equation} 
\label{eq:4.2.5} 
\text{pr}_{1}^{*}\pi_{1}^{-1}\hat{\a} + \text{pr}_{2}^{*}
\pi_{3}^{-1}\hat{\a} = \tilde{m}_{\tilde{P}}^{*}\pi_{2}^{-1}
\hat{\a} + (\pi_{\tilde{P}}\times \pi_{
\tilde{P}})^{*}\b.  
\end{equation} 
Taking $d$ of this equation gives equation~
\ref{eq:4.2.4}.  Also it follows from 
equation~\ref{eq:4.2.5} or   
equation~\ref{eq:4.2.2} 
that $\d(\b) = 0$ in 
$S^{0}(P^{[4]};\Z)$.  Hence 
Proposition~\ref{prop:3.3.2} implies 
that there exists $\c\in S^{0}(P^{[3]};\Z)$ 
such that $\d(\c) = \b$.  Then, from 
equation~\ref{eq:4.2.4}, we get 
$\d(\a) = \d(d\c)$ and hence 
from Proposition~\ref{prop:3.3.2} 
again, we get $\a = d\c + \d(A)$ 
for some $A\in S^{1}(P;\Z)$.  Since 
$d\a=0$ we have that $\d(dA)=0$ and 
hence there exists $F\in S^{2}(X;\Z)$ 
such that $\pi_{P}^{*}F=dA$.  

\begin{lemma} 
\label{lemma:4.2.2} 
$F$ is a 
representative of the Chern class of 
$P\to X$ in $H^{2}(X;\Z)$.    
\end{lemma} 

\begin{proof} 
The problem is to show that $F$ 
is a representative in singular cohomology of the Chern class 
of $P$.  A singular representative 
for the Chern class of $P$ can be constructed  
by observing as above that $P$ gives rise to a $\Z$ bundle gerbe 
$(\tilde{P},P,X)$ on $X$ --- the lifting $\Z$ bundle gerbe 
associated to the short exact sequence of groups 
$\Z \to \C \to \cstar$.  As we saw above, the \v{C}ech 2-cocycle 
$n_{ijk}:U_{ijk}\to \Z$ relative to some 
open cover $\U = \{U_{i}\}_{i\in I}$ of $X$ associated to the 
$\Z$ bundle gerbe $\tilde{P}$ maps to the Chern class 
of $P$ under the isomorphism 
$\check{H}^{2}(X;\underline{\Z}_{M})\simeq H^{2}(X;\Z)$.  

To construct the singular cocycle representing 
the Chern class of $P$ we 
use the fact (see \cite{BotTu} and the note 
following Proposition~\ref{prop:3.3.2}) that 
the singular Mayer-Vietoris sequence 
for the open covering $\{U_{i}\}_{i\in I}$ 
is exact.  Thus we view $n_{ijk}$ as 
living in $S^{0}(U_{ijk};\Z)$ and use 
the exactness of the singular Mayer-Vietoris 
sequence to solve $n_{ijk} = \d(u_{ij})$ 
for some $u_{ij} \in S^{0}(U_{ij};\Z)$.  
Since $n_{ijk}$ is locally constant 
we have $dn_{ijk} = 0$ and  
so $\d(du_{ij}) = 0$.  Applying the exactness of 
the singular Mayer-Vietoris sequence again 
implies that there is an $A_{i} \in S^{1}(U_{i};\Z)$ 
so that $du_{ij} = \d(A_{i})$.  
Applying $d$ to this equation gives us 
$\d(dA_{i}) = 0$ and so there exists $F^{'}$ 
in $S^{2}(X;\Z)$ with $dF^{'} = 0$ and 
$\d(F^{'}) = dA_{i}$.  

We want to show that $F$ and $F^{'}$ represent 
the same class in cohomology.  We use the same 
open cover $\{U_{i}\}_{i\in I}$ as above and 
form the principal $\Z$ bundles $\tilde{P}_{ij} = 
(s_{i},s_{j})^{-1}\tilde{P}$ over $U_{ij}$.  We 
have the following commutative diagram: 
$$
\xymatrix{  
\tilde{P}_{ij} \ar[d] \ar[r]^-{\bar{(s_{i},s_{j}
)}} & \tilde{P} \ar[d] \ar[r]^-{\tilde{\tau}} & 
\C \ar[d]                                       \\ 
U_{ij} \ar[r]^-{(s_{i},s_{j})} & P^{[2]} \ar[r]^-
{\tau} & \cstar,                                } 
$$
which yields by composition the commutative 
diagram 
$$
\xymatrix{ 
\tilde{P}_{ij} \ar[d] \ar[r]^-{\bar{g}_{ij}} 
& \C \ar[d]                                  \\ 
U_{ij} \ar[r]^-{g_{ij}} & \cstar.            } 
$$
We let $\hat{\a}_{ij} \in S^{0}(\tilde{P}_{ij};
\Z)$ denote the pullback $\bar{g}_{ij}^{*}\hat{c}$ 
and we let $\a_{ij} \in S^{1}(U_{ij};\Z)$ denote the 
pullback $g_{ij}^{*}c$.  By pulling back 
equation~\ref{eq:4.2.5} to $S^{0}(\tilde{P}_{jk}
\times_{U_{ijk}}\tilde{P}_{ij};\Z)$ we get 
$$
\text{pr}_{1}^{*}\hat{\a}_{jk} + \text{pr}_{2}^{*}\hat{\a}_{ij} = 
\tilde{m}_{\tilde{P}}^{*}\hat{\a}_{ik} + (\pi_{\tilde{P}_{jk}}
\times_{U_{ijk}} \pi_{\tilde{P}_{ij}})^{*}\b_{ijk}, 
$$
where $\b_{ijk} = (s_{i},s_{j},s_{k})^{*}\b$.        
Since $\b = \d(\c)$, we may rewrite $\b_{ijk}$ 
as $\c_{jk}-\c_{ik}+\c_{ij}$ where $\c_{ij} 
= (s_{i},s_{j})^{*}\c$.  Now the lifts 
$z_{ij}:U_{ij}\to \C$ of the $g_{ij}$ furnish 
us with sections $\sigma_{ij}$ of $\tilde{P}_{ij}$.  
Also it follows that we have $\tilde{m}_{\tilde{P}}(
\sigma_{jk},\sigma_{ij}) = \sigma_{ik} + n_{ijk}$.  
Hence from the equation above we get 
$$
\sigma_{jk}^{*}\hat{\a}_{jk} + \sigma_{ij}^{*}
\hat{\a}_{ij} = (\sigma_{ik} + n_{ijk})^{*}
\hat{\a}_{ik} + \c_{jk} - \c_{ik} + \c_{ij}.  
$$
It follows from the definition of $\hat{\a}_{ik}$ 
that $(\sigma_{ik} + n_{ijk})^{*}\hat{\a}_{ik} = 
\sigma_{ik}^{*}\hat{\a}_{ik} + n_{ijk}$.  This is 
because $(\sigma_{ik}+n_{ijk})^{*}\hat{\a}_{ik} = 
(\sigma_{ik}+n_{ijk})^{*}\bar{g}_{ik}^{*}\hat{c} = 
(z_{ik}+n_{ijk})^{*}\hat{c}$ so if we can show that 
$\hat{c}(\sigma+n) = \hat{c}(\sigma) + n$ for a 
0-simplex $\sigma:\Delta^{0}\to \C$ and $n\in \Z$ 
then we will be done.  By definition $\hat{c}(\sigma+n) = 
(\tau^{*}a)(\sigma+n,s(p(\sigma+n))) = (\tau^{*}a)(
\sigma+n,s(p(\sigma))) = (\tau^{*}a)(\sigma,s(p(\sigma)))+n = 
\hat{c}(\sigma) + n$.      
Hence we have 
$$
\sigma_{jk}^{*}\hat{\a}_{jk} - \c_{jk} + 
\sigma_{ij}^{*}\hat{\a}_{ij} - \c_{ij} = 
\sigma_{ik}^{*}\hat{\a}_{ik} - \c_{ik} + 
n_{ijk}.  
$$
Therefore $\sigma_{ij}^{*}\hat{\a}_{ij} - 
\c_{ij}$ plays the role of $u_{ij}$.  Note 
that $d(\sigma_{ij}^{*}\hat{\a}_{ij} - \c_{ij})$ 
equals $\sigma_{ij}^{*}p^{*}\a_{ij} - d\c_{ij} = 
\a_{ij} - d\c_{ij} = s_{j}^{*}A - s_{i}^{*}A$.  
It follows from this that $F$ and $F^{'}$ 
represent the same class in cohomology.   
\end{proof} 

\section{The singular theory of bundle gerbes} 
\label{sec:4.3} 

In this section we utilise the results of 
Section~\ref{sec:4.2} to show that there are 
singular analogues of the curving and 
three curvature of a bundle gerbe.  So let 
$(Q,Y,M)$ be a bundle gerbe.  We will assume that 
we have constructed $\hat{\a}\in S^{0}(\tilde{Q};\Z)$, 
$\a\in S^{1}(Q^{[2]};\Z)$ and $\b\in S^{0}(Q^{[3]};\Z)$ 
satisfying equations~\ref{eq:4.2.4} and~\ref{eq:4.2.5} 
leading to $A\in S^{1}(Q;\Z)$ and $F\in S^{2}(Y^{[2]};\Z)$.  
   
Form the pullback $\cstar$ 
bundles $\pi_{i}^{-1}Q$ over $Y^{[3]}$ 
for $i =1,2,3$ so that there exist 
$\cstar$ bundle morphisms $\bar{\pi}_{i}
:\pi_{i}^{-1}Q\to Q$ covering $\pi_{i}$ 
as pictured in the following diagram.   
$$ 
\xymatrix{  
\pi_{i}^{-1}Q \ar[r]^{\bar{\pi}_{i}} 
\ar[d]   &  Q \ar[d]                       \\ 
Y^{[3]} \ar[r]^{\pi_{i}} & Y^{[2]}          } 
$$ 
The natural map $\tau:Q^{[2]}\to \cstar$ 
induces maps $\tau^{'}:(\pi_{i}^{-1}Q)^{[2]}
\to \cstar$ such that $\tau^{'} = \tau\circ 
\bar{\pi}_{i}^{[2]}$.  
Form the fibre product 
$\pi_{1}^{-1}Q\times_{Y^{[3]}}\pi_{3}^{-1}Q$.  
Let $\tau_{1}$ denote the map 
$$
(\tau\circ \bar{\pi_{1}}^{[2]}\circ p_{1}^{[2]}, 
\tau\circ \bar{\pi}_{3}^{[2]}\circ p_{2}^{[2]}):
(\pi_{1}^{-1}Q\times_{Y^{[3]}}\pi_{3}^{-1}Q)^{[2]}\to 
\cstar\times \cstar.  
$$ 
where $p_{1}$ and $p_{2}$ denote the 
projections $\pi_{1}^{-1}Q\times_{Y^{[3]}}\pi_{3}^{-1}Q\to \pi_{1}^{-1}Q$ 
and $\pi_{1}^{-1}Q\times_{Y^{[3]}}\pi_{3}^{-1}Q\to \pi_{3}^{-1}Q$ 
respectively.  Also, let 
$\tau_{2}$ denote the map 
$\tau\circ \bar{\pi}_{2}^{[2]}:(\pi_{2}^{-1}Q)^{[2]}\to \cstar$.  
The bundle gerbe product $m_{Q}$ gives 
rise to a map 
$\tilde{m}_{Q}:\pi_{1}^{-1}Q\times_{Y^{[3]}}\pi_{3}^{-1}Q\to \pi_{2}^{-1}Q$ 
and we have a commutative diagram 
$$ 
\xymatrix{   
(\pi_{1}^{-1}Q\times_{Y^{[3]}}\pi_{3}^{-1}Q)^{[2]} 
\ar[r]^-{\tilde{m}^{[2]}_{Q}} \ar[d]_-{\tau_1} & 
(\pi_{2}^{-1}Q)^{[2]} \ar[d]^-{\tau_{2}}              \\ 
\cstar \times \cstar \ar[r]^-{m} & \cstar        } 
$$ 
Hence we have 
$\tau_{2}^{*}(\tilde{m}_{Q}^{[2]})^{*}c = \tau_{1}^{*}m^{*}c$.  
Using equation~\ref{eq:4.2.1} 
we get 
\begin{eqnarray*} 
&  & \tau_{1}^{*}m^{*}c                         \\ 
& = & (p_{1}^{[2]})^{*}(\bar{\pi}_{1}
^{[2]})^{*}\a + 
(p_{2}^{[2]})^{*}(\bar{\pi}_{3}^{[2]})^{*}\a - 
d(\tau\circ \bar{\pi}_{1}^{[2]}\times \tau\circ \bar{\pi}_{3}^{[2]})^{*} 
b                                                      \\ 
& = & (\tilde{m}_{Q}^{[2]})^{*}(\bar{\pi}_{2}^{[2]})^{*}\a.   
\end{eqnarray*} 
We can rewrite this as 
\begin{eqnarray} 
& & \d(p_{1}^{*}\bar{\pi}_{1}^{*}A + p_{2}^{*}\bar{\pi}_{3}^{*}A) - 
d(\tau\circ \bar{\pi}_{1}^{[2]}\circ p_{1}^{[2]}\times \tau\circ\bar{\pi}
_{3}^{[2]}\circ p_{2}^{[2]})^{*}b - d(\tilde{m}_{Q}^{[2]})^{*}
(\bar{\pi}_{2}^{[2]})^{*}\c                  \nonumber           \\ 
& = & \d(\tilde{m}_{Q}^{*}\bar{\pi}_{2}^{*}A) - 
d(p_{1}^{[2]})^{*}(\bar{\pi}_{1}^{[2]})^{*}\c - 
d(p_{2}^{[2]})^{*}(\bar{\pi}_{3}^{[2]})^{*}\c.    \label{eq:4.2.6}      
\end{eqnarray} 
Let $\rho \in S^{0}((\pi_{1}^{-1}Q\times_{Y^{[3]}}\pi_{3}^{-1}Q)^{[2]};\Z)$ 
be defined by  
$$
\rho = (\tau\circ \bar{\pi}_{1}^{[2]},\tau\circ\bar{\pi}_{3}
^{[2]}\circ p_{2}^{[2]})^{*}b + (\tilde{m}_{Q}^{[2]})^{*}
(\bar{\pi}_{2}^{[2]})^{*}\c - (p_{1}^{[2]})^{*}(
\bar{\pi}_{1}^{[2]})^{*}\c - (p_{2}^{[2]})^{*}(
\bar{\pi}_{3}^{[2]})^{*}\c. 
$$
Then equation~\ref{eq:4.2.6} above becomes 
\begin{equation} 
\label{eq:4.2.7} 
\d(p_{1}^{*}\bar{\pi}_{1}^{*}A + p_{2}^{*}\bar{\pi}_{3}^{*}A) 
= \d(\tilde{m}_{Q}^{*}\bar{\pi}_{2}^{*}A) + d\rho. 
\end{equation}   

\begin{lemma} 
\label{lemma:4.2.3}
$\d(\rho) = 0$ in 
$S^{0}((\pi_{1}^{-1}Q\times_{Y^{[3]}}\pi_{3}^{-1}Q)^{[3]};\Z)$. 
\end{lemma} 

We will omit the straightforward but tedious 
proof of this lemma.  

Thus there exists 
$B\in S^{0}(\pi_{1}^{-1}Q\times_{Y^{[3]}}\pi_{3}^{-1}Q;\Z)$ 
such that $\rho = \d(B)$.  Hence 
from equation~\ref{eq:4.2.7} we get 
$\d(p_{1}^{*}\bar{\pi}_{1}^{*}A + 
p_{2}^{*}\bar{\pi}_{3}^{*}A) = 
\d(\tilde{m}_{Q}^{*}\bar{\pi}_{2}^{*}A 
+ dB)$.  
Therefore there 
exists $C\in S^{1}(Y^{[3]};\Z)$ such that 
\begin{equation} 
\label{eq:4.2.8} 
p_{1}^{*}\bar{\pi}_{1}^{*}A + 
p_{2}^{*}\bar{\pi}_{3}^{*}A = 
\tilde{m}_{Q}^{*}\bar{\pi}_{2}^{*}A + dB + \pi_{Q}^{*}C, 
\end{equation} 
where we abuse notation and denote by 
$\pi_{Q}$ the projection 
$\pi_{1}^{-1}Q\times_{Y^{[3]}}\pi_{3}^{-1}Q\to Y^{[3]}$.  
Taking $d$ of equation~\ref{eq:4.2.8}     
implies 
$\pi^{*}(\pi_{1}^{*}F+\pi_{3}^{*}F)=\pi^{*}(\pi_{2}^{*}F +dC)$ 
and hence $\d(F) = dC$.  We would like 
to show that $\d(C) = 0$.  
A singular 1-simplex in 
$S_{1}(Y^{[2]};\Z)$ is a pair 
$(\s_{1},\s_{2})$ of  
1-simplexes $\s_{i}:\Delta^{1}\to Y$ 
such that $\pi(\s_{1})=\pi(\s_{2})$.  
Given such a pair $(\s_{i},\s_{j})$ there 
exists a lift $\s_{ij}:\Delta^{1}\to Q$ 
of the 1-simplex $(\s_{i},\s_{j})$.  
Notice that $(\s_{jk},\s_{ij})$ is then 
a 1-simplex in 
$S_{1}(\pi_{1}^{-1}Q\times_{Y^{[3]}}\pi_{3}^{-1}Q;\Z)$ 
which satisfies 
$\pi_{Q}(\s_{jk},\s_{ij}) = (\s_{i},\s_{j},\s_{k})$. 
From equation~\ref{eq:4.2.8} and from 
the associativity of $m_{Q}$ we see that 
\begin{eqnarray*} 
& & A_{\tilde{m}_{Q}(\s_{34},\s_{23})} + 
A_{\s_{12}} - dB_{(\tilde{m}_{Q}(\s_{34},\s_{23}),\s_{12})} 
-C_{(\s_{1},\s_{2},\s_{4})}                      \\ 
& =& A_{\s_{34}} + A_{\tilde{m}_{Q}(\s_{23},\s_{12})} 
-dB_{(\s_{34},\tilde{m}_{Q}(\s_{23},\s_{12}))} 
-C_{(\s_{1},\s_{3},\s_{4})}                       
\end{eqnarray*} 
which, using equation~\ref{eq:4.2.8} 
again, gives 
\begin{eqnarray*} 
& & A_{\s_{34}} + A_{\s_{23}} + A_{\s_{12}} 
-dB_{(\s_{34},\s_{23})} - dB_{(\tilde{m}_{Q}(\s_{34},\s_{23}),\s_{12})}  
 -C_{(\s_{2},\s_{3},\s_{4})} - 
C_{(\s_{1},\s_{2},\s_{4})}                                       \\ 
& = & A_{\s_{34}} + A_{\s_{23}} + A_{\s_{12}} 
-dB_{(\s_{23},\s_{12})} - dB_{(\s_{34},\tilde{m}_{Q}(\s_{23},\s_{12}))}   
 - C_{(\s_{1},\s_{2},\s_{3})} - C_{(\s_{1},\s_{3},\s_{4})}.    
\end{eqnarray*} 
and hence 
\begin{eqnarray} 
 \d(C)_{(\s_{1},\s_{2},\s_{3},\s_{4})}  
& = & dB_{(\s_{23},\s_{12})} -dB_{(\tilde{m}_{Q}(\s_{34},\s_{23}),\s_{12})} 
+ dB_{(\s_{34},\tilde{m}_{Q}(\s_{23},\s_{12}))} 
-dB_{(\s_{34},\s_{23})}.    \label{eq:4.2.9}                        
\end{eqnarray} 
Notice that if the $\s_{ij}$ are now 
all zero simplices and we can show 
\begin{equation} 
\label{eq:4.2.10} 
B_{(\s_{23},\s_{12})} - B_{(\tilde{m}_{Q}(\s_{34},\s_{23}),\s_{12})} 
+ B_{(\s_{34},\tilde{m}_{Q}(\s_{23},\s_{12}))} - B_{(\s_{34},\s_{23})} = 0, 
\end{equation} 
then we will have $\d(C) = 0$ as we want.  
To establish equation~\ref{eq:4.2.10} 
we first need a more explicit formula 
for $B$.  To achieve this we need to 
digress for a moment.  Notice that if 
we regard $Q$ and $Y^{[2]}$ as sets 
then because $\pi_{Q}:Q\to Y^{[2]}$ is 
onto we can find a section $s$ say, of 
$\pi_{Q}$ (obviously $s$ will not be 
continuous unless $Q$ is trivial).  
Then we can define a set map $c:Y^{[3]}\to \cstar$ 
satisfying 
$$
c(y_{2},y_{3},y_{4})c(y_{1},y_{3},y_{4})^{-1}
c(y_{1},y_{2},y_{4})c(y_{1},y_{2},y_{3})^{-1} = 1,  
$$
for $y_{i}$ all lying in the same fibre of 
$Y$ over $M$, by defining 
$m_{Q}(s(y_{2},y_{3}),s(y_{1},y_{2}))=s(y_{1},y_{3})c(y_{1},y_{2},y_{3})$ 
for $(y_{1},y_{2},y_{3})\in Y^{[3]}$.  
Since $\pi:Y\to M$ is onto, we can 
find a set map $\sigma:M\to Y$ 
which is a section of $\pi$.  Therefore 
we can define a function $c^{'}:Y^{[2]}\to 
\cstar$ by $c^{'}(y_{1},y_{2}) = c(y_{1},
y_{2},\sigma(m))$, where $m = \pi(y_{1}) 
= \pi(y_{2})$.  Then it is easy to check 
that $c^{'}(y_{2},y_{3})c^{'}(y_{1},y_{3})
^{-1}c^{'}(y_{1},y_{2}) = c(y_{1},y_{2},y_{3})$.  
Hence we can rescale our section $s$ so that it 
satisfies 
\begin{equation} 
\label{eq:4.2.11}  
m_{Q}(s(y_{2},y_{3}),s(y_{1},y_{2})) = s(y_{1},y_{3}). 
\end{equation}    
Next, observe from the proof of 
Proposition~\ref{prop:3.3.2} that a 
formula for $B_{(\s_{23},\s_{12})}$ is then 
given by 
$$
B_{(\s_{23},\s_{12})} = b_{(\tau(s(p(\s_{23})),\s_{23}), 
\tau(s(p(\s_{12})),\s_{12}))}. 
$$
From this and equation~\ref{eq:4.2.11} 
we see that $B$ satisfies equation~\ref{eq:4.2.10}.   
Therefore $C$ satisfies $\d(C) = 0$ and 
so we may find $D \in S^{1}(Y^{[2]};\Z)$ such 
that $C = \d(D)$.  Thus if we rename $F$ 
to $F - dD$ then we have $\d(F) = 0$ and $F$ 
is still a representative for the Chern class 
of $Q\to Y^{[2]}$.  Also we can solve $F = \d(f)$ 
for some $f \in S^{2}(Y;\Z)$ which must therefore 
satisfy $\d(df) = 0$.  Hence there exists  
$\omega \in S^{3}(M;\Z)$ satisfying $\pi^{*}\omega = df$ 
and $d\omega = 0$.  We have the following 
proposition. 

\begin{proposition} 
\label{prop:4.2.4} 
The class in $H^{3}(M;\Z)$ defined 
by the singular 3-cocycle $\omega$ 
is a representative for the 
Dixmier-Douady class of the bundle 
gerbe $(Q,Y,M)$.  
\end{proposition} 

We will assume this result for the moment 
and delay a proof until later.   
We use Proposition~\ref{prop:4.2.4} 
to give an alternative proof of the following 
theorem from  
\cite{CaCroMur}. 
\begin{theorem}[\cite{CaCroMur}] 
\label{thm:4.2.5}  
Let $\nu \in S^{2}(G;\Z)$ be a representative 
for the Chern class of $\hat{G}$.   
Then the transgression of $[\nu]$ in the 
principal bundle $P\to M$ is 
the Dixmier-Douady class of the lifting bundle 
gerbe $(Q,P,M)$ associated to 
$P\to M$ via the central extension $\cstar\to 
\hat{G}\to G$.   
\end{theorem} 

The fact that the transgression of $[\nu]$ is the 
Dixmier-Douady class $DD(Q)$ of $Q$ and not minus 
$DD(Q)$, as is the case in \cite{CaCroMur}, is 
due to differing sign conventions.  
  
\begin{proof} 
Suppose in the above construction of the 
singular 2-cocycle $F\in S^{2}(Y^{[2]};\Z)$ 
and the singular 2-cochain $C\in S^{1}(Y^{[3]};
\Z)$ that we had been unable to define a 
set map $s:Y^{[2]}\to Q$ which was a section of 
$\pi_{Q}:Q\to Y^{[2]}$ and which satisfied 
equation~\ref{eq:4.2.11}.  We rewrite equation  
~\ref{eq:4.2.9} as 
\begin{equation} 
\label{eq:4.2.12} 
(\pi_{Q}\times_{Y^{[4]}}\pi_{Q}\times_
{Y^{[4]}}\pi_{Q})^{*}\d(C)_{(\sigma_{34},
\sigma_{23},\sigma_{12})} = dB^{'}_{(\sigma_{34},
\sigma_{23},\sigma_{12})} 
\end{equation} 
where $\pi_{Q}\times_{Y^{[4]}}\pi_{Q}\times_
{Y^{[4]}}\pi_{Q}$ denotes the projection from 
$\pi_{1}^{-1}(\pi_{1}^{-1}Q\times_{Y^{[3]}}
\pi_{3}^{-1}Q)\times_{Y^{[4]}}\pi_{3}^{-1}
\pi_{3}^{-1}Q = \pi_{2}^{-1}\pi_{1}^{-1}Q
\times_{Y^{[4]}}\pi_{4}^{-1}(\pi_{1}^{-1}Q
\times_{Y^{[3]}}\pi_{3}^{-1}Q) \stackrel{\text{def}}
{=} Q_{1234}$ to $Y^{[4]}$ 
(compare with Definition~\ref{def:3.2.1}) and 
where $B^{'}\in S^{0}(Q_{1234};\Z)$ is defined 
by 
$$
B^{'} = 1\times B + (1\times \tilde{m}_{Q})^{*}B - 
B\times 1 - (\tilde{m}_{Q}\times 1)^{*}B. 
$$
If we can show that $\d(B^{'}) = 0$ in 
$S^{0}(Q_{1234}^{[2]};\Z)$ then we will have 
$B^{'} = \d(E)$ for some $E\in S^{0}(Y^{[4]};\Z)$ 
and hence $\d(C) = dE$.  Also it is easy to 
see that we must have $\d(E) = 0$ in 
$S^{0}(Y^{[5]};\Z)$.  We see that $\d(B^{'}) = 0$ 
as follows.  Recall that $\d(B) = \rho\in 
S^{0}((\pi_{1}^{-1}Q\times_{Y^{[3]}}\pi_{3}^{-1}Q);
\Z)$.  Hence we get $\d(B^{'}) = \rho^{'}$ where 
$\rho^{'} \in S^{0}(Q_{1234}^{[2]};\Z)$ is defined 
by 
$$
1\times \rho + (1\times \tilde{m}^{[2]})^{*}
\rho - \rho\times 1 - (\tilde{m}^{[2]}\times 1)^{*}
\rho.  
$$
It is a tedious but straightforward calculation 
to verify, using equation~\ref{eq:4.2.2}, that 
$\rho^{'} = 0$ and we will omit it.     
Notice that we could have considered a 
simplicial line bundle (Definition~\ref{def:3.2.2}) on an arbitrary 
simplicial manifold, the constructions above 
will still work but we will replace the 
$\pi_{i}$ by the face operators $d_{i}$ for 
whatever simplicial manifold we are working with.  
Of course we will not be able to use Proposition~\ref{prop:3.3.1}, 
as that only applies for the special simplicial 
manifold $Y=\{Y_{p}\}$ with $Y_{p} = Y^{[p+1]}$.  So it follows that 
we could consider the central extension 
$\cstar\to \hat{G}\to G$ as defining a simplicial 
line bundle (or simplicial $\cstar$ bundle rather) 
on the simplicial manifold $NG$ 
as per Proposition~\ref{prop:3.2.3} and constructed 
$\nu\in S^{2}(G;\Z)$ representing the Chern class 
of $\hat{G}$ in $H^{2}(G;\Z)$, $\mu\in S^{1}(G^{2};\Z)$ 
such that $d_{0}^{*}\nu - d_{1}^{*}\nu + d_{2}^{*}\nu = 
d\mu$ and $\lambda\in S^{0}(G^{3};\Z)$ such that 
$d_{0}^{*}\mu - d_{1}^{*}\mu + d_{2}^{*}\mu - d_{3}^{*}\mu = 
d\lambda$ and $d_{0}^{*}\lambda - d_{1}^{*}\lambda + d_{2}^{*}\lambda - 
d_{3}^{*}\lambda + d_{4}^{*}\lambda = 0$.  (Thus the triple 
$(\nu,\mu,\lambda)$ represents a class in $H^{3}(S^{\bullet, 
\bullet}(NG);\Z)$ --- the  
cohomology of the total complex of the double complex $S^{\bullet, 
\bullet}(NG)$ of Theorem~\ref{thm:2.2.3} --- 
and hence a class in $H^{3}(BG;\Z)$).       

We have the natural map $\tau:P^{[2]}\to G$ 
defined by $p_{2} = p_{1}\tau(p_{1},p_{2})$ 
for $p_{1}$ and $p_{2}$ in the same fibre 
of $P$.  In fact $\tau$ gives rise to a 
simplicial map $P\to NG$ 
where $P = \{P_{n}\}$ is the simplicial manifold 
with $P_{n} = P^{[n+1]}$, a portion of which is pictured in the 
following diagram: 
$$
\xymatrix{ 
 P^{[3]} \ar[r] \ar[d] \ar@<-2ex>[d] \ar@<2ex>[d]   
 & G\times G \ar[d] 
\ar@<-2ex>[d] \ar@<2ex>[d]             \\ 
P^{[2]} \ar[r]^-{\tau} \ar@<-1ex>[d] 
\ar@<1ex>[d] & G \ar@<-1ex>[d] \ar@<1ex>[d]     \\ 
P \ar[r] &             \{1\}.         } 
$$
By pulling back $\nu$, $\mu$ and $\lambda$ to 
$S^{2}(P^{[2]};\Z)$, $S^{1}(P^{[3]};\Z)$ and 
$S^{0}(P^{[4]};\Z)$ respectively by this simplicial 
map and using the conditions $d_{0}^{*}\nu - 
d_{1}^{*}\nu + d_{2}^{*}\nu = d\mu$, $d_{0}^{*}\mu - 
d_{1}^{*}\mu + d_{2}^{*}\mu - d_{3}^{*}\mu  = d\lambda$ 
and $d_{0}^{*}\lambda - d_{1}^{*}\lambda + 
d_{2}^{*}\lambda - d_{3}^{*}\lambda + d_{4}^{*}\lambda = 0$ 
it follows from Proposition~\ref{prop:3.3.2} that 
we can adjust $\tau^{*}\nu$ by a coboundary $d\a$ 
so that $F = \tau^{*}\nu + d\a$ satisfies $\d(F) = 0$ 
in $S^{2}(P^{[3]};\Z)$.  It follows that there exists 
$f \in S^{2}(P;\Z)$ such that $\d(f) = F$ and hence 
there exists $\omega \in S^{3}(M;\Z)$ with 
$d\omega = 0$ and $\pi^{*}\omega = df$.  We need to 
show that $[\omega]$ is the transgression of $[\nu]$.  
Recall from \cite{Bor} that the transgression $[t(\nu)]$ 
of $[\nu]$ is characterised by the property that there 
exists $e\in S^{2}(P;\Z)$ such that $de = \pi^{*}c$, 
where $c$ is a representative of $[t(\nu)]$, and 
$[i^{*}e|_{P_{m}}] = [\nu]$, where $i:G\to P_{m}$ is the injection 
of the fibre $G$ into $P_{m}$ and where 
$e|_{P_{m}}$ denotes the restriction of $e$ to the 
fibre $P_{m}$.  So we need to show that 
$[i^{*}f|_{P_{m}}] = [\nu]$.  The injection $i$ is 
defined by choosing a point $p$ of $P_{m}$ and 
mapping $g$ to $pg$ (this needs $G$ connected).  If 
we consider $\{p\}\times P_{m}\subset P^{[2]}$ then 
$\d(f) = F = \tau^{*}\nu + d\a$ gives us 
$[i^{*}f|_{P_{m}}] = [\nu]$.  Thus $\omega$ is the 
transgression of $\nu$.  Hence from Proposition~\ref{prop:4.2.4} 
we have that the transgression and the Dixmier-Douady 
class are equal.           
\end{proof} 

In general, given a central extension of 
groups $A\to \hat{G}\to G$, if one gives each of 
the groups the discrete topology then there will 
not exist a homomorphism $s:G\to \hat{G}$ which 
is also a section of $\hat{G}\to G$.  In 
certain circumstances however, such a section can 
be shown to exist.  An example of this occurs if 
$A$ is an abelian topological group and we consider the  
central extension 
of groups $A\to EA\to BA$.  Let $A^{\d}$, $EA^{\d}$ 
and $BA^{\d}$ denote the groups $A$, $EA$ and 
$BA$ equipped with the discrete topology.  Then 
we can define a section $s:BA^{\d}\to EA^{\d}$ 
of $EA^{\d}\to BA^{\d}$ which is also a group 
homomorphism as follows.  A point of $BA$ is of 
the form 
$$
|(x_{0},\ldots,x_{p}),(a_{1},\ldots,a_{p})|, 
$$
or, in the non-homogenous coordinates, 
$$
|(t_{1},\ldots,t_{p}),[a_{1}|\cdots|a_{p}]|, 
$$
where $t_{i} = x_{0} + \cdots + x_{i-1}$.  A 
point of $EA$ is of the form 
$$
|(x_{0},\ldots,x_{p}),(a_{0}^{'},\ldots,a_{p}^{'})|, 
$$
or, expressed in the non-homogenous coordinates, 
$$
|t_{1},\ldots,t_{p},a_{0}^{'}[(a_{0}^{'})^{-1}a_{1}^{'}|
\cdots |(a_{p-1}^{'})^{-1}a_{p}^{'}]|. 
$$
The projection $p:EA\to BA$ is given by 
$$
p(|(x_{0},\ldots,x_{p}),(a_{0}^{'},\ldots, 
a_{p}^{'})|) = |(x_{0},\ldots,x_{p}),(
(a_{0}^{'})^{-1}a_{1}^{'},\ldots, (a_{p-1}^{'})^{-1}
a_{p}^{'})| 
$$
or, alternatively, 
$$
p(|t_{1},\ldots,t_{p},a_{0}^{'}[(a_{0}^{'})^{-1}
a_{1}^{'}|\cdots |(a_{p-1}^{'})^{-1}a_{p}^{'}]|) =  
|t_{1},\ldots,t_{p},[(a_{0}^{'})^{-1}a_{1}^{'}|\cdots 
|(a_{p-1}^{'})^{-1}a_{p}^{'}]|.  
$$
So what we seek is a homomorphism 
$s:BA^{\d}\to EA^{\d}$ such that $p\circ s = 1$.  
Let $|(x_{0},\ldots,x_{p}),(a_{1},\ldots,a_{p})|$ 
be a point of $BA^{\d}$ and assume that there are no $0's$ 
appearing amongst the $x_{i}$.  Define 
$$
s(|(x_{0},\ldots,x_{p}),(a_{1},\ldots,a_{p})|) = 
|(x_{0},\ldots,x_{p}),(1,a_{1},a_{1}a_{2},\ldots, 
a_{1}\cdots a_{p})|.  
$$
We just need to check that $s$ is well defined 
with respect to degeneracies.  We have 
\begin{eqnarray*} 
&    & s((x_{0},\ldots,x_{i}+x_{i+1},\ldots,x_{p}),(a_{1},\ldots,a_{p-1})) \\ 
& = & ((x_{0},\ldots,x_{i}+x_{i+1},\ldots,x_{p}),(1,a_{1},a_{1}a_{2},
\ldots,a_{1}\cdots a_{p-1}))                                                \\ 
& = & ((x_{0},\ldots,x_{p}),(1,a_{1},\ldots,a_{1}\cdots a_{i-1},
a_{1}\cdots a_{i},a_{1}\cdots a_{i},a_{1}\cdots a_{i+1},\ldots, 
a_{1}\cdots a_{p-1}. 
\end{eqnarray*}    
On the other hand we have 
\begin{eqnarray*} 
&   & s((x_{0},\ldots,x_{p}),(a_{1},\ldots,a_{i},1,a_{i+1},\ldots,a_{p-1})) \\ 
& = & ((x_{0},\ldots,x_{p}),(1,a_{1},\ldots,a_{1}\cdots a_{i-1},
a_{1}\cdots a_{i},a_{1}\cdots a_{i},a_{1}\cdots a_{i+1},\ldots, 
a_{1}\cdots a_{p-1}))                                                      \\ 
& = & s((x_{0},\ldots,x_{i}+x_{i+1},\ldots,x_{p}),(a_{1},\ldots,a_{p-1})). 
\end{eqnarray*} 
Hence $s$ is well defined.  To check that 
$s$ is a homomorphism of discrete groups  
we will use the non-homogenous coordinates.  
The condition that there are no $x_{i}$'s equal to zero in 
$(x_{0},\ldots,x_{p})$ translates into no $t_{i}$'s 
repeated in $(t_{1},\ldots,t_{p})$.  $s$ expressed 
in the non-homogenous coordinates is 
$$
s(|t_{1},\ldots,t_{p},[a_{1}|\cdots|a_{p}]|) = 
|t_{1},\ldots,t_{p},1[a_{1}|\cdots|a_{p}]|.  
$$
Let $|t_{p+1},\ldots,t_{p+q},[a_{p+1}|\cdots|
a_{p+q}]|$ be another point of $BA^{\d}$.  We have 
\begin{eqnarray*} 
&   & s(|t_{1},\ldots,t_{p},[a_{1}|\cdots|a_{p}]|)
s(|t_{p+1},\ldots,t_{p+q},[a_{p+1}|\cdots|a_{p+q}]|) \\ 
& = & |t_{1},\ldots,t_{p},1[a_{1}|\cdots|a_{p}]|\cdot 
|t_{p+1},\ldots,t_{p+q},1[a_{p+1}|\cdots|a_{p+q}]|    \\ 
& = & |t_{\s(1)},\ldots,t_{\s(p+q)},1[a_{\s(1)}|\cdots |
a_{\s(p+q)}]|, 
\end{eqnarray*} 
where $\s$ is a permutation of $\{1,\ldots,p+q\}$ 
such that $t_{\s(1)}\leq \cdots\leq t_{\s(p+q)}$.  
This equals 
$$
|t_{\s(1)},t_{\s(i)},\ldots,t_{\s(p+q)},1[a_{\s(1)}|
a^{'}|\cdots |a^{''}]| 
$$
for some $a^{'}$, $a^{''}$ upon removing all repeated 
$t$'s.  This is exactly the result of applying $s$ to 
the point 
$$
|t_{1},\ldots,t_{p},[a_{1}|\cdots|a_{p}]|\cdot |t_{p+1},
\ldots,t_{p+q},[a_{p+1}|\cdots|a_{p+q}]|.  
$$
Therefore $s:BA^{\d}\to EA^{\d}$ is a homomorphism 
and it is straightforward to check that $s$ is also 
a section of $p$.  Note however that $s$ considered as 
a map $s:BA\to EA$ is not continuous.  

It now follows from the remarks made during the course 
of the proof of Theorem~\ref{thm:4.2.5} that if we regard 
the central extension $\cstar\to E\cstar\to B\cstar$ as 
a simplicial $\cstar$ bundle on the simplicial manifold 
$N\cstar$ then we can choose a representative 
$F \in S^{2}(B\cstar;\Z)$ of the Chern class of 
$E\cstar\to B\cstar$ so that there exists 
$C\in S^{1}((B\cstar)^{2};\Z)$ such that the following 
two equations are satisfied: 
\begin{eqnarray} 
&   & d_{0}^{*}F - d_{1}^{*}F + 
d_{2}^{*}F = dC,\label{eq:numu1} \\ 
&   & d_{0}^{*}C - d_{1}^{*}C + 
d_{2}^{*}C- d_{3}^{*}C  = 0.\label{eq:numu2}  
\end{eqnarray} 
It also follows that we can find $A_{2} \in S^{1}(E\cstar;\Z)$ 
such that $p^{*}F = dA_{2}$, $\d(A_{2}) = \tau^{*}c - d\a$, 
and such that there also exists $B_{2} \in S^{0}(E\cstar\times E\cstar;\Z)$ 
so that the following equation in $S^{1}(E\cstar\times E\cstar;\Z)$ 
is valid: 
$$
p_{1}^{*}A_{2} - m^{*}A_{2} + p_{2}^{*}A_{2} = (p\times p)^{*}C + dB_{2}. 
$$
Let us rename $F$ to $\iota_{2}$ and 
$C$ to $\kappa_{2}$.  If one examines the 
construction of $\iota_{2}$ and $\kappa_{2}$ 
then one notices the construction relies upon the 
fact proved in Lemma~\ref{lemma:4.2.1} 
that there exist $c\in S^{1}(\cstar;\Z)$ 
and $b\in S^{0}((\cstar)^{2};\Z)$ satisfying the 
analogues of equations~\ref{eq:numu1} and~\ref{eq:numu2}.  

Let us see what happens if we try to iterate 
this procedure using $\iota_{2}$ in place of 
$c$ and $\kappa_{2}$ in place of $b$.  So we start 
with the principal $B\cstar$ bundle $EB\cstar\to 
BB\cstar$.  We form the canonical map $\tau:(EB\cstar
)^{[2]}\to B\cstar$ and pullback $\iota_{2}$ 
to a class $\tau^{*}\iota_{2}$ on $(EB\cstar)^{[2]}$.  
Then as above we get $\d(\tau^{*}\iota_{2}) = d(
\tau\circ \pi_{1},\tau\circ \pi_{3})^{*}\kappa_{2}$.  
Because $\kappa_{2}$ satisfies the analogue of 
equation~\ref{eq:4.2.2} we get $\d((\tau\circ \pi_{1},
\tau\circ \pi_{3})^{*}\kappa_{2}) = 0$ and hence 
there exists $\a_{3} \in S^{1}((EB\cstar)^{[2]};\Z)$ 
such that $(\tau\circ \pi_{1},\tau\circ \pi_{3})^{*}
\kappa_{2} = \d(\a_{3})$.  Therefore $\tau^{*}\iota_{2} 
-d\a_{3}$ satisfies $\d(\tau^{*}\iota_{2}-d\a_{3}) = 0$ 
and so we can solve $\tau^{*}\iota_{2}-d\a_{3} = \d(A_{3})$ 
for some $A_{3}\in S^{2}(EB\cstar;\Z)$.  Hence there exists 
$\iota_{3} \in S^{3}(BB\cstar;\Z)$ with $p^{*}\iota_{3} 
=dA_{3}$.  Following this procedure through leads in 
a straightforward way to the following analogue    
of equation~\ref{eq:4.2.8}: 
$$
p_{1}^{*}A_{3} + p_{2}^{*}
A_{3} = m^{*}A_{3} + dB_{3} + (p\times p)^{*}\kappa_{3}, 
$$
where $\kappa_{3}\in S^{2}(BB\cstar\times BB\cstar;\Z)$ and 
where $B_{3}\in S^{1}(EB\cstar\times EB\cstar;\Z)$ 
satisfies 
$$
\d(B_{3}) = (\tau\circ p_{1}^{[2]}\times \tau\circ 
p_{2}^{[2]})^{*}\kappa_{2} - (p_{1}^{[2]})^{*}\a_{3} 
- (p_{2}^{[2]})^{*}\a_{3} + (m^{[2]})^{*}\a_{3}. 
$$
This time however there will not exist a section 
$s:S_{1}(BB\cstar;\Z)\to S_{1}(EB\cstar;\Z)$ compatible 
with the group structures on $S_{1}(BB\cstar;\Z)$ and 
$S_{1}(EB\cstar;\Z)$ induced by the group structures 
on $BB\cstar$ and $EB\cstar$ and so we can only conclude, 
following the argument of Theorem~\ref{thm:4.2.5}, with 
$\iota_{2}$ in place of $c$ and $\kappa_{2}$ in place of 
$b$, that there exists $\lambda_{3}\in S^{1}((EB\cstar)^{3};\Z)$ 
such that 
\begin{eqnarray} 
& & d_{0}^{*}\kappa_{3} - d_{1}^{*}\kappa_{3} + d_{2}^{*}\kappa_{3} 
- d_{3}^{*}\kappa_{3} = d\lambda_{3}        \\ \label{eq:lambda2}  
& & d_{0}^{*}\lambda_{3} - d_{1}^{*}\lambda_{3} + d_{2}^{*}\lambda_{3} 
- d_{3}^{*}\lambda_{3} + d_{4}^{*}\lambda_{3} = 0. \label{eq:lambda1}  
\end{eqnarray}  
Therefore it follows that we have 
defined $\iota_{3}\in S^{3}(BB\cstar;\Z)$, $A_{3}\in 
S^{2}(EB\cstar;\Z)$, $\kappa_{3} \in S^{2}((BB\cstar)^{2};\Z)$, 
$B_{3} \in S^{1}((EB\cstar)^{2};\Z)$ and 
$\lambda_{3} \in S^{1}((EB\cstar)^{3};\Z)$ such that 
$d\iota_{3} = 0$, $p^{*}\iota_{3} = dA_{3}$ and 
\begin{eqnarray*} 
& & p_{1}^{*}A_{3} + p_{2}^{*}A_{3} = m^{*}A_{3} + dB_{3} 
+ (p\times p)^{*}\kappa_{3}                                 \\ 
& & d_{1}^{*}\iota_{3} - d_{1}^{*}\iota_{3} + d_{2}^{*}\iota_{3} 
=d\kappa_{3}                                                    \\ 
& & d_{0}^{*}\kappa_{3} - d_{1}^{*}\kappa_{3} +d_{2}^{*}\kappa_{3} 
- d_{3}^{*}\kappa_{3} = d\lambda_{3},            
\end{eqnarray*} 
and $\lambda_{3}$ and 
$\kappa_{3}$ satisfy equations~\ref{eq:lambda1} and~\ref{eq:lambda2} 
above.   
It is clear that $\iota_{3}$ is the transgression 
of $\iota_{2}$ in the principal $B\cstar$ bundle 
$EB\cstar\to BB\cstar$.  It is also clear that 
$\iota_{2}$ is the transgression of $c$ in the $\cstar$ 
bundle $E\cstar\to B\cstar$.  
The fundamental class 
of $H^{p+1}(B^{p}\cstar;\Z)$  
is mapped to the fundamental 
class of $H^{p+2}(B^{p+1}\cstar;\Z)$ under the 
transgression in the principal $B^{p}\cstar$ 
bundle $EB^{p}\cstar\to B^{p+1}\cstar$.  Indeed if we consider the 
long exact sequence in homotopy for the fibering 
$B^{p}\cstar\to EB^{p}\cstar\to B^{p+1}\cstar$ then we 
get $\pi_{p+1}(B^{p}\cstar) 
\simeq \pi_{p+2}(B^{p+1}\cstar)$.  Moreover, since these are 
Eilenberg-MacLane spaces we have the sequence of 
isomorphisms $\pi_{q+1}(B^{q}\cstar)\simeq  
H^{q+1}(B^{q}\cstar;\Z) \simeq  \Z$ and with respect to these 
isomorphisms the isomorphism $\pi_{p+1}(B^{p}\cstar)\simeq 
\pi_{p+2}(B^{p+1}\cstar)$ can easily be seen to be the 
transgression (in the case we are considering the transgression 
$H^{p+1}(B^{p}\cstar;\Z) \to H^{p+2}(B^{p+1}\cstar;\Z)$ is an 
isomorphism).   
It follows that 
$\iota_{2}$ is a representative of 
the fundamental class of $H^{2}(B\cstar;\Z)$ 
and that $\iota_{3}$ represents the fundamental class of 
$H^{3}(BB\cstar;\Z)$.  It would not be too 
hard to go further and prove a similar result   
by induction for the fundamental class 
$\iota_{n} \in S^{n}(B^{n-1}\cstar;\Z)$ although we will not do so here. 

Recall from Section~\ref{sec:4.1} above the recipe for 
constructing from a \v{C}ech cocycle $g_{ijk}:U_{ijk}\to 
\cstar$ maps $g_{ij}:U_{ij}\to B\cstar$ satisfying 
the cocycle condition $g_{ij}g_{jk} = g_{ik}$ and maps 
$\hat{g}_{ij}:U_{ij}\to E\cstar$ lifting $g_{ij}$ such 
that $\hat{g}_{jk}\hat{g}_{ik}^{-1}\hat{g}_{ij} = 
g_{ijk}$.  It is not difficult to go further and 
construct maps $g_{i}:U_{i}\to BB\cstar$ such that 
$g_{i} = g_{j}$ on $U_{ij}$ and lifts $\hat{g}_{i}:U_{i}
\to EB\cstar$ of $g_{i}$ which satisfy $\hat{g}_{j}
\hat{g}_{i}^{-1} = g_{ij}$.  The map $g:M\to BB\cstar$ 
is obtained from the $g_{i}$.  Suppose we are given 
such a cocycle $g_{ijk}$.  We will use the \v{C}ech-singular 
double complex (see \cite{BotTu}) to show that the 
class in $H^{3}(M;\Z)$ represented by $g_{ijk}$ 
and the class represented by $g^{*}\iota_{3}$ are equal.  

Let $i_{U_{i}}:U_{i}\to M$ denote the inclusion map.  
We have $i_{U_{i}}^{*}g^{*}\iota_{3} = g_{i}^{*}\iota_{3} 
= (p\circ \hat{g}_{i})^{*}\iota_{3} = d\hat{g}_{i}^{*}A_{3}$.  
Next, $\hat{g}_{j}^{*}A_{3} - \hat{g}_{i}^{*}A_{3} = 
(\hat{g}_{i},\hat{g}_{j})^{*}\d(A_{3}) = (\hat{g}_{i},
\hat{g}_{j})^{*}(\tau^{*}\iota_{2} - d\a_{3})$.  Since 
$\hat{g}_{j}\hat{g}_{i}^{-1} = g_{ij}$ we get that 
$\hat{g}_{j}^{*}A_{3} - \hat{g}_{i}^{*}A_{3} = g_{ij}^{*}
\iota_{2} - d(\hat{g}_{i},\hat{g}_{j})^{*}\a_{3}$.  
Using $p(\hat{g}_{ij}) = g_{ij}$ we get $g_{ij}^{*}\iota_{2} = 
d\hat{g}_{ij}^{*}A_{2}$ ,where $A_{2}\in S^{1}(E\cstar;\Z)$ is constructed 
as in the paragraph preceding Lemma~\ref{lemma:4.2.2}.  Let 
$\a_{i} = \hat{g}_{i}^{*}A_{2}$ and let $\b_{ij} = \hat{g}_{ij}^{*} 
- (\hat{g}_{i},\hat{g}_{j})^{*}\a_{3}$.  Since 
$\d(\a_{3}) = (\tau\circ \pi_{1},\tau\circ \pi_{3})^{*}\kappa_{2}$, 
we get 
$$
\b_{jk} - \b_{ik} + \b_{ij} = \hat{g}_{jk}^{*}A_{2} - \hat{g}_{ik}^{*}A_{2} 
+ \hat{g}_{ij}^{*}A_{2} - (g_{jk},g_{ij})^{*}\kappa_{2}. 
$$
Using the equation $p_{1}^{*}A_{2} + p_{2}^{*}A_{2} 
= m^{*}A_{2} + dB_{2} + (p\times p)^{*}\kappa_{2}$,    
we get 
$$
\d(\b_{ij}) = (\hat{g}_{ik}\cdot g_{ijk})^{*}A_{2} - \hat{g}_{ik}^{*}A_{2} 
+ d(\hat{g}_{jk},\hat{g}_{ij})^{*}B_{2}. 
$$
Using the same equation again we get 
$$
(i(g_{ijk}))^{*}A_{2} - d(\hat{g}_{ik},i(g_{ijk}))^{*}B_{2} + 
d(\hat{g}_{jk},\hat{g}_{ij})^{*}B_{2} - (g_{ik},1)^{*}\kappa_{2}.  
$$
It is straightforward to check that $(g_{ik},1)^{*}\kappa_{2} = 0$.  
It is also not difficult to check that $i^{*}A_{2} = c$.  Let 
$\hat{g}_{ijk}:U_{ijk}\to \C$ be a lift of $g_{ijk}$.  
We have 
$$
\b_{jk} - \b_{ik} + \b_{ij} = d(\hat{g}_{ijk}^{*}\hat{c} 
- (\hat{g}_{ik},i(g_{ijk}))^{*}B_{2} + (\hat{g}_{jk},\hat{g}_{ij})^{*}B_{2}), 
$$
where $\hat{c} \in S^{0}(\C;\Z)$ is constructed in 
Lemma~\ref{lemma:4.2.2}.  Put $\c_{ijk} = \hat{g}_{ijk}^{*}
\hat{c} - (\hat{g}_{ik},i(g_{ijk}))^{*}B_{2} + (\hat{g}_{jk},
\hat{g}_{ij})^{*}B_{2}$.  Define $z_{ij}$, $z_{jk}$ and 
$z_{ik}$ in $\cstar$ by $z_{ij} = \tau(s(g_{ij}),g_{ij})$, 
$z_{jk} = \tau(s(g_{jk}),g_{jk})$ and $z_{ik} = \tau(s(g_{ik}),
g_{ik})$ respectively.  Then one can check that 
$\c_{ijk}$ is given by 
$$
\hat{g}_{ijk}^{*}\hat{c} - (\sigma(z_{jk}) - \sigma(z_{ik}) + 
\sigma(z_{ij})), 
$$ 
where $\sigma:(\cstar)^{\d}\to (\C)^{\d}$ is the 
section constructed in Lemma~\ref{lemma:4.2.2}.  It now follows that 
$\c_{jkl} - \c_{ikl} + \c_{ijl} - \c_{ijk} = n_{ijkl}$ 
where $n_{ijkl}:U_{ijkl}\to \Z$ is the 
image of the cocycle $g_{ijk}$ under the 
coboundary map $\check{H}^{2}(M;\underline{\C}^{\times}_{M})\to 
\check{H}^{3}(M;\underline{\Z}_{M})$.  Hence $g_{ijk}$ 
and $g^{*}\iota_{3}$ represent the same class in 
$H^{3}(M;\Z)$.  This argument constitutes a proof of 
Proposition~\ref{prop:4.2.4}   

As an application of this, if we are given a bundle gerbe 
$(P,X,M)$ on $M$ then we know that there is a map 
$g:M\to BB\cstar$, unique up to homotopy, such that the 
Dixmier-Douady class $DD(P)$ of $P$ is equal to 
$g^{*}\iota_{3}$.  It follows that the pullback 
$\pi_{X}^{*}g^{*}\iota_{3}$ of the Dixmier-Douady class 
to $X$ is cohomologous to zero.  Hence there is a lift 
$\hat{g}:X\to EB\cstar$ of the map $g\circ \pi_{X}:X\to BB\cstar$, 
ie a map $\hat{g}:X\to EB\cstar$ making the following 
diagram commute: 
$$
\xymatrix{ 
X \ar[r]^-{\hat{g}} \ar[d]_{\pi_{X}} & EB\cstar \ar[d]^{p} \\ 
M \ar[r]_-{g} & BB\cstar              } 
$$
Notice that we can therefore define a map 
$f:X^{[2]}\to B\cstar$ satisfying 
$$
f(x_{2},x_{3})
f(x_{1},x_{2}) = f(x_{1},x_{3})
$$ 
for $x_{1}$, 
$x_{2}$ and $x_{3}$ in the same fibre of $X$, by 
$f(x_{1},x_{2}) = i^{-1}(\hat{g}(x_{2})\hat{g}(x_{1})^{-1})$.  
It is easy to see that the $\cstar$ bundle 
$f^{-1}E\cstar\to X^{[2]}$ has a bundle gerbe 
product and that the resulting bundle gerbe 
$(f^{-1}E\cstar,X,M)$ has Dixmier-Douady class 
equal to $DD(P)$.  Therefore the bundle gerbes 
$f^{-1}E\cstar$ and $P$ are stably isomorphic and so 
there exists a $\cstar$ bundle $L\to X$ such that 
$P \simeq f^{-1}E\cstar \otimes \d(L)$.  Let 
$f_{L}:X\to B\cstar$ be a classifying map for 
$L$ and redefine $\hat{g}$ to $\hat{g}\cdot f_{L}$.  
Denote the rescaled $f$ by the same letter as well.  
It follows that the bundle gerbes $P$ and $f^{-1}E\cstar$ 
are isomorphic.  In effect we have proved the 
following Proposition. 

\begin{proposition} 
Let $(P,X,M)$ be a bundle gerbe.  Then there is a 
bundle gerbe morphism $\bar{g}:(P,X,M)\to 
(\widetilde{EB\cstar},EB\cstar,BB\cstar)$ 
with $\bar{g} = (\tilde{g},\hat{g},g)$.  The map 
$g:M\to BB\cstar$ is unique up to homotopy. 
\end{proposition}

\setcounter{chapter}{4}
\chapter{Open covers}
\label{chapter:5} 
\section{Bundle gerbes via open covers} 
\label{sec:5.1} 
We now turn to some calculations involving 
bundle gerbe connections and curvings on 
open covers.  Similar results for gerbes  
appear in the lecture notes  
\cite{Hit}.   

Suppose $(P,Y,M)$ and $(Q,X,M)$ are bundle gerbes.  
and we have a bundle gerbe morphism 
$\bar{f} = (\hat{f},f,\text{id})$ which 
maps $(P,Y,M)$ to $(Q,X,M)$.  It is easily seen that 
$P$ and $Q$ must have the same Dixmier-Douady class.  We 
will exploit this fact to give a method of calculating 
the Dixmier-Douady class of a bundle 
gerbe $P$ on $M$ from an open cover of $M$.

Let $\{U_{\a}\}_{\a \in \Sigma}$ be an open cover of $M$, and 
suppose we are given line bundles $L_{\a\b} \to U_{\a\b}$.
Suppose further that we have sections $\sigma_{\a\b\c}$ of 
the line bundles $L_{\b\c} \otimes L_{\a\c}^{*} \otimes L_{\a\b}$
over $U_{\a\b\c}$.  These sections are required to satisfy
$$
\sigma_{\b\c\d} \sigma_{\a\c\d}^{-1} 
\sigma_{\a\b\d} \sigma_{\a\b\c}^{-1} = 
\underline{1}, 
$$
where $\underline{1}$ is the non-vanishing 
section induced by the canonical trivialisation 
of 
$$
\d(L_{\a\b\c}) = L_{\b\c\d}\otimes L_{\a\c\d}^{*}\otimes 
L_{\a\b\d}\otimes L_{\a\b\c}^{*} 
$$
over $U_{\a\b\c\d}$ and where we have 
put $L_{\a\b\c} = L_{\b\c} \otimes L_{\a\c}^{*} \otimes L_{\a\b}.$
Assume that we are also given connections 
$\nabla_{\a\b}$ on each $L_{\a\b}$.  Denote 
the induced connection on $L_{\a\b\c}$ by 
$\nabla_{\a\b\c}$.  Hence we can 
find 1-forms $A_{\a\b\c}$ on $U_{\a\b\c}$ such that 
$$
\nabla_{\a\b\c}(\sigma_{\a\b\c}) = A_{\a\b\c} \sigma_{\a\b\c}.  
$$
If we denote the induced connection on 
$\d(L_{\a\b\c})$ by $\d (\nabla_{\a\b\c})$ 
then $\d (\nabla_{\a\b\c})(\underline{1}) = 0$ 
so clearly we will have  
\begin{equation}
A_{\b\c\d} - A_{\a\c\d} + A_{\a\b\d} - A_{\a\b\c} = 0.
\end{equation}
Suppose we have a partition of unity $\{\psi_{\a}\}_{\a\in \Sigma}$
subordinate to the cover $\{U_{\a}\}_{\a\in \Sigma}$ of $M$.
Define a 1-form $B_{\a\b}$ on $U_{\a\b}$ as follows.
Let $B_{\a\b} = \sum_{\c \in \Sigma} \psi_{\c} A_{\a\b\c}$.
Then it is easily checked that we have
\begin{equation}
B_{\b\c} - B_{\a\c} + B_{\a\b} = A_{\a\b\c}
\end{equation}
Define new connections on the $L_{\a\b}$ by 
$\nabla_{\a\b} - B_{\a\b}$.  If $F_{\a\b}$ denotes 
the curvature of the connection $\nabla_{\a\b}$, then 
the new connection will have curvature $\F_{\a\b}$ equal to 
$F_{\a\b} -dB_{\a\b}$.
From the equations above, we see that 
\begin{equation}
F_{\b\c} - F_{\a\c} + F_{\a\b} = dA_{\a\b\c}.
\end{equation}
Hence the new curvature will satisfy
\begin{equation}
\F_{\b\c} - \F_{\a\c} + \F_{\a\b} = 0.
\end{equation}
Using the same trick as before, define 2-forms
$\mu_{\a}$ on $U_{\a}$ satisfying $\mu_{\b}-\mu_{\a} = \F_{\a\b}$ by
\begin{equation}
\mu_{\a} = - \sum_{\b \in \Sigma}\F_{\a\b}\psi_{\b}.
\end{equation}
One can check that on $U_{\a\b}$ we have 
$d\mu_{\a} = d\mu_{\b}$, and so the $d\mu_{\a}$ 
are the restrictions to $U_{\a}$ of a global 3-form $\omega$ 
on $M$.
It will follow from the discussion in the 
next section on Deligne cohomology that $\omega$ lies in the 
image of $H^{3}(M;\Z)$ inside $H^{3}(M;\Reals)$. 

\section{Bundle gerbes and Deligne cohomology} 
\label{sec:5.2} 
We now briefly review the notion of 
\emph{sheaf hypercohomology}.  We refer 
to \cite{Bry} or 
\cite{GriHar}
for more details.  

Suppose we are given a complex of sheaves 
$\underline{K}^{\bullet}$ on our manifold $M$.  Thus 
we have sheaves of abelian groups $\underline{K}^{q}$ on 
$M$ together with morphisms of sheaves 
$d_{q} : \underline{K}^{q} \to \underline{K}^{q+1}$  
which satisfy $d_{q+1}\circ d_{q} = 0$.  
Suppose we are given an open cover 
$\U = \{U_{\a}\}_{\a \in \Sigma}$ of 
$M$.  Then we can form a double complex 
$K_{\U}^{\bullet ,\bullet}$ with 
$K_{\U}^{p,q} = C^{p}(\U ,\underline{K}^{q})$, ie the 
\v{C}ech $p$-cochains with values in the sheaf 
$\underline{K}^{q}$.  The differentials of $K_{\U}^{\bullet ,\bullet}$ 
are $d : K_{\U}^{p,q} \to K_{\U}^{p,q+1}$ and 
$\d : K_{\U}^{p,q} \to K_{\U}^{p+1,q}$ and are induced 
by the morphism of sheaves $d : \underline{K}^{q}\to \underline{K}^{q+1}$ 
and the \v{C}ech coboundary operator 
$\d : C^{p}(\U ,\underline{K}^{q}) \to C^{p+1}(\U ,\underline{K}^{q})$ 
respectively.     
$d$ and $\d$ clearly satisfy $d^{2} = 0$, 
$\d^{2} = 0$ and $d\d + \d d = 0$, thus 
$K_{\U}^{\bullet, \bullet}$ is a double complex.  

Let $Tot(K_{\U}^{\bullet})$ denote the total 
complex of $K_{\U}^{\bullet, \bullet}$.  
A refinement $\V < \U$ will induce 
mappings $Tot(K_{\U}^{\bullet}) \to Tot(K_{\V}^{\bullet})$ 
and mappings 
$H^{\bullet}(Tot(K_{\U}^{\bullet})) 
\to H^{\bullet}(Tot(K_{\V}^{\bullet}))$ 
of the corresponding cohomology groups.  
We define the \emph{sheaf hypercohomology groups} 
$H^{\bullet}(M;\underline{K}^{\bullet})$  
of $M$ with values in the complex of 
sheaves $\underline{K}^{\bullet}$ to be 
$$
H^{\bullet}(M;K^{\bullet}) = \dirlim_{\U} H^{\bullet}(
Tot(K_{\U}^{\bullet})).  
$$
As an example consider the complex 
of sheaves $\underline{\C}^{\times}_{M}
\stackrel{dlog}{\to} \underline{\Omega}^{1}_{M}$ on 
a manifold $M$.  Then, relative to 
an open cover $\U = \{U_{\a}\}_{\a \in \Sigma}$ 
of $M$, a class in the hypercohomology 
group $H^{1}(M;\underline{\C}^{\times}_{M}
\stackrel{dlog}{\to}\underline{\Omega}^{1}_{M})$ 
consists of smooth maps 
$g_{\a\b}: U_{\a\b} \to \cstar$ and 
one forms $A_{\a}$ on $U_{\a}$ 
which together satisfy the conditions 
\begin{eqnarray*}
g_{\b\c} g_{\a\c}^{-1} g_{\a\b} & = & 1   \\
A_{\b} - A_{\a} & = & dg_{\a\b}g_{\a\b}^{-1}   
\end{eqnarray*}
Then it is not hard to see that classes 
in $H^{1}(M;\underline{\C}^{\times}_{M}
\stackrel{dlog}{\to}\underline{\Omega}^{1}_{M})$ are in 
a one to one correspondence with 
isomorphism classes of line bundles 
with connection on $M$.   

Consider the complex of sheaves 
$$
\underline{\C}^{\times}_{M} \stackrel{dlog}{\rightarrow} 
\underline{\Omega^{1}}_{M} \stackrel{d}{\rightarrow} 
\underline{\Omega^{2}}_{M} 
$$
on a manifold $M$.  It is shown in  
\cite{Mur} (see also \cite{MurSte}) that 
a bundle gerbe $(L,Y,M)$ with bundle 
gerbe connection $\nabla$ and curving $f$ gives rise to a 
class $D(L,\nabla,f)$ in the Deligne hypercohomology group 
$H^{2}(M;\underline{\C}^{\times}_{M}\to 
\underline{\Omega}^{1}_{M}\to 
\underline{\Omega}^{2}_{M})$.  Let us briefly recall 
how this class is constructed.  

Choose as usual an open cover 
$\U = \{U_{\a}\}_{\a \in \Sigma}$ of 
$M$ such that there exist local sections 
of $\pi : Y \to M$ above each $U_{\a}$ and 
all of whose non-empty finite intersections 
$U_{\a_{0}}\cap \cdots \cap U_{\a_{p}}$ are 
contractible.  
Form the maps $(s_{\a},s_{\b}) : U_{\a\b} \to Y^{[2]}$ 
and choose sections $\sigma_{\a\b}$ of the 
pullback line bundle $L_{\a\b} = (s_{\a},s_{\b})^{-1}L$ 
over $U_{\a\b}$.  The connection $\nabla$ 
induces a connection $\nabla_{\a\b}$ on 
$L_{\a\b}$ and we get 
$$
\nabla_{\a\b}(\sigma_{\a\b}) = A_{\a\b}\otimes \sigma_{\a\b}
$$
for some one form $A_{\a\b}$ on 
$U_{\a\b}$.  $L_{\a\b}$ has curvature 
$F_{\a\b} = (s_{\a},s_{\b})^{*}F_{\nabla}$ 
with respect to $\nabla_{\a\b}$.  Since 
$F_{\nabla} = \d(f)$ we get 
$F_{\a\b} = f_{\b} - f_{\a}$.  On the other 
hand, we clearly have $F_{\a\b} = dA_{\a\b}$, 
and so we get 
\begin{equation}
\label{eq: one}
dA_{\a\b} = f_{\b} - f_{\a}.
\end{equation}
Finally, the isomorphism of line bundles 
$m : \pi_{1}^{-1}L\otimes \pi_{3}^{-1}L \to \pi_{2}^{-1}L$ 
induces an isomorphism of line bundles 
$m : L_{\b\c}\otimes L_{\a\b} \to L_{\a\c}$ and 
we have $m(\sigma_{\b\c}\otimes \sigma_{\a\b}) = \sigma_{\a\c}g_{\a\b\c}$ 
where $g_{\a\b\c} : U_{\a\b\c} \to \cstar$ 
is a representative of the 
Dixmier-Douady class of $L$ and satisfies 
the cocycle condition 
\begin{equation}
\label{eq: two}
g_{\b\c\d} g_{\a\c\d}^{-1} g_{\a\b\d} g_{\a\b\c}^{-1} = 1.
\end{equation}
Therefore 
we have on the one hand 
$$
(\nabla_{\b\c} + \nabla_{\a\b})(\sigma_{\b\c}\otimes \sigma_{\a\b}) = 
(A_{\b\c} + A_{\a\b})\otimes (\sigma_{\b\c}\otimes \sigma_{\a\b}),
$$
but on the other, from equation~\ref{eq:bundlegerbeconnection},  
\begin{eqnarray*}
(\nabla_{\b\c} + \nabla_{\a\b})(\sigma_{\b\c}\otimes \sigma_{\a\b})    
& = & m^{-1}\circ \nabla_{\a\c}\circ m(\sigma_{\b\c}\otimes \sigma_{\a\b}) \\
 & = & m^{-1}\circ \nabla_{\a\c}(\sigma_{\a\c}g_{\a\b\c})                 \\ 
 & = & m^{-1}((A_{\a\c} + g_{\a\b\c}^{-1}dg_{\a\b\c})\otimes (\sigma_{\a\c} 
g_{\a\b\c})) \\ 
& = & (A_{\a\c} + g_{\a\b\c}^{-1}dg_{\a\b\c})\otimes (\sigma_{\b\c} 
\otimes \sigma_{\a\b}),  
\end{eqnarray*}
and hence we get 
\begin{equation}  
\label{eq: three}
A_{\b\c} - A_{\a\c} + A_{\a\b} = g_{\a\b\c}^{-1}dg_{\a\b\c}.
\end{equation}
From equations~\ref{eq: one},~\ref{eq: two} 
and~\ref{eq: three}, we see that the triple 
$(f_{\a},A_{\a\b},g_{\a\b\c})$ represents 
a class $D(L,\nabla,f)$ in the Deligne hypercohomology 
group $H^{2}(M;\underline{\C}^{\times}_{M}\to 
\underline{\Omega}^{1}_{M}\to \underline{\Omega}
^{2}_{M})$ and one can 
check that this is independent of the 
choices made.  If we consider bundle gerbes 
$(L,Y,M)$ equipped only with a bundle gerbe 
connection $\nabla$ then the discussion above 
shows that every such bundle gerbe $L$ with 
bundle gerbe connection $\nabla$ gives rise to 
a class in the truncated Deligne hypercohomology 
group $H^{2}(M;\underline{\C}^{\times}_{M}\to \underline{\Omega}^{1}_{M})$.   

Suppose we are given a bundle gerbe 
$(L,Y,M)$ with bundle gerbe connection 
$\nabla_{L}$ and curving $f$ and suppose moreover that the 
bundle gerbe $L$ is trivial.  Hence there 
exists a line bundle $J$ on $Y$ such that 
there is an isomorphism of bundle gerbes 
$\phi : L \to \d(J)$ over $Y^{[2]}$.  
Suppose that the line bundle $J$ has a 
connection $\nabla_{J}$.  Then $\d(J)$ 
inherits the tensor product connection 
$\d(\nabla_{J})$.  Therefore there exists 
a complex valued one form $\eta$ on 
$Y^{[2]}$ such that 
$\nabla_{L} = \phi^{-1}\circ \d(\nabla_{J})\circ \phi + \eta$.  
Since $\phi$ commutes with the bundle 
gerbe products, we must have $\d(\eta ) = 0$ 
and hence there exists a complex valued 
one form $\mu$ on $Y$ such that 
$\d(\mu ) = \eta$.  Then if we give 
$J$ the new connection $\nabla_{J} - \mu$ 
and rename this new connection $\nabla_{J}$, 
it is not hard to see that we have an 
isomorphism of bundle gerbes with bundle gerbe  
connection $\phi : L \to \d(J)$ over 
$Y^{[2]}$.  This result is not surprising 
in light of the fact that there is an 
isomorphism $H^{2}(M;\underline{\C}^{\times}_{M})\cong 
H^{2}(M;\underline{\C}^{\times}_{M}\to \underline{\Omega}^{1}_{M})$ 
of cohomology groups. 

Suppose that the class $D(L,\nabla_{L},f)$ of $L$ 
in the Deligne hypercohomology group 
$H^{2}(M;\underline{\C}^{\times}_{M}\to \underline{
\Omega}^{1}_{M}\to \underline{\Omega}_{M}^{2})$ 
is zero, so that (using the notation above) we 
have 
\begin{eqnarray*} 
g_{\a\b\c} & = & h_{\b\c}h_{\a\c}^{-1}h_{\a\b}      \\ 
A_{\a\b} & = & -k_{\b} + k_{\a} + dh_{\a\b}h_{\a\b}^{-1} \\ 
f_{\a} & = & dk_{\a}.                                        
\end{eqnarray*} 
From $f_{\a} = dk_{\a}$ we see that the three curvature 
$\omega$ of the bundle gerbe $L$ equipped with the bundle gerbe 
connection $\nabla_{L}$ and curving $f$ is zero.  In 
other words the map $H^{2}(M;\underline{\C}^{\times}
_{M}\to \underline{\Omega}_{M}^{1}\to \underline{
\Omega}_{M}^{2})\to \Omega^{3}(M)$ applied to 
$D(L,\nabla_{L},f)$ is zero.  So we must have 
$df = \pi^{*}(\omega) = 0$.  The bundle gerbes 
$L$ and $\d(J)$ may not have the same 
curvings, in general the curvings $f$ and 
$F_{\nabla_{J}}$ of $L$ and $\d(J)$ respectively 
will differ by the pullback of a two form $K$ on $M$ 
(here $F_{\nabla_{J}}$ is the 
curvature of $J$ with respect to the connection 
$\nabla_{J}$ defined above).   
Since $F_{\nabla_{J}}$ 
is the curvature two form of a connection on $J$ 
it is closed and we have just seen that $f$ is 
closed, so it follows that $K$ must also be 
closed.  Since $K$ is also an integral two 
form, there exists a line bundle $I$ on 
$M$ with Chern class equal to the class $[K]$ in 
$H^{2}(M;\Z)$.  We may endow $I$ with 
a connection $\nabla_{I}$ whose 
curvature is equal to $K$.  It follows that 
the line bundles $L$ and $\d(J\otimes \pi^{-1}I^{*})$ 
are isomorphic by an isomorphism which 
preserves the bundle gerbe products and 
the connections $\nabla_{L}$ and $\nabla_{J} - 
\pi^{-1}\nabla_{I}$.  What is more, the 
bundle gerbes $L$ and $\d(J\otimes \pi^{-1}I^{*})$ 
have the same curvings.  In effect, we have 
proved the following proposition.  

\begin{proposition}[\cite{MurSte}] 
\label{prop:5.2.1} 
A bundle gerbe $(L,Y,M)$ with bundle gerbe 
connection $\nabla_{L}$ and curving $f$ 
has trivial class $D(L,\nabla_{L},f)$ in 
$H^{2}(M;\underline{\C}^{\times}_{M}\to 
\underline{\Omega}^{1}_{M}\to 
\underline{\Omega}^{2}_{M})$ if and only 
if $\nabla_{L} = \d(\nabla_{J})$, where 
$\nabla_{J}$ is a connection on a line bundle 
$J\to Y$ which trivialises $L$ and $f = F_{
\nabla_{J}}$, where $F_{\nabla_{J}}$ is the 
curvature of the line bundle $J$ with 
respect to the connection $\nabla_{J}$.  
\end{proposition} 

\begin{note} 
One can define a notion of stable 
isomorphism of bundle gerbes with 
bundle gerbe connection and curving --- 
see \cite{MurSte}.  One says that 
two bundle gerbes $(L_{1},Y_{1},M)$ 
and $(L_{2},Y_{2},M)$ with bundle 
gerbe connections $\nabla_{1}$ and 
$\nabla_{2}$ and curvings $f_{1}$ and 
$f_{2}$ respectively are stably 
isomorphic if there exist line bundles 
$S$ and $T$ on $Y_{1}\times_{M}Y_{2}$ 
with connections $\nabla_{S}$ and 
$\nabla_{T}$ respectively and 
there is a line bundle isomorphism 
$(p_{1}^{[2]})^{-1}L_{1}\otimes \d(S)
\to (p_{2}^{[2]})^{-1}L_{2}\otimes 
\d(T)$ which commutes with the 
bundle gerbe products and preserves 
the connections $(p_{1}^{[2]})^{-1}
\nabla_{1}+\d(\nabla_{S})$ and 
$(p_{2}^{[2]})^{-1}\nabla_{2}+\d(
\nabla_{T})$ and such that the 
curvings $p_{1}^{*}f_{1}+F_{\nabla_{S}}$ 
and $p_{2}^{*}f_{2}+F_{\nabla_{T}}$ are 
equal.  One can then prove (\cite{MurSte}) 
a version of the above proposition 
for stable isomorphism classes of 
bundle gerbes with bundle gerbe connections 
and curvings.  One can also show 
(\cite{MurSte}) that there is a 
bijection between stable isomorphism 
classes of bundle gerbes with bundle 
gerbe connections and curvings and 
the Deligne hypercohomology group 
$H^{2}(M;\underline{\C}^{\times}_{M}\to 
\underline{\Omega}^{1}_{M}\to 
\underline{\Omega}^{2}_{M})$.  
\end{note} 

Given a bundle gerbe $L$ with bundle 
gerbe connection $\nabla_{L}$ and 
curving $f$ with trivial Deligne 
class $D(L,\nabla_{L},f)$ we will explicitly 
construct a line bundle $J\to Y$ which 
trivialises $L$ and a connection $\nabla_{J}$ 
on $J$ such that $\nabla_{L} = \d(\nabla_{J})$ 
under the isomorphism $L\to \d(J)$ and 
such that the curvature $F_{\nabla_{J}}$ of 
$\nabla_{J}$ is equal to $f$.   
Since $L$ 
has trivial class in $H^{2}(M;\underline{\C}^{\times}\to 
\underline{\Omega}^{1}_{M}\to \underline{\Omega}^{2}_{M})$ 
then there exists $h_{\a\b}:U_{\a\b} \to \cstar$ 
and $k_{\a} \in \Omega^{1}(U_{\a})$ such that 
\begin{eqnarray}
g_{\a\b\c} & = & h_{\b\c} h_{\a\c}^{-1} h_{\a\b}    \label{eq:six}    \\
A_{\a\b} & = & -k_{\b}+k_{\a}+\frac{dh_{\a\b}}{h_{\a\b}} \label{eq:seven} \\
f_{\a} & = & dk_{\a}    \label{eq:eight}                            
\end{eqnarray}
where we have continued to use 
the notation from above.  
First we construct the bundle $J \to Y$ 
which trivialises $L$.  This is done in 
\cite{Mur} but we need to 
review the construction.  
We define line 
bundles $J_{\a} \to Y_{\a}$ by 
$J_{\a} = \hat s_{\a}^{-1}L$ where 
$\hat{s}_{\a}: Y_{\a} \to Y^{[2]}$ is the 
map $y \mapsto (y,s_{\a}(\pi(y)))$.  
$\nabla_{L}$ induces a connection 
$\nabla_{\a} = \hat{s}_{\a}^{-1}\nabla$  
on $J_{\a}$.  
We use the sections 
$\s_{\a\b} \in \Gamma(U_{\a\b},L_{\a\b})$
to define line bundle isomorphisms 
$\phi_{\a\b}:J_{\a} \to J_{\b}$ by 
$\phi_{\a\b}(u) = m(\s_{\a\b}h_{\a\b}^{-1}\otimes u)$ 
for $u \in J_{\a}$.  The condition~\ref{eq:six} 
above shows that 
$\phi_{\b\c}\circ \phi_{\a\b} = \phi_{\a\c}$ 
on $Y_{\a\b\c}$.  Hence we can 
glue the $J_{\a}$ together using 
the $\phi_{\a\b}$ to get a line 
bundle $J$ on $Y$.  
It is easy to see that if the 
connections $\nabla_{\a}$ on 
the $J_{\a}$ satisfied 
$\nabla_{\a} = \phi_{\a\b}^{-1}\circ \nabla_{\b}\circ \phi_{\a\b}$ 
then we could glue them 
together to form a new 
connection $\nabla_{J}$ on 
$J$.  
We calculate 
$\phi_{\a\b}^{-1}\circ \nabla_{\b}\circ \phi_{\a\b}$ 
for an arbitrary section 
$t_{\a} \in \Gamma (Y_{\a},J_{\a})$ 
restricted to $Y_{\a\b}$.  We have 
\begin{eqnarray*} 
&  & \phi_{\a\b}^{-1}\circ \nabla_{\b} \circ \phi_{\a\b} (t_{\a}) \\ 
& = & \phi_{\a\b}^{-1}\circ \nabla_{\b}\circ 
m(\s_{\a\b}h_{\a\b}^{-1}\otimes t_{\a})                         \\
& = & \phi_{\a\b}^{-1}\circ (
m(\pi^{-1}\nabla_{\a\b}(\s_{\a\b}h_{\a\b}^{-1}) 
\otimes t_{\a}) + m(\s_{\a\b}h_{\a\b}^{-1}
\otimes \nabla_{\a}(t_{\a})))                                    \\
& = & \phi_{\a\b}^{-1}\circ ( 
m(((\pi^{*}A_{\a\b} - \pi^{*}\frac{dh_{\a\b}}{h_{\a\b}}) 
\s_{\a\b}h_{\a\b}^{-1})\otimes t_{\a}) + 
m(\s_{\a\b}h_{\a\b}^{-1}\otimes \nabla_{\a}(t_{\a}))            \\
& = & \pi^{*}A_{\a\b} - \pi^{*}\frac{dh_{\a\b}}{h_{\a\b}} 
+ \nabla_{\a}(t_{\a})                                           \\
& = & -\pi^{*}k_{\b} + \pi^{*}k_{\a} + \nabla_{\a}(t_{\a}).      
\end{eqnarray*} 
Hence if we let $\nabla_{\a}^{'} = \nabla_{\a} + \pi^{*}k_{\a}$ 
then we get 
$\nabla_{\a}^{'} = \phi_{\a\b}^{-1}\circ \nabla_{\b}^{'}\circ \phi_{\a\b}$ 
and so the $\nabla_{\a}^{'}$ 
patch together to define a connection 
$\nabla_{J}$ on $J$.  
Finally we need to show that 
$\nabla_{L} = \d(\nabla_{J})$ and that 
$f = F_{\nabla_{J}}$.  
The map 
$m:\pi_{1}^{-1}L\otimes \pi_{3}^{-1}L \to \pi_{2}^{-1}L$ 
gives rise to isomorphisms 
of bundle gerbes 
$\eta_{\a}:\d(J_{\a}) \to L$ over 
$Y_{\a}^{[2]}$ and furthermore 
we have a commutative diagram 
$$
\xymatrix{ 
\d(J_{\a}) \ar[dd]_{\d(\phi_{\a\b})} 
\ar[dr]^{\eta_{\a}}          &         \\ 
&                L                      \\ 
\d(J_{\b}) \ar[ur]_{\eta_{\b}}      } 
$$
over $Y_{\a\b}^{[2]}$.  Thus we 
can define an isomorphism of 
bundle gerbes 
$\eta : \d(J) \to L$ over 
$Y^{[2]}$.      
We have $\d(\nabla_{\a} + \pi^{*}k_{\a}) = \d(\nabla_{\a})$ 
which equals $\eta_{\a}^{-1}\circ \nabla_{L} \circ \eta_{\a}$ 
and thus we get $\nabla_{J} = \eta^{-1}\circ \nabla_{L} \circ \eta$.   
Finally $\nabla_{\a} + \pi^{*}k_{\a}$ has 
curvature $F_{\a} = f|_{Y_{\a}} +\pi^{*}(dk_{\a} - f_{\a}$ 
which equals $f|_{Y_{\a}}$ by 
equation~\ref{eq:eight}.  Thus $J$ 
has curvature $F_{\nabla_{J}}$ equal to $f$.

%
%
%
%
%

\setcounter{chapter}{5} 
\chapter{Gerbes and bundle gerbes} 
\label{chapter:6}  
 
\section{Torsors} 
\label{sec:6.1} 
In this section we briefly review the 
definition and basic properties of 
\emph{torsors}.  

\begin{definition}[\cite{Bry}] 
\label{def:6.1.1} 
A $\underline{G}$-\emph{torsor}  
on a space $X$ (here $\underline{G}$ is a sheaf 
of groups on $X$) is a 
sheaf $\underline{P}$ on $X$ plus an action 
of $\underline{G}$ on the right of $\underline{P}$, 
ie a morphism of sheaves $\underline{P}\times 
\underline{G}\to \underline{P}$, with the property 
that for any point of $X$ there is an open 
neighbourhood $U$ of that point such that for 
any open set $V\subset U$ the set $\underline{P}(V)$ 
is a right principal homogenous space under 
the group $\underline{G}(V)$.  
\end{definition} 
   
Similarly one can define the notion of a morphism 
$\underline{P}\to \underline{Q}$ between two 
$\underline{G}$-torsors on $X$.  One can show that 
any such morphism is automatically an isomorphism.  
A good example of a torsor is the torsor 
$\underline{P}_{X}$ associated to a principal $G$ bundle $P$ on 
a space $X$.  Here $\underline{P}_{X}$ is the 
sheaf of sections of $P$.  It is then easy to 
see that the sheaf of groups $\underline{G}_{X}$ 
acts on $\underline{P}_{X}$ on the right so as 
to make $\underline{P}_{X}$ into a 
$\underline{G}_{X}$-torsor.   

Suppose we have a morphism of sheaves of 
groups $\underline{G}\to \underline{H}$ 
over a space $X$ and a $\underline{G}$-torsor 
$\underline{P}$ on $X$.  Then one can construct 
(see \cite{Bry}[page 191]) an $\underline{H}$-torsor 
on $X$ denoted by $\underline{P}\times^{\underline{G}}
\underline{H}$.  $\underline{P}\times^{\underline{G}}
\underline{H}$ is the sheaf associated to the 
presheaf $U\mapsto (\underline{P}(U)\times 
\underline{H}(U))/\underline{G}(U)$, where 
$\underline{G}(U)$ acts on $\underline{H}(U)$ 
via the morphism of sheaves $\underline{G}\to 
\underline{H}$.  So a section of this presheaf 
would be an equivalence class $[p,h]$, where 
$p$ is a section of $\underline{P}(U)$ and 
$h$ is a section of $\underline{H}(U)$, under the 
equivalence relation $(p_{1},h_{1})
\equivalencerelation (p_{2},h_{2})$ if and 
only if there is a section $g_{12}$ of 
$\underline{G}(U)$ such that $p_{2} = p_{1}g_{12}$ 
and $h_{2} = g_{12}^{-1}\cdot h_{1}$.  Clearly 
$\underline{H}$ will act on the right of $\underline{P}
\times^{\underline{G}}\underline{H}$ so as to 
make $\underline{P}\times^{\underline{G}}\underline{H}$ 
into a $\underline{H}$-torsor.  From now on 
we will only be interested in torsors under 
the sheaf of groups $\underline{A}_{X}$ of smooth 
$A$ valued functions on 
$X$ associated to an abelian group $A$.    

Given another space $Y$ and a continuous 
mapping $f:Y\to X$ we have the notion of the 
pullback of a $\underline{A}_{X}$-torsor $\underline{P}$ 
on $X$ (see \cite{Bry}[pages 190--191]).  
We can pullback the sheaf $\underline{P}$ 
to a new sheaf $f^{-1}\underline{P}$ on 
$Y$.  This is now a torsor under the pullback 
sheaf of groups $f^{-1}\underline{A}_{Y}$ on 
$Y$.  There is a canonical morphism of sheaves 
of groups $f^{-1}\underline{A}_{X}\to 
\underline{A}_{Y}$ over $Y$.  Hence we can construct,  
using the method outlined above, an $\underline{A}_{Y}$-torsor 
$f^{-1}\underline{P}\times^{f^{-1}\underline{A}_{X}}
\underline{A}_{Y}$ on $Y$.  We will denote this 
$\underline{A}_{Y}$-torsor by $f^{*}\underline{P}$.               
Given another space $Z$ and a continuous mapping 
$g:Z\to Y$, we can form two $\underline{A}_{Z}$-torsors 
on $Z$ from an $\underline{A}_{X}$-torsor 
$\underline{P}$ on $X$; namely $(fg)^{*}\underline{P}$ 
and $g^{*}f^{*}\underline{P}$.  As mentioned in 
\cite{Bry}[page 187] these two torsors are 
not identical but there is an isomorphism 
$g^{*}f^{*}\underline{P}\to (fg)^{*}\underline{P}$ 
between them which satisfies a commutativity 
condition similar to that of diagram (5-4) on 
page 187 of \cite{Bry}.  

Given two $\underline{A}_{X}$ torsors $\underline{P}_{1}$ 
and $\underline{P}_{2}$ on $X$ we have the 
notion of the \emph{contracted product} 
$\underline{A}_{X}$-torsor $\underline{P}_{1}
\otimes \underline{P}_{2}$ on $X$ (see 
\cite{Bry}[page 191]).  To construct 
$\underline{P}_{1}\otimes \underline{P}_{2}$ we first 
take the $\underline{A}_{X}\times \underline{A}_{X}$-torsor 
$\underline{P}_{1}\times\underline{P}_{2}$ on 
$X$ and then using the canonical morphism of 
sheaves $\underline{A}_{X}\times \underline{A}_{X}
\to \underline{A}_{X}$ induced by the product map 
$A\times A\to A$, we form the $\underline{A}_{X}$-
torsor $(\underline{P}_{1}\times\underline{P}_{2})
\times^{\underline{A}_{X}\times \underline{A}_{X}}
\underline{A}_{X} := \underline{P}_{1}\otimes 
\underline{P}_{2}$.   
   
\section{Sheaves of categories and gerbes} 
\label{sec:6.2}
\begin{definition}[\cite{Bre},\cite{Bry}]  
\label{def:6.2.1} 
A \emph{presheaf 
of categories} $\mathcal{C}$ on $X$ 
consists of the following data.  
\begin{enumerate} 

\item for each open set $U\subset X$ a 
category $\mathcal{C}(U)$.  

\item given open sets $V\subset U\subset X$  
there is 
a functor $\rho_{V,U}:\mathcal{C}(U)\to \mathcal{C}(V)$ 
--- the `restriction' functor.  We will sometimes 
denote the image $\rho_{V,U}(C)$ of an 
object $C$ of $\mathcal{C}(U)$ by $C|_{V}$ and the 
image $\rho_{V,U}(\phi)$ of an arrow $\phi$ 
in $\mathcal{C}(U)$ by $\phi|_{V}$.  

\item given open sets $W\subset V\subset U\subset X$ 
we have 
$\rho_{W,V}\circ \rho_{V,U} = \rho_{W,U}$.  

\end{enumerate} 
\end{definition} 

Strictly speaking we should require that 
rather than the functors $\rho_{W,V}\circ \rho_{V,U}$ 
and $\rho_{W,U}$ coinciding, they should differ 
by an invertible natural transformation 
which satisfies a certain coherency condition given a 
fourth open subset of $X$.  In most examples that we 
meet however we will not need this generality.    
A \emph{morphism} or \emph{functor} $\phi:\mathcal{C}\to \mathcal{D}$ 
of presheaves of categories $\mathcal{C}$ and 
$\mathcal{D}$ on $X$ consists of the assignment to every 
open set $U\subset X$ of a functor $\phi_{U}:
\mathcal{C}_{U}\to \mathcal{D}_{U}$ in such a 
way that if $V$ is another open set of $X$ 
with $V\subset U\subset X$ then the following 
diagram commutes: 
$$
\xymatrix{ 
\mathcal{C}_{U} \ar[r]^-{\phi_{U}} \ar[d] & 
\mathcal{D}_{U} \ar[d]                  \\ 
\mathcal{C}_{V} \ar[r]^-{\phi_{V}} & 
\mathcal{D}_{V}                          } 
$$
where the two vertical functors are the 
restriction functors in the 
presheaves of categories $\mathcal{C}$ 
and $\mathcal{D}$.  A \emph{natural transformation} 
$\tau:\phi \Rightarrow \psi$ between morphisms 
of presheaves of categories $\phi,\psi:\mathcal{C}\to 
\mathcal{D}$ consists of a natural 
transformation $\tau_{U}:\phi_{U}\Rightarrow 
\psi_{U}$ for every open set $U\subset X$ 
which is compatible with restrictions to 
smaller open sets.  
  
\begin{definition} 
\label{def:6.2.2} 
A \emph{stack} or \emph{sheaf of categories} 
on $X$ is a presheaf of categories $\mathcal{C}$ 
on $X$ which satisfies the following `gluing' 
laws or `effective descent' conditions on objects 
and arrows.  The effective descent conditions on 
objects are as follows.  Given an open set 
$U\subset X$ and an open cover $\{U_{i}\}_{i\in I}$ 
of $U$ together with objects 
$C_{i} \in Ob \mathcal{C}(U_{i})$ and arrows 
$\phi_{ij}:\rho_{U_{ij},U_{i}}(C_{i})\to \rho_{
U_{ij},U_{j}}(C_{j})$ satisfying 
$\rho_{U_{ijk},U_{jk}}(\phi_{jk})\rho_{U_{ijk},
U_{ij}}(\phi_{ij}) = \rho_{U_{ijk},U_{ik}}(\phi_{ik})$ 
then there is an object $C$ of $\mathcal{C}(U)$ 
and isomorphisms $\psi_{i}:\rho_{U_{i},U}(C)\to C_{i}$ 
compatible with the $\phi_{ij}$.  The effective 
descent conditions for arrows are that given an 
open set $U\subset X$, objects $C_{1}$ and $C_{2}$ 
of $\mathcal{C}(U)$ and an open cover 
$\{U_{i}\}_{i\in I}$ of $U$ together with 
arrows $\phi_{i}:\rho_{U_{i},U}(C_{1})\to
\rho_{U_{i},U}(C_{2})$ satisfying 
$\rho_{U_{ij},U_{i}}(\phi_{i}) = \rho_{U_{ij},U_{j}}
(\phi_{j})$ then there is a unique arrow 
$\phi:C_{1}\to C_{2}$ in $\mathcal{C}(U)$ such that 
$\rho_{U_{i},U}(\phi) = \phi_{i}$ for all 
$i \in I$.  In other words the set of arrows 
from $C_{1}$ to $C_{2}$ in $\mathcal{C}(U)$ forms a 
sheaf $\underline{Hom}(C_{1},C_{2})$.    
\end{definition} 

A \emph{morphism} or \emph{functor} of sheaves of categories is just 
a morphism of the underlying presheaves 
of categories.  A \emph{natural transformation} 
between morphisms of sheaves of categories 
is defined in a similar way.  An equivalence of 
stacks $\phi:\mathcal{C}\to \mathcal{D}$ is 
a morphism of stacks $\phi$ such that there is a 
morphism of stacks $\psi:\mathcal{D}\to 
\mathcal{C}$ together with natural transformations 
$\psi\circ \phi\Rightarrow \text{id}_{\mathcal{C}}$ 
and $\phi\circ \psi\Rightarrow \text{id}_{\mathcal{D}}$.  
If such an equivalence exists then the sheaves of categories 
$\mathcal{C}$ and $\mathcal{D}$ are said to be \emph{equivalent}.  
  
The objects and arrows in an arbitrary 
presheaf of categories may not satisfy the 
effective descent conditions listed above 
however, in analogy with the case for presheaves and 
sheaves of abelian groups, there is a 
functor, the `sheafification' functor or 
`associated sheaf' functor $\mathcal{C}\to 
\tilde{\mathcal{C}}$ which transforms a presheaf of 
categories $\mathcal{C}$ on $X$ into a sheaf of categories 
$\tilde{\mathcal{C}}$.  This is subject to the 
usual universal condition, namely that given a 
sheaf of categories $\mathcal{D}$ together with 
a morphism of presheaves of categories 
$\mathcal{C}\to \mathcal{D}$, there is a unique 
morphism of sheaves of categories $\tilde{\mathcal{C}}\to 
\mathcal{D}$ making the following diagram commute: 
$$
\xymatrix{ 
\mathcal{C} \ar[r] \ar[dr] & 
\mathcal{\tilde{C}} \ar @{.>}[d]    \\ 
& \mathcal{D}.      } 
$$ 
We will briefly discuss the 
construction of $\tilde{C}$ and refer to 
\cite{Bry}[pages 192--194] for a more complete 
discussion.  Given an open set $U\subset X$ 
the category $\tilde{\mathcal{C}}_{U}$ has 
objects consisting of triples $(\{U_{i}\}_
{i\in I},\{P_{i}\}_{i\in I},\{f_{ij}\}_{i,
j\in I})$ where $\{U_{i}\}_{i\in I}$ is an 
open cover of $U$, $P_{i}$ is an object of 
$\mathcal{C}_{U_{i}}$ for each $i\in I$, and 
$f_{ij}:P_{i}|_{U_{ij}}\to P_{j}|_{U_{ij}}$ 
are arrows in $\mathcal{C}_{U_{ij}}$ which 
satisfy the descent condition $f_{jk}|_{U_{ijk}}
f_{ij}|_{U_{ijk}} = f_{ik}|_{U_{ijk}}$ in 
$\mathcal{C}_{U_{ijk}}$.  Given such objects 
$(\{U_{i}\},\{P_{i}\},\{f_{ij}\})$ and 
$(\{V_{\a}\},\{Q_{\a}\},\{g_{\a\b}\})$ of 
$\tilde{\mathcal{C}}_{U}$, an arrow 
$(\{U_{i}\},\{P_{i}\},\{f_{ij}\})\to 
(\{V_{\a}\},\{Q_{\a}\},\{g_{\a\b}\})$ 
is an equivalence class $[\{W_{i,\a}\},\{h_{i,\a}\}]$ 
where $\{W_{i,\a}\}_{i\in I,\a\in \Sigma}$ 
is an open cover of $U$ such that 
$W_{i,\a}\subset U_{i}\cap V_{\a}$ for all 
$i\in I$, $\a\in \Sigma$ and $h_{i,\a}:
P_{i}|_{W_{i,\a}}\to Q_{\a}|_{W_{i,\a}}$ is 
an arrow in $\mathcal{C}_{W_{i,\a}}$ which is 
compatible with the descent data $f_{ij}$ and 
$g_{\a\b}$ in the sense that the following 
diagram of arrows in $\mathcal{C}_{W_{i,\a}
\cap W_{j,\b}}$ commutes: 
$$
\xymatrix{ 
P_{i}|_{W_{i,\a}\cap W_{j,\b}} 
\ar[rr]^-{h_{i,\a}|_{W_{i,\a}\cap W_{j,\b}}} 
\ar[d]_-{f_{ij}|_{W_{i,\a}\cap W_{j,\b}}} & & 
Q_{\a}|_{W_{i,\a}\cap W_{j,\b}} \ar[d]^-{
g_{\a\b}|_{W_{i,\a}\cap W_{j,\b}}}            \\ 
P_{j}|_{W_{i,\a}\cap W_{j,\b}} \ar[rr]_-{
h_{j,\b}|_{W_{i,\a}\cap W_{j,\b}}} & & 
Q_{\b}|_{W_{i,\a}\cap W_{j,\b}},                  } 
$$
and where two such pairs $(\{W_{i,\a}\},\{
h_{i,\a}\})$ and $(\{W_{i,\a}^{'}\},\{
h_{i,\a}^{'}\})$ are declared to be equivalent if 
there is an open cover $\{W_{i,\a}^{''}\}_{i\in I,
\a\in \Sigma}$ of $U$ such that $W_{i,\a}^{''}
\subset W_{i,\a}\cap W_{i,\a}^{'}$ for all $i\in I$ 
and all $\a\in \Sigma$ and such that $h_{i,\a}|_
{W_{i,\a}^{''}} = h_{i,\a}^{'}|_{W_{i,\a}^{''}}$ 
for all $i\in I$ and all $\a\in \Sigma$.  If the 
presheaf of categories $\mathcal{C}$ has the 
property that given any two objects $P$ and $Q$ 
of $\mathcal{C}_{U}$ for some open set $U\subset X$, 
the set of arrows $Hom(P,Q)$ is a sheaf then one 
can show that $\tilde{\mathcal{C}}$ is a sheaf of 
categories.  Otherwise it is not too hard to see 
that the presheaf of categories $\tilde{\mathcal{C}}$ 
has this property and hence $\Tilde{\Tilde{\mathcal{C}}}$ 
is a sheaf of categories (this is analogous to the 
case for the associated sheaf of a presheaf, in 
general one has to apply the sheafification procedure 
twice to yield a sheaf --- see \cite{MaMo} for more 
details on this).  It can also be shown that there is 
a morphism of presheaves of categories $\mathcal{C}\to 
\tilde{\mathcal{C}}$ with the universal property 
referred to above.      

\begin{definition}[\cite{Bre},\cite{Bry}]  
\label{def:6.2.3} 
A \emph{gerbe} $\mathcal{G}$ on $X$ is a sheaf 
of categories $\mathcal{G}$ on $X$ such that the following 
three conditions are satisfied: 
\begin{enumerate} 
\item $\mathcal{G}$ is a sheaf of groupoids on 
$X$, so for each open set $U\subset X$  
the category $\mathcal{G}(U)$ 
is a groupoid.  
\item $\mathcal{G}$ is locally non-empty.  This 
means that there is an open cover $\{U_{i}\}_{i\in I}$ 
of $X$  
so that the groupoid 
$\mathcal{G}(U_{i})$ is non-empty for each $i \in I$.  
\item $\mathcal{G}$ is locally connected.  This 
means that for every open set $U\subset X$  
and all pairs of objects $P_{1}$ and 
$P_{2}$ in $\mathcal{G}(U)$, there 
is an open covering $\{U_{i}\}_{i\in I}$ 
of $U$  
such that $Hom(\rho_{U_{i},U}(P_{1}),\rho_{U_{i},U}(P_{2}))$ is 
non-empty.  
\end{enumerate} 
\end{definition} 

\begin{example} 
\label{ex:6.2.4} 
Let $\underline{G}$ be a sheaf of groups on $X$.  
For each open subset $U\subset X$, let 
$Tors(\underline{G})_{U}$ be the groupoid 
with objects the $\underline{G}$ torsors 
$\underline{P}$ on $U$ and morphisms the 
isomorphisms between $\underline{G}$ torsors 
on $U$.  Given another open set  
$V\subset U\subset X$ there is a natural 
pullback functor 
$\rho_{V,U}:\mathcal{C}(U)\to 
\mathcal{C}(V)$ induced by the restriction 
functors of the $\underline{G}$ torsors.    
By definition of a sheaf these 
satisfy $\rho_{W,V}\rho_{V,U} = 
\rho_{W,U}$ given a third open set $W\subset 
V\subset U\subset X$.  
It is straightforward to show that the 
assignment $U\mapsto Tors(\underline{G})_{U}$ 
defines a sheaf of groupoids on $X$.   
$Tors(\underline{G})$ is locally non-empty 
since to any open covering $\{U_{i}\}_{i\in I}$ 
of $X$ there is an object of each groupoid 
$Tors(\underline{G})_{U_{i}}$ --- namely the 
trivial $\underline{G}$ torsor on $U_{i}$.  
$Tors(\underline{G})$ is also locally connected: 
given an open set $U\subset X$ and $\underline{G}$ 
torsors $\underline{P}_{1}$ and 
$\underline{P}_{2}$ then locally $\underline{P}_{1}$ 
and $\underline{P}_{2}$ are isomorphic to the 
trivial $\underline{G}$ torsor and hence to each 
other.  Therefore $Tors(\underline{G})$ is a 
gerbe on $X$.    
\end{example}          

\begin{note} 
In several examples that we will encounter we will be 
given a presheaf of groupoids $\mathcal{G}$ on 
$X$ which is locally non-empty and locally 
connected.  One can show that the associated sheaf 
of groupoids 
$\tilde{\mathcal{G}}$ on $X$ will be locally 
non-empty and locally connected as well and hence 
will define a gerbe on $X$.   
\end{note}  

\begin{definition} 
\label{def:6.2.5} 
We say that a gerbe $\mathcal{G}$ on $X$ is 
\emph{neutral} if there is a global object 
of $\mathcal{G}$; ie if $\mathcal{G}_{X}$ is 
non-empty.  
\end{definition} 
 
\begin{example} 
\label{ex:6.2.6}
The gerbe $Tors(\underline{G})$ is neutral 
since it has a global object --- namely the 
trivial $\underline{G}$ torsor on $X$.  
\end{example} 

\begin{definition} 
\label{def:6.2.7} 
Let $\mathcal{G}$ be a gerbe on $X$ and let 
$\underline{A}$ be a sheaf of abelian 
groups on $X$.  We say that $\mathcal{G}$ 
is \emph{bound} by $\underline{A}$, or 
has \emph{band} equal to $\underline{A}$ 
if for any object $P$ of $\mathcal{G}(U)$, 
where $U\subset X$ is an open set, there 
is an isomorphism of sheaves of groups 
$$
\underline{Aut}(P) \stackrel{\eta_{P}}{\to} 
\underline{A} 
$$
over $U$ which is compatible with 
restriction to smaller open sets.  
If $P_{1}$ and $P_{2}$ are two objects 
of $\mathcal{G}(U)$ and $f:P_{1}\to P_{2}$ 
is an isomorphism between them then we also require 
that $\eta_{P_{2}}\underline{Aut}(f) = \eta_{P_{1}}$.  
Here $\underline{Aut}(f):\underline{Aut}(P_{1})
\to \underline{Aut}(P_{2})$ 
is the isomorphism of sheaves of groups 
induced by $f$.  
\end{definition} 

We will only be interested in gerbes with 
band equal to $\underline{\C}^{\times}_{X}$. 
A morphism of gerbes $\mathcal{G}$ to 
$\mathcal{H}$ is a morphism $\phi:\mathcal{G}
\to \mathcal{H}$ of the underlying sheaves of 
categories which satisfies the extra condition that 
given an open set $U\subset X$ and an object $P$ of 
$\mathcal{G}_{U}$, then the following diagram 
commutes: 
$$
\xymatrix{ 
\underline{Aut}_{\mathcal{G}_{U}}(P) 
\ar[r]^{\phi_{U}} \ar[dr]_{\eta_{P}} & 
\underline{Aut}_{\mathcal{H}_{U}}(\phi_{U}(P)) 
\ar[d]^{\eta_{\phi_{U}(P)}}                   \\ 
& \underline{\C}^{\times}_{U}.                 } 
$$
A natural transformation between morphisms of gerbes 
is defined similarly (notice that any natural 
transformation between morphisms of gerbes must 
be invertible).  There is also the notion of an 
\emph{equivalence} of gerbes and hence the notion 
of \emph{equivalent} gerbes.  Equivalence of gerbes 
is clearly an equivalence relation.  

Recall from Definition~\ref{def:6.2.2} above 
that the set of arrows between two objects 
$P_{1}$ and $P_{2}$ in a sheaf of categories is a 
sheaf $\underline{Hom}(P_{1},P_{2})$.  If 
$P_{1}$ and $P_{2}$ are objects of a gerbe then 
more is true.  The sheaf of groups $\underline{Aut}
(P_{2})$ acts on $\underline{Hom}(P_{1},P_{2})$ 
on the right and in fact makes $\underline{Hom}(
P_{1},P_{2})$ into a $\underline{\C}^{\times}_{X}$ 
torsor.  We have the following Proposition from 
\cite{Bry}. 
 
\begin{proposition}[\cite{Bry}]  
\label{prop:6.2.8} 
Let $\mathcal{G}$ be a gerbe on $X$ with band 
equal to $\underline{A}$ for some sheaf of 
commutative groups $\underline{A}$ on $X$.  
Then we have 
\begin{enumerate} 
\item Let $U\subset X$ be an open subset 
and let $P_{1}$ and $P_{2}$ be two objects 
of $\mathcal{G}_{U}$.  Then the sheaf of 
arrows $\underline{Hom}(P_{1},P_{2})$ is 
an $\underline{A}$-torsor.  
\item Let $P$ be an object of $\mathcal{G}_{U}$ 
for some open set $U\subset X$ and let $\underline{I}$ 
be an $\underline{A}$-torsor on $U$.  Then there 
is an object $Q$ of $\mathcal{G}_{U}$ such 
that there is an isomorphism $\underline{Hom}(
P,Q)\simeq \underline{I}$ of $\underline{A}$-torsors 
on $U$.  The object $Q$ is unique up to 
isomorphism and is denoted by $P\times^{\underline{A}}
\underline{I}$.  
\end{enumerate} 
\end{proposition} 
 
\section{Properties of gerbes} 
\label{sec:6.3}

We will briefly review some ways of 
constructing new gerbes from existing gerbes 
from \cite{Bry}[pages 198--200].  First of all 
suppose we have a gerbe $\mathcal{G}$ on $X$ with 
band equal to the sheaf of abelian groups $\underline{A}$ 
on $X$.  Suppose that there is a morphism of sheaves 
of abelian groups $\underline{A}\to \underline{B}$ on 
$X$.  Then we can construct a new gerbe on $X$, 
bound by $\underline{B}$ and denoted $\mathcal{G}\times^{
\underline{A}}\underline{B}$ as follows.  First we 
construct a presheaf of groupoids on $X$ by letting 
$(\mathcal{G}\times^{\underline{A}}\underline{B})_{U}$ 
have objects the objects of $\mathcal{G}_{U}$ but 
given two objects $P_{1}$ and $P_{2}$ of 
$(\mathcal{G}\times^{\underline{A}}\underline{B})_{U}$ 
we let $Hom_({\mathcal{G}\times^{\underline{A}}
\underline{B})_{U}}(P_{1},P_{2}) = Hom_{\mathcal{G}_{U}}
(P_{1},P_{2})\times^{\underline{A}}\underline{B}$.  
This presheaf of groupoids is clearly locally non-empty 
and locally connected.  We define $\mathcal{G}\times^{
\underline{A}}\underline{B}$ to be the sheaf of 
groupoids associated to this presheaf.     

Given two gerbes $\mathcal{G}_{1}$ and $\mathcal{G}_{2}$ 
on $X$ bound by $\underline{\C}^{\times}_{X}$ we have 
the notion of the \emph{contracted product} 
gerbe $\mathcal{G}_{1}\otimes \mathcal{G}_{2}$ on $X$.  
This is also a gerbe bound by $\underline{\C}^{\times}_{X}$.  
It is constructed as follows.  We first form the gerbe 
$\mathcal{G}_{1}\times \mathcal{G}_{2}$ on $X$, ie 
the sheaf of groupoids $U\mapsto \mathcal{G}_{1}(U)\times 
\mathcal{G}_{2}(U)$.  This is gerbe bound by $\underline{\C}
^{\times}_{X}\times \underline{\C}^{\times}_{X}$.  
Ordinary multiplication in $\cstar$ induces a morphism of 
sheaves of abelian groups $\underline{\C}^{\times}_{X}
\times \underline{\C}^{\times}_{X}\to \underline{\C}^{\times}_{X}$.  
We define $\mathcal{G}_{1}\otimes \mathcal{G}_{2} = 
(\mathcal{G}_{1}\times \mathcal{G}_{2})\times^{\underline{\C}
^{\times}_{X}\times \underline{\C}^{\times}_{X}}
\underline{\C}^{\times}_{X}$ --- see \cite{Bry}[page 200] 
for a direct description of this gerbe.  
 
Suppose $\mathcal{G}$ is a gerbe on $X$.  Then we have 
the notion of the \emph{dual} gerbe $\mathcal{G}^{*}$ 
on $X$ associated to $\mathcal{G}$.  $\mathcal{G}^{*}$ is 
defined as follows.  We let $\mathcal{G}^{*}$ have fibre 
$\mathcal{G}^{*}(U)$ at open set $U\subset X$ 
equal to the opposite groupoid $\mathcal{G}(U)^{\circ}$ 
of $\mathcal{G}(U)$.  Given open sets 
$V\subset U\subset X$, $\mathcal{G}^{*}$ inherits 
restriction functors $\mathcal{G}^{*}(U)\to 
\mathcal{G}^{*}(V)$ from the $\rho_{U,V}:\mathcal{G}
(U)\to \mathcal{G}(V)$.  These restriction functors 
satisfy the compatibility condition given a third 
open set $W\subset V\subset U\subset X$.  It is not 
hard to see that the assignment $U\mapsto \mathcal{G}
^{*}(U)$ defines a sheaf of categories and it is also easy to see 
that this stack is locally non-empty and locally 
connected and hence defines a gerbe.   

Given a smooth map $f:Y\to X$ between smooth manifolds 
$X$ and $Y$ and a gerbe $\mathcal{C}$ on $X$ 
bound by $\underline{\C}^{\times}_{X}$, we have 
the notion of the \emph{pullback gerbe} $f^{-1}\mathcal{C}$ 
on $Y$ --- see \cite{Bry}[page 198].  
This is a gerbe bound by the sheaf 
of groups $f^{-1}\underline{\C}^{\times}_{X}$ on $X$.  
$f^{-1}\mathcal{C}$ is obtained by 
sheafifying the presheaf of categories 
$U\mapsto f^{-1}\mathcal{C}_{U}$, where 
$f^{-1}\mathcal{C}_{U}$ is the groupoid whose 
objects consist of pairs $(V,P)$, where 
$V\subset X$ is an open set of $X$ such that 
$f(U)\subset V$ and $P$ is an object of 
$\mathcal{C}_{V}$.  An arrow $(V_{1},P_{1})
\to (V_{2},P_{2})$ is an equivalence class 
$[W_{12},\a_{12}]$, where $W_{12}\subset X$ is an open subset 
of $X$ such that $f(U)\subset W_{12}$ and $W_{12}\subset 
V_{1}\cap V_{2}$ and $\a_{12}:P_{1}|_{W_{12}}\to P_{2}
|_{W_{12}}$ is an arrow in $\mathcal{C}_{W_{12}}$, under the 
equivalence relation $(W_{12},\a_{12})\equivalencerelation 
(W_{12}^{'},\a_{12}^{'})$ if there exists an open set 
$W_{12}^{''}\subset X$ with $W_{12}^{''}\subset W_{12}\cap 
W_{12}^{'}$, $f(U)\subset W_{12}^{''}$ and $\a_{12}|_{W_{12}^{''}} = 
\a_{12}^{'}|_{W_{12}^{''}}$. If we have a third object 
$(V_{3},P_{3})$ and an arrow $[W_{23},\a_{23}]$ 
from $(V_{2},P_{2})$ to $(V_{3},P_{3})$, then the 
composite arrow $[W_{23},\a_{23}][W_{12},\a_{12}]:
(V_{1},P_{1})\to (V_{2},P_{2})$ is defined by 
$[W_{23},\a_{23}][W_{12},\a_{12}] = [W_{23}\cap W_{12},
\a_{23}|_{W_{23}\cap W_{12}}\a_{12}|_{W_{23}\cap W_{12}}]$.  
One can check that this respects the equivalence 
relation $\equivalencerelation$ defined above.  
One can check (see \cite{Bry}[Proposition 5.2.6]) 
that the presheaf of groupoids $f^{-1}\mathcal{C}$ 
is locally non-empty and locally connected and hence the same will 
be true of the sheaf of groupoids associated to 
$f^{-1}\mathcal{C}$.  

We will usually want $f^{-1}\mathcal{C}$ to 
be a gerbe bound by $\underline{\C}^{\times}_{Y}$.  
We can arrange this as follows (see \cite{Bry}[page 
203]).  Observe that we have a morphism of sheaves 
of groups $f^{-1}\underline{\C}^{\times}_{X}\to 
\underline{\C}^{\times}_{Y}$ and hence we can 
form a new gerbe bound by $\underline{\C}^{\times}
_{Y}$ in the manner described above.  We will 
use the notation of \cite{Bry} and denote this 
new gerbe on $Y$ bound by $\underline{\C}^{\times}_{Y}$ 
by $f^{*}\mathcal{C}$.   

We also have the following Lemma.  
\begin{lemma}[\cite{Bry1}]   
\label{lemma:6.3.2} 
Suppose we have smooth manifolds $X$, $Y$ and 
$Z$ and smooth maps $f:Y\to X$ and $g:Y\to Z$.  
Suppose we are also given a gerbe $\mathcal{C}$ 
on $X$.  Then there is an equivalence of gerbes 
$\phi_{f,g}:g^{-1}f^{-1}\mathcal{C}\stackrel{\simeq}{\to} 
(fg)^{-1}\mathcal{C}$ on $Z$ such that, given a 
fourth smooth manifold $W$ and a smooth map 
$h:W\to Z$ then there is a natural transformation 
between the equivalences of gerbes as pictured 
in the following diagram: 
$$
\xymatrix{ 
h^{-1}g^{-1}f^{-1}\mathcal{C} 
\ar[r]^-{h^{-1}\phi_{f,g}} 
\ar[d]_-{\phi_{g,h}} & h^{-1}(
fg)^{-1}\mathcal{C} \ar[d]^-{
\phi_{fg,h}} \ar @2{->}[dl]^-{
\theta_{f,g,h}}                 \\ 
(gh)^{-1}f^{-1}\mathcal{C} 
\ar[r]_-{\phi_{f,gh}} & (fgh)^{-1}
\mathcal{C}.                         } 
$$
The natural transformation $\theta_{f,g,h}$ 
satisfies the following coherency condition:  
given a fifth smooth manifold $V$ and a 
smooth map $k:V\to W$ we require the following 
diagram of natural transformations to commute. 
$$
\xymatrix@!C=40pt{ 
& k^{-1}h^{-1}g^{-1}f^{-1}\mathcal{C} \ar[dl] 
\ar[dr] \ar[rr] & & k^{-1}h^{-1}(fg)^{-1}
\mathcal{C} \ar[dr] \ar @2{->}[dl]^-{k^{-1}
\theta_{f,g,h}} &                               \\ 
(hk)^{-1}g^{-1}f^{-1}\mathcal{C} \ar[dr] & &  
k^{-1}(gh)^{-1}f^{-1}\mathcal{C} 
\ar @2{->}[ll]^-{\theta_{g,h,k}} \ar[dl] 
\ar[rr] & & k^{-1}(fgh)^{-1}\mathcal{C} \ar[dl] 
\ar @2{->}[dlll]^-{\theta_{f,gh,k}}                   \\ 
& (ghk)^{-1}f^{-1}\mathcal{C} \ar[rr] & 
\ar @2{-}[d] & (fghk)^{-1}\mathcal{C}  &            \\ 
& k^{-1}h^{-1}g^{-1}f^{-1}\mathcal{C} \ar[dl] 
\ar[rr] & & k^{-1}h^{-1}(fg)^{-1}\mathcal{C} \ar[dr] 
\ar[dl] \ar @2{-}[dlll]                 &                \\  
(hk)^{-1}g^{-1}f^{-1}\mathcal{C} \ar[dr] \ar[rr] & & 
(hk)^{-1}(fg)^{-1}\mathcal{C} \ar[dr] \ar @2{->}[dl]^
-{\theta_{f,g,hk}} & & k^{-1}(fgh)^{-1}\mathcal{C} 
\ar[dl] \ar @2{->}[ll]^-{\theta_{fg,h,k}}            \\ 
& (ghk)^{-1}f^{-1}\mathcal{C} \ar[rr] & & 
(fghk)^{-1}\mathcal{C} . &                              } 
$$
There is also a version of this Lemma with 
$f^{-1}$ replaced by $f^{*}$.  
\end{lemma} 
 
We now turn to discuss the characteristic class 
associated to gerbes, the \emph{Dixmier-Douady} 
class.  This is a class in the sheaf cohomology 
group $H^{2}(X;\underline{\C}^{\times}_{X})$ and hence 
we can regard it as a class in $H^{3}(X;\Z)$ 
via the long exact sequence in sheaf cohomology 
induced by the short exact sequence of sheaves 
of abelian groups 
$$
1\to \underline{\Z}_{X}\to \underline{\C}_{X}\to 
\underline{\C}^{\times}_{X}\to 1.  
$$
We briefly recall how to construct this class 
and refer to \cite{Bry} and \cite{Bre} for more 
details --- particularly \cite{Bre} where there 
is no assumption made about gerbes being bound by 
sheaves of abelian groups.  We also assume 
throughout the remainder of the chapter that 
all the gerbes we encounter are gerbes over a 
smooth paracompact manifold $X$.  Let $\mathcal{C}$ 
be such a gerbe and choose an open 
covering $\{U_{i}\}_{i\in I}$ of $X$ such that 
each non-empty intersection $U_{i_{0}}\cap \cdots 
\cap U_{i_{p}}$ is contractible.  Since 
$U_{i}$ is contractible we can choose objects 
$C_{i}$ of $\mathcal{C}(U_{i})$.  
Choose arrows 
$f_{ij}:\rho_{U_{ij},U_{i}}(C_{i})\to \rho_{U_{ij},
U_{j}}(C_{j})$ in $\mathcal{C}(U_{ij})$.  It 
is not immediately obvious that such  
arrows exist --- indeed the axioms for a gerbe 
only imply that given a choice of the objects 
$C_{i}$ there is an open cover $\{U_{ij}^{\a}\}_{
\a\in \Sigma_{ij}}$ of $U_{ij}$ such that there exist 
arrows $f_{ij}^{\a}:\rho_{U_{ij}^{\a},U_{i}}(C_{i})
\to \rho_{U_{ij}^{\a},U_{j}}(C_{j})$ for each 
$\a \in \Sigma_{ij}$.  This is the viewpoint taken in 
\cite{Bre}.  To see that in fact the arrows 
$f_{ij}$ exist in the present case we use  
Proposition~\ref{prop:6.2.8} to conclude that  
the sheaf $\underline{Isom}(\rho_{U_{ij},U_{i}}
(C_{i}),\rho_{U_{ij},U_{j}}(C_{j}))$ is a 
$\underline{\C}^{\times}_{U_{ij}}$ torsor.  Such 
objects are classified by the sheaf cohomology 
group $H^{1}(U_{ij},\underline{\C}^{\times}_{U_{ij}})$ 
which is zero, since $U_{ij}$ is assumed contractible.  
This means that there is an isomorphism 
$f_{ij}:\rho_{U_{ij},U_{i}}(C_{i})\to \rho_{U_{ij},U_{j}}
(C_{j})$.  Now, in $\mathcal{C}_{U_{ijk}}$, we 
form the automorphism $\rho_{U_{ijk},U_{ki}}(f_{ki})
\rho_{U_{ijk},U_{jk}}(f_{jk})\rho_{U_{ijk},U_{ij}}(
f_{ij})$ of $\rho_{U_{ijk},U_{i}}(C_{i})$ and 
identify it with a map $\epsilon_{ijk}:U_{ijk}\to 
\cstar$ under the isomorphism $\underline{Aut}(C_{i})\to 
\underline{\C}^{\times}_{U_{i}}$.  It is easy to see 
that in fact $\epsilon_{ijk}$ defines a $\cstar$ valued \v{C}ech 
two cocycle which is independent of all choices 
involved.  The class $\epsilon_{ijk}$ defines in 
$H^{2}(X;\underline{\C}^{\times}_{X})$ is called the 
Dixmier-Douady class of the gerbe $\mathcal{C}$,         
and denoted by $DD(\mathcal{C})$.  
We summarise the main properties of the 
Dixmier-Douady  
class in the following proposition.   

\begin{proposition}[\cite{Bry}] 
\label{prop:6.3.3} 
The Dixmier-Douady class has the following 
properties:
\begin{enumerate} 
\item The Dixmier-Douady class of the trivial 
gerbe $Tors(\underline{\C}^{\times}_{X})$ is zero.  
\item The contracted product $\mathcal{C}\otimes 
\mathcal{D}$ of two gerbes $\mathcal{C}$ and 
$\mathcal{D}$ on $X$ has $DD(\mathcal{C}\otimes 
\mathcal{D}) = DD(\mathcal{C}) + DD(\mathcal{D})$.  
\item The dual $\mathcal{C}^{*}$ of a gerbe 
$\mathcal{C}$ has $DD(\mathcal{C}^{*}) = - DD(
\mathcal{C})$.  
\item If $f:Y\to X$ is a smooth map and 
$\mathcal{C}$ is a gerbe on $X$ then $f^{*}\mathcal{C}$ 
has Dixmier-Douady class $DD(f^{*}\mathcal{C}) = 
f^{*}DD(\mathcal{C})$.  
\end{enumerate}  

Equivalence between gerbes defines an equivalence       
relation on the class of all gerbes on $X$ and 
the operations of contracted product and duals 
make the set of equivalence classes of gerbes on 
$X$ into an abelian group with identity element 
the trivial gerbe of torsors on $X$.  The Dixmier-
Douady class supplies an homomorphism from the 
group of all equivalence classes of gerbes to 
$H^{2}(X;\underline{\C}^{\times}_{X})$.  This homomorphism 
is an isomorphism.   
\end{proposition}  

It is worth noting that given gerbes $\mathcal{C}$ 
and $\mathcal{D}$ on $X$ with a morphism 
$\phi:\mathcal{C}\to \mathcal{D}$ between them 
which is not necessarily an equivalence, then 
$\mathcal{C}$ and $\mathcal{D}$ have the same 
Dixmier-Douady class.  Indeed, given an open cover 
$\{U_{i}\}_{i\in I}$ of $X$ and objects $C_{i}$ 
of $\mathcal{C}(U_{i})$ and arrows $f_{ij}$ as above,
then we get objects $\phi_{U_{i}}(C_{i})$ of 
$\mathcal{D}(U_{i})$ and arrows $\phi_{U_{ij}}(f_{ij})$
of $\mathcal{D}(U_{ij})$ which give rise to the 
same \v{C}ech two cocycle $\epsilon_{ijk}$ as for 
$\mathcal{C}$.  

\section{The gerbe $HOM(\P,\Q)$} 
\label{sec:6.4}  
\begin{proposition}[\cite{Gir}] 
\label{prop:6.4.1}
Given two gerbes $\P$ and $\Q$  
on $M$, both bound by $\underline{\C}^{\times}_{M}$,  
we can form a new gerbe $HOM(\P,\Q)$ 
on $M$ also bound by $\underline{\C}_{M}^{\times}$.  
The fibre of  
$HOM(\P,\Q)$ over an open set $U$ of 
$M$ is $HOM(\P,\Q)(U) = Hom(\P(U),\Q(U))$ --- 
that is to say the category with objects 
equal to the morphisms between the gerbes 
$\P(U)$ and $\Q(U)$ and whose arrows are 
the natural isomorphisms between these 
morphisms.  Furthermore there is an 
equivalence of the gerbes $HOM(\P,\Q)$ and 
$\P^{*}\otimes \Q$ on $M$.  There is also 
a `subgerbe' of $HOM(\P,\Q)$, namely the 
gerbe $EQ(\P,\Q)$ whose fibre at an open set 
$U\subset M$ is $EQ(\P,\Q)(U) = Eq(\P(U),\Q(U))$, 
ie the category whose objects are the equivalences 
of gerbes $\P_{U}\to \Q_{U}$ and whose arrows 
are the natural transformations between these 
equivalences.  Clearly there is a morphism of 
gerbes $EQ(\P,\Q)\to HOM(\P,\Q)$ and hence the 
gerbes $EQ(\P,\Q)$, $HOM(\P,\Q)$ and $\Q\otimes 
\P^{*}$ are all equivalent.           
\end{proposition} 
Note that this does not require that $\P$ 
and $\Q$ are equivalent gerbes.   
\begin{proof} 
To show that this does indeed define a gerbe 
we first of all need to show that the 
assignment $U\to HOM(\P,\Q)(U)$ 
defines a stack in groupoids on $M$ --- that is the 
gluing laws for objects and arrows are 
satisfied --- and that $HOM(\P,\Q)$ 
is locally non-empty and locally connected.
Note that by definition $HOM(\P,\Q)_{U}$ 
is a groupoid.  Also $HOM(\P,\Q)$ defines a 
presheaf of categories by definition of a 
morphism of gerbes.  Explicitly if $V\subset U$ 
are open subsets of $M$ then there are natural 
restriction functors $HOM(\P,\Q)_{U}\to 
HOM(\P,\Q)_{V}$ obtained by restricting a morphism 
$\phi_{U}:\P_{U}\to \Q_{U}$ to $V$ and similarly 
for natural transformations.  We have to show that 
in fact $HOM(\P,\Q)$ defines a sheaf of categories.  

Suppose then that we are given an open set $U\subset M$ 
together with an open cover $\{U_{i}\}_{i\in I}$ 
of $U$ by open subsets $U_{i}$ of $M$ such that 
there exist objects (ie morphisms of gerbes) 
$\phi_{i}$ of $HOM(\P,\Q)_{U_{i}}$ 
for every $i\in I$ and arrows (ie natural transformations)  
$\theta_{ij}:\phi_{i}|_{U_{ij}}\Rightarrow 
\phi_{j}|_{U_{ij}}$ in $HOM(\P,\Q)_{U_{ij}}$ which 
satisfy the gluing condition $\theta_{jk}|_{U_{ijk}}
\theta_{ij}|_{U_{ijk}} = \theta_{ik}|_{U_{ijk}}$ 
in $HOM(\P,\Q)_{U_{ijk}}$.  We need to construct a 
morphism of gerbes $\phi:\P_{U}\to \Q_{U}$ which 
is isomorphic to $\phi_{i}$ when restricted to 
$U_{i}$.  Let $V\subset U$ and suppose $P$ 
is an object of $\P_{V}$.  Consider the objects 
$\phi_{i}|_{V\cap U_{i}}(P|_{V\cap U_{i}})$ of $\Q_{V\cap U_{i}}$.        
The natural transformations $\theta_{ij}$ 
provide arrows $\theta_{ij}|_{V\cap U_{ij}}(P|_{V\cap U_{ij}})
:\phi_{i}|_{V\cap U_{i}}(P|_{V\cap U_{i}})|_{V\cap U_{ij}}\to 
\phi_{j}|_{V\cap U_{j}}(P|_{V\cap U_{j}})|_{V\cap U_{ij}}$ 
which are descent data in the gerbe $\Q$ because of 
the condition $\theta_{jk}|_{U_{ijk}}\theta_{ij}|_{U_{ijk}} 
= \theta_{ik}|_{U_{ijk}}$.  Let the descended object 
of $\Q_{V}$ be $Q$ say.  We put $\phi_{V}(P) = Q$.  
In a similar way one can define the action of 
$\phi_{V}$ on arrows of $\P_{V}$ and it follows from 
the uniqueness of gluing in sheaves that this 
defines a functor.  It is also not hard to see that 
that this functor is compatible with the restriction 
functors in $\P$ and $\Q$ and induces the identity 
morphism on the bands of the gerbes $\P$ and $\Q$.  
Finally, note that since $Q$ is the descended object 
for the descent data $\phi_{i}|_{V\cap U_{i}}(P|_{V\cap 
U_{i}})$ and $\theta_{ij}|_{V\cap U_{ij}}(P|_{V\cap 
U_{ij}})$ there are isomorphisms $\chi_{i}:Q|_{V\cap U_{i}}\to 
\phi_{i}|_{V\cap U_{i}}(P|_{V\cap U_{i}})$ compatible 
with the $\theta_{ij}|_{V\cap U_{ij}}(P|_{V\cap U_{ij}})$.  
These isomorphisms induce natural transformations, 
also denoted $\chi_{i}$ from $\phi_{U_{i}}$ to 
$\phi_{i}$.  

We now have to show that arrows glue properly.     
Suppose we are given morphisms of gerbes $\phi, \psi:
\P_{U}\to \Q_{U}$, ie objects of $HOM(\P,\Q)_{U}$, and 
suppose also that we are given an open cover $\{U_{i}\}_{i\in I}$ 
of $U$ by open subsets $U_{i}$ of $M$ such that there 
exist natural transformations $\tau_{i}:\phi_{U_{i}}
\Rightarrow \psi_{U_{i}}$ for each $i\in I$ such that 
$\tau_{i}|_{U_{ij}} = \tau_{j}|_{U_{ij}}$.  Let $P$ 
be an object of $\P_{U}$ and let $P_{i} = P|_{U_{i}}$.  
The $\tau_{i}$ provide arrows $\tau_{i}(P_{i}):
\phi_{U_{i}}(P_{i})\to \psi_{U_{i}}(P_{i})$ in $\Q_{U_{i}}$.  
Since $\phi$ and $\psi$ are morphisms of gerbes we can 
think of the $\tau_{i}(P_{i})$ as arrows $\phi_{U}(P)|_{U_{i}}
\to \psi_{U}(P)|_{U_{i}}$.  We have $\tau_{i}(P_{i})|_{U_{ij}} 
= \tau_{i}|_{U_{ij}}(P|_{U_{ij}}) = \tau_{j}|_{U_{ij}}(
P|_{U_{ij}}) = \tau_{j}(P_{j})|_{U_{ij}}$.  Since $\Q$ is a 
gerbe and hence a sheaf of groupoids, there is a unique 
arrow $\tau_{U}(P):\phi_{U}(P)\to \psi_{U}(P)$ such 
that $\tau_{U}(P)|_{U_{i}} = \tau_{i}(P|_{U_{i}})$.  
The assignment $P\mapsto \tau_{U}(P)$ defines a natural 
transformation.  Hence $HOM(\P,\Q)$ is a sheaf of 
groupoids.   

We will omit the proof that $HOM(\P,\Q)$ is 
locally non-empty and locally connected and 
refer instead to \cite{Gir} where it is also 
shown that $HOM(\P,\Q)$ has band equal to 
$\underline{\C}^{\times}_{M}$.  

To show that there is an equivalence of gerbes 
between $HOM(\P,\Q)$ and $\P^{*}\otimes \Q$, we calculate 
the Dixmier-Douady class of $HOM(\P,\Q)$.  Choose an open 
cover $\U = \{U_{i}\}_{i \in I}$ of $M$ such that there 
exist morphisms $\phi_{i}:\P_{i}\to \Q_{i}$ where 
$\P_{i}$ and $\Q_{i}$ denote the restriction of 
$\P$ and $\Q$ to $U_{i}$ respectively.  We can choose 
the open cover $\U$ so that  
there exist natural isomorphisms 
$\a_{ij}:\phi_{i}|_{U_{ij}}\Rightarrow \phi_{j}|_{U_{ij}}$ 
in $HOM(\P,\Q)(U_{ij})$.  Then the Dixmier-Douady class 
of $HOM(\P,\Q)$ will be given by choosing an object 
$P \in Ob(\P(U_{ijk}))$ and forming the cocycle 
$\a_{ki}|_{U_{ijk}}(P)\a_{jk}|_{U_{ijk}}(P)\a_{ij}|_{U_{ijk}}(P)$.  
We want to relate this cocycle to the cocycles representing 
the Dixmier-Douady classes of $\P$ and $\Q$.  

To do this we first need to choose objects 
$P_{i} \in Ob\ (\P(U_{i}))$ for each $i \in I$.  Again 
we can assume that there exist 
arrows $f_{ij}:P_{i}|_{U_{ij}}\to P_{j}|_{U_{ij}}$ 
in $\P(U_{ij})$.  Since $\a_{ij}$ is a natural 
transformation we get the commutative diagram 
$$
\xymatrix{  
\phi_{i}|_{U_{ij}}(P_{i}|_{U_{ij}}) \ar[r]^{\phi_{i}|_{U_{ij}}(f_{ij})} 
\ar[d]_{\a_{ij}(P_{i}|_{U_{ij}})} & \phi_{i}|_{U_{ij}}(P_{j}|_{U_{ij}}) 
\ar[d]^{\a_{ij}(P_{j}|_{U_{ij}})}                                      \\ 
\phi_{j}|_{U_{ij}}(P_{i}|_{U_{ij}}) \ar[r]^{\phi_{j}|_{U_{ij}}(f_{ij})} 
& \phi_{j}|_{U_{ij}}(P_{j}|_{U_{ij}}).                   } 
$$
This induces the following diagram of arrows in 
$\Q(U_{ijk})$ (We have omitted some restriction 
functors for convenience of notation).  
$$ 
\xymatrix{ 
\phi_{i}(P_{i})|_{U_{ijk}}  \ar[r]^-{\phi_{i}
(f_{ij})|_{U_{ijk}}} \ar[d]_{\a_{ij}|_{U_{ijk}}(
P_{i}|_{U_{ijk}})} \ar[dr] & \phi_{i}(P_{j})|_{U_{ijk}} 
\ar[r]^-{\phi_{i}(f_{jk})|_{U_{ijk}}} \ar[d] & 
\phi_{i}(P_{k})|_{U_{ijk}} \ar[d] \ar[r]^-{
\phi_{i}(f_{ik})|_{U_{ijk}}} & \phi_{i}
(P_{i})|_{U_{ijk}}   \ar[d]                                  \\ 
\phi_{j}(P_{i})|_{U_{ijk}} \ar[r] \ar[d]_{\a_{jk}
(P_{i})|_{U_{ijk}}} & \phi_{j}(P_{j})|_{U_{ijk}} 
\ar[r] \ar[d] \ar[dr] & \phi_{j}(P_{k})|_{U_{ijk}} \ar[r] 
\ar[d] & \phi_{j}(P_{i})|_{U_{ijk}} \ar[d]        \\ 
\phi_{k}(P_{i})|_{U_{ijk}}) \ar[d]_{\a_{ki}
(P_{i})|_{U_{ijk}}} \ar[r] & \phi_{k}(P_{j})|_{U_{ijk}} 
\ar[r] \ar[d] & \phi_{k}(P_{k})|_{U_{ijk}} \ar[r] 
\ar[d] \ar[dr] & \phi_{k}(P_{i})|_{U_{ijk}} \ar[d] \\ 
\phi_{i}(P_{i})|_{U_{ijk}} \ar[r] & 
\phi_{i}(P_{j})|_{U_{ijk}} \ar[r] & \phi_{i}
(P_{k})|_{U_{ijk}} \ar[r] & \phi_{i}(P_{i})|_{U_{ijk}}        } 
$$
The diagonal arrow is a cocycle representing the 
Dixmier-Douady class of the gerbe $\Q$, while 
the bottom horizontal arrow is a cocycle representing the 
class of $\P$.  Since the left vertical arrow 
represents the Dixmier-Douady class of $HOM(\P,\Q)$, we 
must have $DD(HOM(\P,\Q)) = DD(\Q) - DD(\P)$.  
\end{proof}

For us the main point of this Proposition 
is that one can glue a family of locally defined morphisms 
between two gerbes together to form a global 
morphism between the two gerbes.  This is also 
the case for bundle gerbes and stable morphisms 
between bundle gerbes.  More precisely, suppose 
we have bundle gerbes $(P,X,M)$ and $(Q,Y,M)$ 
together with stable morphisms $\phi_{i}:P|_{U_{i}}\to 
Q_{U_{i}}$ where $\U = \{U_{i}\}_{i\in I}$ is an 
open cover of $M$.  Suppose we also have 
transformations $\theta_{ij}:\phi_{i}\Rightarrow 
\phi_{j}$ between the restrictions of the stable morphisms 
$\phi_{i}$ and $\phi_{j}$ to $U_{ij}$ which satisfy 
the gluing condition $\theta_{jk}\theta_{ij} 
= \theta_{ik}$ on $U_{ijk}$.  Then there is a 
stable morphism $\phi:P\to Q$ together with 
transformations $\xi_{i}:\phi|_{U_{i}}\Rightarrow 
\phi_{i}$ which are compatible with $\theta_{ij}$.  

To see this suppose the stable morphisms 
$\phi_{i}$ correspond to trivialisations 
$L_{\phi_{i}} = L_{i}$ of $P^{*}\otimes Q$ 
on $X\times_{U_{i}}Y$.  Then the transformations 
$\theta_{ij}$ provide isomorphisms 
$\tilde{\theta}_{ij}:L_{i}\to L_{j}$ over 
$X\times_{U_{ij}}Y$.  The gluing condition 
$\theta_{jk}\theta_{ij} = \theta_{ik}$ 
translates into the condition $\tilde{\theta}
_{jk}\circ \tilde{\theta}_{ij} = \tilde{
\theta}_{ik}$.  Hence we can use the 
standard clutching construction to glue 
the various principal $\cstar$ bundles 
$L_{i}$ together to form a trivialisation 
$L$ defined on the whole of $X\times_{M}Y$.  

Unfortunately there is no general method 
for gluing together bundle gerbe morphisms; 
ie given bundle gerbes $(P,X,M)$ and $(Q,Y,M)$ 
together with locally defined bundle gerbe 
morphisms $\bar{f}_{i}:P|_{U_{i}}\to Q|_{U_{i}}$ 
and transformations of bundle gerbe morphisms 
$\theta_{ij}:\bar{f}_{i}|_{U_{ij}}\Rightarrow 
\bar{f}_{j}|_{U_{ij}}$ satisfying the 
gluing condition $\theta_{jk}\theta_{ij} = 
\theta_{ik}$ over $U_{ijk}$ then there is 
no general way to construct a bundle gerbe 
morphism $\bar{f}:P\to Q$ which restricts 
to $\bar{f}_{i}$.  
   
\section{The gerbe associated to a bundle gerbe} 
\label{sec:6.5} 
We will briefly review a construction given 
in \cite{MurSte} which associates a gerbe 
$\mathcal{G}(P)$ to a bundle gerbe $(P,X,M)$ 
in such a way that $DD(P) = DD(\mathcal{G}(P))$.  
Given an open set $U\subset M$ we define 
a groupoid $\mathcal{G}(P)_{U}$ as follows.  The 
objects of $\mathcal{G}(P)_{U}$ consist of 
line bundles $L\to X_{U} = \pi^{-1}(U)$ 
together with an isomorphism $\eta:\d(L)\to P_{X_{U}^{[2]}}$ 
of line bundles covering the identity 
on $X_{U}^{[2]} = (\pi^{[2]})^{-1}(U)$ 
and commuting with the bundle gerbe products.  
(Here we regard $\d(L)$ as a trivial bundle gerbe). 
An arrow in $\mathcal{G}(P)_{U}$ from an object 
$(L_{1},\eta_{1})$ to an object $(L_{2},\eta_{2})$ 
is an isomorphism of line bundles $\a:L_{1}\to L_{2}$ 
which is compatible with $\eta_{1}$ and $\eta_{2}$ 
in the sense that $\eta_{2}\d(\a) = \eta_{1}$.   
Given another open set $V$ with 
$V\subset U\subset M$ then are canonical pullback 
functors $\rho_{V,U}:\mathcal{G}(P)_{U}\to \mathcal{G}(P)_{V}$ 
satisfying $\rho_{W,V}\rho_{V,U} = \rho_{W,U}$.  
In \cite{MurSte} it is shown that the assignment 
$U\mapsto \mathcal{G}(P)_{U}$ defines a gerbe.  
To see that the band of this 
gerbe is equal to $\underline{\C}^{\times}_{M}$, 
note that if $(L,\eta)$ is an object of 
$\mathcal{G}(P)_{U}$ then, since we may identify 
an automorphism $\a$ of $(L,\eta)$ with 
a function $f:X_{U}\to \cstar$, the condition 
$\eta\d(\a) = \eta$ forces $\d(f) = 1$ and 
so $f$ must descend to a function $\tilde{f}:U\to \cstar$.  
In this way we can construct an isomorphism 
$\underline{Aut}(L,\eta)\cong \underline{\C}^{\times}_{U}$ 
which is compatible with restriction to smaller open sets 
and with isomorphisms $Aut(L_{1},\eta_{1})\to 
Aut(L_{2},\eta_{2})$ induced by an isomorphism 
$(L_{1},\eta_{1})\to (L_{2},\eta_{2})$.  
The advantage of this construction is that one 
does not need to apply a `sheafification' 
procedure to force objects and arrows to 
satisfy the effective descent conditions.  
To see that the objects of $\mathcal{G}(P)_{U}$ 
satisfy effective descent note that if 
$\{U_{i}\}_{i\in I}$ is an open cover of $U$,  
$(L_{i},\eta_{i})$ are objects of 
$\mathcal{G}(P)_{U_{i}}$ and $\a_{ij}:(L_{i},\eta_{i})|_{U_{ij}}
\to (L_{j},\eta_{j})|_{U_{ij}}$ are 
arrows of $\mathcal{G}(P)_{U_{ij}}$ satisfying 
$\a_{jk}\a_{ij} = \a_{ik}$ then we can glue 
the line bundles $L_{i}$ together with the 
$\a_{ij}$ using the standard clutching 
construction to form a line bundle $L$ on 
$X_{U}$.  The $\eta_{i}$ glue together to 
form an isomorphism $\eta:\d(L)\to P_{X_{U}^{[2]}}$ 
which trivialises $P_{X_{U}^{[2]}}$.  Similarly 
one can show that the arrows of $\mathcal{G}(P)_{U}$ 
satisfy the effective descent condition. 

To see that $\mathcal{G}(P)$ had Dixmier-Douady 
class equal to that of the bundle gerbe $P$, 
choose an open cover $\{U_{i}\}_{i\in I}$ of 
$M$, all of whose non-empty intersections 
$U_{i_{0}}\cap \cdots \cap U_{i_{p}}$ are 
contractible, such that there exist local sections 
$s_{i}:U_{i}\to X$ of $\pi$.  We can define objects 
$(L_{i},\eta_{i})$ of $\mathcal{G}(P)_{U_{i}}$ by 
setting $L_{i} = \hat{s}_{i}^{-1}P_{X_{U_{i}}^{[2]}}$, where 
$\hat{s}_{i}:X_{U_{i}}\to X_{U_{i}}^{[2]}$ is 
the map which sends $x\in X_{U_{i}}$ to 
$(x,s_{i}(\pi(x)))\in X_{U_{i}}^{[2]}$.   
$\eta_{i}:\d(L_{i})\to P_{X_{U_{i}}^{[2]}}$ 
is induced by the bundle gerbe product on $P$.  
We can define arrows $\a_{ij}:(L_{i},\eta_{i})|_{U_{ij}}
\to (L_{j},\eta_{j})|_{U_{ij}}$ by first choosing 
sections $\sigma_{ij}:U_{ij}\to (s_{i},s_{j})^{-1}P$ 
and then defining $\a_{ij}(u_{i}) = m_{P}(\sigma(\pi_{L_{i}}
(u_{i}))\otimes u_{i})$ for $u_{i}\in L_{i}$, where 
$\pi_{L_{i}}:L_{i}\to X_{U_{i}}$ denotes the projection.  
If we now calculate the \v{C}ech cocycle representative 
for the Dixmier-Douady class of $\mathcal{G}(P)$ from 
this data, then we get it exactly equal to the 
cocycle for the bundle gerbe $P$.  
We summarise this discussion in the following 
proposition. 
\begin{proposition}[\cite{MurSte}] 
\label{prop:6.5.1} 
Given a bundle gerbe $(P,X,M)$ the assignment 
$U\mapsto \mathcal{G}(P)_{U}$ defined above is a  
gerbe with Dixmier-Douady class equal to the 
Dixmier-Douady class of the bundle gerbe $P$.  
\end{proposition} 
 
\begin{note} 
\label{note:6.5.2} 
One can show that a bundle gerbe morphism 
$\bar{f}:P\to Q$ induces a morphism of 
gerbes $\mathcal{G}(f):\mathcal{G}(P)\to 
\mathcal{G}(Q)$ and that a transformation 
$\theta:\bar{f}\Rightarrow \bar{g}$ induces 
a transformation $\mathcal{G}(\theta):\mathcal{G}(f)
\Rightarrow \mathcal{G}(g)$ between the induced 
morphisms of gerbes.  Thus we can think of $\mathcal{G}$ 
as furnishing us with a 2-functor from the 2-category 
of bundle gerbes and bundle gerbe morphisms 
to the 2-category of gerbes.  Given two bundle gerbes 
$P_{1}$ and $P_{2}$ on $M$ one can show that there 
is an equivalence of gerbes $\mathcal{G}(P_{1}\otimes 
P_{2})\simeq \mathcal{G}(P_{1})\otimes \mathcal{G}(P_{2})$ 
and one can show as well that given a bundle gerbe 
$P$ on $M$ and a smooth map $f:N\to M$ there is an 
equivalence of gerbes $\mathcal{G}(f^{-1}P)\simeq 
f^{*}\mathcal{G}(P)$.  Unfortunately I have been unable 
to show that these equivalences behave as one would like; 
namely that the equivalence $\mathcal{G}(f^{-1}P)\simeq 
f^{*}\mathcal{G}(P)$ commutes with the equivalences 
$\phi_{f,g}$ of Lemma~\ref{lemma:6.3.2} and the equivalence 
$\mathcal{G}(P_{1}\otimes P_{2})\simeq \mathcal{G}(P_{1})
\otimes \mathcal{G}(P_{2})$ commutes with the equivalence 
$\mathcal{G}_{1}\otimes (\mathcal{G}_{2}\otimes \mathcal{G}_{3})
\simeq (\mathcal{G}_{1}\otimes \mathcal{G}_{2})\otimes 
\mathcal{G}_{3}$ of gerbes $\mathcal{G}_{1}$, $\mathcal{G}_{2}$ 
and $\mathcal{G}_{3}$ on $M$.      
\end{note} 

We now wish to show how a bundle gerbe $(P,X,M)$ 
with bundle gerbe connection $\nabla$ and curving 
$f$ gives rise to a connective structure and curving on 
the associated gerbe $\mathcal{G}(P)$.  We first need 
to recall the definition of connective structure 
and curving from \cite{Bry} and \cite{BryMcL1}.  

\begin{definition}[\cite{Bry}, \cite{BryMcL1}] 
\label{def:6.5.3} 
Let $\mathcal{G}$ be a gerbe on a smooth manifold 
$M$.  A \emph{connective structure} on $\mathcal{G}$ 
s given by the following data.  

\begin{enumerate} 
\item For every open set $U\subset M$, an 
assignment to each object $P$ of $\mathcal{G}(U)$ 
of an $\underline{\Omega}^{1}_{M}$ torsor 
$\underline{Co}(P)$ on $U$ which is compatible with the 
pullback functors $\rho_{V,U}$.  Sections of 
$\underline{Co}(P)$ are usually denoted by $\nabla$.    
\item Given objects $P_{1}$ and $P_{2}$ of 
$\mathcal{G}_{U}$ and an isomorphism 
$f:P_{1}\to P_{2}$, there is an isomorphism 
$f_{*}:\underline{Co}(P_{1})\to \underline{Co}(P_{2})$ of 
$\underline{\Omega}^{1}_{M}$ 
torsors on $U$, which is compatible with 
composition of morphisms in $\mathcal{G}(U)$ and 
with the pullback functors $\rho_{V,U}$.  
\item Given an object $P$ of $\mathcal{G}(U)$ and 
an automorphism $f:P\to P$ then we require that 
$f_{*}:\underline{Co}(P)\to \underline{Co}(P)$ is the automorphism 
defined by $\nabla\mapsto \nabla + df/f$, where 
we regard $f$ as an element of $\underline{\C}^{\times}_{U}$.  
\end{enumerate} 

A \emph{curving} for a connective structure on a 
gerbe $\mathcal{G}$ on $M$ consists of an assignment 
to each open set $U\subset M$, each object $P$ of 
$\mathcal{G}(U)$ and each section $\nabla$ of 
$\underline{Co}(P)$ of a two form $f(P,\nabla)$ on 
$U$ satisfying the following conditions: 

\begin{enumerate} 
\item For objects $P_{1}$ and $P_{2}$ of $\mathcal{G}(U)$ 
and isomorphisms $\phi:P_{1}\to P_{2}$ in $\mathcal{G}(U)$ 
we have $f(P_{1},\nabla) = f(P_{2},\phi_{*}(\nabla))$ 
for all sections $\nabla$ of $\underline{Co}(P_{1})$.  
\item If $A$ is a one form on $U$ then for any object 
$P$ of $\mathcal{G}(U)$ and any section 
$\nabla$ of $\underline{Co}(P)$ we have 
$f(P,\nabla + A) = f(P,\nabla) + dA$.  
\item Given open sets $V\subset U\subset M$, an 
object $P$ of $\mathcal{G}(U)$ and a section $\nabla$ 
of $\underline{Co}(P)$, the assignment $(P,\nabla)\mapsto 
f(P,\nabla)$ is required to be compatible with 
restrictions to smaller open sets.   
\end{enumerate}
\end{definition}        
 
We now show how a bundle gerbe connection $\nabla$ 
with a curving $f$ on a bundle gerbe $(P,X,M)$ 
gives rise to a connective structure and a curving 
on the associated gerbe $\mathcal{G}(P)$ on $M$.  Let 
$U\subset M$ be an open subset of $M$.  Let $(L,\eta)$ 
be an object of $\mathcal{G}(P)_{U}$.  Let $\nabla_{L}$ be 
any connection on $L$.  $\nabla_{L}$ induces a 
connection $\d(\nabla_{L})$ on $\d(L)$.  We also have 
the connection $\eta^{-1}\circ \nabla\circ \eta$ on 
$\d(L)$ induced by $\nabla$.  These two connections 
must differ by the pullback of a one form $A$ on 
$X_{U}^{[2]}$.  Since $\nabla$ is a bundle gerbe 
connection, we must have $\d(A) = 0$.  Hence we 
can solve $A = \d(B)$ for some one form $B$ on 
$X_{U}$.  Therefore if we had originally defined 
$\nabla_{L}$ to be equal to $\nabla_{L} - B$ then 
we would have $\d(\nabla_{L}) = \eta^{-1}\circ \nabla\circ \eta$.  
Therefore we define $\underline{Co}(L,\eta)$ to 
be the $\underline{\Omega}^{1}_{U}$ torsor 
consisting of the connections $\nabla_{L}$ on 
$L$ satisfying $\d(\nabla_{L}) = \eta^{-1}\circ \nabla\circ \eta$.  
It is clear from the previous discussion that 
this does indeed define an $\underline{\Omega}^{1}_{U}$ 
torsor.  It is also clear that the conditions (1) 
and (2) above are satisfied.  To see that condition 
(3) is satisfied, let $f:U\to \cstar$ be an 
automorphism of $(L,\eta)$ under the isomorphism 
$\underline{Aut}(L,\eta)\cong \underline{\C}^{\times}_{M}$.  
Then if $\nabla_{L}$ is a section of 
$\underline{Co}(L,\eta)$ then we have 
$f_{*}(\nabla) = \nabla + \pi^{*}(df/f)$, as 
required.   
To define a curving for this connective 
structure, simply $f(L,\eta,\nabla_{L})$ be 
the two form $F_{\nabla_{L}} - f$ where $F_{
\nabla_{L}}$ is the curvature of the connection 
$\nabla_{L}$ on $L$.  Since $\d(F_{\nabla_{L}} 
- f) = 0$ we may identify $f(L,\eta,\nabla_{L})$ 
with a two form on $U$.   
This clearly defines a curving for the connective 
structure on $\mathcal{G}(P)$.  
 
\section{Higher gluing laws} 
\label{sec:6.6} 

In this section we discuss how to glue 
together a family of gerbes $\mathcal{G}_{i}$ defined on an 
open cover $\{U_{i}\}_{i\in I}$ of $M$ 
to get a gerbe defined on the whole of $M$.  
This has an extra level of complication in that 
the gluing morphisms $\phi_{ij}:\mathcal{G}_{j}\to 
\mathcal{G}_{i}$ need not satisfy $\phi_{jk}\phi_{ij} 
= \phi_{ik}$ but rather there is a natural 
transformation $\psi_{ijk}:\phi_{jk}\phi_{ij}\Rightarrow 
\phi_{ik}$ satisfying a certain cocycle condition 
over $U_{ijkl}$.  The triple $(\mathcal{G}_{i},
\phi_{ij},\psi_{ijk})$ is called \emph{2-descent data}.  
We will need the results of this section in 
Chapter~\ref{chapter:13}.   

\begin{proposition}[\cite{Bre}]   
\label{prop:6.6.1} 
Suppose we are given an open cover $\{U_{i}\}_{i\in I}$ of 
$X$ and gerbes $\mathcal{G}_{i}$ on $U_{i}$ for 
each $i\in I$ as well as morphisms of gerbes 
$\varphi_{ij}:\mathcal{G}_{i}\to \mathcal{G}_{j}$ 
between the restrictions of the gerbes $\mathcal{G}_{i}$ 
and $\mathcal{G}_{j}$ to $U_{ij}$, which satisfy 
$\varphi_{ii} = \text{id}_{\mathcal{G}_{i}}$.  Suppose we are 
also given natural transformations $\psi_{ijk}:
\varphi_{jk}|_{U_{ijk}}\varphi_{ij}|_{U_{ijk}}
\Rightarrow \varphi_{ik}|_{U_{ijk}}$ over $U_{ijk}$ 
which satisfy the following condition (the non-abelian 
2-cocycle condition): the diagram of natural 
transformations below commutes: 
$$
\xymatrix{ 
\mathcal{G}_{i} \ar[rr]^-{\varphi_{ij}} 
\ar[dd]_-{\varphi_{il}} \ar[ddrr] & & 
\mathcal{G}_{j} \ar[dd]_-{\varphi_{jk}} 
\ar @2{->}[dl]_-{\psi_{ijk}}     & 
\mathcal{G}_{i} \ar[dd]^<<<<<{\varphi_{il}} 
\ar[rr]^-{\varphi_{ij}} & & 
\mathcal{G}_{j} \ar[dd]^-{\varphi_{jk}} \ar[ddll]    \\ 
& \ar @2{->}[l]_-{\psi_{ikl}} &  
\ar @2{-}[r] & & \ar @2{->}[l]^-{\psi_{ijl}}  &    \\ 
\mathcal{G}_{l} & & \mathcal{G}_{k} 
\ar[ll]_-{\varphi_{kl}}  & \mathcal{G}_{l} 
&   & \mathcal{G}_{k} 
\ar[ll]_-{\varphi_{kl}} \ar @2{->}[ul]_-{
\psi_{jkl}}.                                      } 
$$
What this means is that after pasting together 
the natural transformations in each diagram, 
the two natural transformations between the morphisms 
bounding the above diagrams are equal; so in 
other words we have 
$$
\psi_{ikl}(1_{\varphi_{kl}}\circ \psi_{ijk}) = 
\psi_{ijl}(\psi_{jkl}\circ 1_{\varphi_{ij}}). 
$$
Then there exists a gerbe $\mathcal{G}$ on $X$ 
and an equivalence of gerbes $\chi_{i}:\mathcal{G}|_{U_{i}}\to 
\mathcal{G}_{i}$ over $U_{i}$ together with 
natural transformations $\xi_{ij}:\varphi_{ij}
\chi_{i}|_{U_{ij}}\Rightarrow \chi_{j}|_{U_{ij}}$ 
which are compatible with the natural transformations 
$\psi_{ijk}$ in the sense that the following diagram of 
natural transformations commutes: 
$$
\xymatrix{ 
\mathcal{G}|_{U_{ijk}} \ar[rr]^-{\chi_{i}} 
\ar[dd]_-{\chi_{k}} & & \mathcal{G}_{i} 
\ar[ddll] \ar[dd]_-{\varphi_{ij}} & 
\mathcal{G}|_{U_{ijk}} \ar[dd]_<<<<{\chi_{k}} 
\ar[rr]^-{\chi_{i}} \ar[ddrr] & & 
\mathcal{G}_{i} \ar @2{->}[dl]_-{\xi_{ij}} 
\ar[dd]^-{\varphi_{ij}}                    \\ 
& \ar @2{->}[l]^-{\xi_{ik}} & \ar @2{-}[r] & 
& \ar @2{->}[l]_-{\xi_{jk}} &               \\ 
\mathcal{G}_{k} & & \mathcal{G}_{j} 
\ar[ll]^-{\varphi_{jk}} \ar @2{->}[ul]^-{
\psi_{ijk}} & \mathcal{G}_{k} & & 
\mathcal{G}_{j} \ar[ll]^-{\varphi_{jk}}.   } 
$$
In other words we have 
$$
\xi_{ik}(\psi_{ijk}\circ 1_{\chi_{i}}) = 
\xi_{jk}(1_{\varphi_{jk}}\circ \xi_{ij}). 
$$
\end{proposition} 

We have an analogue 
of the above proposition for bundle gerbes 
and stable morphisms.  We have to exercise a 
bit more care here though because composition 
of stable morphisms is not strictly associative but 
rather associative up to a coherent transformation 
(see Proposition~\ref{prop:stablebigrpd}).   

\begin{proposition} 
\label{prop:6.6.2} 
Let $\{U_{i}\}_{i\in I}$ be an open cover of 
$M$ by open subsets $U_{i}$ of $M$.  Suppose 
for each $i\in I$ we have a bundle gerbe 
$(Q_{i},Y_{i},U_{i})$ over $U_{i}$ such that for 
each pair $i, j\in I$ with the intersection 
$U_{ij}$ non-empty there is a stable morphism 
$f_{ij}:Q_{i}\to Q_{j}$ corresponding to a 
trivialisation $(L_{f_{ij}},\phi_{f_{ij}})$ of $Q_{i}^{*}\otimes 
Q_{j}$.  We require that the stable morphisms 
$f_{ij}$ satisfy $f_{ji} = \tilde{f}_{ij}$ 
in the notation of Section~\ref{sec:3.4} so that 
$\tilde{f}_{ij}$ corresponds to the trivialisation 
of $Q_{j}^{*}\otimes Q_{i}$ given by 
$(L_{f_{ij}}^{*},\phi_{f_{ij}}^{*})$.  
Suppose also there is a transformation 
$\theta_{ijk}:f_{jk}\circ f_{ij}\Rightarrow f_{ik}$ 
of stable morphisms over $U_{ijk}$ which satisfies 
the condition 
$$
\theta_{ijl}(\theta_{jkl}\circ 1_{f_{ij}})\theta_{
f_{ij},f_{jk},f_{kl}} = 
\theta_{ikl}(1_{f_{kl}}\circ \theta_{ijk}) 
$$
over $U_{ijkl}$.  Then there is a bundle gerbe 
$(Q,Y,M)$ on $M$ together with stable morphisms 
$g_{i}:Q|_{U_{i}}\to Q_{i}$ over $U_{i}$ such 
that there exist transformations of stable 
morphisms $\xi_{ij}:f_{ij}\circ g_{i}\Rightarrow 
g_{j}$ over $U_{ij}$ which satisfy the compatibility 
condition 
$$
\xi_{jk}(1_{f_{jk}}\circ \xi_{ij}) = 
\xi_{ik}(\theta_{ijk}\circ 1_{g_{i}}) 
$$
with $\theta_{ijk}$ over $U_{ijk}$.  
\end{proposition} 

\begin{proof} 
Recall that the stable morphism 
$f_{jk}\circ f_{ij}$ corresponds to 
the trivialisation 
$$
(L_{f_{jk}\circ f_{ij}},\phi_{f_{jk}\circ 
f_{ij}})
$$ 
of $Q_{i}^{*}\otimes 
Q_{k}$ which is obtained by descending the 
$\cstar$ bundle $L_{f_{jk}}\otimes 
L_{f_{ij}}\otimes Q_{j}^{*}$ on 
$Y_{i}\times_{M}
Y_{k}\times_{M}Y_{j}^{[2]}$  
to $Y_{i}\times_{M}
Y_{k}$.  In particular, this implies that the 
following isomorphism is true on 
$Y_{i}\times_{M}Y_{k}\times_{M}Y_{j}^{[2]}$: 
\begin{equation} 
\label{eq:stableone} 
L_{f_{jk}\circ f_{ij}} 
\simeq L_{f_{jk}}\otimes L_{f_{ij}}\otimes Q_{j}^{*}.   
\end{equation} 
The two trivialisations $L_{f_{ik}}$ and 
$L_{f_{jk}\circ f_{ij}}$ of $Q_{i}^{*}\otimes 
Q_{k}$ differ by the pullback to 
$Y_{i}\times_{M}Y_{k}$ of a $\cstar$ 
bundle $J_{ijk}$ on $U_{ijk}$.  So we have 
\begin{equation} 
\label{eq:stabletwo}  
L_{f_{jk}\circ f_{ij}} \simeq L_{f_{ik}}
\otimes \pi_{Y_{i}\times_{M}Y_{k}}^{-1}
J_{ijk}.   
\end{equation} 
$\theta_{ijk}$ is then a section 
of the bundle $J_{ijk}$ on $U_{ijk}$.    
From equations~\ref{eq:stableone} and~\ref{eq:stabletwo}  
we have 
\begin{eqnarray*} 
L_{(f_{kl}\circ f_{jk})\circ f_{ij}} & 
\simeq & L_{f_{kl}\circ f_{jk}}\otimes L_{f_{ij}}\otimes Q_{k}^{*}   \\ 
& \simeq & L_{f_{jl}}\otimes L_{f_{ij}}\otimes Q_{k}^{*}\otimes  
\pi^{-1}J_{jkl}                               \\ 
& \simeq & L_{f_{jl}\circ f_{ij}}\otimes Q_{j}^{*}\otimes Q_{j}  
\otimes \pi^{-1}J_{jkl}                          \\ 
& \simeq & L_{f_{il}}\otimes   
\pi^{-1}(J_{jkl}\otimes J_{ijl}),     
\end{eqnarray*} 
where we have denoted the various projections  
to $M$ by $\pi$.  Using equations~\ref{eq:stableone} 
and~\ref{eq:stabletwo} again we get 
\begin{eqnarray*} 
L_{f_{kl}\circ (f_{jk}\circ f_{ij})} 
 & \simeq & L_{f_{kl}}\otimes 
L_{f_{jk}\circ f_{ij}}\otimes Q_{k}^{*}       \\ 
& \simeq & L_{f_{kl}}\otimes L_{f_{ik}}\otimes Q_{k}^{*}\otimes  
\pi^{-1}J_{ijk}                        \\ 
& \simeq & L_{f_{kl}\circ f_{ik}}\otimes Q_{k}\otimes Q_{k}^{*} 
\otimes \pi^{-1}J_{ijk}                  \\ 
& \simeq & L_{f_{il}}\otimes  
\pi^{-1}(J_{ikl}\otimes J_{ijk}),  
\end{eqnarray*} 
where again we have denoted the projections  
to $M$ by $\pi$.  The transformation 
$\theta_{f_{ij},f_{jk},f_{kl}}$ provides 
an isomorphism $L_{f_{kl}\circ (f_{jk}\circ f_{ij})} 
\simeq L_{(f_{kl}\circ f_{jk})\circ f_{ij}}$ which 
is compatible with the trivialisations  
$\phi_{f_{kl}\circ (f_{jk}\circ f_{ij})}$ and 
$\phi_{(f_{kl}\circ f_{jk})\circ f_{ij}}$.  It 
follows from the series of equations above and the property 
of $\theta_{f_{ij},f_{jk},f_{kl}}$ just mentioned 
that there is an isomorphism  
$$
J_{jkl}\otimes J_{ijl}\simeq J_{ikl}\otimes J_{ijk} 
$$
over $U_{ijkl}$.  The cocycle condition on 
$\theta_{ijk}$ is that the section $\theta_{jkl}
\otimes \theta_{ijl}$ of $J_{jkl}\otimes J_{ijl}$ 
is mapped to the section $\theta_{ikl}\otimes 
\theta_{ijk}$ of $J_{ikl}\otimes J_{ijk}$ under 
the this isomorphism.  We now define the bundle 
gerbe $(Q,Y,M)$ as follows.  We put $Y$ equal to 
the disjoint union $\coprod_{i\in I}Y_{i}$ with 
the obvious projection $Y\to M$.  The fibre product 
$Y^{[2]}$ is equal to the disjoint union 
$\coprod_{i,j\in I}Y_{i}\times_{U_{ij}}Y_{j}$.  
Define a $\cstar$ bundle $Q$ on $Y^{[2]}$ by 
putting $Q|_{Y_{i}\times_{M}Y_{j}} = L_{f_{ij}}$.  
We now need to define a bundle gerbe product 
on $Q$.  This will be a $\cstar$ bundle isomorphism 
$L_{f_{kl}}\otimes L_{f_{jk}}\to L_{f_{ik}}$ 
covering the identity on $Y_{i}\times_{M}Y_{j}
\times_{M}Y_{k}$ which satisfies the 
associativity condition.  Note that from 
equations~\ref{eq:stableone} and~\ref{eq:stabletwo} 
we have the isomorphism 
$$
L_{f_{jk}}\otimes L_{f_{ij}} \simeq 
L_{f_{ik}}\otimes Q_{j}\otimes 
\pi^{-1}J_{ijk} 
$$
of $\cstar$ bundles on $Y_{i}\times_{M}
Y_{j}^{[2]}\times_{M}Y_{k}$.  If we embed 
$Y_{i}\times_{M}Y_{j}\times_{M}Y_{k}$ 
inside $Y_{i}\times_{M}Y_{j}^{[2]}\times_{
M}Y_{k}$ by sending $(y_{i},y_{j},y_{k})$ 
to $(y_{i},y_{j},y_{j},y_{k})$, then using 
the section $\theta_{ijk}$ of $J_{ijk}$ and the 
identity section of $Q_{j}$ over the diagonal, 
we can construct an isomorphism $L_{f_{jk}}\otimes 
L_{f_{ij}}\to L_{f_{ik}}$.  The coherency condition 
on $\theta_{ijk}$ shows that this product is associative 
and hence $(Q,Y,M)$ is a bundle gerbe.    

We now need to show that there is a 
stable morphism $g_{i}:Q|_{U_{i}}\to Q_{i}$, 
ie a trivialisation $(L_{g_{i}},\phi_{g_{i}})$ of $Q^{*}
\otimes Q_{i}$ over $U_{i}$.  The $\cstar$ 
bundle $L_{g_{i}}$ lives over $Y\times_{M}Y_{i}= 
\coprod_{j\in I}Y_{j}\times_{M}Y_{i}$.  
Let $L_{g_{i}}$ be the disjoint union of the $\cstar$ 
bundles $L_{f_{ji}}$.  Since $f_{ji} = 
\tilde{f}_{ij}$ the isomorphism 
$L_{f_{kj}}\otimes Q_{i}\simeq 
L_{f_{ki}}\otimes L_{f_{ij}}$ becomes 
$L_{f_{jk}}^{*}\otimes Q_{i}\simeq 
L_{f_{ki}}\otimes L_{f_{ji}}^{*}$ which shows that 
$L_{g_{i}}$ trivialises $Q^{*}\otimes Q_{i}$.  

Finally we have to define transformations 
$\xi_{ij}:f_{ij}\circ g_{i}\Rightarrow g_{j}$, 
compatible with $\theta_{ijk}$ in the 
appropriate sense.  We have two trivialisations 
$L_{f_{ij}\circ g_{i}}$ and $L_{g_{j}}$ of 
$Q^{*}\otimes Q_{j}$ and hence there is a $\cstar$ 
bundle $D_{ij}$ on $U_{ij}$ such that 
$L_{f_{ij}\circ g_{i}} = L_{g_{j}}\otimes \pi^{-1}
D_{ij}$ where $\pi$ denotes the projection 
$Y\times_{M}Y_{j}\to U_{j}$.  Therefore we 
have $L_{f_{ij}\circ g_{i}}\otimes Q_{i}\simeq 
L_{g_{j}}\otimes Q_{i}\otimes \pi^{-1}D_{ij}$, an 
isomorphism of $\cstar$ bundles on $\coprod_{k\in I}
Y_{k}\times_{M}Y_{j}\times_{M}Y_{i}^{[2]}$.  
On each constituent $Y_{k}\times_{M}Y_{j}\times_{M}
Y_{i}^{[2]}$ of the above disjoint union this 
isomorphism takes the form 
$$
L_{f_{ij}}\otimes L_{f_{ki}}\simeq L_{f_{kj}}\otimes 
Q_{i}\otimes \pi^{-1}D_{ij}.  
$$
From equations~\ref{eq:stableone} and~\ref{eq:stabletwo} 
we have 
\begin{eqnarray*} 
L_{f_{ki}}\otimes L_{f_{ij}}& \simeq &  
L_{f_{ki}\circ f_{ij}}\otimes Q_{i} \\    
& \simeq & L_{f_{kj}}\otimes \pi^{-1}J_{kij}
\otimes Q_{i}. 
\end{eqnarray*} 
This shows that there is an isomorphism 
$\pi^{-1}D_{ij} \simeq \pi^{-1}J_{kij}$ on 
$Y_{k}\times_{M}Y_{j}\times_{M}
Y_{i}^{[2]}$.  Therefore $\pi^{-1}\theta_{kij}$ 
induces a section $\hat{\xi}_{ij}$ of $\pi^{-1}D_{ij}$ 
which descends to a section $\xi_{ij}$ of 
$D_{ij}$.  Because of the coherency condition 
satisfied by $\theta_{kij}$, it follows that 
$\xi_{ij}$ will be compatible with $\theta_{ijk}$.     
\end{proof} 
We have the following generalisation of the 
above proposition which we will need in 
Chapter~\ref{chapter:13}.  
\begin{proposition} 
\label{prop:6.6.3} 
Suppose $\pi:X\to M$ is a surjection 
admitting local sections and that $(P,Y,X)$ 
is a bundle gerbe on $X$ such that there is 
a stable morphism $f:\pi_{1}^{-1}P\to 
\pi_{2}^{-1}P$ corresponding to a trivialisation 
$(L_{f},\phi_{f})$ of $\pi_{1}^{-1}P^{*}\otimes 
\pi_{2}^{-1}P$.  Suppose also that there is a 
transformation $\psi:\pi_{1}^{-1}f\circ 
\pi_{3}^{-1}f\Rightarrow \pi_{2}^{-1}f$ 
as pictured in the following diagram: 
$$
\xymatrix{ 
\pi_{1}^{-1}\pi_{1}^{-1}P \ar[r]^-{\pi_{1}^{-1}
f} \ar @2{-}[dd] & \pi_{1}^{-1}\pi_{2}^{-1}P 
\ar @2{-}[r] & \pi_{3}^{-1}\pi_{1}^{-1}P 
\ar @2{->}[ddll]^-{\psi} \ar[dd]^-{\pi_{3}^{-1}f}   \\ 
& &                                                \\ 
\pi_{2}^{-1}\pi_{1}^{-1}P \ar[r]_-{\pi_{2}^{-1}
f} & \pi_{2}^{-1}\pi_{2}^{-1}P \ar @2{-}[r] & 
\pi_{3}^{-1}\pi_{2}^{-1}P.                          } 
$$
$\psi$ is required to satisfy the coherency 
condition that the following diagram of 
transformations commutes: 
$$
\xymatrix{ 
\pi_{3}^{-1}\pi_{3}^{-1}f\circ \pi_{1}^{-1}
\pi_{3}^{-1}f\circ \pi_{1}^{-1}\pi_{1}^{-1}
f \ar @2{->}[d]_-{1_{\pi_{3}^{-1}\pi_{3}^{-1}
f}\circ \pi_{1}^{-1}\psi} \ar @2{->}[r]^{\theta} & 
\pi_{4}^{-1}\pi_{3}^{-1}f\circ \pi_{4}^{-1}
\pi_{1}^{-1}f\circ \pi_{2}^{-1}\pi_{1}^{-1}
\phi \ar @2{->}[d]^-{\pi_{4}^{-1}\psi\circ 
1_{\pi_{2}^{-1}\pi_{1}^{-1}f}}                  \\ 
\pi_{3}^{-1}\pi_{3}^{-1}f\circ \pi_{1}^{-1}
\pi_{2}^{-1}f \ar @2{-}[d] & \pi_{4}^{-1}
\pi_{2}^{-1}f\circ \pi_{2}^{-1}\pi_{1}^{-1}
f \ar @2{-}[d]                                   \\ 
\pi_{3}^{-1}\pi_{3}^{-1}f\circ \pi_{3}^{-1}
\pi_{1}^{-1}f \ar @2{->}[d]_-{\pi_{3}^{-1}\psi} 
& \pi_{2}^{-1}\pi_{3}^{-1}f\circ \pi_{2}^{-1}
\pi_{1}^{-1}f \ar @2{->}[d]^-{\pi_{2}^{-1}\psi}  \\ 
\pi_{3}^{-1}\pi_{2}^{-1}f \ar @2{-}[r] & 
\pi_{2}^{-1}\pi_{2}^{-1}f.                         } 
$$
We will write this condition figuratively as 
$\pi_{1}^{-1}\psi\otimes \pi_{2}^{-1}\psi^{*}
\otimes \pi_{3}^{-1}\psi\otimes \pi_{4}^{-1}
\psi^{*} =1$.  Then there exists a bundle 
gerbe $(Q,Z,M)$ on $M$ and a stable morphism 
$\eta:P\to \pi^{-1}Q$ together with 
transformations $\xi:\pi_{2}^{-1}\eta\circ \phi
\Rightarrow \pi_{1}^{-1}\eta$ which are compatible 
with the transformations $\psi$ in the 
appropriate sense.  The converse is also true.  
\end{proposition} 

\begin{proof} 
Choose an open cover $\{U_{i}\}_{i\in I}$ of $M$ 
such that there exist local sections $s_{i}:U_{i}\to X$ 
of $\pi:X\to M$.  Form the pullback bundle gerbes 
$P_{i} = s_{i}^{-1}P$.  The stable morphism $f:
\pi_{1}^{-1}P\to \pi_{2}^{-1}P$ induces a stable 
morphism $f_{ij}:P_{j}\to P_{i}$ by pullback: 
$f_{ij} = (s_{i},s_{j})^{-1}f$, and the 
transformation $\psi$ induces a transformation 
$\psi_{ijk}:f_{ij}\circ f_{jk}\Rightarrow 
f_{ik}$ by pullback as well: $\psi_{ijk} = 
(s_{i},s_{j},s_{k})^{-1}\psi$.  The transformation 
$\psi_{ijk}$ satisfies the  
condition 
$$
\psi_{ijl}(1_{\phi_{ij}}\circ \psi_{jkl})\theta_{
f_{ij},f_{jk},f_{kl}} = 
\psi_{ikl}(\psi_{ijk}\circ 1_{\phi_{kl}}). 
$$
Therefore we can apply Proposition~\ref{prop:6.6.2} 
to find a bundle gerbe $(Q,Z,M)$ on $M$ and a 
stable morphism $g_{i}:Q|_{U_{i}}\to P_{i}$ and 
transformations $\xi_{ij}:\phi_{ij}\circ g_{j}
\Rightarrow g_{i}$ which are compatible with the 
transformations $\psi_{ijk}$ (the fact that the 
stable morphisms $f_{ij}$ go from $P_{j}$ to $P_{i}$ 
rather than from $P_{i}$ to $P_{j}$ is of no 
consequence).  We need to construct a stable 
morphism $\eta:\pi^{-1}Q\to P$.  Let $\hat{s}_{i}:
X_{i}\to X^{[2]}$ denote the map $x\mapsto (x,
s_{i}(\pi(x)))$.  $f$ induces a map $(s_{i}\circ 
\pi)^{-1}P\to P|_{X_{i}}$ by $\hat{s}_{i}^{-1}f$.  
Therefore we can define a stable morphism $\eta_{i}:(\pi^{-1}Q)|
_{X_{i}}\to P|_{X_{i}}$ by composition: 
$$
(\pi^{-1}Q)|_{X_{i}} = \pi^{-1}(Q|_{U_{i}}) 
\stackrel{\pi^{-1}g_{i}}{\to} \pi^{-1}s_{i}^{-1}
P \stackrel{\hat{s}_{i}^{-1}f}{\to} P|_{X_{i}}.  
$$
Next we define transformations $\theta_{ij}:\eta_{i}
\Rightarrow \eta_{j}$ by composition as indicated 
in the following diagram: 
$$
\xymatrix{ 
& & \pi^{-1}s_{i}^{-1}P \ar[ddrr]^-{\hat{s}_{i}
^{-1}f} & &                                   \\ 
& \ar @2{->}[dr]^-{\pi^{-1}\xi_{ij}^{-1}} & & &   \\ 
\pi^{-1}(Q|_{U_{ij}}) \ar[uurr]^-{\pi^{-1}g_{i}} 
\ar[ddrr]_-{\pi^{-1}g_{j}} & & \ar @2{->}[dr]^-{
\hat{s}_{ij}^{-1}\psi} & &     P|_{X_{ij}}         \\ 
& & & &                                            \\ 
& & \pi^{-1}s_{j}^{-1}P \ar[uurr]_-{\hat{s}_{j}^{-1}
f} \ar[uuuu]^-{\pi^{-1}f_{ij}}               } 
$$
ie $\theta_{ij} = (\hat{s}_{ij}^{-1}\psi\circ 
1_{\pi^{-1}g_{j}})(1_{\hat{s}_{i}^{-1}f}\circ 
\pi^{-1}\xi_{ij}^{-1})$.  Here $\hat{s}_{ij}:X_{ij}\to 
X^{[3]}$ is the map $x\mapsto (x,s_{i}(\pi(x)),s_{j}(
\pi(x)))$.  Because of the compatibility of $\xi_{ij}$ 
with $\psi_{ijk}$ we have $\theta_{jk}\theta_{ij} 
= \theta_{ik}$.  Hence we can glue the stable morphisms 
$\eta_{i}:(\pi^{-1}Q)|_{X_{i}}\to P|_{X_{i}}$ 
together to form a stable morphism $\eta:\pi^{-1}Q\to P$.  
We now need to define the transformations $\chi:f
\circ \pi_{1}^{-1}\eta\Rightarrow \pi_{2}^{-1}\eta$ 
and prove they are compatible with $\psi$.  It is 
sufficient to define transformations $\chi_{i}:f\circ 
\pi_{1}^{-1}\eta_{i}\Rightarrow \pi_{2}^{-1}\eta_{i}$ 
which are compatible with the gluing transformations 
$\theta_{ij}$.  $\chi_{i}$ is induced from $\psi$ by 
pullback with the map $X_{i}^{[2]}\to X^{[3]}$ 
which sends $(x_{1},x_{2})\in X_{i}^{[2]}$ to 
$(x_{1},x_{2},s_{i}(m))$ where $\pi(x_{1}) = 
\pi(x_{2}) = m$.  We try to indicate this in the 
following diagram: 
$$
\xymatrix{ 
\pi_{1}^{-1}\pi^{-1}(Q|_{U_{i}}) 
\ar @2{-}[d] \ar[rr]^-{\pi_{1}^{-1}
\pi^{-1}g_{i}} & \ar @2{-}[d] & 
\pi_{1}^{-1}\pi^{-1}s_{i}^{-1}P \ar @2{-}[d] 
\ar[rr]^-{\pi_{1}^{-1}\hat{s}_{i}^{-1}
\phi} & \ar @2{->}[d]^-{\chi_{i}} & 
\pi_{1}^{-1}P \ar[d]^-{f}            \\ 
\pi_{2}^{-1}\pi^{-1}(Q|_{U_{i}}) 
\ar[rr]_-{\pi_{2}^{-1}\pi^{-1}g_{i}} & & 
\pi_{2}^{-1}\pi^{-1}s_{i}^{-1}P 
\ar[rr]_-{\pi_{2}^{-1}\hat{s}_{i}
^{-1}f} & & \pi_{2}^{-1}P.    }
$$
One can check that the $\chi_{i}$ are 
compatible with the $\theta_{ij}$ and also 
that the glued together transformation 
$\chi$ is compatible with $\psi$.       
\end{proof}    

We conjecture that there is also  
a version of this proposition for 
gerbes: 
\begin{conjecture} 
\label{conj:6.6.4} 
Let $\pi:X\to M$ be a surjection admitting local 
sections.  Suppose there is a gerbe $\mathcal{G}$ 
on $X$ together with a morphism of gerbes 
$\phi:\pi_{1}^{*}\mathcal{G}\to \pi_{2}^{*}\mathcal{G}$ 
over $X^{[2]}$ which makes the following 
diagram commute: 
$$ 
\xymatrix{ 
\Delta^{*}\pi_{1}^{*}\mathcal{G} 
\ar[d]_-{\simeq} \ar[r]^-{\Delta^{*}\phi} & 
\Delta^{*}\pi_{2}^{*}\mathcal{G} 
\ar[d]^-{\simeq}                              \\ 
\mathcal{G} \ar[r]^-{\text{id}_{\mathcal{G}}} 
& \mathcal{G},                                    } 
$$ 
where the two vertical equivalences are induced 
by the equivalence of Lemma~\ref{lemma:6.3.2} 
and $\Delta:X\to X^{[2]}$ is the inclusion of the 
diagonal.  
Suppose also there is a natural transformation 
$\psi$ between morphisms of gerbes over $X^{[3]}$ 
as pictured in the following diagram: 
$$
\xymatrix{ 
\pi_{1}^{*}\pi_{1}^{*}\mathcal{G} 
\ar[r]^-{\pi_{1}^{*}\phi} \ar[d]_-{\simeq} 
& \pi_{1}^{*}\pi_{2}^{*}\mathcal{G} 
\ar[r]^-{\simeq} & (\pi_{2}\circ \pi_{1})^{-1})^{*}
\mathcal{G} \ar @2{-}[r] & (\pi_{1}\circ \pi_{3})^{*}
\mathcal{G} \ar @2{->}[dddlll]^-{\psi} \ar[d]^-{\simeq} \\ 
(\pi_{1}\circ \pi_{1})^{*}\mathcal{G} \ar @2{-}[d] & & & \pi_{3}^{*}
\pi_{1}^{*}\mathcal{G} \ar[d]^-{\pi_{3}^{*}\phi}          \\   
(\pi_{1}\circ \pi_{2})^{*}\mathcal{G} \ar[d]_-{\simeq} 
& & & \pi_{3}^{*}\pi_{2}^{*}\mathcal{G} \ar[d]^-{\simeq}  \\ 
\pi_{2}^{*}\pi_{1}^{*}\mathcal{G} \ar[r]_-{\pi_{2}^{*}\phi} & 
\pi_{2}^{*}\pi_{2}^{*}\mathcal{G} \ar[r]^-{\simeq} & 
(\pi_{2}\circ \pi_{2})^{*}\mathcal{G} \ar @2{-}[r] & 
(\pi_{2}\circ \pi_{3})^{*}\mathcal{G}.                       } 
$$
We also require the natural transformation 
$\psi$ to satisfy a coherency condition over 
$X^{[4]}$.  This is basically the non-abelian 
2-cocycle condition referred to in the last 
proposition but with the extra complication 
of having to insert equivalences of gerbes 
between double pullbacks and natural 
transformations between diagrams of such equivalences 
as in Lemma~\ref{lemma:6.3.2}.  
Observe that we may pullback the natural 
transformations $\psi$ to form new natural 
transformations $\pi_{1}^{*}\psi$, $\pi_{2}^{*}\psi$, 
and so on between the appropriate morphisms 
of gerbes.  We can then form the diagrams 
$$ 
\xymatrix{ 
\pi_{1}^{*}\pi_{1}^{*}\pi_{2}^{*}\mathcal{G} 
\ar[r]^-{\simeq} & \pi_{1}^{*}\pi_{3}^{*}\pi_{1}^{*}
\mathcal{G} \ar[rr]^-{\pi_{1}^{*}\pi_{3}^{*}\phi} 
\ar @2{->}[dr]^-{\pi_{1}^{*}\psi} & & 
\pi_{1}^{*}\pi_{3}^{*}\pi_{2}^{*}\mathcal{G} \ar[d]^-{\simeq} \\ 
& & \pi_{1}^{*}\pi_{2}^{*}\pi_{2}^{*}\mathcal{G} 
\ar[d]_-{\simeq} \ar @2{->}[r] 
\ar[ur]^-{\simeq} & \pi_{3}^{*}\pi_{3}^{*}\pi_{1}^{*}
\mathcal{G} \ar[ddd]^-{\pi_{3}^{*}\pi_{3}^{*}\phi}              \\ 
& \pi_{1}^{*}\pi_{2}^{*}\pi_{1}^{*}\mathcal{G} \ar[d]_{\simeq} 
\ar[ur]^-{\pi_{1}^{*}\pi_{2}^{*}\phi} \ar @2{-}[r] & \pi_{3}^{*}
\pi_{1}^{*}\pi_{2}^{*}\mathcal{G} \ar[ur]^-{\simeq} \ar @2{->}[dr]^{\pi_{3}
^{*}\psi} &  \\ 
\pi_{1}^{*}\pi_{1}^{*}\pi_{1}^{*}\mathcal{G} \ar @2{->}[r]  
\ar[uuu]^-{\pi_{1}^{*}\pi_{1}^{*}\phi} \ar[ur]^-{\simeq} & 
\pi_{3}^{*}\pi_{1}^{*}\pi_{1}^{*}\mathcal{G} 
\ar[ur]^-{\pi_{3}^{*}\pi_{1}^{*}\phi} & &                       \\   
\pi_{3}^{*}\pi_{2}^{*}\pi_{1}^{*}\mathcal{G} 
\ar[u]^-{\simeq} \ar[rr]_-{\pi_{3}^{*}\pi_{2}^{*}\phi} \ar[ur]^-{\simeq} & & 
\pi_{3}^{*}\pi_{2}^{*}\pi_{2}^{*}\mathcal{G} \ar[r]_-{\simeq} 
& \pi_{3}^{*}\pi_{3}^{*}\pi_{2}^{*}\mathcal{G}                   } 
$$
and 
$$
\xymatrix{ 
\pi_{2}^{*}\pi_{1}^{*}\pi_{2}^{*}\mathcal{G} 
\ar[r]^-{\simeq} \ar[dr]_-{\simeq} & \pi_{4}^{*}
\pi_{1}^{*}\pi_{1}^{*}\mathcal{G} \ar[dr]_-{\simeq} 
\ar @2{->}[d] \ar[rrr]^-{\pi_{4}^{*}\pi_{1}^{*}\phi} & & &  
\pi_{4}^{*}\pi_{1}^{*}\pi_{2}^{*}\mathcal{G} \ar[d]^-{\simeq} \\ 
& \pi_{2}^{*}\pi_{3}^{*}\pi_{1}^{*}\mathcal{G} 
\ar[r]^-{\simeq} \ar[dr]_{\pi_{2}^{*}\pi_{3}^{*}\phi} & 
\pi_{4}^{*}\pi_{2}^{*}\pi_{1}^{*}\mathcal{G} \ar @2{-}[d] 
\ar[dr]^-{\pi_{4}^{*}\pi_{2}^{*}\phi} & & \pi_{4}^{*}
\pi_{3}^{*}\pi_{1}^{*}\mathcal{G} \ar[dd]^-{\pi_{4}^{*}
\pi_{3}^{*}\phi} \ar @2{->}[dl]^-{\pi_{4}^{*}\psi}                \\ 
\pi_{2}^{*}\pi_{1}^{*}\pi_{1}^{*}
\mathcal{G} \ar[uu]^-{\pi_{2}^{*}\pi_{1}^{*}
\phi}& & \pi_{2}^{*}\pi_{3}^{*}\pi_{2}^{*}\mathcal{G} 
\ar @2{->}[dl]^-{\pi_{2}^{*}\psi} \ar[r]^-{\simeq} 
\ar[dr]^-{\simeq} & \pi_{4}^{*}\pi_{2}^{*}\pi_{2}^{*}
\mathcal{G} \ar @2{->}[d] \ar[dr]^-{\simeq} &                       \\ 
\pi_{2}^{*}\pi_{2}^{*}\pi_{1}^{*}\mathcal{G} \ar[u]^-{\simeq} 
\ar[rrr]_{\pi_{2}^{*}\pi_{2}^{*}\phi} & & & \pi_{2}^{*}
\pi_{2}^{*}\pi_{2}^{*}\mathcal{G} \ar[r]_{\simeq} & 
\pi_{4}^{*}\pi_{3}^{*}\pi_{2}^{*}\mathcal{G}.                     } 
$$ 
Paste together the natural transformations in the first 
of the above two diagrams to yield a single natural 
transformation between the morphisms of gerbes bounding the 
diagram.  Finally we can state the coherency condition 
on the natural transformation $\psi$: the following 
diagram of natural transformations, after being pasted together, 
is the same natural transformation as the pasted-together 
natural transformation in the second diagram above: 
$$ 
\xymatrix{ 
\pi_{2}^{*}\pi_{1}^{*}\pi_{2}^{*}\mathcal{G} 
\ar[dr]^-{\simeq} \ar[rr]^-{\simeq} & & 
\pi_{4}^{*}\pi_{1}^{*}\pi_{1}^{*}\mathcal{G} 
\ar[rr]^-{\pi_{4}^{*}\pi_{1}^{*}\phi} \ar[d]_-{\simeq} 
\ar @2{->}[dl] & \ar @2{-}[dl] & \pi_{4}^{*}
\pi_{1}^{*}\pi_{2}^{*}\mathcal{G} \ar[dd]^-{\simeq} 
\ar[dl]^-{\simeq}                                      \\ 
& \pi_{1}^{*}\pi_{1}^{*}\pi_{2}^{*}\mathcal{G} 
\ar[r]^-{\simeq} \ar @2{-}[dl] & \pi_{1}^{*}\pi_{3}^{*}
\pi_{1}^{*}\mathcal{G} \ar[r] & \pi_{1}^{*}\pi_{3}^{*}
\pi_{2}^{*}\mathcal{G} \ar[d]^-{\simeq} \ar @2{->}[ddll] &  \\ 
\pi_{2}^{*}\pi_{1}^{*}\pi_{1}^{*}\mathcal{G} 
\ar[uu]_{\pi_{2}^{*}\pi_{1}^{*}\phi} \ar[r]^-{\simeq} & 
\pi_{1}^{*}\pi_{1}^{*}\pi_{1}^{*}\mathcal{G} \ar[u]  
\ar @2{->}[dl] &     & \pi_{3}^{*}\pi_{3}^{*}\pi_{1}^{*}
\mathcal{G} \ar[d] & \pi_{4}^{*}\pi_{3}^{*}\pi_{1}^{*}
\mathcal{G} \ar[l]^-{\simeq} \ar[dd]^-{\pi_{4}^{*}
\pi_{3}^{*}\phi} \ar @2{-}[dl]                           \\ 
& \pi_{3}^{*}\pi_{2}^{*}\pi_{1}^{*}\mathcal{G} 
\ar[u]^-{\simeq} \ar[r] & \pi_{3}^{*}\pi_{2}^{*}
\pi_{2}^{*}\mathcal{G} \ar[r]^-{\simeq} \ar @2{-}[dl] & 
\pi_{3}^{*}\pi_{3}^{*}\pi_{2}^{*}\mathcal{G} \ar @2{->}[dl] & \\ 
\pi_{2}^{*}\pi_{2}^{*}\pi_{1}^{*}\mathcal{G} 
\ar[uu]^-{\simeq} \ar[ur]^-{\simeq} \ar[rr]_-{\pi_{2}^{*}
\pi_{2}^{*}\phi} & & \pi_{2}^{*}\pi_{2}^{*}\pi_{2}^{*}
\mathcal{G} \ar[u]^-{\simeq} \ar[rr]^-{\simeq} & & 
\pi_{4}^{*}\pi_{3}^{*}\pi_{2}^{*}\mathcal{G} \ar[ul]^-{\simeq}.  }  
$$
Then there exists a gerbe $\mathcal{D}$ on $M$ 
together with an equivalence of gerbes $\chi:\pi^{*}
\mathcal{D}\to \mathcal{G}$ and natural transformations 
$\xi$ between the morphisms of gerbes as pictured in 
the following diagram: 
$$
\xymatrix{ 
\pi_{1}^{*}\mathcal{G} \ar[r]^-{\phi} & 
\pi_{2}^{*}\mathcal{G}                      \\ 
\pi_{1}^{*}\pi^{*}\mathcal{D} \ar[u]^-{\pi_{1}^{*}\chi} 
\ar @2{->}[r]^{\xi} & \pi_{2}^{*}\pi^{*}\mathcal{D} 
\ar[u]^-{\pi_{2}^{*}\chi}                     \\ 
(\pi\circ \pi_{1})^{*}\mathcal{D} \ar[u]^-{\simeq} 
\ar @2{-}[r] & (\pi\circ \pi_{2})^{*}\mathcal{D} 
\ar[u]^-{\simeq}.                                 } 
$$ 
We then require that the natural transformations 
$\pi_{1}^{*}\xi$, $\pi_{2}^{*}\xi$ and $\pi_{3}^{*}\xi$ 
are compatible with $\psi$ in the sense that using the 
natural transformations $\pi_{1}^{*}\xi$, $\pi_{2}^{*}\xi$ 
and $\pi_{3}^{*}\xi$ and the natural transformations 
$\theta_{f,g,h}$ of Lemma~\ref{lemma:6.3.2} 
we can construct a natural transformation between the 
morphisms of gerbes bounding diagram (1) above.  The compatibility 
condition is that this natural transformation and $\psi$ 
coincide.  The converse is also true.          
\end{conjecture} 

\begin{note} 
\label{note:6.6.6} 
\begin{enumerate} 
\item One should be able to prove this conjecture 
in a similar manner to the proof of Proposition~\ref{prop:6.6.3} 
but with the additional complication of 
having to insert the equivalences $\phi_{f,g}$ 
and natural transformations $\theta_{f,g,h}$ 
of Lemma~\ref{lemma:6.3.2}.  
\item If the equivalences $\phi_{f,g}:g^{*}f^{*}
\mathcal{G}\to (fg)^{*}\mathcal{G}$ of Lemma~\ref{lemma:6.3.2} 
were the 
identity, we could write the coherency condition 
on the natural transformation $\psi$ more 
succinctly as 
$$
\pi_{3}^{*}\psi(1_{\pi_{3}^{*}\pi_{3}^{*}
\phi}\circ \pi_{1}^{*}\psi) = \pi_{2}^{*}\psi 
(\pi_{4}^{*}\psi\circ 1_{\pi_{2}^{*}\pi_{1}
^{*}\phi}). 
$$
Unfortunately we have to cope with having 
the equivalences $\phi_{f,g}$ and natural 
transformations $\theta_{f,g,h}$ and so, 
in order to introduce a new notation to 
avoid the diagrams above, we will abbreviate 
the coherency condition on $\psi$ to $\d(\psi) =1$.  
We will encounter such diagrams again in 
Chapter~\ref{chapter:13}.   
\end{enumerate} 
\end{note}

\setcounter{chapter}{6}
\chapter{Definition of a bundle 2-gerbe} 
\label{chapter:7} 
\section{Simplicial bundle gerbes and bundle 2-gerbes}
\label{sec:7.1} 
Recall in Chapter~\ref{chapter:3}, 
Section~\ref{sec:3.2}, 
that we saw how 
bundle gerbes $(P,X,M)$ could be 
regarded as simplicial line bundles 
on the simplicial manifold 
$X^{[p]}$.  We will use this 
observation to motivate our 
definition of a bundle two gerbe.  

Our first step is to see if there 
is a workable notion of a 
\emph{simplicial bundle gerbe}.  
Note that Brylinski and 
McLaughlin define a notion of 
\emph{simplicial gerbe} (see 
\cite{BryMcL1}, \cite{BryMcL2} 
and also Definition~\ref{def:13.1.1})  
so this will also serve to guide 
us.  We will work on the 
simplicial manifold $X=\{X_{p}\}$ 
with $X_{p} = X^{[p+1]}$ associated to a surjective map 
$\pi:X\to M$ admitting local 
sections.  Recall, see 
Definition~\ref{def:3.2.2} 
and also \cite{BryMcL1}, 
that a simplicial line bundle 
on a simplicial manifold 
$X = \{X_{p}\}$ consists of a line 
bundle $L$ on $X_{1}$, together 
with a non-vanishing section $s$ 
of $\d(L)$ over $X_{2}$ which 
satisfies $\d(s) = 1$ on $X_{3}$.  
We attempt to generalise this to 
the case of bundle gerbes as 
follows.  Let $(Q,Y,X^{[2]})$ 
be a bundle gerbe on $X^{[2]}$.  
We suppose there is a bundle gerbe 
morphism 
$$
\bar{m}:\pi_{1}^{-1}Q\otimes \pi_{3}^{-1}Q 
\to \pi_{2}^{-1}Q 
$$
over $X^{[3]}$.  We think of $\bar{m}$ 
as defining a `product' on the bundle gerbe $Q$.      
Next observe that we have 
$$
\pi_{1}^{-1}(\pi_{1}^{-1}Q\otimes 
\pi_{3}^{-1}Q)\otimes \pi_{3}^{-1}
\pi_{3}^{-1}Q = \pi_{2}^{-1}\pi_{1}^{-1}
Q\otimes\pi_{4}^{-1}(\pi_{1}^{-1}Q
\otimes \pi_{3}^{-1}Q), 
$$
by virtue of the simplicial identities.  
The crucial observation is that 
$\bar{m}$ gives rise to two bundle gerbe 
morphisms over $X^{[4]}$ as indicated in the 
following diagram. 
$$  
\xymatrix{ 
\pi_{1}^{-1}(\pi_{1}^{-1}Q\otimes \pi_{3}^{-1}Q)
\otimes \pi_{3}^{-1}\pi_{3}^{-1}Q \ar @2{-}[r] 
\ar[d]_-{\pi_{1}^{-1}\bar{m}\otimes 1} & 
\pi_{2}^{-1}\pi_{1}^{-1}Q\otimes \pi_{4}^{-1}(
\pi_{1}^{-1}Q\otimes \pi_{3}^{-1}Q) \ar[d]^-{
1\otimes \pi_{4}^{-1}\bar{m}}                     \\ 
\pi_{1}^{-1}\pi_{2}^{-1}Q\otimes \pi_{3}^{-1}
\pi_{3}^{-1}Q \ar @2{-}[d] & \pi_{2}^{-1}\pi_{1}^{-1}
Q\otimes \pi_{4}^{-1}\pi_{2}^{-1}Q \ar @2{-}[d]     \\   
\pi_{3}^{-1}(\pi_{1}^{-1}Q\otimes \pi_{3}^{-1}Q) 
\ar[d]_-{\pi_{3}^{-1}\bar{m}} & \pi_{2}^{-1}(
\pi_{1}^{-1}Q\otimes \pi_{3}^{-1}Q) \ar[d]^-{\pi_{2}^{-1}\bar{m}} \\ 
\pi_{3}^{-1}\pi_{2}^{-1}Q \ar @2{-}[r] & 
\pi_{2}^{-1}\pi_{2}^{-1}Q.                          } 
$$
We define bundle gerbe morphisms 
$\bar{m}_{1}$ and $\bar{m}_{2}$ as follows.  
\begin{eqnarray*}  
\bar{m}_{1} & = & \pi_{3}^{-1}\bar{m}\circ 
(\pi_{1}^{-1}\bar{m}\otimes 1)           \\ 
\bar{m}_{2} & = & \pi_{2}^{-1}\bar{m}\circ 
(1\otimes \pi_{4}^{-1}\bar{m}).            
\end{eqnarray*} 
At this point it is convenient to 
introduce a different notation which 
avoid cumbersome diagrams such as the one 
above.  We will denote $\pi_{1}^{-1}Y$ 
by $Y_{23}$, $\pi_{2}^{-1}Y$ by $Y_{13}$, 
$\pi_{3}^{-1}Y$ by $Y_{12}$, $\pi_{1}^{-1}
\pi_{3}^{-1}Y$ by $Y_{24}$ and so on.  
Similarly we will denote the $\cstar$ 
bundle $\pi_{1}^{-1}Q$ on $Y_{23}^{[2]}$ 
by $Q_{23}$ and so forth.  More generally 
we shall denote $\pi_{1}^{-1}Y\times_{X^{[3]}}
\pi_{3}^{-1}Y = Y_{23}\times_{X^{[3]}}Y_{12}$ 
by $Y_{123}$, $\pi_{1}^{-1}Q\otimes \pi_{3}^{-1}Q$ 
by $Q_{123}$ and similarly for other combinations.  
Hence $\bar{m}:Q_{123}\to Q_{13}$ and 
$\bar{m}_{1}, \bar{m}_{2}:Q_{1234}\to Q_{14}$.  
Therefore the diagram above becomes in the 
new notation 
$$
\xymatrix{ 
Q_{1234} \ar @2{-}[r] \ar [d]_{\bar{m}\otimes 
1} & Q_{1234} \ar[d]^{1
\otimes \bar{m}}                               \\ 
Q_{24}\otimes Q_{12} \ar[d]_{\bar{m}} & 
Q_{34}\otimes Q_{13} \ar[d]^{\bar{m}}    \\ 
Q_{14} \ar @2{-}[r] & Q_{14} } 
$$
Rather than demand the two morphisms 
$\bar{m}_{1}$ and $\bar{m}_{2}$ are equal, we settle 
for the weaker condition that there is 
a transformation $a:\bar{m}_{1}\Rightarrow \bar{m}_{2}$ 
of bundle gerbe morphisms.  Recall, see 
Definition~\ref{def:3.4.9}, that this 
means that $a$ is a section of the 
descended $\cstar$ bundle 
$D_{\bar{m}_{1},\bar{m}_{2}}$ on 
$X^{[4]}$ 
of Lemma~\ref{lemma:3.4.6}.  The 
transformation $a$ gives, in some sense, 
a measure of the `non-associativity' 
of the `product' $\bar{m}$.      
Accordingly, we will call $a$ the 
`associator' transformation or 
`associator' section and we will 
call the (trivial) $\cstar$ bundle 
$A = D_{\bar{m}_{1},\bar{m}_{2}}$ on $X^{[4]}$ 
the `associator' bundle.  

The final condition is a coherency 
condition on $a$.  To write it down 
first observe that $Q$ induces a 
bundle gerbe $(Q_{12345},Y_{12345},X^{[5]})$ 
over $X^{[5]}$.    
The bundle gerbe morphism 
$\bar{m}$ induces five bundle gerbe 
morphisms $\bar{M}_{i}:Q_{12345}\to Q_{15}$ 
given as follows. 
\begin{eqnarray*} 
& \bar{M}_{1} = \bar{m}\circ (\bar{m}\otimes 1)\circ 
(\bar{m}\otimes 1\otimes 1) &                           \\ 
& \bar{M}_{2} = \bar{m}\circ (\bar{m}\otimes 1)\circ 
(1\otimes \bar{m}\otimes 1) &                            \\ 
& \bar{M}_{3} = \bar{m}\circ (1\otimes \bar{m})\circ 
(1\otimes \bar{m}\otimes 1) &                            \\ 
& \bar{M}_{4} = \bar{m}\circ (1\otimes \bar{m})\circ 
(1\otimes 1\otimes \bar{m}) &                            \\ 
& \bar{M}_{5} = \bar{m}\circ (\bar{m}\otimes 1)\circ 
(1\otimes 1\otimes \bar{m}) &      
\end{eqnarray*} 
Notice that $\bar{M}_{5}$ can also be written as 
$$
\bar{M}_{5} = \bar{m}\circ (1\otimes \bar{m})
\circ (\bar{m}\otimes 1\otimes 1). 
$$
It is not too hard to see that we have the 
following isomorphisms of $\cstar$ bundles 
on $X^{[5]}$.  We have 
$D_{\bar{M}_{1},\bar{M}_{2}}\simeq \pi_{1}^{-1}A$, 
$D_{\bar{M}_{2},\bar{M}_{3}}\simeq \pi_{3}^{-1}A$, 
$D_{\bar{M}_{3},\bar{M}_{4}}\simeq \pi_{5}^{-1}A$, 
$D_{\bar{M}_{4},\bar{M}_{5}}\simeq \pi_{2}^{-1}A^{*}$ 
and $D_{\bar{M}_{5},\bar{M}_{1}}\simeq \pi_{4}^{-1}A^{*}$.  
Since $D_{\bar{M}_{1},\bar{M}_{1}}$ is canonically 
trivialised, and the bundle gerbe product 
on $Q$ gives an isomorphism 
$$
D_{\bar{M}_{1},\bar{M}_{2}}\otimes D_{\bar{M}_{2},\bar{M}_{3}}\otimes 
D_{\bar{M}_{3},\bar{M}_{4}}\otimes D_{\bar{M}_{4},\bar{M}_{5}}\otimes 
D_{\bar{M}_{5},\bar{M}_{1}}\simeq D_{\bar{M}_{1},\bar{M}_{1}},
$$
we have that $\d(A)$ on $X^{[5]}$ is 
canonically trivialised.  The final condition that 
we require is that the section $\d(a)$ 
of $\d(A)$ matches this canonical 
trivialisation.  

Notice that the thing that makes all this 
possible is the fact that the $\pi_{i}$ 
are the face operators for the simplicial 
manifold $X_{p} = X^{[p+1]}$.  Therefore all that we 
have said above applies equally well to 
an arbitrary simplicial manifold 
$X = \{X_{p}\}$.  
Thus we make the following definition. 

\begin{definition} 
\label{def:7.1.1} 
A \emph{simplicial bundle gerbe} on a simplicial 
manifold $X = \{X_{p}\}$ consists of the 
following data. 
\begin{enumerate} 
\item A bundle gerbe $(Q,Y,X_{1})$ on 
$X_{1}$.  
\item A bundle gerbe morphism 
$$
\bar{m}:d_{0}^{-1}Q\otimes d_{2}^{-1}Q 
\to d_{1}^{-1}Q 
$$
over $X_{2}$.  
\item A transformation 
$a:\bar{m}_{1}\Rightarrow \bar{m}_{2}$ 
between the two induced bundle gerbe 
morphisms $\bar{m}_{1}$ and 
$\bar{m}_{2}$ over $X_{3}$.  $\bar{m}_{1}$ 
and $\bar{m}_{2}$ are defined in the 
following diagram.  
$$
\xymatrix{ 
d_{0}^{-1}(d_{0}^{-1}Q\otimes d_{2}^{-1}Q)\otimes 
d_{2}^{-1}d_{2}^{-1}Q \ar[d]_-{d_{0}^{-1}\bar{m}\otimes 
1} \ar @2{-}[r] & d_{1}^{-1}d_{0}^{-1}Q\otimes 
d_{3}^{-1}(d_{0}^{-1}Q\otimes d_{2}^{-1}Q) 
\ar[d]^-{1\otimes d_{3}^{-1}\bar{m}}             \\ 
d_{0}^{-1}d_{1}^{-1}Q\otimes d_{2}^{-1}d_{2}^{-1}Q 
\ar @2{-}[d] \ar @<-4.2ex> @2{->}[r]^-{a} & 
d_{1}^{-1}d_{0}^{-1}Q\otimes d_{3}^{-1}d_{1}^{-1}Q 
\ar @2{-}[d]                                      \\ 
d_{2}^{-1}(d_{0}^{-1}Q\otimes d_{2}^{-1}Q) \ar[d]_-{
d_{2}^{-1}\bar{m}} & d_{1}^{-1}(d_{0}^{-1}Q\otimes 
d_{2}^{-1}Q) \ar[d]^-{d_{1}^{-1}\bar{m}}            \\ 
d_{2}^{-1}d_{1}^{-1}Q \ar @2{-}[r] & d_{1}^{-1}d_{1}^{-1}Q.  } 
$$
So $\bar{m}_{1} = d_{2}^{-1}\bar{m}\circ 
(d_{0}^{-1}\bar{m}\otimes 1)$ and $\bar{m}_{2} = 
d_{1}^{-1}\bar{m}\circ (1\otimes d_{3}^{-1}\bar{m})$. 
Thus $a$ 
is a section of the $\cstar$ bundle 
$A = D_{\bar{m}_{1},\bar{m}_{2}}$ over $X_{3}$.  
\item The transformation $a$ satisfies the 
coherency condition 
$$
d_{0}^{-1}a\otimes d_{1}^{-1}a^{*}\otimes 
d_{2}^{-1}a\otimes d_{3}^{-1}a^{*}\otimes 
d_{4}^{-1}a = 1,
$$
where $1$ is the canonical section of the 
$\cstar$ bundle $\d(A)$ over $X_{4}$. 
\end{enumerate} 
\end{definition} 

This definition should be compared with 
the definition of simplicial gerbe given 
in \cite{BryMcL1} and \cite{BryMcL2}.  
With this definition in hand, we now 
turn to the problem of defining a bundle 
2-gerbe.  We make two definitions, the 
first of which includes the following one  
as a special case.   

\begin{definition} 
\label{def:7.1.2} 
A \emph{bundle 2-gerbe} is a quadruple 
of manifolds $(Q,Y,X,M)$ where $\pi:X\to M$ 
is a smooth surjection admitting local sections, 
and $(Q,Y,X^{[2]})$ is a simplicial bundle 
gerbe on the simplicial manifold $X=\{X_{p}\}$ 
associated to $\pi:X\to M$ with $X_{p} = X^{[p+1]}$.  
\end{definition} 

\begin{definition} 
\label{def:7.1.4} 
A \emph{strict bundle 2-gerbe} is a 
quadruple of manifolds $(Q,Y,X,M)$ where  
$\pi:X\to M$ is a smooth surjection 
admitting local sections and where $(Q,Y,X^{[2]})$ 
is a bundle gerbe on $X^{[2]}$ such that 
there is a bundle gerbe morphism 
$\bar{m}:\pi_{1}^{-1}Q\otimes \pi_{3}^{-1}Q\to 
\pi_{2}^{-1}Q$ such that the two induced 
bundle gerbe morphisms $\bar{m}_{1}$ and 
$\bar{m}_{2}$ over $X^{[4]}$ are equal.  
\end{definition} 
  
The following diagram shows what a bundle 
2-gerbe looks like. 
$$
\begin{array}{ccccc}
Q  &       &     &     &          \\
\downarrow &   &    &     &           \\
Y^{[2]} & \dra       & Y   &    &        \\
  &       &  \downarrow &     &          \\
 &     &   X^{[2]}  & \dra     &  X         \\
   &     &      &       &   \downarrow  \\
    &    &       &       &   M.         \\
\end{array}
$$
For future reference, we will record 
here some properties of the associator 
transformation $a$ and its lift $\hat{a}$.  
So $\hat{a}$ is a section of the 
$\cstar$ bundle $\hat{A} = (m_{1},m_{2})^{-1}Q$ 
on $Y_{1234}$.  Since $\hat{a}$ 
descends to the section $a$ of $A$, it 
must commute with the descent isomorphism 
of $\hat{A}$.  This implies that 
$\hat{a}$ enjoys the following 
property: if $u_{34}\in Q_{(y_{34},y_{34}^{'})}$, 
$u_{23} \in Q_{(y_{23},y_{23}^{'})}$ and 
$u_{12} \in Q_{(y_{12},y_{12}^{'})}$ where 
$y_{34},y_{34}^{'} \in Y_{(x_{3},x_{4})}$, 
$y_{23},y_{23}^{'} \in Y_{(x_{2},x_{3})}$ 
and $y_{12},y_{12}^{'} \in Y_{(x_{1},x_{2})}$, 
then 
\begin{eqnarray*} 
&   &  m_{Q}(\hat{m}_{2}(u_{34}\otimes 
u_{23}\otimes u_{12})\otimes \hat{a}(y_{34},
y_{23},y_{12}))                                  \\ 
& = & m_{Q}(\hat{a}(y_{34}^{'},y_{23}^{'},
y_{12}^{'})\otimes \hat{m}_{1}(u_{34}\otimes 
u_{23}\otimes u_{12})).                           
\end{eqnarray*} 
The coherency condition on $a$ is a 
consequence of the following 
property of $\hat{a}$: suppose 
$y_{45} \in Y_{(x_4,x_5)}$, 
$y_{34} \in Y_{(x_3,x_4)}$, 
$y_{23} \in Y_{(x_2,x_3)}$ and 
$y_{12} \in Y_{(x_1,x_2)}$.  Then 
we have the following equation 
(we have denoted the bundle gerbe 
product $m_{Q}$ by $\cdot$ --- since 
the bundle gerbe product is associative there is no danger of 
confusion) 
\begin{multline*} 
\hat{a}(m(y_{45},y_{34}),y_{23},y_{12})^{-1}\cdot 
\hat{a}(y_{45},y_{34},m(y_{23},y_{12}))^{-1}\cdot 
\hat{m}(e(y_{45})\otimes \hat{a}(y_{34},y_{23},
y_{12}))                                          \\ 
\cdot \hat{a}(y_{45},m(y_{34},y_{23}),y_{12})\cdot 
\hat{m}(\hat{a}(y_{45},y_{34},y_{23})\otimes 
e(y_{12}))                                         \\ 
= e(m(m(m(y_{45},y_{34}),y_{23}),y_{12}))       
\end{multline*} 

\section{Stable bundle 2-gerbes} 
\label{sec:7.2} 

We also define a notion of bundle 
2-gerbe with bundle gerbe morphisms 
replaced by stable morphisms.  It 
will turn out that this definition 
of a `stable bundle 2-gerbe' will 
turn out to be the most versatile, 
for example in Chapter~\ref{chapter:13} 
we will show that a stable bundle 
2-gerbe on $M$ gives rise to a 2-gerbe on $M$.   

\begin{definition} 
\label{def:7.2.1} 
A \emph{stable bundle 2-gerbe} consists of a 
quadruple of smooth manifolds $(Q,Y,X,M)$ 
where $\pi:X\to M$ is a surjection admitting 
local sections and where $(Q,Y,X^{[2]})$ is a 
bundle gerbe on $X^{[2]}$.  We suppose there is a 
stable morphism $m:\pi_{1}^{-1}Q\otimes \pi_{3}^{-1}Q
\to \pi_{2}^{-1}Q$ corresponding to a trivialisation 
$\d(L_{m})$ of $\pi_{1}^{-1}Q\otimes\pi_{2}^{-1}Q^{*}
\otimes \pi_{3}^{-1}Q$.  We also suppose there is a 
transformation $a:m\circ (m\otimes 1)\Rightarrow 
m\circ (1\otimes m)$ as pictured in the following 
diagram: 
$$
\xymatrix{ 
Q_{1234} \ar[r]^-{m\otimes 1} 
\ar[d]_-{1\otimes m} & Q_{124} 
\ar[d]^-{m} \ar @2{->}[dl]^-{a}  \\ 
Q_{134} \ar[r]_-{m} & Q_{14}.    } 
$$
$a$ is required to satisfy the coherency 
condition that the following two diagrams of 
transformations are equal:   
$$
\xymatrix{   
& Q_{12345} \ar[dl] \ar[dr] \ar[rr] & & 
Q_{1235} \ar[dr] \ar @2{->}[dl]^-{1\otimes 
a_{2345}} &                                \\ 
Q_{1345} \ar[dr] & & Q_{1245} \ar[dl] 
\ar @2{->}[ll]_-{a_{1234}\otimes 1} \ar[rr] 
& & Q_{125} \ar[dl] \ar @2{->}[dlll]^-{a_{1245}} \\ 
& Q_{145} \ar[rr] & \ar @2{-}[d] & Q_{15} &     \\ 
& Q_{12345} \ar[dl] \ar[rr] &  & Q_{1235} 
\ar[dl] \ar[dr] \ar @2{-}[dlll] &                 \\ 
Q_{1345} \ar[dr] \ar[rr] & & Q_{135} \ar[dr] 
\ar @2{->}[dl]_-{a_{1345}} & & Q_{125} \ar @2{->}[ll]^-
{a_{1235}} \ar[dl]                                \\ 
& Q_{145} \ar[rr] & & Q_{15}.                     }   
$$
\end{definition} 
 
\begin{note} 
\label{note:7.2.2} 
\begin{enumerate} 
\item Notice that by (1) of Note~\ref{note:3.4.13}, 
every bundle 2-gerbe on $M$ gives 
rise to a stable bundle 2-gerbe on $M$.  

\item As in Section~\ref{sec:7.1} 
let $m_{1}:Q_{1234}\to Q_{14}$ be the 
stable morphism $m\circ (m\otimes 1)$ 
and let $m_{2}:Q_{1234}\to Q_{14}$ be the 
stable morphism $m\circ (1\otimes m)$.  
Then the transformation $a:m_{1}\Rightarrow 
m_{2}$ is a section of the $\cstar$ 
bundle $A$ on $X^{[4]}$ such that 
$L_{m_{1}} \simeq L_{m_{2}}\otimes \pi_{1234}^{-1}A$ 
on $Y_{1234}$.  As before, the $\cstar$ 
bundle $\d(A)$ on $X^{[4]}$ is canonically 
trivialised and the coherency condition 
on $a$ is that the induced section 
$\d(a)$ of $\d(A)$ matches this canonical 
trivialisation.   
\end{enumerate} 
\end{note}

\setcounter{chapter}{7}
\chapter{Relationship of bundle 2-gerbes with 
bicategories} 
\label{chapter:10} 
\section{Bicategories}
\label{sec:10.1} 
We first recall the definition of a 
bicategory (see for instance \cite{Ben})   
\begin{definition}[\cite{Ben}] 
\label{def:10.1.1} 
A \emph{bicategory} $\B$ consists of the following data: 
\begin{enumerate} 
\item  A set $Ob(\B)$ whose elements are 
called the \emph{objects} of $\B$.  

\item  Given two objects $A$ and $B$ of $\B$ 
there is a category $\text{Hom}(A,B)$.  
The objects of $\text{Hom}(A,B)$ 
are called \emph{1-arrows} of $\B$ while the 
arrows of $\text{Hom}(A,B)$ are called 
\emph{2-arrows} of $\B$.  If $\a$ and $\b$ are 
1-arrows of $\text{Hom}(A,B)$ and 
$\phi$ is a 2-arrow with source $\a$ and 
target $\b$ then this is denoted by  
$\phi:\a \Rightarrow \b$.  The following diagram gives  
a useful way to picture 2-arrows: 
$$ 
\diagram 
A \rtwocell^{\a}_{\b}{\phi} & B \\ 
\enddiagram 
$$

\item  Given three objects $A,B$ and 
$C$ of $\B$, there is a functor, called 
the \emph{composition} functor, 
$$
m(A,B,C):\text{Hom}(B,C)\times 
\text{Hom}(A,B) \to \text{Hom}(A,C).
$$
If $\a$ and $\b$ are 1-arrows in 
$\text{Hom}(B,C)$ and $\text{Hom}(A,B)$ 
respectively then we usually write $m(A,B,C)(\a,\b)$ 
as $\a \circ \b$ and  
$m(A,B,C)(\phi,\psi)$ is usually  
written as $\phi \circ \psi$ where  
$\phi : \a \Rightarrow \a^{\prime}$ and 
$\psi : \b \Rightarrow \b^{\prime}$ are 2-arrows 
in $\text{Hom}(B,C)$ and $\text{Hom}(A,B)$ 
respectively.   

\item  For each object $A$ of $\B$ there is a 
1-arrow $\text{Id}_{A}$ in $\text{Hom}(A,A)$ called the 
\emph{identity 1-arrow} of $A$.  The identity map 
of $\text{Id}_{A}$ in $\text{Hom}(A,A)$ is denoted by 
$\text{id}_{A}:\text{Id}_{A} \Rightarrow \text{Id}_{A}$ and is called 
the \emph{identity 2-arrow} of $A$.  

\item  Let $A,B,C$ and 
$D$ be objects of $\B$ and let $\text{Hom}(A,B,C,D)$ 
denote the category $\text{Hom}(C,D)\times \text{Hom}(B,C)
\times \text{Hom}(A,B)$.  There is a natural isomorphism, 
called the \emph{associativity isomorphism}  
and denoted $a(A,B,C,D)$, as pictured in the following 
diagram: 
$$
\xymatrix{ 
\text{Hom}(A,B,C,D) \ar[r]^-{m(B,C,D)\times \text{id}} 
\ar[d]_-{\text{id}\times m(A,B,C)} & 
\text{Hom}(B,D)\times \text{Hom}(A,B) \ar @2{->}[dl]^-{a(A,B,C,D)} 
\ar[d]^-{m(A,B,D)}                                      \\ 
\text{Hom}(C,D)\times \text{Hom}(A,C) \ar[r]_-{m(A,C,D)} & 
\text{Hom}(A,D).                                                } 
$$
So $a$ consists of an arrow $a(\c,\b,\a)
=a(A,B,C,D)(\c,\b,\a)$ 
in $\text{Hom}(A,D)$ from $(\c\circ 
\b)\circ \a$ to $\c\circ (\b\circ 
\a)$ for all 1-arrows $\c$ in $\text{Hom}(
C,D)$, $\b$ in $\text{Hom}(B,C)$ 
and $\a$ in $\text{Hom}(A,B)$ such that 
given 1-arrows $\c^{'}$, $\b^{'}$ and 
$\a^{'}$ together with 2-arrows $\phi:\c
\Rightarrow \c^{'}$, $\psi:\b\Rightarrow 
\b^{'}$ and $\chi:\a\Rightarrow \a^{'}$ 
then the following diagram commutes: 
$$ 
\xymatrix{ 
(\c\circ \b)\circ \a 
\ar @2{->}[rr]^{(\phi\circ \psi)
\circ \chi} \ar @2{->}[d]_{a(\c,\b,
\a)} & & (\c^{'}\circ \b^{'})
\circ \a^{'} \ar @2{->}[d]^{a(\c^{'},
\b^{'},\a^{'})}                           \\ 
\c\circ (\b\circ \a) \ar @2{->}[rr]^{
\phi\circ (\psi\circ \chi)} 
& & \c^{'}\circ (\b^{'}\circ \a^{'}).} 
$$

\item  Suppose $A$ and $B$ are objects of $\B$ and that 
$\a$ is a 1-arrow of $\text{Hom}(A,B)$.  Then there are 
natural isomorphisms $\text{id}_{L}(A,B)(\a):\a\Rightarrow 
\a\circ \text{Id}_{A}$ and $\text{id}_{R}(A,B)(\a):
\a\Rightarrow \text{Id}_{B}\circ \a$ in $\text{Hom}(A,B)$ 
called \emph{left} and \emph{right identity isomorphisms} 
respectively.  The isomorphisms are natural in 
the sense that if $\phi:\a\Rightarrow \b$ is a 2-arrow 
in $\text{Hom}(A,B)$ then the following diagrams of 2-arrows 
commute: 
$$
\xymatrix{ 
\a \ar @2{->}[d]_-{\phi} \ar @2{->}[rr]^-{\text{id}_{L}(A,B)(\a)} & 
& \a\circ \text{Id}_{A} \ar @2{->}[d]^-{\phi\circ \text{id}_{A}} \\   
\b \ar @2{->}[rr]^-{\text{id}_{L}(A,B)(\b)} & & \b\circ 
\text{Id}_{A},                                                  } 
$$ 
and 
$$
\xymatrix{ 
\a \ar @2{->}[d]_-{\phi} \ar @2{->}[rr]^-{\text{id}_{R}
(A,B)(\a)} & & \text{Id}_{B}\circ \a \ar @2{->}[d]^-{
\text{id}_{B}\circ \phi}                               \\ 
\b \ar @2{->}[rr]^-{\text{id}_{R}(A,B)(\b)} & & \text{Id}_{B}\circ \b. } 
$$ 
 
\item  This axiom forces a coherency condition 
on the associativity isomorphism $a$.  So 
suppose we have objects $A,B,C,D$ 
and $E$ of $\B$.  Suppose we have 1-arrows 
$\a$ of $\text{Hom}(A,B)$, $\b$ of $\text{Hom}(B,C)$, 
$\c$ of $\text{Hom}(C,D)$ and $\d$ of $\text{Hom}(D,E)$.  
Then we require that the following 
pentagonal diagram of 2-arrows in $\text{Hom}(A,E)$ commutes: 
$$
\xymatrix@!C=40pt{ 
& & (\d\circ (\c\circ \b))\circ \a 
\ar @2{->}[drr]^-{\ \ \ \ a(\d,\c\circ \b,\a)} & & \\ 
((\d\circ \c)\circ \b)\circ \a 
\ar @2{->}[urr]^-{a(\d,\c,\b)\circ 1_{\a}\ \ \ \ } 
\ar @2{->}[dr]^-{a(\d\circ \c,\b,\a)} & 
& & & \d\circ ((\c\circ \b)\circ \a) 
\ar @2{->}[dl]^-{\ \ \ \ 1_{\d}\circ a(\c,\b,\a)}   \\ 
& (\d\circ \c)\circ (\b\circ \a) 
\ar @2{->}[rr]_-{a(\d,\c,\b\circ \a)} & & 
\d\circ (\c\circ (\b\circ \a)) &            } 
$$
This condition is 
called \emph{associativity coherence}.   

\item  Suppose $A,B$ and $C$ are objects 
of $\B$ and $\a$ and $\b$ are 
1-arrows in $\text{Hom}(A,B)$ and 
$\text{Hom}(B,C)$ respectively.  We require 
that the following diagram commutes:

$$ 
\xymatrix{ 
\b\circ (\text{Id}_{B}\circ \a) 
\ar[rr]^{a} \ar[dr]_{1_{\b}\circ   
\text{id}_{R}(A,B)\ \ \ } &     & (\b\circ \text{Id}_{B})
\circ \a \ar[dl]^{\ \ \ \text{id}_{L}(B,C)
\circ 1_{\a}}                               \\ 
& \b\circ \a.                         } 
$$ 

\end{enumerate} 
\end{definition} 

One can also define a notion of \emph{morphism} of 
bicategories, \emph{transformations} between morphisms 
of bicategories and finally a notion of 
\emph{modification} between transformations 
of morphisms of bicategories, however we will 
not spell out these definitions and refer instead to 
\cite{Str}.  We will also need the notion of 
a bigroupoid.  

\begin{definition}[\cite{Bre}] 
\label{def:10.1.2} 
A \emph{bigroupoid} consists of a bicategory 
$\mathcal{B}$ which satisfies the 
following two additional axioms: 
\begin{enumerate} 
\item 1-arrows are coherently invertible.  This means 
that if $\a$ is a 1-arrow of $\text{Hom}(A,
B)$ there is a 1-arrow $\b$ of $\text{Hom}(
B,A)$ together with 2-arrows 
$\phi:\b\circ \a\Rightarrow \text{Id}_{A}$ 
in $\text{Hom}(A,A)$ and 
$\psi:\a\circ \b\Rightarrow \text{Id}_{B}$ 
in $\text{Hom}(B,B)$,  
which satisfy the following compatibility 
condition.  We require that 
the following diagram of 2-arrows in 
$\text{Hom}(A,B)$ commutes: 
$$
\xymatrix{ 
& & \a\circ (\b\circ \a) 
\ar @2{->}[drr]^-{1_{\a}\circ \phi} & & \\  
(\a\circ \b)\circ \a \ar @2{->}[urr]^-{
a(\a,\b,\a)} \ar @2{->}[dr]_-{\psi 
\circ 1_{\a}} & & & & \a\circ \text{Id}_{A} 
\ar @2{->}[dl]^-{\text{id}_{L}(A,B)(\a)}   \\ 
& \text{Id}_{B}\circ \a \ar @2{->}[rr]_-{
\text{id}_{R}(A,B)(\a)} & &  \a.             } 
$$
\item We also require that all 2-arrows are invertible.  
\end{enumerate} 
\end{definition} 

One can show as a consequence of the definition 
of a bigroupoid that the following 
diagram of 2-arrows in $\text{Hom}(B,
A)$ commutes: 
$$
\xymatrix{ 
& & \b\circ (\a\circ \b) \ar @2{->}[drr]^-{
1_{\b}\circ \psi} &      &                      \\  
(\b\circ \a)\circ \b \ar @2{->}[urr]^-{
a(\b,\a,\b)} \ar @2{->}[dr]_-{\phi\circ 
1_{\b}} & & & & \b\circ \text{Id}_{B} 
\ar @2{->}[dl]^-{\text{id}_{L}(B,A)(\b)}        \\ 
& \text{Id}_{A}\circ \b \ar @2{->}[rr]_-{
\text{id}_{R}(B,A)(\b)} & & \b.                  } 
$$

\begin{example} 
\label{ex:10.1.3} 
Let $X$ be a topological space.  Define a 
bicategory $\Pi_{2}(X)$ as follows.  Take 
as the objects of $\Pi_{2}(X)$ the points of 
$X$.  Given two points $x_{1}$ and $x_{2}$ define 
a category $\text{Hom}(x_{1},x_{2})$ whose 
objects are the paths $\c:I\to X$ with 
$\c(0) = x_{1}$ and $\c(1) = x_{2}$, where 
$I$ is the unit interval $[0,1]$.  Given two 
such paths $\c_{1}$ and $\c_{2}$ we define the 
set of arrows from $\c_{1}$ to $\c_{2}$ in 
$\text{Hom}(x_{1},x_{2})$ to be the homotopy classes  
$[\mu]$ of maps $\mu:I\times I\to X$ such that $\mu(0,t) = \c_{1}
(t)$, $\mu(1,t) = \c_{2}(t)$, $\mu(s,0) = x_{1}$ 
and $\mu(s,1) = x_{2}$ where it is 
understood that two such maps $\mu$ and $\mu^{'}$ 
belong to the same homotopy class if and only if 
there is a map $H:I\times I\times I\to X$ 
such that $H(0,s,t) = \mu(s,t)$, $H(1,s,t) = 
\mu^{'}(s,t)$, $H(r,0,t) = \c_{1}(t)$, 
$H(r,1,t) = \c_{2}(t)$, $H(r,s,0) = x_{1}$ 
and $H(r,s,1) = x_{2}$.  Notice that 
every arrow in $\text{Hom}(x_{1},x_{2})$ is invertible.   
If we have a third 
path $\c_{3}$ from $x_{1}$ to $x_{2}$ and a 
homotopy class $[\lambda]$ of maps 
$\lambda:I\times I\to X$ from $\c_{2}$ 
to $\c_{3}$ leaving the endpoints fixed, then 
we define the composite arrow $[\lambda \mu]$ 
from $\c_{1}$ to $\c_{3}$ to be the homotopy 
class of the map $\lambda \mu:I\times I\to X$ given by 
$$
(\lambda \mu)(s,t) = \begin{cases} 
                     \mu(2s,t),\ s\in [0,\frac{1}{2}],\  
t\in [0,1],                                     \\ 
                     \lambda(2s-1,t),\ s\in [\frac{1}{2},1],\  
t\in [0,1].                                         \\ 
                     \end{cases} 
$$
where $\mu$ and $\lambda$ are 
representatives of the homotopy 
classes $[\mu]$ and $[\lambda]$ under the 
homotopies $H$ above.   
It is straightforward to check that $[\lambda \mu]$ 
is well defined  
and that this law of 
composition makes $\text{Hom}(x_{1},x_{2})$ 
into a category. 
To recap, the objects of 
$\Pi_{2}(X)$ are the points of $X$, the  
1-arrows of $\Pi_{2}(X)$ are the paths in 
$X$ and the 2-arrows of $\Pi_{2}(X)$ are 
the homotopy classes of homotopies with endpoints fixed between 
such paths.  We need to define a composition 
functor 
$$
m(x_{1},x_{2},x_{3}):\text{Hom}(x_{2},x_{3})\times \text{
Hom}(x_{1},x_{2})\to \text{Hom}(x_{1},x_{3}).  
$$ 
If $\c_{23}$ is a 1-arrow of $\text{Hom}(x_{2},x_{3})$ 
and $\c_{12}$ is a 1-arrow of $\text{Hom}(x_{1},x_{2})$, 
then we define $m(x_{1},x_{2},x_{3})(\c_{23},\c_{12}) = 
\c_{23}\circ \c_{12}$, where $\c_{23}\circ \c_{12}:I\to X$ is the path 
$$ 
(\c_{23}\circ \c_{12})(t) = \begin{cases} 
                  \c_{12}(2t),\ t\in [0,\frac{1}{2}], \\ 
                  \c_{23}(2t-1),\ t\in [\frac{1}{2},1]. \\ 
                  \end{cases} 
$$
If we are also given 1-arrows $\c_{23}^{'}$ of 
$\text{Hom}(x_{2},x_{3})$ and $\c_{12}^{'}$ of 
$\text{Hom}(x_{1},x_{2})$ together with 
2-arrows $[\mu_{23}]:\c_{23}\Rightarrow \c_{23}^{'}$ 
and $[\mu_{12}]:\c_{12}\Rightarrow \c_{12}^{'}$ in 
$\text{Hom}(x_{2},x_{3})$ and $\text{Hom}(x_{1},
x_{2})$ respectively, then we define 
$m(x_{1},x_{2},x_{3})([\mu_{23}],[\mu_{12}]) = [\mu_{23}
\circ \mu_{12}]$ to be the homotopy 
class of the map 
$$
(\mu_{23}\circ \mu_{12})(s,t) = \begin{cases} 
                      \mu_{12}(s,2t),\ s\in [0,1],\  
t\in [0,\frac{1}{2}],                                       \\ 
                   \mu_{23}(s,2t-1),\ s\in [0,1],\  
t\in [\frac{1}{2},1].                                       \\ 
                            \end{cases} 
$$
where $\mu_{23}$ is a representative 
of $[\mu_{23}]$ and $\mu_{12}$ is a 
representative of $[\mu_{12}]$.  
Again, it is straightforward to check that 
this defines a functor.  Next, we need to define 
identity 1-arrows and identity 2-arrows.  Given 
an object $x$ of $\Pi_{2}(X)$, we define $\text{Id}_{x}$  
to be the constant path at $x$ and we define the 
identity 2-arrow $\text{id}_{x}$ to be the homotopy class of the constant 
homotopy from the constant path to itself.    
Now we need to define the associator isomorphism.  
Given 1-arrows $\c_{34}$ in $\text{Hom}(x_{3},x_{4})$, 
$\c_{23}$ in $\text{Hom}(x_{2},x_{3})$ and 
$\c_{12}$ in $\text{Hom}(x_{1},x_{2})$, we need to 
define a 2-arrow $a(\c_{34},\c_{23},\c_{12}):(\c_{34}
\circ \c_{23})\circ \c_{12}\Rightarrow \c_{34}\circ 
(\c_{23}\circ \c_{12})$ in such a way that given 
1-arrows $\c_{34}^{'}$ in $\text{Hom}(x_{3},x_{4})$, 
$\c_{23}^{'}$ in $\text{Hom}(x_{2},x_{3})$ and 
$\c_{12}^{'}$ in $\text{Hom}(x_{1},x_{2})$, together 
with 2-arrows $[\mu_{34}]:\c_{34}\Rightarrow \c_{34}^{'}$, 
$[\mu_{23}]:\c_{23}\Rightarrow \c_{23}^{'}$ and 
$[\mu_{12}]:\c_{12}\Rightarrow \c_{12}^{'}$ then the 
following diagram commutes: 
$$ 
\xymatrix{  
(\c_{34}\circ \c_{23})\circ \c_{12} 
\ar @2{->}[rr]^{(\mu_{34}\circ \mu_{23})\circ \mu_{12}} 
\ar @2{->}[d]_{a(\c_{34},\c_{23},\c_{12})} &   &  
(\c_{34}^{'}\circ \c_{23}^{'})\circ \c_{12}^{'} 
\ar @2{->}[d]^{a(\c_{34}^{'},\c_{23}^{'},\c_{12}^{'})} \\ 
\c_{34}\circ (\c_{23}\circ \c_{12}) \ar @2{->}[rr]_{\mu_{34}
\circ (\mu_{23}\circ \mu_{12})} &   &  \c_{34}^{'}\circ 
(\c_{23}^{'}\circ \c_{12}^{'}),            } 
$$ 
and that the pentagonal diagram of axiom (7) above 
commutes.  There is a standard choice for 
$a(\c_{34},\c_{23},\c_{12})$ --- see for example 
\cite{Spa}.  We set $a(\c_{34},\c_{23},\c_{12})$ 
to be the homotopy class of the map $\bar{a}(
\c_{34},\c_{23},\c_{12}):I\times I\to X$ 
defined by 
$$
\bar{a}(\c_{34},\c_{23},\c_{12})(s,t) = \begin{cases} 
            \c_{12}(\frac{4t}{2-s}),\ s\in [0,1],\  
t\in [0,\frac{2-s}{4}],                            \\ 
            \c_{23}(4t-2+s),\ s\in [0,1],\  
t\in [\frac{2-s}{4},\frac{3-s}{4}],                               \\ 
            \c_{34}(\frac{4t-3+s}{1+s}),\ s\in [0,1],\  
t\in [\frac{3-s}{4},1].                                 \\ 
                              \end{cases} 
$$
Next suppose we are given 2-arrows 
$[\mu_{34}]:\c_{34}\Rightarrow \c_{34}^{'}$, 
$[\mu_{23}]:\c_{23}\Rightarrow \c_{23}^{'}$ and 
$[\mu_{12}]:\c_{12}\Rightarrow \c_{12}^{'}$.  
We need to show that $\bar{a}(\c_{34},\c_{23},\c_{12})
((\mu_{12}\circ \mu_{23})\circ \mu_{34}) \simeq  
(\mu_{12}\circ (\mu_{23}\circ \mu_{34}))\bar{a}(
\c_{34},\c_{23},\c_{12})$, where $\mu_{ij}$ is a 
representative of the homotopy class $[\mu_{ij}]$.  
One can check that a homotopy with endpoints 
fixed between these two maps is 
$$
\begin{cases} 
\mu_{12}(2-2rs,\frac{4t}{2(1-r-s+rs)+1)},
\ r\in [0,1],\ s\in [0,\frac{1}{2}],\  
t\in [0,\frac{2(1-r-s+rs)+1}{4}],                           \\ 
\mu_{23}(2-2rs,4t-2(1-r-s+rs)-1),
\ r\in [0,1],\ s\in [0,\frac{1}{2}],                         \\  
t\in [\frac{2(1-r-s+rs)+1}{4},\frac{2(1-r-s+rs)+2}{4}],      \\  
\mu_{34}(2-2rs,\frac{4t-2(1-r-s+rs)-2}{2r+2s-2rs}),\ r\in [0,1],\ s\in [0,
\frac{1}{2}],\ t\in [\frac{2(1-r-s+rs)+2}{4},1],                           \\  
\mu_{12}(2rs-2s+1,\frac{4t}{r(1-2s)+1}),\ r\in [0,1],\ s\in [
\frac{1}{2},1],\ t\in [0,\frac{r(1-2s)+1}{4}],                         \\  
\mu_{23}(2rs-2s+1,4t-r(1-2s)-1),\ r\in [0,1],\ s\in [\frac{1}{2},1],\\   
t\in [\frac{r(1-2s)+1}{4},\frac{r(1-2s)+2}{4}],                \\ 
\mu_{34}(2rs-2s+1,\frac{4t-r(1-2s)-2}{2-r(1-2s)}),\ 
r\in [0,1],\ s\in [\frac{1}{2},1], t\in [\frac{r(1-2s)+2}{4},1].     \\ 
\end{cases} 
$$
This proves that $a$ is a natural transformation.  
We need to check that the associativity coherence condition 
is satisfied.  So we need to check that the 
following diagram of 2-arrows is the 
identity 2-arrow from $((\c_{45}\circ \c_{34})\circ 
\c_{23})\circ \c_{12}$ to itself.  
$$
\xymatrix@!C=50pt{ 
& & (\c_{45}\circ (\c_{34}\circ \c_{23}))\circ \c_{12}  
\ar @2{->}[drr]^-{\ \ \ \ a(\c_{45},\c_{34}\circ \c_{23},\c_{12})} & & \\ 
((\c_{45}\circ \c_{34})\circ \c_{23})\circ \c_{12} 
\ar @2{->}[urr]^-{a(\c_{45},\c_{34},\c_{23})\circ 
1_{\c_{12}}\ \ \ \ \ } & & & & \c_{45}\circ ((\c_{34}
\circ \c_{23})\circ \c_{12}) \ar @2{->}[dl]^-{\ \ \ 1_{\c_{45}}
\circ a(\c_{34},\c_{23},\c_{12})}                                 \\ 
& (\c_{45}\circ \c_{34})\circ (\c_{23}\circ \c_{12}) 
\ar @2{->}[ul]^-{a(\c_{45}\circ \c_{34},\c_{23},\c_{12})^{-1}\ \ \ \ } 
& & \c_{45}\circ (\c_{34}\circ (\c_{23}\circ \c_{12})) 
\ar @2{->}[ll]_-{a(\c_{45},\c_{34},\c_{23}\circ \c_{12})^{-1}}.    } 
$$
We have the following expressions for the  
2-arrows (here $s\in [0,1]$): 
\begin{eqnarray*} 
(\bar{a}(\c_{45},\c_{34},\c_{23})\circ \text{id}_{\c_{12}})(s,t)  & = & 
\begin{cases} 
\c_{12}(2t),\ t\in [0,\frac{1}{2}],      \\ 
\c_{23}(\frac{8t-4}{2-s}),\ t\in [\frac{1}{2},\frac{6-s}{8}], \\ 
\c_{34}(8t+s-6),\ t\in [\frac{6-s}{8},\frac{7-s}{8}],          \\ 
\c_{45}(\frac{8t+s-7}{1+s}),\ t\in [\frac{7-s}{8},1],          \\ 
\end{cases}                                                    \\ 
\bar{a}(\c_{45},\c_{34}\circ \c_{23},\c_{12})(s,t) & = & 
\begin{cases} 
\c_{12}(\frac{4t}{2-s}),\ t\in [0,\frac{2-s}{4}],                 \\ 
\c_{23}(8t+2s-4),\ t\in [\frac{2-s}{4},\frac{5-2s}{8}],          \\ 
\c_{34}(8t+2s-5),\ t\in [\frac{5-2s}{8},\frac{3-s}{4}],           \\ 
\c_{45}(\frac{4t+s-3}{1+s}),\ t\in [\frac{3-s}{4},1],               \\ 
\end{cases}                                                      \\ 
(\text{id}_{\c_{45}}\circ \bar{a}(\c_{34},\c_{23},\c_{12}))(s,t) & = & 
\begin{cases} 
\c_{12}(\frac{8t}{2-s}),\ t\in [0,\frac{2-s}{8}],               \\ 
\c_{23}(8t+s-2),\ t\in [\frac{2-s}{8},\frac{3-s}{8}],            \\ 
\c_{34}(\frac{8t+s-3}{1+s}),\ t\in [\frac{3-s}{8},\frac{1}{2}],  \\ 
\c_{45}(2t-1),\ t\in [\frac{1}{2},1],                            \\ 
\end{cases}                                                         \\  
\bar{a}(\c_{45},\c_{34},\c_{23}\circ \c_{12})^{-1}(s,t) & = & 
\begin{cases} 
\c_{12}(\frac{8t}{1+s}),\ t\in [0,\frac{1+s}{8}],                   \\ 
\c_{23}(\frac{8t-s-1}{1+s}),\ t\in [\frac{1+s}{8},\frac{1+s}{4}],  \\ 
\c_{34}(4t-s-1),\ t\in [\frac{1+s}{4},\frac{2+s}{4}],               \\ 
\c_{45}(\frac{4t-s-2}{2-s}),\ t\in [\frac{2+s}{4},1],               \\ 
\end{cases}                                                        \\ 
\bar{a}(\c_{45}\circ \c_{34},\c_{23},\c_{12})^{-1}(s,t) & = & 
\begin{cases} 
\c_{12}(\frac{4t}{1+s}),\ t\in [0,\frac{1+s}{4}],                 \\ 
\c_{23}(4t-s-1),\ t\in [\frac{1+s}{4},\frac{2+s}{4}],            \\ 
\c_{34}(\frac{8t-2s-4}{2-s}),\ t\in [\frac{2+s}{4},\frac{6+s}{8}], \\ 
\c_{45}(\frac{8t-s-6}{2-s}),\ t\in [\frac{6+s}{8},1].                \\ 
\end{cases}                                                         
\end{eqnarray*} 
Composing these 2-arrows gives the 2-arrow 
from $((\c_{45}\circ \c_{34})\circ \c_{23})
\circ \c_{12}$ to itself.  We want to show 
that this 2-arrow is the identity 2-arrow at 
$((\c_{45}\circ \c_{34})\circ \c_{23})\circ \c_{12}$.  
The plan is to find homotopies from each 
of the individual 2-arrows above to the 
identity 2-arrow at $((\c_{45}\circ \c_{34})\circ 
\c_{23})\circ \c_{12}$.  Of course such 
homotopies will not fix endpoints but we will see  
that the homotopy obtained by composing each of the 
individual homotopies will.  The homotopies are 
given as follows (here $r\in [0,1]$ and $s\in [0,1]$): 
\begin{eqnarray*} 
F_{1}(r,s,t) & = & \begin{cases} 
\c_{12}(2t),\ t\in [0,\frac{1}{2}],               \\ 
\c_{23}(\frac{8t-4}{(r-1)s+2}),\ t\in [\frac{1}{2},\frac{(r-1)s+6}{8}], \\ 
\c_{34}(8t+s(1-r)-6),\ t\in [\frac{(r-1)s+6}{8},\frac{7-(1-r)s}{8}],   \\ 
\c_{45}(\frac{8t-7+(1-r)s}{(1-r)s+1}),\ t\in 
[\frac{7-(1-r)s}{8},1],                              \\ 
\end{cases}                                         \\ 
F_{2}(r,s,t) & = & \begin{cases}                         
\c_{12}(\frac{4t}{(r-1)s+2}),\ t\in [0,\frac{(r-1)s+2}{4}],    \\ 
\c_{23}(\frac{8t-2(r-1)s-4}{1+r}),\ t\in 
[\frac{(r-1)s+2}{4},\frac{2(r-1)s+r+5}{8}],                     \\ 
\c_{34}(8t-2(r-1)s-r-5),\ t\in 
[\frac{2(r-1)s+r+5}{8},\frac{2(r-1)s+r+6}{8}],                  \\ 
\c_{45}(\frac{8t-2(r-1)s-r-6}{2-r-2(r-1)s}),\ t\in 
[\frac{2(r-1)s+r+6}{8},1],                                      \\ 
\end{cases}                                                     \\ 
F_{3}(r,s,t) & = & \begin{cases}                           
\c_{12}(\frac{8t}{(r-1)s+2r+2}),\ t\in [0,\frac{(r-1)s+2+2r}{8}], \\ 
\c_{23}(\frac{8t-(r-1)s-2-2r}{1+r}),\ t\in [\frac{(r-1)s+2+2r}{8},
\frac{(r-1)s+3r+3}{8}],                                              \\ 
\c_{34}(\frac{8t-(r-1)s-3-3r}{1-(r-1)s}),\ t\in [\frac{(r-1)s+3r+3}{8}, 
\frac{3r+4}{8}],                                                   \\ 
\c_{45}(\frac{8t-3r+4}{4-3r}),\ t\in [\frac{3r+4}{8},1],                \\ 
\end{cases}                                                       \\ 
F_{4}(r,s,t) & = & \begin{cases}                                   
\c_{12}(\frac{8t}{(1-r)s+3r+1}),\ t\in [0,\frac{(1-r)s+3r+1}{8}],     \\ 
\c_{23}(\frac{8t-(1-r)s-3r-1}{(1-r)s+r+1}),\ t\in [
\frac{(1-r)s+3r+1}{8},\frac{(1-r)s+2r+1}{4}],                          \\ 
\c_{34}(\frac{8t-2(1-r)s-4r-2}{2-r}),\ t\in [
\frac{(1-r)s+2r+1}{4},\frac{2(1-r)s+3r+4}{8}],                        \\ 
\c_{45}(\frac{8t-2(1-r)s-3r-4}{4-3r-2(1-r)s}),\ tin [
\frac{2(1-r)s+3r+4}{8},1],                                             \\ 
\end{cases}                                                          \\ 
F_{5}(r,s,t) & = & \begin{cases}                                   
\c_{12}(\frac{4t}{s(1-r)+r+1}),\ t\in [0,\frac{s(1-r)+r+1}{4}],        \\ 
\c_{23}(4t-s(1-r)-r-1),\ t\in [\frac{s(1-r)+r+1}{4},\frac{s(1-r)+r+2}{4}], \\ 
\c_{34}(\frac{8t-2s(1-r)-2r-4}{2-r-s(1-r)}),\ t\in [
\frac{s(1-r)+r+2}{4},\frac{s(1-r)+r+6}],                                  \\ 
\c_{45}(\frac{8t-s(1-r)-r-6}{2-r-s(1-r)},\ t\in [\frac{s(1-r)+r+6}{8},1]. \\ 
\end{cases}
\end{eqnarray*} 
It is not too hard to check that 
if one composes these homotopies 
then one has a homotopy with endpoints 
fixed between a representative of 
the 2-arrow above and the identity 
2-arrow at $((\c_{45}\circ \c_{34})\circ 
\c_{23})\circ \c_{12}$.  This proves the 
associativity coherence of $a$.     

We now need to define left and 
right identity transformations 
$\text{id}_{L}$ and $\text{id}_{R}$ 
respectively.  
If $\c \in \text{Hom}(x_{1},x_{2})$ 
is a 1-arrow then $\text{id}_{L}(x_{1},
x_{2})(\c)$ is a 2-arrow $\c\Rightarrow 
\c\circ \text{Id}_{x_{1}}$.  We define 
$\text{id}_{L}(x_{1},x_{2})(\c)$ to 
be the homotopy class of the map 
$$
(s,t)\mapsto \begin{cases} 
             x_{1},\ t\in [0,\frac{s}{2}],    \\ 
             \c(\frac{2t-s}{2-s}),\ t\in [\frac{s}{2},1].  \\  
             \end{cases} 
$$
We need to prove that the assignment 
$\c \mapsto \text{id}_{L}(x_{1},x_{2})(\c)$ 
defines a natural transformation, so we 
need to show that if we are given 1-arrows 
$\c$ and $\c^{'}$ of $\text{Hom}(x_{1},x_{2})$ 
together with a 2-arrow $[\mu]:\c\Rightarrow 
\c^{'}$ in $\text{Hom}(x_{1},x_{2})$ then the 
following diagram of 2-arrows commutes.  
$$
\xymatrix{ 
\c \ar @2{->}[r]^{[\mu]} \ar @2{->}[d]_{
\text{id}_{L}(x_{1},x_{2})(\c)} & 
\c^{'} \ar @2{->}[d]^{\text{id}_{L}(
x_{1},x_{2})(\c^{'})}                       \\ 
\c\circ \text{Id}_{x_{1}} \ar @2{->}[r]_{
[\mu]\circ \text{id}_{x_{1}}} & 
\c^{'}\circ \text{Id}_{x_{1}} } 
$$
We show that the homotopy 
classes $([\mu]\circ \text{Id}_{x_{1}}
\text{id}_{L}(x_{1},x_{2})(\c)$ and 
$\text{id}_{L}(x_{1},x_{2})(\c^{'})[\mu]$ 
are equal.  Firstly, 
$\text{id}_{L}(x_{1},x_{2})(\c^{'})[\mu]$ is 
represented by the map 
$$
(s,t)\mapsto \begin{cases} 
       \mu(2s,t),\ s\in [0,\frac{1}{2}], t\in [0,1],     \\ 
       x_{1},\ s\in [\frac{1}{2},1],\ t\in [0,\frac{2s-1}{2}], \\ 
       \c(\frac{2t-2s+1}{3-2s}),\ s\in [\frac{1}{2},1],\  
t\in [\frac{2s-1}{2},1],               \\ 
              \end{cases} 
$$
where $\mu:I^{2}\to X$ is a representative of 
the homotopy class $[\mu]$, and the homotopy class  
$([\mu]\circ \text{Id}_{x_{1}})\text{id}_{L}
(x_{1},x_{2})(\c)$ is represented by the map 
$$
(s,t)\mapsto \begin{cases} 
         x_{1},\ s\in [0,\frac{1}{2}],\ t\in [0,s],        \\ 
         \c(\frac{t-s}{1-s}),\ s\in [0,\frac{1}{2}],\  
t\in [s,1],              \\ 
         x_{1},\ s\in [\frac{1}{2},1],\ t\in [0,\frac{2s-1}{2}], \\   
         \mu(2s-1,\frac{2t-2s+1}{3-2s}),\ s\in [\frac{1}{2},1],\  
t\in [\frac{2s-1}{2},1].                       \\ 
               \end{cases} 
$$
Then 
$$
(r,s,t)\mapsto \begin{cases} 
         x_{1},\ r\in [0,1],\ s\in [0,\frac{1}{2}],\  
t\in [0,rs],                                              \\ 
\mu(2(1-r)s,\frac{t-rs}{1-rs}),\ r\in [0,1],\  
s\in [0,\frac{1}{2}],\ t\in [rs,1],                    \\ 
x_{1},\ r\in [0,1],\ s\in [\frac{1}{2},1],\  
t\in [0,\frac{2s-1}{2}],                 \\ 
\mu((2s-2)r+1,\frac{2t-2s+1}{3-2s}),\ r\in [0,1],\  
s\in [\frac{1}{2},1],\ t\in [\frac{2s-1}{2},1],            \\ 
            \end{cases} 
$$
is a homotopy with endpoints fixed between 
$(\mu\circ \text{Id}_{x_{1}})\text{id}_{L}
(x_{1},x_{2})(\c)$ and 
$\text{id}_{L}(x_{1},x_{2})(\c^{'})\mu$.   
Similarly, if $\c$ is a 
1-arrow of $\text{Hom}(x_{1},x_{2})$ then we define 
$\text{id}_{R}(x_{1},x_{2})(\c):\c\Rightarrow 
\text{Id}_{x_{2}}\circ \c$ by setting 
$\text{id}_{R}(x_{1},x_{2})(\c)$ equal 
to the homotopy class of the map  
$$
(s,t)\mapsto \begin{cases} 
\c((s+1)t),\ t\in [0,\frac{1}{s+1}],          \\ 
x_{2},\ t\in [\frac{1}{s+1},1].               \\ 
\end{cases} 
$$
Again one can show that given another 
1-arrow $\c^{'}$ in $\text{Hom}(x_{1},x_{2})$ 
and a 2-arrow $[\mu]:\c\Rightarrow \c^{'}$ 
in $\text{Hom}(x_{1},x_{2})$ then the 
diagram of 2-arrows below commutes.  
$$
\xymatrix{ 
\c \ar @2{->}[r]^{[\mu]} \ar @2{->}[d]_{
\text{id}_{R}(x_{1},x_{2})(\c)} & 
\c^{'} \ar @2{->}[d]^{\text{id}_{R}(
x_{1},x_{2})(\c^{'})}                           \\ 
\text{Id}_{x_{2}}\circ \c \ar @2{->}[r]_{
\text{id}_{x_{2}}\circ [\mu]} & \text{Id}_{x_{2}}
\circ \c^{'}                                        } 
$$
Finally, we need to show that axiom (8) 
of Definition~\ref{def:10.1.1} is 
satisfied.  So we need to show that, given 
objects $x_{i}$, $x_{j}$ and $x_{k}$ of 
$\Pi_{2}(X)$ and 1-arrows $\c_{ij}$ and 
$\c_{jk}$ in $\text{Hom}(x_{i},x_{j})$ 
and $\text{Hom}(x_{j},x_{k})$ respectively, 
that the following diagram of 2-arrows 
in $\text{Hom}(x_{i},x_{k})$ commutes. 
$$
\xymatrix{ 
\c_{jk}\circ (\text{Id}_{x_{j}}\circ \c_{ij}) 
 & & (\c_{jk}\circ \text{Id}_{x_{j}})
\circ \c_{ij} \ar @2{->}[ll]^-{a(\c_{jk},\text{Id}_{x_{j}},\c_{ij})}  \\ 
&  \c_{jk}\circ \c_{ij} \ar @2{->}[ul]^{  
\text{id}_{\c_{jk}}\circ \text{id}_{R}(x_{i},
x_{j})(\c_{ij})\ \ \ \ \ } \ar @2{->}[ur]_{\ \ \ \ \ \text{id}_{L}(
x_{i},x_{j})(\c_{jk})\circ \text{id}_{\c_{ij}}} } 
$$
The 2-arrow $a(\c_{jk},\text{Id}_{x_{j}},\c_{ij})
(\text{id}_{L}(x_{i},x_{j})(\c_{jk})\circ \text{id}_{\c_{ij}})$ 
is represented by the map 
$$
(s,t)\mapsto \begin{cases} 
\c_{ij}(2t),\ s\in [0,\frac{1}{2}],\ t\in [0,\frac{1}{2}], \\ 
x_{j},\ s\in [0,\frac{1}{2}],\ t\in [\frac{1}{2},\frac{s+1}{2}], \\ 
\c_{jk}(\frac{2t-s-1}{1-s}),\ s\in [0,\frac{1}{2}],\ 
t\in [\frac{s+1}{2},1],                                           \\  
\c_{ij}(\frac{4t}{1-2s}),\ s\in [\frac{1}{2},1],\ 
t\in [0,\frac{1-2s}{4}],                                            \\ 
x_{j},\ s\in [\frac{1}{2},1],\ t\in [\frac{1-2s}{4},\frac{2-s}{2}], \\ 
\c_{jk}(\frac{2t+s-2}{s}),\ s\in [\frac{1}{2},1],\ 
t\in [\frac{2-s}{2},1],                                               \\ 
\end{cases} 
$$
while the 2-arrow $\text{id}_{\c_{jk}}\circ 
\text{id}_{R}(x_{i},x_{j})(\c_{ij})$   
is represented 
by the map 
$$
(s,t)\mapsto \begin{cases} 
\c_{ij}(2(2s+1)t),\ s\in [0,\frac{1}{2}],\ 
t\in [0,\frac{1}{2(2s+1)}],                     \\ 
x_{j},\ s\in [0,\frac{1}{2}],\ 
t\in [\frac{1}{2(2s+1)},1],                     \\ 
\c_{jk}(2t-1),\ s\in [\frac{1}{2},1],\ 
t\in [\frac{1}{2},1].                             \\ 
\end{cases} 
$$
We want to show that these two maps are 
homotopic by a homotopy fixing endpoints.  
One can check that such a homotopy is 
given by 
$$
(r,s,t)\mapsto \begin{cases} 
\c_{ij}(2(2rs+1)),\ s\in [0,\frac{1}{2}],\ 
t\in [0,\frac{1}{2(2rs+1)}],                      \\ 
x_{j},\ s\in [0,\frac{1}{2}],\ 
t\in [\frac{1}{2(2rs+1)},\frac{1+(1+r)s}{2}],     \\ 
\c_{jk}(\frac{2t-(1+r)s-1}{1-(1+r)s}),\ 
s\in [0,\frac{1}{2}],\ 
t\in [\frac{1+(1+r)s}{2},1],                       \\ 
\c_{ij}(\frac{2(2rs+2-r)}{1-2(1-r)s}),\ 
s\in [\frac{1}{2},1],\ 
t\in [0,\frac{1-2(1-r)s}{2(2rs+2-r)}],               \\ 
x_{j},\ s\in [\frac{1}{2},1],\ 
t\in [\frac{1-2(1-r)s}{2(2rs+2-r)},\frac{2-r-(1-r)s}{2}], \\ 
\c_{jk}(\frac{2t+(1-r)s-2+r}{r+(1-r)s}),\ 
s\in [\frac{1}{2},1],\ 
t\in [\frac{2-r-(1-r)s}{2},1],                          \\ 
\end{cases} 
$$
Thus we have shown that $\Pi_{2}(X)$ is 
an example of a bicategory.  In fact, 
$\Pi_{2}(X)$ is a bigroupoid as we will 
now show.  We will call $\Pi_{2}(X)$ 
the \emph{homotopy bigroupoid} of $X$.  
Firstly we have to show that the 1-arrows 
of $\Pi_{2}(X)$ are coherently invertible.  
So let $x_{1}$ and $x_{2}$ be points of 
$X$ and let $\c$ be a 1-arrow of $\text{Hom}
(x_{1},x_{2})$.  There is an obvious candidate 
for the 1-arrow $\c^{-1}$ which is to be a weak coherent 
inverse of $\c$, namely we set $\c^{-1}$ to 
be the path opposite to $\c$, so $\c^{-1}(t) = 
\c(1-t)$.  A homotopy $\mu$ joining $\c^{-1}\circ 
\c$ to the constant path at $x_{1}$ is given by 
$$
\mu(s,t) = \begin{cases} 
\c(2t(1-s)),\ 0\leq s\leq 1,\ 0\leq t\leq \frac{1}{2},   \\ 
\c((1-s)(2-2t)),\ 0\leq s\leq 1,\ \frac{1}{2}\leq t\leq 1, \\ 
\end{cases} 
$$
while a homotopy $\nu$ joining $\c\circ \c^{-1}$ 
to the constant path at $x_{2}$ is 
given by 
$$
\nu(s,t) = \begin{cases} 
\c(1-2t+2ts),\ 0\leq s\leq 1,\ 0\leq t\leq \frac{1}{2}, \\ 
\c(2t-1-2st+2s),\ 0\leq s\leq 1, \frac{1}{2}\leq t\leq 1. \\ 
\end{cases} 
$$
One can check that these satisfy the 
conditions in (1) of Definition~\ref{def:10.1.2}. 
Clearly all 2-arrows in $\Pi_{2}(X)$ are 
invertible so $\Pi_{2}(X)$ is a bigroupoid.  
We will 
need this example when we come to discuss 
the tautological bundle 2-gerbe.       
\end{example} 
 
Another important example of a bicategory 
is a \emph{2-category}.  This is a 
bicategory in which the composition 
functor is strictly associative rather 
than associative up to a coherent 
natural transformation and the left and 
right identity isomorphisms are the 
identity.  We will write the definition 
of a 2-category out in full below, as we 
will refer to it  
several times.  

\begin{definition}[\cite{KelStr},]    
\label{def:10.1.4} 
A \emph{2-category} $\mathcal{C}$ consists of the 
following data: 
\begin{enumerate} 
\item A set $Ob(\mathcal{C})$ whose elements 
are called objects of $\mathcal{C}$.  

\item Given two objects $x_{1}$ and $x_{2}$ 
of $\mathcal{C}$ there is a category 
$\text{Hom}(x_{1},x_{2})$ whose objects are 
called 1-arrows of $\mathcal{C}$ and whose 
arrows are called 2-arrows of $\mathcal{C}$.  
As before if $\a$ and $\b$ are 1-arrows of 
$\text{Hom}(x_{1},x_{2})$ and $\phi$ is a 
2-arrow of $\text{Hom}(x_{1},x_{2})$ with 
source $\a$ and target $\b$ then we 
indicate this as $\phi:\a\Rightarrow \b$.  

\item Given three objects $x_{1}$, $x_{2}$ 
and $x_{3}$ of $\mathcal{C}$ there is a 
composition functor 
$$
m(x_{1},x_{2},x_{3}):\text{Hom}(x_{2},x_{3})
\times \text{Hom}(x_{1},x_{2})\to \text{Hom}
(x_{1},x_{3}) 
$$
such that, if $x_{4}$ is another object of 
$\mathcal{C}$, then the two functors 
\begin{eqnarray*} 
& & m(x_{1},x_{3},x_{4})(1\times m(x_{1},x_{2},x_{3})) \\  
& & m(x_{1},x_{2},x_{4})(m(x_{2},x_{3},x_{4})
\times 1) 
\end{eqnarray*} 
from $\text{Hom}(x_{1},x_{2},x_{3},x_{4})$  
to $\text{Hom}(x_{1},x_{4})$ coincide (here 
$\text{Hom}(x_{1},x_{2},x_{3},x_{4})$ denotes the 
category $\text{Hom}(x_{3},x_{4})\times \text{Hom}(x_{2},x_{3})
\times \text{Hom}(x_{1},x_{2})$).  
As before, if $\a$ and $\b$ are 1-arrows of 
$\text{Hom}(x_{2},x_{3})$ and $\text{Hom}(x_{1},
x_{2})$, then we usually write $\a\circ \b$ 
for $m(x_{1},x_{2},x_{3})(\a,\b)$.  Also if 
$\phi:\a\Rightarrow \a^{'}$ and $\psi:\b
\Rightarrow \b^{'}$ are 2-arrows in 
$\text{Hom}(x_{2},x_{3})$ and $\text{Hom}(x_{1},x_{2})$ 
respectively, then we write $\phi\circ \psi$ 
for the 2-arrow $m(x_{1},x_{2},x_{3})(\phi,\psi)$ 
in $\text{Hom}(x_{1},x_{3})$.  

\item For each object $x$ of $\mathcal{C}$ there 
is a 1-arrow $\text{Id}_{x}$ of $\text{Hom}(x,x)$ called 
the identity 1-arrow of $x$.  The identity map  
of $\text{Id}_{x}$ is denoted by $\text{id}_{x}$ and is called 
the identity 2-arrow of $x$.  If $y$ is another 
object of $\mathcal{C}$ and $\a$ is a 1-arrow 
of $\text{Hom}(x,y)$ then we have $\text{Id}_{y}\circ \a = 
\a = \a\circ \text{Id}_{x}$.  If $\b$ is another 1-arrow 
of $\text{Hom}(x,y)$ and $\phi:\a\Rightarrow \b$ 
is a 2-arrow of $\text{Hom}(x,y)$ then 
$\text{id}_{y}\circ \phi = \phi = \phi\circ \text{id}_{x}$.  
\end{enumerate} 
\end{definition} 

\section{Bundle 2-gerbes over a point and bicategories} 
\label{sec:10.2} 

We first review a construction of \cite{Mur} which 
shows how to associate a $\cstar$ groupoid to a 
bundle gerbe $(P,X,M)$ restricted to a single point 
$m_{0}$ of $M$.  Let $X_{m_{0}} = \pi^{-1}(m_{0})$.     
Take the 
objects of the $\cstar$ groupoid to be the points of 
$X_{m_{0}}$.  Given two such points $x_{1}$ and $x_{2}$, 
let the set of arrows from $x_{1}$ to $x_{2}$ be the 
points of the fibre $P_{(x_{1},x_{2})}$.  We have to be able 
to compose arrows and we also have to show that there 
exist identity arrows.  Given points $x_{1}$, $x_{2}$ and 
$x_{3}$ of $X_{m_{0}}$ as well as  
arrows $u_{12}:
x_{1}\to x_{2}$ and $u_{23}:x_{2}\to x_{3}$ (ie points 
$u_{12}\in P_{(x_{1},x_{2})}$ and $u_{23}\in P_{(x_{2},x_{3})}$) define the 
composed arrow $u_{23}u_{12}:x_{1}\to x_{3}$ to be the point 
$m_{P}(u_{23}\otimes u_{12})\in P_{(x_{1},x_{3})}$.  Here 
$m_{P}$ denotes the bundle gerbe product in $P$.  Since the 
bundle gerbe product is associative, this process of 
composing arrows is associative.  The identity arrow 
$1:x\to x$ is given by the identity section $e(y)\in P_{(x,x)}$.  
Note that the category so defined is a groupoid (ie every 
arrow is invertible) and that the set of automorphisms of 
an object is isomorphic to $\cstar$.  

Now suppose that we have two bundle gerbes $(P,X,M)$ and 
$(Q,Y,M)$.  Denote by $P_{m_{0}}$ and $Q_{m_{0}}$ the $\cstar$ 
groupoids constructed from $P$ and $Q$ respectively.  
Suppose we have a bundle gerbe morphism $\bar{f} = (\hat{f},f,\text{id}_{M})
:P\to Q$.  Then it is easy to see that because $\hat{f}$ 
commutes with the respective bundle gerbe products, we can 
use $\bar{f}$ to define a functor $f:P_{m_{0}}\to Q_{m_{0}}$.  
If there is a second bundle gerbe morphism $\bar{g} = (\hat{g},
g,\text{id}_{M}):P\to Q$ then we can form the $\cstar$ bundle 
$(f,g)^{-1}Q\to X_{m_{0}}$ as in Lemma~\ref{lemma:3.4.6} and, 
by the results of that Lemma, we know that $(f,g)^{-1}Q$ 
`descends to $m_{0}$'.  Hence it is possible to choose 
a section of $(f,g)^{-1}Q$ which is compatible with the 
descent isomorphism.  Such a section then gives rise to 
a natural transformation $f\Rightarrow g$.   

Bicategories are related to bundle 2-gerbes 
over a point in the following way.  Let 
$(Q,Y,X,M)$ be a bundle 2-gerbe and let $m$ be a point of $M$.  We define 
a bicategory $\mathcal{Q}_{m}$ as follows.  
The objects of $\mathcal{Q}_{m}$ the points 
of $X_{m} = \pi^{-1}(m)$.  Given two such objects 
$x_{i}$ and $x_{j}$, define a category $\text{Hom}(x_{i},x_{j})$ 
as above.  So we let the objects of 
$\text{Hom}(x_{i},x_{j})$ (the 1-arrows of 
$\mathcal{Q}_{m}$) be the points $y$ lying 
in the fibre $Y_{(x_{i},x_{j})}$.  Given 
two such points $y$ and $y^{'}$ we define 
the set of arrows from $y$ to $y^{'}$ to be 
points $u\in Q_{(y,y^{'})}$.  Composition of arrows 
is via the bundle gerbe product in $Q$ as explained above.   

The bundle gerbe morphism $\bar{m}:\pi_{1}^{-1}Q
\otimes \pi_{3}^{-1}Q\to \pi_{2}^{-1}Q$ can 
then be used to define the composition functor 
$$
m:\text{Hom}(x_{j},x_{k})\times \text{Hom}(x_{i},x_{j}) 
\to \text{Hom}(x_{i},x_{k}) 
$$
by mapping an object $(y_{jk},y_{ij})$ of 
$\text{Hom}(x_{j},x_{k})\times \text{Hom}(x_{i},x_{j})$ 
to the object $m(y_{jk},y_{ij})$ of 
$\text{Hom}(x_{i},x_{k})$ and mapping an arrow 
$(u_{jk},u_{ij})$ in $\text{Hom}(x_{j},x_{k})\times 
\text{Hom}(x_{i},x_{j})$ from $(y_{jk},y_{ij})$ 
to $(y_{jk}^{'},y_{ij}^{'})$ to the arrow 
$\hat{m}(u_{jk}\otimes u_{ij})$ in $\text{Hom}(x_{i},x_{k})$ 
from $m(y_{jk},y_{ij})$ to $m(y_{jk}^{'},y_{ij}^{'})$.  
The fact that $\bar{m}$ is a bundle gerbe morphism and 
hence commutes with the bundle gerbe products 
ensures that $m$ is a functor.  We now need to define 
the associator natural transformation $a(x_{i},x_{j},
x_{k},x_{l})$ between the functors bounding the 
diagram below: 
$$
\xymatrix{ 
\text{Hom}(x_{i},x_{j},x_{k},x_{l}) \ar[r]^-{m\times 1} 
\ar[d]_-{1\times m} & \text{Hom}(x_{j},x_{k})\times 
\text{Hom}(x_{i},x_{j}) \ar[d]^-{m} \ar @2{->}[dl]^-{a
(x_{i},x_{j},x_{k},x_{l})}                              \\ 
\text{Hom}(x_{k},x_{l})\times \text{Hom}(x_{i},x_{k}) 
\ar[r]_-{m} & \text{Hom}(x_{i},x_{k}).                   } 
$$
Given an object $(y_{kl},y_{jk},y_{ij})$ of 
$\text{Hom}(x_{i},x_{j},x_{k},x_{l})$, we set 
$a(x_{i},x_{j},x_{k},x_{l})(y_{kl},y_{jk},y_{ij})$ 
to be the arrow $\hat{a}(y_{kl},y_{jk},y_{ij})$.  
The fact that the section $\hat{a}$ of $\hat{A}$ 
descends to a section $a$ of $A$ ensures that 
$a(x_{i},x_{j},x_{k},x_{l})$ is a natural 
transformation and the coherency condition satisfied 
by the section $a$ guarantees that the natural 
transformation $a(x_{i},x_{j},x_{k},x_{l})$ 
satisfies the associativity coherency condition.  

To define an identity 1-arrow at an object $x$ of 
$\mathcal{Q}_{m}$, first choose a set map $I:\Delta(X)\to 
Y$ which is a section of $\pi_{Y}:Y\to X^{[2]}$ (here 
$\Delta(X)$ denotes the image of $X$ under the diagonal 
map $\Delta:X\to X^{[2]}$, $x\mapsto (x,x)$).  So 
$I$ associates to each point $x$ of $X$ an element 
of the fibre $Y_{(x,x)}$.  Let $\text{Id}_{x} = I(x)$ 
and let $\text{id}_{x}:\text{Id}_{x}\Rightarrow \text{Id}_{x}$ 
be the identity section of the bundle gerbe $Q$ evaluated 
at the point $I(x)$ --- ie $e(I(x))$.  We now need to 
define left and right identity transformations.  
Ignoring such concepts as continuity and smoothness, notice 
that we can define bundle gerbe morphisms $\bar{\text{Id}}_{L}:
Q\to Q$ and $\bar{\text{Id}}_{R}:Q\to Q$ as follows: 
\begin{eqnarray*} 
& & \text{Id}_{L}:Y\to Y               \\ 
& & y\mapsto m(y,I(\pi_{2}(\pi_{Y}(y)))), \\ 
& & \text{Id}_{R}:Y\ to Y                  \\ 
& & y\mapsto m(I(\pi_{1}(\pi_{Y}(y))),y),    
\end{eqnarray*} 
and 
\begin{eqnarray*} 
& & \hat{\text{Id}}_{L}:Q\to Q              \\ 
& & u\mapsto \hat{m}(u\otimes e(I(\pi_{2}(\pi_{Y}(\pi_{Q}(u)))))), \\ 
& & \hat{\text{Id}}_{R}:Q\to Q               \\ 
& & u\mapsto \hat{m}(e(I(\pi_{1}(\pi_{Y}(\pi_{Q}(u))))),u).    
\end{eqnarray*} 
Therefore we obtain $\cstar$ bundles (in the set theoretic 
sense) $I_{L} = D_{\bar{1}_{Q},\bar{\text{Id}}_{L}}$ 
and $I_{R} = D_{(\bar{1}_{Q},\bar{\text{Id}}_{R}}$ 
on $X^{[2]}$ and hence by restriction set theoretic 
$\cstar$ bundles $I_{L}$ and $I_{R}$ on $X_{m}\times X_{m}$.  
Choose set maps $\text{id}_{L}:X_{m}^{2}\to I_{L}$ 
and $\text{id}_{R}:X_{m}^{2}\to I_{R}$ which are sections 
of the projections $\pi_{I_{L}}:I_{L}\to X_{m}^{2}$ and 
$\pi_{I_{R}}:I_{R}\to X^{2}_{m}$ respectively.  We denote 
also by $\text{id}_{L}$ and $\text{id}_{R}$ the lifts 
of the sections to the set theoretic $\cstar$ bundles 
$(1,\text{Id}_{L})^{-1}Q$ and $(1,\text{Id}_{R})^{-1}Q$ 
respectively.  It is easy to see that these sections 
provide 2-arrows $\text{id}_{L}(x_{i},x_{j})(y):y\Rightarrow 
m(y,I_{x_{i}})$ and $\text{id}_{R}(x_{i},x_{j})(y):y\Rightarrow 
m(I_{x_{j}},y)$ for a 1-arrow $y\in Y_{(x_{i},x_{j})}$ 
which are natural in the sense of (6) of Definition~\ref{def:10.1.1}.  
We need to check that the left and right identity 
natural transformations are compatible with the associator 
natural transformation $a(x_{i},x_{j},x_{k},x_{l})$.              
It is easy to see that we have an isomorphism 
of set theoretic $\cstar$ bundles over $X_{m}^{3}$ 
given fibrewise as follows: 
$$
I_{R}(x_{1},x_{2})\otimes I^{*}_{L}(x_{2},x_{3})\simeq 
A(x_{1},x_{2},x_{2},x_{3}). 
$$
We can define a set map 
$f:X_{m}^{3}\to \cstar$ such that 
$$
a(x_{1},x_{2},x_{2},x_{3})\text{id}_{L}(x_{2},x_{3}) 
= \text{id}_{R}(x_{1},x_{2})\cdot f(x_{1},x_{2},x_{3}).  
$$
Therefore if $y_{12} \in Y_{(x_{1},x_{2})}$ and 
$y_{23}\in Y_{(x_{2},x_{3})}$ then $f$ satisfies 
$$
\hat{a}(y_{12},I(x_{2}),y_{23})(\text{id}_{L}(x_{2},x_{3})(y_{23})
\circ 1_{y_{12}}) = 1_{y_{23}}\circ \text{id}_{R}(x_{1},x_{2})(y_{12})\cdot 
f(x_{1},x_{2},x_{3}). 
$$
Similarly we can define a function 
$g:X_{m}^{3}\to \cstar$ by 
$$
a(x_{1},x_{1},x_{2},x_{3})\text{id}_{L}(x_{1},x_{3}) = 
\text{id}_{L}(x_{1},x_{2})\cdot g(x_{1},x_{2},x_{3}). 
$$
Again, if $y_{12}\in Y_{(x_{1},x_{2})}$ and 
$y_{23} \in Y_{(x_{2},x_{3})}$ then $g$ satisfies 
$$
\hat{a}(I(x_{1}),y_{12},y_{23})\text{id}_{L}(x_{1},x_{3})(y_{23}\circ y_{12}) 
= 1_{y_{23}}\circ \text{id}_{L}(x_{1},x_{2})(y_{12})\cdot 
g(x_{1},x_{2},x_{3}). 
$$
By considering the diagram of Axiom TA3 of \cite{GorPowStr} we have 
the following equation: 
\begin{equation} 
f(x_{1},x_{2},x_{3}) = f(x_{1},x_{2},x_{4})
g(x_{2},x_{3},x_{4})^{-1}.\label{eq:axiom}  
\end{equation} 
We now need to be more specific about our 
choices of $\text{id}_{L}$ and $\text{id}_{R}$.  
First of all, let us agree to put $\text{id}_{L}(x,x) = 
\text{id}_{R}(x,x)$ so that the two 2-arrows 
$\text{id}_{L}(x,x)(I_{x}):I_{x}\Rightarrow I_{x}\circ I_{x}$ 
and $\text{id}_{R}(x,x)(I_{x}):I_{x}\Rightarrow I_{x}\circ I_{x}$ 
are equal.  Notice also that we can choose 
$\text{id}_{L}$ and $\text{id}_{R}$ so that we have 
$\hat{a}(I_{x},I_{x},I_{x})(\text{id}_{L}(x,x)(I_{x})\circ 
1_{I_{x}}) = 1_{I_{x}}\circ \text{id}_{R}(x,x)(I_{x})$.  
Therefore we have made an explicit choice of 
$\text{id}_{L}$ and $\text{id}_{R}$ over the diagonal 
$\Delta(X)\subset X_{m}^{2}$.  For $y_{12} \in 
Y_{(x_{1},x_{2})}$ with $x_{1}\neq x_{2}$ let us agree 
to choose $\text{id}_{R}(x_{1},x_{2})$ so that the 
2-arrow $\text{id}_{R}(x_{1},x_{2})(y_{12}):y_{12}\Rightarrow 
I_{x_{2}}\circ y_{12}$ satisfies 
$$
\hat{a}(y_{12},I_{x_{2}},I_{x_{2}})(\text{id}_{L}(x_{2},x_{2})
(I_{x_{2}})\circ 1_{y_{12}}) = 1_{I_{x_{2}}}\circ \text{id}_{R}
(x_{1},x_{2})(y_{12}). 
$$
For $y_{21}\in Y_{(x_{2},x_{1})}$ with $x_{1}\neq x_{2}$ let us 
agree to choose $\text{id}_{L}(x_{2},x_{1})$ so that the 
2-arrow $\text{id}_{L}(x_{2},x_{1})(y_{21}):y_{21}\Rightarrow 
y_{21}\circ I_{x_{2}}$ satisfies: 
$$
a(I_{x_{2}},I_{x_{2}},y_{21})(\text{id}_{L}(x_{2},x_{1})(y_{21})\circ 
1_{I_{x_{2}}}) = 1_{y_{12}}\circ \text{id}_{R}(x_{2},x_{2})(I_{x_{2}}). 
$$
Notice that this has now fixed a choice of $\text{id}_{R}$ and 
$\text{id}_{L}$.  Also notice that we now have 
\begin{eqnarray*} 
f(x,x,x) & = & 1 \\ 
f(x_{1},x_{2},x_{2}) & = & 1 \\ 
f(x_{2},x_{2},x_{1}) & = & 1 
\end{eqnarray*} 
Put $x_{4} = x_{2}$ and $x_{1} = x_{2}$ in 
equation~\ref{eq:axiom}.  Then we get 
$$
f(x_{2},x_{2},x_{3}) = f(x_{2},x_{2},x_{2})g(x_{2},x_{3},x_{2})^{-1}. 
$$
This implies that $g(x_{2},x_{3},x_{2}) = 1$.  Now put 
$x_{4} = x_{2}$ in equation~\ref{eq:axiom}.  We get 
$$
f(x_{1},x_{2},x_{3}) = f(x_{1},x_{2},x_{2})g(x_{2},x_{3},x_{2})^{-1} = 1. 
$$
This shows that $\text{id}_{L}$ and $\text{id}_{R}$ 
satisfy the compatibility condition with $a$.  
We record the above discussion in the following 
Proposition. 

\begin{proposition} 
Given a bundle 2-gerbe $(Q,Y,X,M)$ on $M$ we 
can associate to each point $m$ of $M$ a 
bicategory $\mathcal{Q}_{m}$ where $\mathcal{Q}_{m}$ 
is defined as above. 
\end{proposition} 

\begin{note} 
\begin{enumerate} 
\item  Recall from \cite{Mur} that a $\cstar$ groupoid 
is a principal $\cstar$ bundle $P$ over $X^{2}$ 
with an associative product.  The point  
of this result was that the existence of 
an identity and of inverses was a consequence 
of the product being associative and one did not 
need to impose further axioms.  Let us agree to 
call a `$\cstar$ bicategory' a bicategory $\mathcal{Q}$ 
such that for any pair of objects $A$ and $B$ of 
$\mathcal{Q}$, the category $\text{Hom}(A,B)$ 
is a $\cstar$ groupoid.  It then follows that 
every $\cstar$ bicategory arises from a bundle 
2-gerbe over a point, or alternatively, as a 
simplicial bundle gerbe on the simplicial 
manifold $X = \{X_{p}\}$ with $X_{p} = X^{p}$ 
where $X$ is the set of objects of the $\cstar$ 
bicategory.  In fact one can show that  
the $\cstar$ bicategory arising 
from a bundle 2-gerbe is a $\cstar$ 
bigroupoid and hence any $\cstar$ bigroupoid 
would arise from a bundle 2-gerbe over a point.   
  
\item  Notice that we could weaken the definition 
of a bundle 2-gerbe by removing the requirement 
that there exist an `associator' section and still 
obtain a $\cstar$ bicategory (in fact a $\cstar$ 
bigroupoid) when one restricts the bundle 2-gerbe 
to a point of the base manifold.  This is because we 
could choose a (set theoretic) section $a$ of the `associator 
bundle' $A\to X_{m}^{4}$ and define a $\cstar$ valued 
function $f:X_{m}^{5}\to \cstar$ satisfying a cocycle condition 
on $X_{m}^{6}$ by comparing $\d(a)$ with the canonical 
trivialisation of $\d(A)$.  It then follows that we 
could rescale $a$ to obtain the associativity coherence 
condition $\d(a) = 1$.  However, in order to 
associate a class in $H^{4}(M;\Z)$ to a bundle 2-gerbe 
we need to be able to choose $a$ smoothly.  

\item  We could require further axioms to hold 
in the definition of a bundle 2-gerbe so that the 
construction of the $\cstar$ bicategory above was 
`smooth' in a certain sense, ie one had a smooth choice 
of the section $I:X\to \Delta^{-1}Y$ and smooth choices 
of the left and right identity isomorphisms $\text{id}_{L}$ 
and $\text{id}_{R}$.  These conditions will in fact hold in 
the main examples of bundle 2-gerbes we consider in the 
next chapter.    
\end{enumerate} 
\end{note}

\setcounter{chapter}{8}
\chapter{Some examples of bundle 2-gerbes} 
\label{chapter:11} 

\section{The lifting bundle 2-gerbe} 
\label{sec:11.1} 
Suppose we have an extension of groups 
$B\cstar\to \hat{G}\to G$, 
with $B\cstar$ central in $\hat{G}$.  
Suppose also that $P\to M$ is a 
principal $G$ bundle.  Since $B\cstar$ 
is an abelian group, we can form the 
lifting $B\cstar$ bundle gerbe $(\tilde{P},
P,M)$.  So the bundle gerbe product on 
$\tilde{P}$ is a $B\cstar$ bundle 
isomorphism $m:\pi_{1}^{-1}\tilde{P}\otimes 
\pi_{3}^{-1}\tilde{P}\to \pi_{2}^{-1}\tilde{P}$ 
which is associative in the usual sense.  
(As we have previously remarked, the fact that 
$B\cstar$ is an abelian group enables us 
to form the contracted product 
$\pi_{1}^{-1}\tilde{P}\otimes \pi_{3}^{-1}
\tilde{P}$).  Now we can form the lifting 
bundle gerbe $(Q,\tilde{P},P^{[2]})$.  The 
$B\cstar$ bundle gerbe product $m$ then 
induces a bundle gerbe morphism 
$\bar{m}:\pi_{1}^{-1}Q\otimes \pi_{3}^{-1}Q
\to \pi_{2}^{-1}Q$ which is strictly associative 
in the sense that the two induced bundle 
gerbe morphisms $\bar{m}_{1}$ and 
$\bar{m}_{2}$ are equal.  It follows that 
the quadruple $(Q,\tilde{P},P,M)$ 
defines a bundle 2-gerbe in the strict 
sense (Definition~\ref{def:7.1.4}).  
An example of this situation occurs with the 
short exact sequence of abelian groups 
$B\cstar\to EB\cstar\to B^{2}\cstar$ and the 
universal $B^{2}\cstar$ bundle $EB^{2}\cstar\to 
B^{3}\cstar$.  We record this discussion in 
the following proposition. 

\begin{proposition} 
\label{prop:11.1.1} 
Let $B\cstar\to \hat{G}\to G$ be a short 
exact sequence of groups with $B\cstar$ 
central in $\hat{G}$.  Let $P\to M$ be 
a principal $G$ bundle.  Then the 
quadruple $(Q,\tilde{P},P,M)$ defined 
above is a strict bundle 2-gerbe.  
\end{proposition} 

Another situation in which bundle 
2-gerbes arise is the following.  
Suppose we have a surjection $\pi:X\to M$ 
admitting local sections and a $B\cstar$ 
bundle $P\to X^{[2]}$ together with a 
$B\cstar$ bundle isomorphism $m:\pi_{1}^{-1}
P\otimes \pi_{3}^{-1}P\to \pi_{2}^{-1}P$ 
covering the identity on $X^{[3]}$.  
However we do not assume that the two 
induced $B\cstar$ bundle isomorphisms 
$m_{1}$ and $m_{2}$ are equal.  
Here $m_{1}$ and $m_{2}$ are defined as 
follows 
\begin{eqnarray*} 
& m_{1}:\pi_{2}^{-1}\pi_{1}^{-1}P\otimes 
\pi_{4}^{-1}(\pi_{1}^{-1}P\otimes\pi_{3}
^{-1}P) \to \pi_{2}^{-1}\pi_{2}^{-1}P       \\ 
& m_{1} = \pi_{2}^{-1}m\circ (1\otimes 
\pi_{4}^{-1}m)                               \\ 
& m_{2}:\pi_{1}^{-1}(\pi_{1}^{-1}P\otimes 
\pi_{3}^{-1}P)\otimes \pi_{3}^{-1}\pi_{3}
^{-1}P \to \pi_{3}^{-1}\pi_{2}^{-1}P          \\ 
& m_{2} = \pi_{3}^{-1}m\circ (\pi_{1}^{-1}m
\otimes 1),                                    
\end{eqnarray*} 
where we have identified $\pi_{2}^{-1}\pi_{1}^{-1}
P\otimes \pi_{4}^{-1}(\pi_{1}^{-1}P\otimes 
\pi_{3}^{-1}P) = \pi_{1}^{-1}(\pi_{1}^{-1}P
\otimes \pi_{3}^{-1}P)\otimes \pi_{3}^{-1}\pi_{3}
^{-1}P$ and $\pi_{2}^{-1}\pi_{2}^{-1}P = 
\pi_{3}^{-1}\pi_{2}^{-1}P$ by virtue of the 
simplicial identities satisfied by the $\pi_{i}$ 
on the simplicial 
manifold $X=\{X_{p}\}$ with $X_{p} = X^{[p+1]}$.  
Therefore $m_{1}$ and $m_{2}$ differ by a map 
$\phi:X^{[4]}\to B\cstar$ which must satisfy 
$\d(\phi) = 1$ on $X^{[5]}$.  Suppose finally 
that $\phi$ has a lift $\hat{\phi}:X^{[4]}\to E\cstar$ 
which satisfies $\d(\hat{\phi}) = 1$ on $X^{[5]}$.  
Then if we let $(\tilde{P},P,X^{[2]})$ denote 
the lifting bundle gerbe associated to the 
$B\cstar$ bundle $P\to X^{[2]}$ then the quadruple 
$(\tilde{P},P,X,M)$ is a bundle 2-gerbe.     
As before the $B\cstar$ bundle isomorphism 
$m$ provides a bundle gerbe morphism 
$\bar{m}:\pi_{1}^{-1}\tilde{P}\otimes 
\pi_{3}^{-1}\tilde{P}\to \pi_{2}^{-1}\tilde{P}$ 
covering the identity on $X^{[3]}$.  One can check 
that the associator bundle $A\to X^{[4]}$ 
has classifying map equal to $\phi$ and so 
$\hat{\phi}$ provides a section $a$ of $A$ --- 
the associator section --- and the condition 
$\d(\hat{\phi}) = 1$ is the requirement 
that $\d(a)$ matches the canonical trivialisation 
of $\d(A)$ on $X^{[5]}$. 

\section{The tautological bundle 2-gerbe}
\label{sec:11.2} 
Here we verify that the tautological 
bundle 2-gerbe introduced in 
\cite{CaMuWa} is still a bundle 2-gerbe 
as defined in these notes.  To be more 
precise we only verify this for the 
special case in which the base manifold is 
3-connected or at least has 
$\pi_{3}(M) = 0 $.   

To define this tautological object 
recall that we start with a manifold $M$, 
which we assume to be 3-connected,
together with a closed four form $\Theta$ 
on $M$ with integral periods representing 
a class in $H^{4}(M;\Z)$.  Next we take 
the path fibration $\pi:\mathcal{P}M \to M$ consisting 
of piecewise smooth paths $\c:[0,1]\to M$, $\c(0) = m_{0}$ where 
$m_{0}$ is a basepoint of $M$ and with $\pi$ the 
map which sends a path $\c$ to its 
end point $\c(1)$.  

The fibre product $\mathcal{P}M^{[2]}$ consists 
of all pairs of piecewise smooth paths $(\c_{1},\c_{2})$ 
with the same end point, that is  
$\c_{1}(1) = \c_{2}(1)$.  Denote by 
$\c^{\circ}$ the path $\c$ traversed in 
the opposite direction and denote by 
$\c_{1}*\c_{2}$ the path formed from 
a pair of paths $\c_{1},\c_{2}$ by 
traversing $\c_{1}$ at double speed and 
then $\c_{2}$ at double speed.  Of 
course, this is only a piecewise smooth path if 
$\c_{1}(1) = \c_{2}(0)$.  Thus if 
$(\c_{1},\c_{2}) \in \mathcal{P}M^{[2]}$ then 
$\c_{1}*\c_{2}^{\circ} \in \Omega M$, where 
$\Omega M$ is the based piecewise smooth loop space of $M$. 

Therefore we can define a map 
$\text{ev}:\mathcal{P}M^{[2]} \times S^{1} \to M$ by 
evaluating the loop $\c_{1}*\c_{2}^{\circ}$ 
against the angle $\theta$.  Therefore 
we may define a closed three form 
$\omega$ on $\mathcal{P}M^{[2]}$ by pulling back 
$\Theta$ to $\mathcal{P}M^{[2]}$ with $\text{ev}$ and then 
integrating $\text{ev}^{*}(\Theta)$ over the 
circle.  One can check that $\omega$ 
so defined is a closed three form with 
integral periods.  The next step is to 
perform the tautological bundle gerbe 
construction on $\mathcal{P}M^{[2]}$ using the 
integral three form $\omega$.   

Bundle 2-gerbes enter the picture because 
we can define a bundle gerbe morphism 
covering the map $\mathcal{P}M^{[3]} \to \mathcal{P}M^{[2]}$ 
defined by $(\c_{1},\c_{2},\c_{3}) \mapsto (\c_{1},\c_{3})$.
For the details see \cite{CaMuWa}.

To make the calculation simpler assume 
for the moment that the base manifold $M$ is a point.  Then 
we start with a manifold $X$ and a 
closed 3-from $\omega$ on $X$ with 
integral periods.  We construct a bundle gerbe 
on $X^{2} = X \times X$ in the usual way.
We define a space $Y$ above $X^{2}$ 
with fibre at $(x_{1},x_{2}) \in X^{2}$ 
equal to 
$Y_{(x_{1},x_{2})} = \{\a:I \to X, 
\ \a \ \text{piecewise smooth}: \a(0) = x_{1},\ \a(1) = x_{2}\}$.  
Next define a $\cstar$-bundle 
$Q \to Y^{[2]}$ whose fibre at $(\a,\b) \in Y^{[2]}$ 
is all equivalence classes $[\mu,z]$, 
where $z$ is a non-zero complex number 
and $\mu:\a \Rightarrow \b$ is a 
homotopy with end points fixed,
that is $\mu:I\times I\to X$ and $\mu(0,t) = \a(t)$, 
$\mu(1,t) = \b(t)$, $\mu(s,0) = \a(0) = \b(0)$ 
and $\mu(s,1) = \a(1) = \b(1)$.   
The equivalence relation is defined by 
$(\mu_{1},z_{1}) \equivalencerelation (\mu_{2},z_{2})$ 
if for any homotopy $F:\mu_{1} \Rightarrow \mu_{2}$ 
with endpoints fixed (so $F:I \times I \times I \to X$, 
$F(0,s,t) = \mu_{1}(s,t)$, $F(1,s,t) = \mu_{2}(s,t)$ 
and for each $r \in [0,1]$, $F(r,-,-)$ is 
a homotopy with endpoints fixed between 
$\a$ and $\b$) we have 
$$
z_{2} = \exp (\int_{I^{3}} F^{*}(\omega)) z_{1}
$$
If $(x_{1},x_{2},x_{3}) \in X^{3} = X \times X \times X$ 
note that we can define a map 
$m:Y_{(x_{2},x_{3})}\times Y_{(x_{1},x_{2})}\to Y_{(x_{1},x_{3})}$  
by $(\a_{23},\a_{12}) \mapsto \a_{23}\circ\a_{12}$,
where 
$$
(\a_{23}\circ \a_{12})(t) = \begin{cases}
                   \a_{12}(2t),\ 0 \leq t \leq 1/2,   \\
                   \a_{23}(2t-1),\ 1/2 \leq t \leq 1, \\
                   \end{cases}
$$
and where $\a_{23} \in Y_{(x_{2},x_{3})}$ 
and $\a_{12}\in Y_{(x_{1},x_{2})}$.   
This map $m$ clearly extends to a map 
$m:\pi_{1}^{-1}Y\times_{X^{3}}\pi_{3}^{-1}
Y\to \pi_{2}^{-1}Y$ covering the identity 
on $X^{3}$.  We can also define a $\cstar$ 
bundle morphism $\hat{m}:\pi_{1}^{-1}Q\otimes 
\pi_{3}^{-1}Q\to \pi_{2}^{-1}Q$ covering the 
map $m^{[2]}:(\pi_{1}^{-1}Y\times_{X^{3}}
\pi_{3}^{-1}Y)^{[2]}\to (\pi_{2}^{-1}Y)^{[2]}$ 
by $\hat{m}([\mu_{23},z_{23}]\otimes [\mu_{12},
z_{12}]) = [\mu_{23}\circ \mu_{12},z_{23}z_{12}]$, 
where $\mu_{23}$ is a homotopy with endpoints 
fixed between $\a_{23}$ and $\b_{23}$, $\mu_{12}$ 
is a homotopy with endpoints fixed between 
$\a_{12}$ and $\b_{12}$ and $\mu_{23}\circ \mu_{12}$ 
is the homotopy with endpoints fixed between 
$\a_{23}\circ \a_{12}$ and $\b_{23}\circ \b_{12}$ 
given by 
$$
(\mu_{23}\circ \mu_{12})(s,t) = \begin{cases} 
\mu_{12}(s,2t),\ 0\leq s\leq 1,\ 0\leq t\leq 1/2,   \\ 
\mu_{23}(s,2t-1),\ 0\leq s\leq 1,\ 1/2\leq t\leq 1.  \\ 
\end{cases} 
$$
One can check, see \cite{CaMuWa}, that 
$\hat{m}$ is well defined and commutes with the bundle gerbe 
products.  Clearly $\hat{m}$ is $\cstar$ 
equivariant.  Therefore $\bar{m} = (\hat{m},
m,\text{id})$ is a bundle gerbe morphism.   

As usual, $\bar{m}$ defines two bundle 
gerbe morphisms $\bar{m}_{1} = (\hat{m}_{1},
m_{1},\text{id})$ and $\bar{m}_{2} = 
(\hat{m}_{2},m_{2},\text{id})$ between 
the appropriately defined bundle gerbes 
over $X^{4} = X\times X\times X\times X$.  
The maps $m_{1}$ and $m_{2}$ are defined 
fibrewise by 
\begin{eqnarray*}
m_{1}:Y_{(x_{3},x_{4})}\times Y_{(x_{2},x_{3})}
\times Y_{(x_{1},x_{2})} 
\to Y_{(x_{1},x_{4})}             \\
(\a_{34},\a_{23},\a_{12})\mapsto 
(\a_{34}\circ \a_{23})\circ \a_{12},   \\
m_{2}:Y_{(x_{3},x_{4})}\times Y_{(x_{2},x_{3})}
\times Y_{(x_{1},x_{2})} 
\to Y_{(x_{1},x_{4})}              \\
(\a_{34},\a_{23},\a_{12})\mapsto 
\a_{34}\circ (\a_{23}\circ \a_{12})   
\end{eqnarray*}
for $\a_{34}\in Y_{(x_{3},x_{4})}$, 
$\a_{23}\in Y_{(x_{2},x_{3})}$ and 
$\a_{12}\in Y_{(x_{1},x_{2})}$.   

It is a standard result (see \cite{Spa}) 
that $m_{1}\simeq m_{2}$. 
This means that the $\cstar$-bundle 
$(m_{1},m_{2})^{-1}Q \to Y_{1234}$ has 
a section $\hat{a}$, where 
we write for notational convenience 
$Y_{1234}$ for either of the spaces 
$\pi_{2}^{-1}\pi_{1}^{-1}Y\times_{X^{4}}
\pi_{4}^{-1}(\pi_{1}^{-1}Y\times_{X^{3}}
\pi_{3}^{-1}Y$ or $\pi_{1}^{-1}(\pi_{1}
^{-1}Y\times_{X^{3}}\pi_{3}^{-1}Y)\times_{X
^{4}}\pi_{3}^{-1}\pi_{3}^{-1}Y$ over 
$X^{4}$.  From the general theory 
of simplicial bundle gerbes (see Section~\ref{sec:7.1})
we know that the $\cstar$ bundle 
$\hat{A} = (m_{1},m_{2})^{-1}Q$ 
descends to a $\cstar$ bundle $A$ on 
$X^{4}$ such that $\d(A)$ is canonically 
trivialised on $X^{5} = X\times X\times X
\times X\times X$.  To show that the 
section $\hat{a}$ of $\hat{A}$ descends 
to a section $a$ of $A$, we first need 
an explicit formula for $\hat{a}$.      
We get this from the 
homotopy $m_{1}\simeq m_{2}$ (see Example~
\ref{ex:10.1.3} and \cite{Spa}).  
Set $\hat{a}(\a_{34},\a_{23},\a_{12}) = [
\bar{a}(\a_{34},\a_{23},\a_{12}),1]$, 
where $\bar{a}(\a_{34},\a_{23},\a_{12})
:I\times I\to X$ is given by 
$$
\bar{a}(\a_{34},\a_{23},\a_{12})(s,t) = \begin{cases}
             \a_{12}(\frac{4t}{2-s}),\ 
 0\leq t \leq\frac{2-s}{4},\ 0\leq s\leq 1,               \\
             \a_{23}(4t+s-2),\ \frac{2-s}{4}\leq 
t \leq\frac{3-s}{4},\ 0\leq s\leq 1,                         \\
             \a_{34}(\frac{4t+s-3}{1+s}),\ \frac{3-s}{4}
\leq t\leq 1,\ 0\leq s\leq 1.  \\
            \end{cases}
$$
From Example~\ref{ex:10.1.3}, 
we know that $\bar{a}(\a_{34},\a_{23},\a_{12})$ 
is the associator natural transformation 
for the bicategory $\Pi_{2}(X)$ evaluated at 
$(\a_{34},\a_{23},\a_{12})$.  The fact that 
this is a natural transformation is exactly 
the requirement that the section $\hat{a}$ 
descend to $X^{4}$ --- that is, commute with 
the descent isomorphism for the $\cstar$ 
bundle $\hat{A}$.  

We now need to show that the descended 
section $a$ satisfies the required 
coherency condition on $X^{5}$, that is 
that $\d(a)$ matches the canonical trivialisation 
of $\d(A)$.  Again this follows from 
the results of Example~\ref{ex:10.1.3}. 
One can check that the homotopy 
$F = F_{5}F_{4}F_{3}F_{2}F_{1}$ joining 
\begin{eqnarray*} 
&   & \bar{a}(\a_{45}\circ \a_{34},\a_{23},\a_{12})^{-1}
\bar{a}(\a_{45},\a_{34},\a_{23}\circ \a_{12})^{-1}(\text{id}_{\a_{45}}
\circ \bar{a}(\a_{34},\a_{23},\a_{12}))          \\ 
&   & \bar{a}(\a_{45},\a_{34}\circ 
\a_{23},\a_{12})(\bar{a}(\a_{45},\a_{34},\a_{23})\circ 
\text{id}_{\a_{12}}) 
\end{eqnarray*} 
to the constant homotopy from $((\a_{45}\circ \a_{34})  
\circ \a_{23})\circ \a_{12}$ to itself has 
image that is at most two dimensional.  Hence 
$F^{*}(\omega) = 0$ which shows that $a$ satisfies 
the required coherency condition.  

In the general case, where $M$ is not a point, 
we can apply the construction given above to 
each fibre of the path fibration $\pi:\mathcal{P}M\to M$ 
to obtain a bundle 2-gerbe.   

\section{The bundle 2-gerbe of a principal $G$ bundle} 
\label{sec:11.3} 
Suppose we are given a principal $G$ 
bundle $P \to  M$, where $G$ is a 
compact, simply connected, simple Lie group.  
Then it is known \cite{Car} that  
$\pi_{2}(G) = 0$ and $H^{3}(G;\Z) = \Z$.  
It is shown in \cite{Bry} that there 
is a closed, bi-invariant three form $\nu$ on 
$G$ with integral periods which 
represents the canonical generator of 
$H^{3}(G;\Z)=\Z$.  If $G = SU(N)$, then 
$\nu$ is the three form $\frac{1}{24\pi^{2}}
\text{tr}(dgg^{-1})^{3}$ --- see \cite{BryMcL1}.  

Recall from \cite{CaMuWa} that 
we can define a bundle gerbe 
$(Q,PG,G)$ on $G$ with three 
curvature equal to $\nu$.  
The fibre 
of $Q\to PG^{[2]}$ at a 
point $(\a,\b) \in PG^{[2]}$ 
is all pairs $(\phi,z)$, where $z \in \cstar$ 
and $\phi : \a \Rightarrow \b$ is a homotopy 
$\phi:I^{2}\to G$ satisfying 
$\phi(0,t) = \a(t), \phi(1,t) = \b(t), \phi(s,0) = e$ 
and $\phi(s,1) = \a(1) = \b(1)$, under the 
equivalence relation $\equivalencerelation$ defined by 
$(\phi_{1},z_{1}) \equivalencerelation (\phi_{2},z_{2})$ if and 
only if for all homotopies $F:\phi_{1}\Rightarrow \phi_{2}$ 
with endpoints fixed between $\phi_{1}$ and $\phi_{2}$ we have 
$z_{2} = z_{1}\cdot \exp(\int_{I^{3}}F^{*}\nu)$.  

The bundle gerbe product is defined by 
$$
[\phi_{1},z_{1}]\otimes [\phi_{2},z_{2}] \mapsto 
[\phi_{1}\phi_{2},z_{1}z_{2}],
$$   
where $\phi_{1}\phi_{2}$ denotes 
the homotopy defined by 
$$
(\phi_{1}\phi_{2})(s,t) = \begin{cases} 
                        \phi_{1}(2s,t) & \text{for $0 \leq s \leq 1/2$} \\
                        \phi_{2}(2s-1,t) &\text{for $1/2 \leq s \leq 1$} \\
                        \end{cases}     
$$
It is shown in \cite{CaMuWa} that 
this is well defined, associative, etc.  

\begin{proposition} 
\label{prop:11.3.1} 
The bundle gerbe $(Q,PG,G)$ is a 
simplicial bundle gerbe on the 
simplicial manifold $NG$.
\end{proposition} 

\begin{proof} 
We first need to define the bundle gerbe 
morphism $\bar{m} = (\hat{m},m,\text{id})$ 
which maps 
$$
\bar{m}:d_{0}^{-1}Q\otimes d_{2}^{-1}Q 
\to d_{1}^{-1}Q.
$$
Define 
$m:d_{0}^{-1}PG\times_{G^{2}}d_{2}^{-1}PG\to d_{1}^{-1}PG$ 
covering the identity on $G^{2} = G\times G$ 
by sending $(\a,\b)$ to the path 
$\a\circ \a(1)\b$ given by 
$$
(\a\circ \a(1)\b)(t) = \begin{cases} 
                    \a(2t), & 0 \leq t \leq 1/2   \\ 
                    \a(1)\b(t), & 1/2\leq t \leq 1. \\ 
                    \end{cases}  
$$
Next, we need to define a $\cstar$ 
equivariant map 
$\hat{m}:d_{0}^{-1}Q\otimes d_{2}^{-1}Q\to d_{1}^{-1}Q$ 
covering 
$$
m^{[2]}:(d_{0}^{-1}PG\times_{G^{2}}
d_{1}^{-1}PG)^{[2]}\to d_{1}^{-1}PG^{[2]}
$$ 
and check that it commutes with the 
bundle gerbe product.  
So take pairs $(\phi, z)$ and 
$(\psi, w)$ where $z,w \in \cstar$ 
and $\phi:I^{2}\to G$ and 
$\psi:I^{2}\to G$ are homotopies with 
endpoints fixed between paths 
$\a_{1},\a_{2}$ and $\b_{1},\b_{2}$ 
respectively.  Then we put 
$$
\hat{m}((\phi,z),(\psi,w)) = 
(\phi\circ \phi(0,1)\psi ,zw)    
$$
where $\phi\circ \phi(0,1)\psi :I^{2}\to G$ 
is the homotopy with endpoints fixed 
between $\a_{1}\circ \a_{1}(1)\b_{1}$ and 
$\a_{2}\circ \a_{2}(1)\b_{2}$ given by 
$$
(\phi\circ \phi(0,1)\psi)(s,t) = 
                     \begin{cases}  
                      \phi(s,2t), & 0\leq t \leq 1/2    \\ 
                      \phi(0,1)\psi(s,2t-1), & 1/2 \leq t \leq 1. \\ 
                      \end{cases} 
$$
We need to check firstly that this 
map is well defined --- that is it respects 
the equivalence relation 
$\equivalencerelation$ --- and secondly that 
$\hat{m}$ commutes with the bundle 
gerbe products.  

So suppose 
$(\phi,z)\equivalencerelation (\phi^{'},z^{'})$ and 
$(\psi,w)\equivalencerelation (\psi^{'},w^{'})$, 
where $\phi$ and $\phi^{'}$ are homotopies 
with endpoints fixed between paths $\a_{1}$ and 
$\a_{2}$ and where $\psi$ and $\psi^{'}$ are 
homotopies with endpoints fixed between 
paths $\b_{1}$ and $\b_{2}$.  We want to show 
that 
\begin{equation} 
(\phi\circ \phi(0,1)\psi,zw)
\equivalencerelation 
(\phi^{'}\circ \phi^{'}(0,1)\psi^{'},z^{'}w^{'}).   
\end{equation} 
Therefore we want to show that for all 
homotopies $H:I^{3}\to G$ with endpoints fixed 
between $\phi\circ \phi(0,1)\psi$ and 
$\phi^{'}\circ \phi^{'}(0,1)\psi^{'}$ we have 
$$
z^{'}w^{'}
 = zw
\exp (\int_{I^{3}}H^{*}\nu).
$$
Note that if $\Phi:I^{3}\to G$ is a 
homotopy with endpoints fixed between 
$\phi$ and $\phi^{'}$ and  
$\Psi:I^{3}\to G$ is a homotopy with 
endpoints fixed between $\psi$ and 
$\psi^{'}$, then by integrality 
of $\nu$ we have  
$$
\exp (\int_{I^{3}}H^{*}\nu) = 
\exp (\int_{I^{3}}(\Phi\circ \Phi(0,0,1)\Psi)^{*}\nu). 
$$
Therefore we are reduced to showing 
that 
\begin{equation} 
z^{'}w^{'} 
= zw \exp (\int_{I^{3}}(\Phi\circ \Phi(0,0,1)\Psi)^{*}\nu).   
\end{equation} 
We have 
$$
\exp (\int_{I^{3}}(\Phi\circ \Phi(0,0,1)\Psi)^{*}\nu) = 
\exp (\int_{I^{3}}\Phi^{*}\nu) \exp (\int_{I^{3}}
(\Phi(0,0,1)\Psi)^{*}\nu). 
$$
By the bi-invariantness of $\nu$, we get 
$(\Phi(0,0,1)\Psi)^{*}\nu = \Psi^{*}\nu$, hence 
$$
\exp (\int_{I^{3}}(\Phi * \Phi(0,0,1)\Psi)^{*}\nu ) 
= \exp (\int_{I^{3}}\Phi^{*}\nu)\exp (\int_{I^{3}}\Psi^{*}\nu), 
$$
which implies the result.  
Hence $\hat{m}$ is well defined.  It is 
a straightforward matter to verify that 
$\hat{m}$ respects the bundle gerbe products.  

It remains to show that there is a 
transformation of the bundle gerbe 
morphisms $\bar{m}_{1}$ and 
$\bar{m}_{2}$ over $G\times G\times G$ 
which satisfies the compatibility 
criterion over $G\times G\times G\times G$.  
This has already been done above for 
the tautological bundle 2-gerbe and the 
proof given there carries over to this case.  
\end{proof} 

Suppose that we have a principal 
$G$ bundle $\pi:P\to M$.   
Pulling back the 
simplicial bundle gerbe 
$(Q,PG,G)$ to $P^{[2]}$ via 
the canonical map 
$\tau:P^{[2]}\to G$ defines 
a quadruple of manifolds 
$(\tilde{Q},\tilde{P},P,M)$  
which is in fact a bundle 2-gerbe.  
We have the following proposition.  

\begin{proposition} 
\label{prop:11.3.2} 
The quadruple of manifolds $(\tilde{Q},\tilde{P},P,M)$ 
is a bundle 2-gerbe.  
\end{proposition} 

\begin{note} 
In the next section we will define a class in 
$H^{4}(M;\Z)$ associated to a bundle 2-gerbe.  
It will turn out (see Section~\ref{sec:9.3}) 
that the class in $H^{4}(M;\Z)$ associated to 
the bundle 2-gerbe $\tilde{Q}$ is the 
Pontryagin class of the principal $G$ bundle 
$P$.  
\end{note}

\setcounter{chapter}{9}
\chapter{Bundle 2-gerbe connections and 2-curvings}
\label{chapter:8} 
\section{Definitions and some preliminary lemmas} 
\label{sec:8.1} 
Let $(Q,Y,X,M)$ be a bundle  
2-gerbe.  So there is a surjection $\pi:X\to M$  
admitting local sections.  Recall 
from Section~\ref{sec:3.3} the 
extended Mayer-Vietoris sequence 
$$
0\to \Omega^{p}(M) \stackrel{\pi^{*}}{\to} 
\Omega^{p}(Y) \stackrel{\d}{\to} 
\Omega^{p}(Y^{[2]}) \stackrel{\d}{\to} 
\cdots \stackrel{\d}{\to} 
\Omega^{p}(Y^{[q]}) 
\stackrel{\d}{\to} \cdots 
$$ 
It is shown in \cite{Mur} that this 
sequence is exact (see 
Proposition~\ref{prop:3.3.1}).  
We will exploit this fact to define the 
notion of a bundle 2-gerbe connection and 
a 2-curving for a bundle 2-gerbe connection and   
construct a closed integral four form 
on $M$ associated to the bundle 
2-gerbe $(Q,Y,X,M)$.  

Recall from Section~\ref{sec:3.3} that a bundle
gerbe connection $\nabla$ on the bundle gerbe
$(Q,Y,X^{[2]})$ has a curving $f_{1}$ which is a two form
$f_{1} \in \Omega^{2}(Y)$.  Recall from \cite{Mur}  
(see also Section~\ref{sec:3.3}) that 
there is a closed three form $\omega$ on $X^{[2]}$
with integral periods which satisfies
$\pi_{Y}^{*}\omega = df$.  $\omega$ is called the 
Dixmier-Douady three form of the bundle gerbe $(Q,Y,X^{[2]})$
or sometimes the three-curvature.

\begin{definition}
\label{def:8.1.1}  
We use the notation of the previous paragraph.
A \emph{bundle 2-gerbe connection} on a bundle 2-gerbe
$(Q,Y,X,M)$ is a bundle gerbe connection $\nabla$ on
the bundle gerbe $(Q,Y,X^{[2]})$ together with a choice of curving
$f_{1} \in \Omega^{2}(Y)$ for $\nabla$, such that
$\delta(\omega) = 0$ in $\Omega^{3}(X^{[3]})$, 
where $\omega \in \Omega^{3}(X^{[2]})$ is the 
Dixmier-Douady three form with $\pi_{Y}^{*}\omega = df_{1}$.
\end{definition}

\begin{proposition}
\label{prop:8.1.2}  
Let $(Q,Y,X,M)$ be a bundle 2-gerbe.
Then (i) bundle 2-gerbe connections on 
$(Q,Y,X,M)$ exist and (ii) the three curvature
$\omega$ defines a closed four form $\Theta$
on $M$ with integral periods.
\end{proposition}

Assuming the existence of bundle 2-gerbe 
connections for now, suppose we have 
chosen one on the bundle 2-gerbe 
$(Q,Y,X,M)$.  Let $\omega \in \Omega^{3}(X^{[2]})$ 
be the three curvature, then $\d(\omega) =0$ and 
we can solve $\omega = \d(f_{2})$ 
for some $f_{2} \in \Omega^{3}(X)$.  We call 
a choice of $f_{2}$ a \emph{2-curving} for the bundle 2-gerbe 
connection $(\nabla , f_{1})$.  Since 
$\omega$ is closed, we have $\d(df_{2}) = 0$, 
and hence there exists $\Theta \in \Omega^{4}(M)$ 
such that $df_{2} = \pi^{*}\Theta$ and 
$\Theta$ is closed.  We will see later that 
$\Theta$ is in fact an integral four form 
and hence is a de Rham representative 
of a class in $H^{4}(M;\Z)$.  
$\Theta$ is called the four class of the 
bundle 2-gerbe $(Q,Y,X,M)$. 

\begin{example} 
As an example of this we consider the 
tautological bundle 2-gerbe on a 3-connected 
manifold $M$ associated to a closed integral 
four form $\Theta$ on $M$ (see Section~\ref{sec:11.2}).  
Recall that we first pull $\Theta$ back to $\mathcal{P}M$, 
the piecewise smooth path space of $M$, using the 
evaluation map $\text{ev}:\mathcal{P}M\times I\to M$ 
sending $(\c,t)$ to $\c(t)$.  We then define an 
integral three form $f_{2}$ on $\mathcal{P}M$ by 
$f_{2} = \int_{I}\text{ev}^{*}\Theta$.  It is straightforward 
to check that $df_{2} = \pi^{*}\Theta$ where 
$\pi:\mathcal{P}M\to M$ is the projection.  It then 
follows that the three form $\d(f_{2})$ on 
$\mathcal{P}M^{[2]}$ defines a class in the image 
of $H^{3}(\mathcal{P}M^{[2]};\Z)$ inside 
$H^{3}(\mathcal{P}M^{[2]};\Reals)$.  Alternatively, 
we can think of $\d(f_{2})$ as defining a closed 
integral three form on $\Omega M$, the space of based,  
piecewise smooth loops on $M$.  Recall from 
Section~\ref{sec:11.2} that we then construct the 
tautological bundle gerbe on each fibre of $\pi:\mathcal{P}M 
\to M$.  This amounts to constructing the 
tautological bundle gerbe with bundle gerbe 
connection $\nabla$ and curving $f$ (see example~\ref{ex:tautbgconn}) 
whose three curvature 
equals $\d(f_{2})$ on $\mathcal{P}M^{[2]}$.  Therefore 
$(\nabla,f)$ is an example of a bundle 2-gerbe 
connection on the tautological bundle 2-gerbe 
and $f_{2}$ is a 2-curving of this bundle 2-gerbe 
connection.  
\end{example} 
   
We need the following series of Lemmas.

\begin{lemma}
\label{lemma:8.1.3}  
Let $(P,X,M)$ and $(Q,Y,M)$ be bundle gerbes with
bundle gerbe connections $\nabla_{P}$ and
$\nabla_{Q}$ respectively.  Suppose
$(\hat{f},f,\text{id}):(P,X,M) \to (Q,Y,M)$ is a bundle 
gerbe morphism.  $\hat{f}$ induces an isomorphism
of line bundles
$\tilde{f}: P \to (f^{[2]})^{-1}Q$ over 
$X^{[2]}$ which commutes with the bundle
gerbe products.  Then we have
$$
\nabla_{P} = \tilde{f}^{-1} \circ (f^{[2]})^{-1} \nabla_{Q}
\circ \tilde{f} + \delta(A),
$$
for some $A \in \Omega^{1}(X)$.
\end{lemma}
\begin{proof}  Because $\hat{f}$ commutes with the bundle
gerbe products on $P$ and $Q$, $\tilde{f}$
commutes with the bundle gerbe products on 
$P$ and $(f^{[2]})^{-1}Q$.  Since $\nabla_{Q}$
is a bundle gerbe connection, so is
$(f^{[2]})^{-1}\nabla_{Q}$ and therefore $\nabla_{P}$
and $\tilde{f}^{-1}\circ (f^{[2]})^{-1}\nabla_{Q} \circ \tilde{f}$
are bundle gerbe connections on $P$.  Since 
any two connections differ by the pullback of
a complex valued one form on the base, we have 
$\nabla_{P}=\tilde{f}^{-1}\circ (f^{[2]})^{-1}\nabla_{Q}\circ \tilde{f} + \a$, 
for some $\a \in \Omega^{1}(X^{[2]})$.  Because
$\nabla_{P}$ and
$\tilde{f}^{-1}\circ (f^{[2]})^{-1}\nabla_{Q}\circ \tilde{f}$
are both bundle gerbe connections, we must
have $\delta(\a) = 0$ and so $\a = \delta(A)$
for some $A \in \Omega^{1}(X)$.
\end{proof} 

Suppose now that $(P,X,M)$ and $(Q,Y,M)$
are the bundle gerbes of the previous lemma but we now
have a pair of bundle gerbe morphisms
$(\hat{f}_{i},f_{i},\text{id}):(P,X,M) 
\to (Q,Y,M)$, $i = 1,2.$
By Lemma~\ref{lemma:8.1.3} there exist 
$A_{1},A_{2} \in \Omega^{1}(X)$
such that
$$
\nabla_{P} = \tilde{f}_{i}^{-1} \circ 
(f^{[2]})^{-1}\nabla_{Q} \circ
\tilde{f}_{i} + \delta(A_{i}).
$$
We have the pullback bundle $\abfQ \to X$
of Lemma~\ref{lemma:3.4.6}.  Let
$\phi:\pi_{1}^{-1}\abfQ \simeq \pi_{2}^{-1}\abfQ$
denote the descent isomorphism of Lemma~\ref{lemma:3.4.6} 
and let $(f_{1},f_{2})^{-1}\nabla_{Q}$
denote the pullback connection on
$\abfQ \to X$.

\begin{lemma}
\label{lemma:8.1.4}   
The pullback connection $\abfdQ$ satisfies
$$
\pi_{1}^{-1}\abfdQ = \phi^{-1}\circ \pi_{2}^{-1}\abfdQ 
\circ \phi + \d (A_{1} - A_{2}).
$$
\end{lemma}
\begin{proof}  Define a map $<f_{1},f_{2}>:X^{[2]} \to Y^{[2]}$
by sending the point $(x_{1},x_{2})$ of $X^{[2]}$ to
the point $(f_{1}(x_{1}),f_{2}(x_{2}))$ of $Y^{[2]}$.  The bundle gerbe
product on $Q$ gives the following isomorphisms
of line bundles with connection 
over $X^{[2]}$.
\begin{eqnarray} 
(f_{2}^{[2]})^{-1}Q \otimes \pi_{2}^{-1}(f_{1},f_{2})^{-1}Q  
\stackrel{\simeq}{\rightarrow} 
<f_{1},f_{2}>^{-1}Q              \label{eq:8.1.1}                        \\  
\pi_{1}^{-1}(f_{1},f_{2})^{-1}Q \otimes (f_{1}^{[2]})^{-1}Q 
\stackrel{\simeq}{\rightarrow}   
<f_{1},f_{2}>^{-1}Q            \label{eq:8.1.2}                    
\end{eqnarray} 
Here the line bundle $\abFQ$ is equipped 
with the connection $\abFdQ$, while the 
line bundles $(f_{2}^{[2]})^{-1}Q \otimes \pi_{2}^{-1}
(f_{1},f_{2})^{-1}Q$ 
and $\pi_{1}^{-1}(f_{1},f_{2})
^{-1}Q\otimes (f_{1}^{[2]})^{-1}Q$ 
are endowed with the connections
$(f_{2}^{[2]})^{-1}\nabla_{Q} + 
\pi_{2}^{-1}(f_{1},f_{2})^{-1}\nabla_{Q}$ 
and $\pi_{1}^{-1}(f_{1},f_{2})^{-1}
\nabla_{Q} + (f_{1}^{[2]})^{-1}\nabla_{Q}$ 
respectively.   
Hence we can combine the isomorphisms~\ref{eq:8.1.1} 
and~\ref{eq:8.1.2}  
to get an isomorphism
of line bundles  
$$
\pi_{1}^{-1}\abfQ \otimes (f_{1}^{[2]})^{-1}Q 
\stackrel{\simeq}{\rightarrow} 
(f_{2}^{[2]})^{-1}Q \otimes \pi_{2}^{-1}\abfQ.
$$
Note that this is an isomorphism 
of line bundles with connection.  
We have the following sequence of 
isomorphisms of line bundles
\begin{eqnarray} 
& & \pi_{1}^{-1}(f_{1},f_{2})^{-1}Q 
\rightarrow  \pi_{1}^{-1}
(f_{1},f_{2})^{-1}Q\otimes P\otimes 
P^{*}         \label{eq:8.1.3}                        \\  
& & \pi_{1}^{-1}(f_{1},f_{2})^{-1}Q\otimes P
\otimes P^{*}           
\rightarrow  \pi_{1}^{-1} 
(f_{1},f_{2})^{-1}Q   
\otimes  
(f_{1}^{[2]})^{-1}
Q\otimes (f_{2}^{[2]})^{-1}Q^{*} \label{eq:8.1.4}   \\    
& & \pi_{1}^{-1}(f_{1},f_{2})^{-1}Q\otimes 
(f_{1}^{[2]})^{-1}Q\otimes (f_{2}^{[2]})^{-1}
Q^{*} \stackrel{m_{Q}}{\rightarrow}  \pi_{2}^{-1}
(f_{1},f_{2})^{-1}Q            \label{eq:8.1.5}    \\                           
\end{eqnarray}  
The descent isomorphism $\phi$ of 
Lemma~\ref{lemma:3.4.6} is the composition of these 
line bundle isomorphisms.  Here the 
isomorphism~\ref{eq:8.1.3} is induced by 
the canonical trivialisation of $P\otimes P^{*}$, 
the isomorphism~\ref{eq:8.1.4} is induced by the 
isomorphisms $\tilde{f}_{1}$ and $\tilde{f}_{2}^{*}$
and the isomorphism~\ref{eq:8.1.5} is 
induced by the bundle gerbe product $m_{Q}$ 
(in other words by the isomorphisms~\ref{eq:8.1.1} 
and~\ref{eq:8.1.2}).   
From this, and Lemma~\ref{lemma:8.1.3} 
we see that the difference
$\pi_{1}^{-1}\abfdQ - \phi^{-1}\circ\pi_{2}^{-1}\abfdQ\circ \phi$
is equal to $\d(A_{1}-A_{2})$. 
\end{proof} 
 
Suppose that we have bundle gerbes
$(P,X,M)$, $(Q,Y,M)$ and $(R,Z,M)$ equipped
with bundle gerbe connections
$\nabla_{P}$, $\nabla_{Q}$ and
$\nabla_{R}$ respectively.  Suppose also that
there exist bundle gerbe morphisms
$\bar{f}:P\to Q$ and
$\bar{g}:Q\to R$ with $\bar{f} = (\hat{f},f,\text{id})$ 
and $\bar{g} = (\hat{g},g,\text{id})$.  
With the notation of Lemma~\ref{lemma:8.1.3} we have
$$
\nabla_{P} = \tilde{f}^{-1}\circ (f^{[2]})^{-1}\nabla_{Q}\circ 
\tilde{f} + \delta (A),
$$
$$
\nabla_{Q} = \tilde{g}^{-1}\circ (g^{[2]})^{-1}\nabla_{R}\circ 
\tilde{g} + \delta (B),
$$
for some $A \in \Omega^{1}(X)$,
$B \in \Omega^{1}(Y)$.  Let $\widetilde{g\circ f}$
denote the isomorphism 
$\widetilde{g\circ f}:P \simeq (g^{[2]}\circ f^{[2]})^{-1}R$
induced by the bundle map
$\hat{g}\circ \hat{f}:P \to R$ covering
$g^{[2]}\circ f^{[2]}:X^{[2]} \to Z^{[2]}$.

\begin{lemma}
\label{lemma:8.1.5}  
With the notation of the above 
paragraph, we have the following expression:
$$
\nabla_{P} = (\widetilde{g\circ f})^{-1}\circ ((g^{[2]}\circ f^{[2]})^{-1}
\nabla_{R})\circ (\widetilde{g\circ f}) + \delta(A + f^{*}B).
$$
\end{lemma}
\begin{proof} The isomorphism
$\tilde{g}:Q \to (g^{[2]})^{-1}R$ gives rise
to an isomorphism over $X^{[2]}$:
$$
(f^{[2]})^{-1}\tilde{g}:(f^{[2]})^{-1}Q\simeq (f^{[2]})^{-1}(g^{[2]})^{-1}R.
$$
Also by definition we have an 
isomorphism of line bundles with
connection
$$ 
(f^{[2]})^{-1}(g^{[2]})^{-1}R\simeq (g^{[2]}\circ f^{[2]})^{-1}R,
$$ 
where each bundle is equipped with the 
obvious pullback connection.  Next, the
diagram
$$
\xymatrix{ 
(f^{[2]})^{-1}(g^{[2]})^{-1}R 
\ar[rr]^-{\text{canonical}} & &  
(g^{[2]}\circ f^{[2]})^{-1}R   \\ 
& P \ar[ul]^{((f^{[2]})^{-1}
\tilde{g})\circ \tilde{f}\ \ \ } 
\ar[ur]_{\widetilde{g\circ f}} & } 
$$
commutes.  This is really the 
definition of $(f^{[2]})^{-1}\tilde{g}$.
Therefore we will have
\begin{eqnarray*}
& & (\widetilde{g\circ f})^{-1}\circ     
((g^{[2]}\circ f^{[2]})^{-1}\nabla_{R})\circ  
(\widetilde{g\circ f})              \\
& = & (((f^{[2]})^{-1}\tilde{g})\circ \tilde{f})^{-1}\circ
((f^{[2]})^{-1}(g^{[2]})^{-1}\nabla_{R})\circ (((f^{[2]})^{-1}
\tilde{g})\circ \tilde{f})  \\
& = & \tilde{f}^{-1}\circ ((f^{[2]})^{-1}\tilde{g})^{-1}\circ 
((f^{[2]})^{-1}(g^{[2]})^{-1}\nabla_{R})\circ ((f^{[2]})^{-1}
\tilde{g})\circ \tilde{f}  \\
& = & \tilde{f}^{-1}\circ (f^{[2]})^{-1}(\nabla_{Q} - \delta
(B))\circ \tilde{f}  \\
& = & \tilde{f}^{-1}\circ (f^{[2]})^{-1}\nabla_{Q}\circ \tilde{f} 
 - \delta(f^{*}B)  \\
& = & \nabla_{P} - \delta(A + f^{*}B)  \\
\end{eqnarray*}
\end{proof} 

Suppose, as in Lemma~\ref{lemma:8.1.5}, that we have
bundle gerbes $(P,X,M),\ (Q,Y,M)$ and
$(R,Z,M)$, but that we have pairs of bundle
gerbe morphisms
$(\hat{f}_{i},f_{i},\text{id}):(P,X,M)\to (Q,Y,M)$ and
$(\hat{g}_{i},g_{i},\text{id}):(Q,Y,M)\to (R,Z,M)$ for
$i = 1,2$.  Using Lemma~\ref{lemma:8.1.3} again we get
$$
\nabla_{P} = \tilde{f}_{i}^{-1}\circ (f_{i}^{[2]})^{-1}\nabla_{Q}
\circ \tilde{f}_{i} + \delta(A_{i}),
$$
$$
\nabla_{Q} = \tilde{g}_{i}^{-1}\circ (g_{i}^{[2]})^{-1}\nabla_{R}
\circ \tilde{g}_{i} + \delta(B_{i}),
$$
for some $A_{i} \in \Omega^{1}(X)$,
$B_{i} \in \Omega^{1}(Y)$ for $i = 1,2$.

\begin{lemma}
\label{lemma:8.1.6}  
With the above notation, the pullback connection
$(g_{1}\circ f_{1},g_{2}\circ f_{2})^{-1}\nabla_{R}$
on $(g_{1}\circ f_{1},g_{2}\circ f_{2})^{-1}R \to X$
satisfies
$$
\pi_{1}^{-1}(g_{1}\circ f_{1},g_{2}\circ f_{2})^{-1}\nabla_{R} = 
\phi^{-1}\circ \pi_{2}^{-1}(g_{1}\circ f_{1},g_{2}\circ f_{2})^{-1}
\nabla_{R}\circ \phi + \delta(A_{1} - A_{2} + f_{1}^{*}B_{1} - 
f_{2}^{*}B_{2}),
$$
where $\phi:\pi_{1}^{-1}\abfg ^{-1}R\simeq \pi_{2}^{-1}\abfg ^{-1}R$
is the descent isomorphism of~\ref{lemma:3.4.6}.    
\end{lemma}
\begin{proof}  Same as for Lemma~\ref{lemma:8.1.4}. 
\end{proof}

\section{Proof of the main result} 
\label{sec:8.2} 
We now turn to the proof of Proposition~\ref{prop:8.1.2} 
\begin{proof}  Suppose we are
given a bundle 2-gerbe $(Q,Y,X,M)$.  Let
$\nabla$ be a bundle gerbe connection
on $(Q,Y,X^{[2]})$ with curving $f \in \Omega^{2}(Y)$
such that $\delta(f) = F_{\nabla}$, 
where $F_{\nabla}$ is the curvature of $Q$
with respect to $\nabla$.  The three curvature
$\omega$ then satisfies $\pi_{Y}^{*}\omega = df$.

Let us look at what happens over $X^{[3]}$.
We have induced connections $\pi_{i}^{-1}\nabla$
on the bundle gerbes $(\pi_{i}^{-1}Q,\pi_{i}^{-1}Y,X^{[3]})$.
Let $\pi_{i}^{-1}f$ denote the curving 
for $\pi_{i}^{-1}\nabla$ induced by $f$. 
Let $(p_{1}^{[2]})^{-1}\pi_{1}^{-1}\nabla + 
(p_{2}^{[2]})^{-1}\pi_{3}^{-1}\nabla$
denote the tensor product bundle gerbe 
connection on $\tensorQ \to \abYtwo$.  The 
bundle gerbe morphism 
$$
\bar{m}:(\tensorQ ,\abYtwo ,X^{[3]}) 
\to (\pi_{2}^{-1}Q,\pi_{2}^{-1}Y,X^{[3]})
$$
gives rise to a $\cstar$-bundle isomorphism
$$
\xymatrix{ 
\pi_{1}^{-1}Q\otimes \pi_{3}^{-1}Q 
\ar[rr]^{\tilde{m}} \ar[dr] & & 
(m^{[2]})^{-1}\pi_{2}^{-1}Q \ar[dl]               \\ 
& (\pi_{1}^{-1}Y\times_{X^{[3]}}\pi_{3}^{-1}Y)
^{[2]}                                   } 
$$ 
$\tilde{m}$ clearly commutes with the
respective bundle gerbe products.  By Lemma~\ref{lemma:8.1.3}   
there is a one form 
$\rho \in \Omega^{1}(\abY )$ such that
\begin{equation}
\label{eq:8.2.1.1}  
(p_{1}^{[2]})^{-1}\pi_{1}^{-1}\nabla + 
(p_{2}^{[2]})^{-1}\pi_{3}^{-1}\nabla = 
\tilde{m}^{-1}\circ (m^{[2]})^{-1}\pi_{2}^{-1}\nabla 
\circ \tilde{m} + \delta(\rho).
\end{equation}
Hence we get
\begin{equation}
\label{eq:8.2.1.2}  
(p_{1}^{[2]})^{*}F_{\pi_{1}^{-1}\nabla} + 
(p_{2}^{[2]})^{*}F_{\pi_{3}^{-1}\nabla} = 
(m^{[2]})^{*}F_{\pi_{2}^{-1}\nabla} + \d (d\rho),
\end{equation}
and hence 
\begin{equation}
\label{eq:8.2.1.3}  
\d (p_{1}^{*}\pi_{1}^{-1}f + p_{2}^{*}\pi_{3}^{-1}f - 
m^{*}\pi_{2}^{-1}f - d\rho ) = 0.
\end{equation}

Therefore, there exists 
$\mu \in \Omega^{2}(X^{[3]})$ such 
that
\begin{equation}
\label{eq:8.2.1.4}  
p_{1}^{*}\pi_{1}^{-1}f + p_{2}^{*}\pi_{3}^{-1}f = 
m^{*}\pi_{2}^{-1}f + d\rho + \pi_{Y_{123}}^{*}\mu,
\end{equation}
where $\pi_{Y_{123}}$ denotes the projection 
$Y_{123} = \pi_{1}^{-1}Y\times_{X^{[3]}}\pi_{3}^{-1}Y
\to X^{[3]}$.  It is easily checked that
\begin{equation}
\label{eq:8.2.1.5}  
\pi_{1}^{*}\omega - \pi_{2}^{*}\omega 
+ \pi_{3}^{*}\omega = d\mu .
\end{equation}
Now let us look at the situation
over $X^{[4]}$.  We will use the notation 
introduced in Section~\ref{sec:7.1} so 
that we will denote $\pi_{1}^{-1}Y$ by 
$Y_{23}$, $\pi_{1}^{-1}Q\otimes \pi_{3}^{-1}
Q$ by $Q_{123}$, $\pi_{1}^{-1}\pi_{1}^{-1}Q$ 
by $Q_{34}$ and so on.  As mentioned in 
Section~\ref{sec:7.1}, the bundle gerbe 
morphisms $\bar{m}_{1}$ and $\bar{m}_{2}$ 
are defined as in the following diagram: 
$$
\xymatrix{ 
Q_{1234} \ar[r]^-{\bar{m}\otimes 1} 
\ar[d]_-{1\otimes \bar{m}} & Q_{124} 
\ar[d]^-{\bar{m}}                    \\ 
Q_{134} \ar[r]^-{\bar{m}} & Q_{14}.   } 
$$
So $\bar{m}_{1} = \bar{m}\circ 
(\bar{m}\otimes 1)$ and $\bar{m}_{2} 
= \bar{m}\circ (1\otimes \bar{m})$.  Recall from 
Definition~\ref{def:7.1.2} that 
$\hat{A} = (m_{1},m_{2})^{-1}Q_{14}$ is the 
line bundle on $Y_{1234}$ which descends 
to the associator line bundle $A$ on 
$X^{[4]}$.  Furthermore $\hat{A}$ has a
section $\hat{a}$ which descends to a 
section $a$ of $A$ --- the associator 
section.  $\hat{A}$ inherits the 
pullback connection $\nabla_{\hat{A}} = 
(m_{1},m_{2})^{-1}\nabla_{14}$ (here 
$\nabla_{14}$ denotes the induced 
connection on $Q_{14}$) relative to 
which $\hat{A}$ has curvature    
equal to 
$$
(m_{1},m_{2})^{*}\d(\pi_{2}^{-1}
\pi_{2}^{-1}f) = m_{2}^{*}\pi_{2}^{-1}
\pi_{2}^{-1}f - m_{1}^{*}\pi_{3}^{-1}
\pi_{2}^{-1}f.
$$
Now 
\begin{eqnarray*} 
&   & m_{2}^{*}\pi_{2}^{-1}\pi_{2}^{-1}f        \\ 
& = & (\pi_{2}^{-1}m\circ (1\times_{X^{[4]}}
\pi_{4}^{-1}m))^{*}\pi_{2}^{-1}\pi_{2}^{-1}f     \\ 
& = & (1\times_{X^{[4]}}\pi_{4}^{-1}m)^{*} 
(\pi_{2}^{-1}m)^{*}\pi_{2}^{-1}\pi_{2}^{-1}f     \\ 
& = & (1\times_{X^{[4]}}\pi_{4}^{-1}m)^{*}
\pi_{2}^{-1}(p_{1}^{*}\pi_{1}^{-1}f + 
p_{2}^{*}\pi_{3}^{-1}f - d\rho - \pi_{Y_{134}}^{*}
\pi_{2}^{*}\mu)                                  \\ 
& = & p_{1}^{*}\pi_{2}^{-1}\pi_{1}^{-1}f + 
(1\times_{X^{[4]}}\pi_{4}^{-1}m)^{*}p_{2}^{*}
\pi_{4}^{-1}\pi_{2}^{-1}f - d(1\times_{X^{[4]}}
\pi_{4}^{-1}m)^{*}\pi_{2}^{-1}\rho - \pi_{Y_{1234}}^{*}
\pi_{2}^{*}\mu                                   \\ 
& = & p_{1}^{*}\pi_{2}^{-1}\pi_{1}^{-1}f + 
p_{2}^{*}\pi_{4}^{-1}\pi_{1}^{-1}f + p_{3}^{*}
\pi_{4}^{-1}\pi_{3}^{-1}f - d(1\times_{X^{[4]}}
\pi_{4}^{-1}m)^{*}\pi_{2}^{-1}\rho                \\ 
&  & - 1\times d\pi_{4}^{-1}\rho - \pi_{Y_{1234}}^{*}(
\pi_{2}^{*}\mu + \pi_{4}^{*}\mu).                  
\end{eqnarray*} 
Similarly we have 
\begin{eqnarray*} 
m_{1}^{*}\pi_{3}^{-1}\pi_{2}^{-1}f & = & 
p_{1}^{*}\pi_{1}^{-1}\pi_{1}^{-1}f + 
p_{2}^{*}\pi_{1}^{-1}\pi_{3}^{-1}f + 
p_{3}^{*}\pi_{3}^{-1}\pi_{3}^{-1}f - 
d\pi_{1}^{-1}\rho \times 1                 \\ 
&  & - d(\pi_{1}^{-1}m\times_{X^{[4]}}1)^{*}
\pi_{3}^{-1}\rho - \pi_{Y_{1234}}^{*}(\pi_{1}^{*}\mu 
+ \pi_{3}^{*}\mu).                          
\end{eqnarray*} 
Hence we have 
\begin{eqnarray*}   
m_{2}^{*}\pi_{2}^{-1}\pi_{2}^{-1}f - 
m_{1}^{*}\pi_{3}^{-1}\pi_{2}^{-1}f   
& = & \pi_{Y_{1234}}^{*}\d(\mu) + d(\pi_{1}^{-1}\rho
\times 1) + (\pi_{1}^{-1}m\times_{X^{[4]}}1)^{*}
\pi_{3}^{-1}\rho                                   \\  
&   & - 1\times \pi_{4}^{-1}\rho 
- (1\times_{X^{[4]}}\pi_{4}^{-1}m)^{*}\pi_{2}
^{-1}\rho).                                         
\end{eqnarray*}   
On the other hand we can write 
$(\nabla_{\hat{A}})(\hat{a})=\hat{\a}\otimes \hat{a}$ 
for some one form $\hat{\a}$ on 
$Y_{1234}$.  Thus we must 
have $d\hat{\a} = \pi_{Y_{1234}}^{*}\d(\mu) + d\rho_{1}$ 
where $\rho_{1}\in \Omega^{1}(Y_{1234})$ 
is the one form 
$$
\pi_{1}^{-1}\rho\times 1 + 
(\pi_{1}^{-1}m\times_{X^{[4]}}1)^{*}\pi_{3}^{-1}
\rho - 1\times \pi_{4}^{-1}\rho - 
(1\times_{X^{[4]}}\pi_{4}^{-1}m)^{*}\pi_{2}^{-1}
\rho. 
$$
Let $\phi:\pi_{1}^{-1}\hat{A}\to \pi_{2}^{-1}\hat{A}$ 
denote the descent isomorphism 
of Lemma~\ref{lemma:3.4.6}.  By Lemma's~\ref{lemma:8.1.3} 
and~\ref{lemma:8.1.6}  we have 
$$ 
\pi_{1}^{-1}\nabla_{\hat{A}} = 
\phi^{-1}\circ \pi_{2}^{-1} 
\nabla_{\hat{A}}\circ \phi + 
\d(\rho_{1}).    
$$
Therefore the connection $\nabla_{\hat{A}}-\rho_{1}$ on 
$\hat{A}$ descends to a connection $\nabla_{A}$ on 
$A$.  Let $\nabla_{A}(a) = \a\otimes a$ for some 
$\a\in \Omega^{1}(X^{[4]})$.  Then we have $\hat{a}-\rho_{1}
= \pi_{Y_{1234}}^{*}\a$.  
Hence we get 
\begin{equation} 
\label{eq:8.2.1.6} 
\pi_{1}^{*}\mu - \pi_{2}^{*}\mu  + 
\pi_{3}^{*}\mu - \pi_{4}^{*}\mu = 
d\a.
\end{equation}   
If we can show $\d(\a) = 0$ then we are done.  
To do this we need to examine the situation 
over $X^{[5]}$.   
$\bar{m}$ induces five bundle gerbe morphisms 
$\bar{M}_{i}:Q_{12345}\to Q_{15}$ for 
$i = 1,\ldots,5$.  We can form five $\cstar$ 
bundles $\hat{D}_{M_1 , M_2}=(M_{1},M_{2})^{-1}Q_{15}$, 
$\hat{D}_{M_{2},M_{3}}=(M_{2},M_{3})^{-1}Q_{15}$,\ldots, 
$\hat{D}_{M_{5},M_{1}}=(M_{5},M_{1})^{-1}Q_{15}$ 
on $Y_{12345}$, all of which descend to $\cstar$ 
bundles $D_{\bar{M}_{1},\bar{M}_{2}}$,$D_{\bar{M}_{2},\bar{M}_{3}}$,\ldots, 
$D_{\bar{M}_{5},\bar{M}_{1}}$ on $X^{[5]}$.  Since each pair of 
maps $(M_{i},M_{i+1})$ has $M_{i}$ factoring through 
$m_{1}$ and $M_{i+1}$ factoring through $m_{2}$ or 
vice versa, we can identify $D_{\bar{M}_{1},\bar{M}_{2}}$ 
with $\pi_{5}^{-1}A$, $D_{\bar{M}_{2},\bar{M}_{3}}$ with 
$\pi_{3}^{-1}A$, $D_{\bar{M}_{3},\bar{M}_{4}}$ with 
$\pi_{1}^{-1}A$, $D_{\bar{M}_{4},\bar{M}_{5}}$ with $\pi_{4}^{-1}A^{*}$ 
and $D_{\bar{M}_{5},\bar{M}_{1}}$ with $\pi_{2}^{-1}A^{*}$ as 
explained in Section~\ref{sec:7.1}.  Also note that the 
associator section $\hat{a}$ induces sections 
$\hat{a}_{ij}$ of $\hat{D}_{M_{i},M_{j}}$.  

Note that the bundle gerbe product in $Q$ supplies 
an isomorphism 
\begin{equation}
\label{eq:8.2.1.7} 
\hat{D}_{M_{1},M_{2}}\otimes \hat{D}_{M_{2},M_{3}}
\otimes \cdots \otimes \hat{D}_{M_{5},M_{1}} 
\to \hat{D}_{M_{1},M_{1}} 
\end{equation}   
and the coherency condition on $\hat{a}$ is that 
$\hat{a}_{12}\otimes \hat{a}_{23}\otimes \cdots \otimes \hat{a}_{51}$ 
is mapped to the identity section of 
$\hat{D}_{M_{1},M_{1}} = (M_{1},M_{1})^{-1}Q_{15}$.  

Let $\nabla_{15}$ denote the bundle gerbe 
connection on $Q_{15}$ induced by $\nabla$.  
Then $\nabla_{15}$ induces a connection 
$\hat{\nabla}_{ij}$ on the (line) bundle 
$\hat{D}_{M_{i},M_{j}}$ by pullback: 
$\hat{\nabla}_{ij} = (M_{i},M_{j})^{-1}\nabla_{15}$.  
Since $\nabla_{15}$ is a bundle gerbe connection, 
under the isomorphism~\ref{eq:8.2.1.7} above, the connection 
$\hat{\nabla}_{12}+\cdots +\hat{\nabla}_{51}$ on 
$\hat{D}_{M_{1},M_{2}}\otimes \cdots \otimes \hat{D}_{M_{5},M_{1}}$ 
is mapped to the connection 
$\hat{\nabla}_{11} = (M_{1},M_{1})^{-1}\nabla_{15}$ 
on $\hat{D}_{M_{1},M_{1}}$.  Hence we have that 
\begin{equation} 
\label{eq:8.2.1.8} 
(\hat{\nabla}_{12} + \cdots + \hat{\nabla}_{51})
(\hat{a}_{12}\otimes \cdots \otimes \hat{a}_{51}) = 
0.
\end{equation} 
We now need to calculate $\hat{\nabla}_{ij}(\hat{a}_{ij})$. 
Using Lemmas~\ref{lemma:8.1.3}--~\ref{lemma:8.1.6} 
one finds 
\begin{eqnarray*} 
\hat{\nabla}_{12}(\hat{a}_{12}) & = & (\pi^{*}\pi_{1}^{*}\a + 
1\times \rho_{1} + (m_{1}\times 1)^{*}\rho - (m_{2}\times 1)^{*}
\rho )\otimes \hat{a}_{12}                                       \\ 
\hat{\nabla}_{23}(\hat{a}_{23}) & = & (\pi^{*}\pi_{3}^{*}\a + 
(1\times m\times 1)^{*}\rho_{1})\otimes \hat{a}_{23}             \\ 
\hat{\nabla}_{34}(\hat{a}_{34}) & = & (\pi^{*}\pi_{5}^{*}\a + 
(1\times m_{1})^{*}\rho - (1\times m_{2})^{*}\rho)\otimes \hat{a}_{34} \\ 
\hat{\nabla}_{45}(\hat{a}_{45}) & = & (-\pi^{*}\pi_{2}^{*}\a - 
(1\times 1\times m)^{*}\rho_{1})\otimes \hat{a}_{45}              \\ 
\hat{\nabla}_{51}(\hat{a}_{51}) & = & (-\pi^{*}\pi_{4}^{*}\a - 
(m\times 1\times 1)^{*}\rho_{1})\otimes \hat{a}_{51}.            \\ 
\end{eqnarray*} 
Here we have abused notation and denoted by $\pi$ 
the projection $\pi_{Y_{12345}}:Y_{12345}\to X^{[5]}$.  If one does 
the calculation then one finds that 
\begin{equation} 
\label{eq:8.2.1.9}  
(\hat{\nabla}_{12} + \cdots + \hat{\nabla}_{51})
(\hat{a}_{12}\otimes \cdots \otimes \hat{a}_{51}) = 
\d(\a)\otimes (\hat{a}_{12}\otimes \cdots \otimes \hat{a}_{51}) 
\end{equation} 
and so we conclude that $\d(\a) = 0$, completing 
the proof.     
\end{proof}       

Let us explicitly write down the 
bundle 2-gerbe connection $(\nabla,f)$ on $(Q,Y,X,M)$ 
constructed in the above proof.  Since $\d(\a)= 0$, we can 
solve $\a=\d(\b)$ for some $\b\in \Omega^{1}(X^{[3]})$.  
Therefore we have $\d(\mu) = \d(d\b)$ and so we can 
find $\nu\in \Omega^{2}(X^{[2]})$ such that $\mu = 
d\b + \d(\nu)$.  Now give the bundle gerbe $(Q,Y,X^{[2]})$ 
the bundle gerbe connection $\nabla$ with curving 
$f^{'} = f-\pi_{Y}^{*}\nu$.  From equation~\ref{eq:8.2.1.4} we have 
\begin{equation} 
\label{eq:8.2.1.10} 
p_{1}^{*}\pi_{1}^{-1}f^{'} + 
p_{2}^{*}\pi_{3}^{-1}f^{'} = 
m^{*}\pi_{2}^{-1}f^{'} + d\rho 
+ \pi_{Y_{123}}^{*}(\mu - \d(\nu)) = 
m^{*}\pi_{2}^{-1}f^{'} + \pi_{Y_{123}}^{*}d\b. 
\end{equation} 
This equation implies that the three curvature 
$\omega^{'} = \omega -d\nu$ of the bundle gerbe 
$(Q,Y,X^{[2]})$ with bundle gerbe connection $\nabla$ 
and curving $f^{'}$ satisfies $\d(\omega^{'}) = 0$.  
Therefore $(\nabla,f^{'})$ is a bundle 2-gerbe 
connection for $Q$.

\setcounter{chapter}{10}
\chapter{Comparison of the \v{C}ech, Deligne and de Rham classes 
associated to a bundle 2-gerbe}
\label{chapter:9}
\section{A \v{C}ech three class}
\label{sec:9.1} 
Let $(Q,Y,X,M)$ be a bundle 2-gerbe 
(Definition~\ref{def:7.1.2}).
We will show how to construct a 
$\cstar$-valued \v{C}ech 3-cocycle associated
to $Q$.

Choose an open covering $\{U_{i}\}_{i \in I}$
of $M$, all of whose finite intersections
are empty or contractible, and such that
there exist local sections 
$s_{i}:U_{i}\to X|_{U_{i}}$ of $\pi:X\to M$.  We have 
the usual maps $(s_{i},s_{j}):U_{ij} \to X^{[2]}$,
$m \mapsto (s_{i}(m),s_{j}(m))$.  Let 
$(Q_{ij},Y_{ij},U_{ij})$ denote the pullback
of the bundle gerbe $(Q,Y,X^{[2]})$ to $U_{ij}$ by
$(s_{i},s_{j})$.  Thus $Y_{ij}$ is the space
over $U_{ij}$ whose fibre at $m \in U_{ij}$
is $Y_{ij_{m}} = Y_{(s_{i}(m),s_{j}(m))}$.

Since $U_{ij}$ is contractible, the bundle 
gerbe $(Q_{ij},Y_{ij},U_{ij})$ is trivial.
Hence there is a $\cstar$-bundle 
$L_{ij} \to Y_{ij}$ such that 
$\delta(L_{ij}) = \pi_1^{-1}L_{ij}\otimes 
\pi_2^{-1}L_{ij}^{*}\simeq Q_{ij}$
on $Y_{ij}^{[2]}$.

Over $U_{ijk}$ the bundle gerbe morphism 
$\bar{m} = (\hat{m},m,\text{id})$ induces 
a map $m:Y_{ijk} = Y_{jk}\times_{M}Y_{ij}\to Y_{ik}$.  
Let $\hat L_{ijk} = L_{jk}\otimes m^{-1}
L_{ik}^{*}\otimes L_{ij}$.
Notice that we have 
$\delta(\hat L_{ijk}) \simeq Q_{jk}\otimes 
(m^{[2]})^{-1}Q_{ik}^{*}\otimes Q_{ij}$,
which is trivial.  In fact, the $\cstar$-bundle
$\hat L_{ijk}$ descends to a bundle 
$L_{ijk} \to U_{ijk}$.  To see this we 
need the following lemma.

\begin{lemma}
Suppose $(P,X,M)$ and $(Q,Y,M)$ are 
bundle gerbes with a bundle gerbe 
morphism $(\hat{f},f,\text{id}):P \to Q$.
If $P$ and $Q$ are both trivial, so there
exist $L \to X$, $J \to Y$ such that
$P \simeq \delta(L)$, $Q \simeq \delta(J)$,
then the bundle $L\otimes f^{-1}J^{*}$
over $X$ descends to $M$.
\end{lemma}

\begin{proof}  $\hat{f}:P \to Q$ induces 
an isomorphism $\delta(L) \isom \delta(f^{-1}J)$
as pictured in the following diagram:
$$
\xymatrix{   
P \ar[r]^-{\simeq} \ar[d] & (f^{[2]})^{-1}Q 
\ar[d]                                     \\ 
\d(L) \ar[r]^-{\simeq} & \d(f^{-1}J)           } 
$$ 
This diagram lives over $X^{[2]}$.  Clearly
over $X^{[3]}$ we get the following 
diagram, induced from the one above by
tensor product and pullback,
$$ 
\xymatrix{  
\d(P) \ar[r]^-{\simeq} \ar[d]_{\simeq} 
& \d((f^{[2]})^{-1}Q) \ar[d]^{\simeq}   \\ 
\d\d(L) \ar[r]^-{\simeq} & \d\d(f^{-1}J) } 
$$ 
Let $s_{P}$ and $s_{Q}$ denote the 
sections of $\d(P)$ and $\d(Q)$ induced 
by the bundle gerbe products respectively.  
In the second diagram above, the 
isomorphism $\d(P)\to \d((f^{[2]})^{-1}Q)$ 
maps $s_{P}$ to $s_{(f^{[2]})^{-1}Q} = 
(f^{[3]})^{-1}s_{Q}$, since this isomorphism 
is induced by a bundle gerbe morphism.
Similarly, the vertical isomorphism 
$\d(P)\to \d\d(L)$ maps $s_{P}$ to 
$1_{L}$ and the vertical 
isomorphism $\d((f^{[2]})^{-1}Q)\to \d\d(f^{-1}J)$ 
maps $(f^{[3]})^{-1}s_{Q}$ to 
$(f^{[3]})^{-1}1_{J}$.  
Here $1_{L}$ and 
$1_{J}$ are the sections 
induced by the canonical trivialisations of 
$\d\d(L)$ and $\d\d(J)$ respectively.  
It follows now that if we define a section 
$t$ of $\d(L\otimes f^{-1}J^{*})$ using 
the isomorphism $\d(L)\to \d(f^{-1}J)$ 
from the first diagram above, then $t$ 
satisfies $\d(t) = 1$ and hence 
$L\otimes f^{-1}J^{*}$ descends to $M$.  
\end{proof} 

It is easy to see from this that 
$\hat L_{ijk}$ descends to some $\cstar$-bundle 
$L_{ijk} \to U_{ijk}$.
Next, over $U_{ijkl}$ we have two bundle 
gerbe morphisms $\bar{m}_{1}$ and 
$\bar{m}_{2}$ obtained from the bundle 
gerbe morphism 
$\bar{m}=(\hat{m},m,\text{id}):\pi_{1}^{-1}
Q\otimes \pi_{3}^{-1}Q\to \pi_{2}^{-1}Q$ over $X^{[3]}$:
$$
\bar{m}_{1},\bar{m}_{2}: 
Q_{kl}\otimes Q_{jk}\otimes Q_{ij} \to Q_{ik} 
$$
The $\cstar$ bundle $(m_{1},m_{2})^{-1}Q_{il}$ 
on $Y_{ijkl} = Y_{kl}\times_{M}Y_{jk}\times_{M}Y_{ij}$ 
descends to a $\cstar$-bundle $A_{ijkl}$
on $U_{ijkl}$.  It is clear from  
Definition~\ref{def:7.1.2}    
that $A_{ijkl}$ is canonically trivial with 
a section $a_{ijkl} = a(s_{i},s_{j},s_{k},s_{l})$
which clearly satisfies $\d(a)_{ijklm} = 1$,
where $1$ is the section of $\d(A_{ijkl})$ 
on $U_{ijklm}$ induced by the canonical section 
of $\d(A) \to X^{[5]}$. We need 
to relate $(m_{1},m_{2})^{-1}Q_{il}$ to the 
principal $\cstar$ bundles $L_{ijk}$.   
The map $m_{1}$ is defined by 
composition 
as in the following diagram.
$$ 
\xymatrix{ 
Y_{ijkl} 
\ar[r]^-{m\times 1} 
& Y_{ijl} 
\ar[r]^-{m}  & 
Y_{il}                                 }  
$$ 
Therefore
\begin{eqnarray*}
&   &  L_{kl}\otimes L_{jk}\otimes L_{ij}\otimes m_{2}^{-1}L_{il}^{*} \\ 
& = &  L_{kl}\otimes L_{jk}\otimes L_{ij}\otimes 
       (m\circ (m\times 1))^{-1}L_{il}^{*}         \\ 
& = &  L_{kl}\otimes L_{jk}\otimes L_{ij}\otimes 
       (m\times 1)^{-1}m^{-1}L_{il}^{*}             \\ 
& \simeq &  L_{kl}\otimes L_{jk}\otimes L_{ij}\otimes 
       (m\times 1)^{-1}[\pi_{Y_{ijl}}^{-1}L_{ijl}\otimes L_{jl}^{*}
       \otimes L_{ij}^{*}]                            \\ 
& = &  L_{kl}\otimes m^{-1}L_{jl}^{*}\otimes L_{jk}
       \otimes L_{ij}\otimes L_{ij}^{*}\otimes 
       \pi_{Y_{ijkl}}^{-1}L_{ijl}                              \\ 
& \simeq &  \pi_{Y_{ijkl}}^{-1}(L_{ijl}\otimes L_{jkl}).            
\end{eqnarray*}  
The map $m_{2}$ is induced by composition as 
indicated in the following diagram. 
$$
\xymatrix{ 
Y_{ijkl} 
\ar[r]^-{1\times m} & Y_{ikl}
 \ar[r]^-{m} & Y_{il}          }  
$$ 
Then 
\begin{eqnarray*} 
&    & L_{kl}\otimes L_{jk}\otimes L_{ij}\otimes 
m_{2}^{-1}L_{il}^{*}                                \\ 
& = & L_{kl}\otimes L_{jk}\otimes L_{ij}\otimes 
(m\circ (1\times m))^{-1}L_{il}^{*}                  \\ 
& = & L_{kl}\otimes L_{jk}\otimes L_{ij}\otimes 
(1\times m)^{-1}m^{-1}L_{il}^{*}                     \\ 
& \simeq & L_{kl}\otimes L_{jk}\otimes (1\times m)^{-1}
[\pi_{Y_{ikl}}^{-1}L_{ikl}\otimes L_{kl}^{*}L_{ik}^{*}]         \\ 
& = & L_{kl}\otimes L_{kl}^{*}\otimes L_{jk}\otimes 
m^{-1}L_{ik}^{*}\otimes L_{ij}\otimes \pi_{Y_{ijkl}}^{-1}L_{ikl} \\ 
& \simeq & \pi_{Y_{ijkl}}^{-1}(L_{ijk}\otimes 
L_{ikl}).                                              
\end{eqnarray*} 
Therefore we get 
\begin{eqnarray*} 
&   &  (L_{kl}\otimes L_{jk}\otimes L_{ij}\otimes 
       m_{2}^{-1}L_{il}^{*})^{*}\otimes (L_{kl}\otimes 
       L_{jk}\otimes L_{ij}\otimes m_{1}^{-1}L_{il}^{*}
       )                                           \\ 
& \simeq &  m_{2}^{-1}L_{il}\otimes m_{1}^{-1}L_{il}^{*}    \\ 
& = &  (m_{1},m_{2})^{-1}\d(L_{il})                    \\ 
& \simeq &  (m_{1},m_{2})^{-1}Q_{il}                          \\ 
& \simeq &  \pi_{Y_{ijkl}}^{-1}A_{ijkl}.                                   
\end{eqnarray*} 
Therefore we must have 
$$
L_{jkl}\otimes L_{ikl}^{*}\otimes L_{ijl}\otimes 
L_{ijk}^{*} \simeq A_{ijkl} 
$$
over $U_{ijkl}$.  
Now choose sections $\s_{ijk} \in \Gamma(U_{ijk},L_{ijk})$ 
and define $g_{ijkl}:U_{ijkl} \to \cstar$ 
by 
$$
\s_{jkl}\otimes \s_{ikl}^{*} \otimes 
\s_{ijl}\otimes \s_{ijk}^{*}\cdot g_{ijkl}  = 
a_{ijkl}.
$$
We have the following proposition. 
\begin{proposition} 
\label{prop:9.1} 
$g_{ijkl}$ satisfies the \v{C}ech 
3-cocycle condition 
$$ 
g_{jklm}g_{iklm}^{-1}g_{ijlm}g_{ijkm}^{-1}
g_{ijkl} = 1, 
$$
and hence is a representative   
of a class in $\check{H}^{3}(M;\underline{\C}^{\times}_{M}) = H^{4}(M;\Z)$.
\end{proposition} 

\begin{proof} 
All that one really needs to do 
is to check that $g_{ijkl}$ does 
indeed define a cocycle.  This 
follows from the coherency 
condition on $a$.  One can also 
check that another choice of the 
sections $s_{i}$, or the line 
bundles $L_{ij}$, or the sections 
$\sigma_{ijk}$ changes $g_{ijkl}$ 
by a coboundary.    
\end{proof} 

There is another method of calculating 
the \v{C}ech cocycle $g_{ijkl}$ which is 
similar in spirit to the method used 
to calculate the \v{C}ech representative 
of the Dixmier-Douady class of a 
bundle gerbe.    
Let $(Q,Y,X,M)$ be a bundle 
2-gerbe.  Choose an open cover 
$\{U_{i}\}_{i \in I}$ of $M$ such that there 
exist local sections $s_{i}:U_{i}\to X$ of 
$\pi:X\to M$ and such that each non-empty 
finite intersection $U_{i_{1}}\cap \cdots \cap U_{i_{k}}$ 
is contractible.  As above, form the maps 
$(s_{i},s_{j}):U_{ij} \to X^{[2]}$ and 
let $(Q_{ij},Y_{ij},U_{ij})$ denote the 
pullback of the bundle gerbe $(Q,Y,X^{[2]})$ 
with the map $(s_{i},s_{j})$.   
So $Y_{ij}$ 
has fibre $Y_{ij_{m}}$ at 
$m\in U_{ij}$ equal to 
$Y_{(s_{i}(m),s_{j}(m))}$ and there 
is an induced projection 
$\pi_{Y_{ij}}:Y_{ij}\to U_{ij}$.                           

In certain circumstances, for 
instance if $\pi_{Y}:Y\to X^{[2]}$ 
is a fibration, then one can choose 
sections $\sigma_{ij}:U_{ij}\to Y_{ij}$ 
of $\pi_{Y_{ij}}$.  We will assume 
that this is the case (in general 
one would only be able to choose an 
open cover $\{U_{ij}^{\a}\}_{\a\in 
\Sigma_{ij}}$ of $U_{ij}$ such that there 
exist local sections $\sigma_{ij}^{\a}$ 
of $\pi_{Y_{ij}}:Y_{ij}\to U_{ij}$).  
Denote, for ease of notation, 
$m(\sigma_{jk},\sigma_{ij})$ by 
$\sigma_{jk}\circ \sigma_{ij}$.  
Then we have a map $U_{ijk} \to Y_{ik}^{[2]}$ 
which sends $m\in U_{ijk}$ to 
$(\s_{ik}(m),(\s_{jk}\circ \s_{ij})(m))$.  Let 
$Q_{ijk} = (\s_{ik},\s_{jk}\circ \s_{ij})^{-1}Q_{ik}$ 
denote the pullback bundle on 
$U_{ijk}$.  Choose a non-vanishing section 
$\rho_{ijk}$ of 
$Q_{ijk}\to U_{ijk}$. 
Define a non-zero section $s_{ijkl}$ of the line bundle 
$(\s_{il},(\s_{kl}\circ \s_{jk})\circ \s_{ij})^{-1}Q_{il}$ 
over $U_{ijkl}$ by    
\begin{equation} 
\label{eq:9.1} 
s_{ijkl} = m_{Q}(\hat{m}(\rho_{jkl}\otimes e(\s_{ij}))\otimes 
\rho_{ijl}). 
\end{equation} 
Here $m_{Q}$ denotes the bundle gerbe 
multiplication on $Q$.  We can also define 
a non-zero section $t_{ijkl}$ of the line bundle 
$(\s_{il},\s_{kl}\circ(\s_{jk}\circ \s_{ij}))^{-1}Q_{il}$ 
by putting $t_{ijkl}$ equal to 
\begin{equation} 
\label{eq:9.2} 
m_{Q}(\hat{m}(e(\s_{kl})\otimes \rho_{ijk})\otimes 
\rho_{ikl}).
\end{equation} 
Therefore we can use the associator 
section to define another section of 
$(\s_{il},(\s_{kl}\circ \s_{jk})\circ \s_{ij})^{-1}Q_{il}$ by 
$$
m_{Q}(a(\s_{kl},\s_{jk},\s_{ij})\otimes s_{ijkl}).
$$
Thus we can define a function 
$\epsilon_{ijkl}:U_{ijkl}\to \cstar$ which 
satisfies 
\begin{equation} 
\label{eq:9.3} 
t_{ijkl} = m_{Q}(a(\s_{kl},\s_{jk},\s_{ij})\otimes 
s_{ijkl})\cdot \epsilon_{ijkl}.
\end{equation} 
One can show that $\epsilon_{ijkl}$ 
satisfies the cocycle condition
$$
\epsilon_{jklm}\epsilon_{iklm}^{-1}\epsilon_{ijlm}
\epsilon_{ijkm}^{-1}\epsilon_{ijkl} = 1.
$$
We have the following lemma.  
\begin{lemma} 
\label{lemma:9.2} 
The cocycle $\epsilon_{ijkl}$ defined 
above gives rise to the same class 
in the \v{C}ech cohomology group 
$\check{H}^{3}(M;\underline{\C}^{\times}_{M})$ 
as the cocycle $g_{ijkl}$ 
of Proposition~\ref{prop:9.1} 
\end{lemma} 

\begin{proof} 
We have $(\s_{jk},\s_{ij})^{-1}\hat{L}
_{ijk} \simeq L_{ijk}$.  But 
$\hat{L}_{ijk} = L_{jk}\otimes m^{-1}
L_{ik}^{*}\otimes L_{ij}$ and so 
$(\s_{jk},\s_{ij})^{-1}\hat{L}_{ijk} = 
\s_{jk}^{-1}L_{jk}\otimes (\s_{jk}\circ 
\s_{ij})^{-1}L_{ik}^{*}\otimes \s_{ij}^{-1}
L_{ij}$.  Hence we get 
$(\s_{jk},\s_{ij})^{-1}L_{ik} = 
L_{ijk}\otimes \s_{jk}^{-1}L_{jk}^{*}\otimes 
\s_{ij}^{-1}L_{ij}^{*}$.  Also we have 
\begin{eqnarray*} 
&    &  (\s_{jk}\circ \s_{ij},\s_{ik})^{-1}Q_{ik}        \\ 
& \simeq & (\s_{jk}\circ \s_{ij},\s_{ik})^{-1}\d(L_{ik}) \\ 
& =  & \s_{ik}^{-1}L_{ik}\otimes (\s_{jk}\circ 
\s_{ij})^{-1}L_{ik}^{*}                                    \\ 
& \simeq & L_{ijk}\otimes (\s_{jk}^{-1}L_{jk}\otimes 
\s_{ik}^{-1}L_{ik}^{*}\otimes \s_{ij}^{-1}L_{ij})^{*}.     
\end{eqnarray*} 
Therefore we can construct a section 
$\rho_{ijk}$ of $(\s_{jk}\circ \s_{ij},
\s_{ik})^{-1}Q_{ik}$ by taking the section 
$\s_{ijk}$ of $L_{ijk}$ and choosing sections 
$s_{ij}$ of $\s_{ij}^{-1}L_{ij}$ and forming the 
tensor product section 
$\s_{ijk}\otimes (s_{jk}\otimes s_{ik}^{*}
\otimes s_{ij})^{*}$.   
From here it is not too hard to show 
that the cocycle $\epsilon_{ijkl}$ is 
cohomologous to $g_{ijkl}$.  
\end{proof} 
  
\begin{example} 
Consider the lifting bundle 2-gerbe $(Q,\tilde{P},P,M)$  
of Proposition~\ref{prop:11.1.1} associated 
to a principal $G$ bundle $P\to M$ where $G$ 
is part of a central extension of groups 
$B\cstar\to \hat{G}\to G$.  We associate a 
$\cstar$ valued \v{C}ech 3-cocycle $g_{ijkl}:U_{ijkl}\to \cstar$ 
to $Q$ as follows.  Choose an open covering 
$\{U_{i}\}_{i\in I}$ of $M$ all of whose 
finite intersections are contractible or empty and 
such that $P$ is locally trivialised over $U_{i}$ --- 
say there exist sections $s_{i}$ of 
$P\to M$ over $U_{i}$.  Now let $\sigma_{ij}$ be 
sections of the pullback $B\cstar$ bundle $\tilde{P}_{ij} = 
(s_{i},s_{j})^{-1}\tilde{P}$ over $U_{ij}$.  We 
can now define a \v{C}ech 2-cocycle $g_{ijk}:U_{ijk}\to B\cstar$ 
representing the class in $\check{H}^{2}(M;\underline{B\C}^{\times}_{M})$ 
associated to the $B\cstar$ bundle gerbe $(\tilde{P},P,M)$ 
by $m_{\tilde{P}}(\sigma_{jk}\otimes_{B\cstar}\sigma_{ij}) = 
\sigma_{ik}\cdot g_{ijk}$.  It follows that the 
image of the class $[g_{ijk}]$ defined by 
$g_{ijk}$ under the coboundary map 
$H^{2}(M;\underline{B\C}^{\times}_{M})\to 
H^{3}(M;\underline{\C}_{M}^{\times})$ is the 
class defined by the \v{C}ech 3-cocycle 
$g_{ijkl}:U_{ijkl}\to \cstar$.  $g_{ijkl}$ would 
then be defined by choosing lifts 
$\hat{g}_{ijk}:U_{ijk}\to E\cstar$ of the $g_{ijk}$ and 
then setting $g_{ijkl} = \hat{g}_{jk}\hat{g}_{ik}^{-1}
\hat{g}_{ij}$.  
\end{example} 

We can also construct a \v{C}ech 3-cocycle 
$g_{ijkl}$ representing a class in 
$\check{H}^{3}(M;\underline{\C}^{\times}_{M})= H^{4}(M;\Z)$ 
associated to a stable bundle 2-gerbe.  
We record this in the following 
proposition. 

\begin{proposition} 
\label{prop:9.3} 
Let $(Q,Y,X,M)$ be a stable bundle 2-gerbe.  
Then there is a \v{C}ech 3-cocycle $g_{ijkl}$ 
representing a class in $\check{H}^{3}(M;
\underline{\C}^{\times}_{M})$ associated to the stable bundle 
2-gerbe $Q$.   
\end{proposition} 
  
The class $g_{ijkl}$ is constructed as 
follows.  Choose an open cover $\U = \{U_{i}
\}_{i\in I}$ of $M$ all of whose 
intersections are contractible and such that 
there exist local sections $s_{i}:U_{i}\to 
X$ of $\pi:X\to M$.  Form the maps 
$(s_{i},s_{j}):U_{ij}\to X^{[2]}$ and 
let $(Q_{ij},Y_{ij},U_{ij})$ denote 
the pullback bundle gerbe $(s_{i},s_{j})
^{-1}Q$.  Since $U_{ij}$ is contractible 
there is a trivialisation $L_{ij}$ of 
$Q_{ij}$.  Also the stable morphism 
$m:\pi_{1}^{-1}Q\otimes \pi_{3}^{-1}Q
\to \pi_{2}^{-1}Q$ induces a stable 
morphism $m_{ijk}:Q_{jk}\otimes Q_{ij}
\to Q_{ik}$ corresponding to a 
trivialisation $(L_{m_{ijk}},\phi_{m_{ijk}})$ of 
$Q_{jk}\otimes Q_{ik}^{*}\otimes Q_{ij}$.  
Therefore there is a $\cstar$ bundle 
$L_{ijk}$ on $U_{ijk}$ such that 
$L_{m_{ijk}} \simeq L_{jk}\otimes 
L_{ik}^{*}\otimes L_{ij}\otimes 
\pi^{-1}L_{ijk}$, where $\pi:Y_{jk}
\times_{M}Y_{ik}\times_{M}Y_{ij}\to 
U_{ijk}$ denotes the projection.  It is 
not hard to see that there is an 
isomorphism $L_{jkl}\otimes L_{ikl}^{*}
\otimes L_{ijl}\otimes L_{ijk}^{*}\simeq 
A_{ijkl}$ of $\cstar$ bundles on $U_{ijkl}$ 
(here $A_{ijkl}$ denotes the pullback 
$\cstar$ bundle $(s_{i},s_{j},s_{k},s_{l})
^{-1}A$).  Therefore if we choose 
sections $\sigma_{ijk}$ of the 
bundles $L_{ijk}$, then we can define 
a \v{C}ech 3-cocycle $g_{ijkl}:U_{ijkl}
\to \cstar$ by 
$$
\sigma_{jkl}\otimes \sigma_{ikl}^{*}
\otimes \sigma_{ijl}\otimes \sigma_{ijk}^{*} 
= a_{ijkl}\cdot g_{ijkl}, 
$$
where $a_{ijkl}$ denotes the pullback 
$(s_{i},s_{j},s_{k},s_{l})^{-1}a$ of 
the associator section $a$ of $A$.  
 
\section{A Deligne hypercohomology class}
\label{sec:9.2} 
Let $(Q,Y,X,M)$ be a bundle 2-gerbe.  Suppose 
the bundle gerbe $(Q,Y,X^{[2]})$ comes 
equipped with a bundle 2-gerbe connection $(\nabla,f)$ 
so that the three curvature form $\omega$ 
of the bundle gerbe $(Q,Y,X^{[2]})$ 
satisfies $\d(\omega) = 0$ in $\Omega^{3}(X^{[3]})$, 
as in Definition~\ref{def:8.1.1}.  
Let $K\in \Omega^{3}(X)$ be a 2-curving for 
$\omega$ so that $\d(K) = \omega$.  
We will construct a class $D(Q,\nabla,f,K)$ 
associated to the bundle 2-gerbe $Q$ 
with bundle 2-gerbe connection 
$(\nabla,f)$ and 2-curving $K$ 
in the Deligne hypercohomology group 
$H^{3}(M;\underline{\C}^{\times}_{M}\to 
\underline{\Omega}^{1}_{M}\to 
\underline{\Omega}^{2}_{M}\to 
\underline{\Omega}^{3}_{M})$.   

Choose an open covering $\{U_{i}\}_{i \in I}$ 
of $M$ such that each finite intersection 
$U_{i_{1}}\cap \cdots \cap U_{i_{p}}$ is either 
empty or contractible and such that there exist local 
sections $s_{i}:U_{i} \to X$ of $\pi$ for all $i \in I$.  
As in Section~\ref{sec:9.1}, 
let $(Q_{ij},Y_{ij},U_{ij})$ denote 
the bundle gerbe on $U_{ij}$ obtained by pulling back 
$Q$ with the map $(s_{i},s_{j}):U_{ij} \to X^{[2]}$.
Therefore $\nabla$ induces a bundle gerbe connection 
$\nabla_{ij}$ on $Q_{ij}$ by pullback.  Let 
$F_{\nabla_{ij}}$ denote the curvature of 
$Q_{ij}$ with respect to $\nabla_{ij}$.  The  
curving $f$ on $Y$ induces a curving $f_{ij}$ 
for $\nabla_{ij}$ by pullback with the  
map 
$Y_{ij} \to Y$.  Similarly the three curvature  
of $Q_{ij}$ with respect to the bundle gerbe 
connection $\nabla_{ij}$ and curving $f_{ij}$ is 
$\omega_{ij} = (s_{i},s_{j})^{*}\omega \in \Omega^{3}(U_{ij})$.   
Since $U_{ij}$ is contractible, $\omega_{ij}$ 
must be exact and so, using the same notation 
as in Section~\ref{sec:9.1}, let $L_{ij} \to Y_{ij}$ 
denote the $\cstar$-bundle which trivialises $Q_{ij}$.
As noted in Section~\ref{sec:5.1}, the 
isomorphism $H^{2}(M;\underline{\C}^{\times}_{M})
\cong H^{2}(M;\underline{\C}^{\times}_{M}\to 
\underline{\Omega}^{1}_{M})$ of truncated Deligne 
hypercohomology groups implies that a bundle 
gerbe with trivial Dixmier-Douady class 
can be endowed with a trivial bundle gerbe 
connection --- so that if $Q$ is a trivial bundle 
gerbe with a trivialisation $Q\simeq \d(L)$, 
then, given a bundle gerbe connection $\nabla_{Q}$ 
on $Q$, one may choose a connection $\nabla_{L}$ 
on $L$ so that the isomorphism $Q\simeq \d(L)$ 
is an isomorphism of bundle gerbes with 
connection: $(Q,\nabla_{Q})\simeq (\d(L),\d(\nabla_{L}))$.  
Therefore we may choose a connection 
$\tilde{\nabla}_{ij}$ on 
$L_{ij} \to Y_{ij}$ such that 
$(Q_{ij},\nabla_{ij}) \simeq (\d(L_{ij}),\d(\tilde{\nabla_{ij}}))$ 
is an isomorphism of bundle gerbes with 
connection on $Y_{ij}^{[2]}$.  Let 
$\tilde{F}_{ij}$ denote the curvature of 
$L_{ij}$ with respect to $\tilde{\nabla}_{ij}$.
Then we must have 
$F_{\nabla_{ij}} = \d(f_{ij}) = \d(\tilde{F}_{ij})$.  
Therefore there exists $\c_{ij} \in \Omega^{2}(U_{ij})$ 
such that 
\begin{equation}
\label{eq:9.2.1} 
\tilde{F}_{ij} = f_{ij} - \pi^{*}\c_{ij}.
\end{equation}
From this we deduce 
$0 = d\tilde{F}_{ij} = df_{ij} - \pi^{*}d\c_{ij}$ 
and hence 
\begin{equation}
\label{eq:9.2.2} 
\omega_{ij} = d\c_{ij}.
\end{equation} 
Next we examine the situation over $U_{ijk}$.  
Recall equation~\ref{eq:8.2.1.1} 
from Section~\ref{sec:8.2}, 
\begin{equation} 
\label{eq:9.2.3} 
(p_{1}^{[2]})^{-1}\pi_{1}^{-1}\nabla 
+ (p_{2}^{[2]})^{-1}\pi_{3}^{-1}\nabla 
= \tilde{m}^{-1}\circ (m^{[2]})^{-1}
\pi_{2}^{-1}\nabla\circ \tilde{m} + 
\d(\rho).  
\end{equation} 
Since $(\nabla,f)$ is a bundle  
2-gerbe connection and the associated 
three curvature $\omega$ of the bundle 
gerbe $(Q,Y,X^{[2]})$ satisfies 
$\d(\omega) = 0$, we must have 
\begin{equation} 
\label{eq:9.2.4} 
p_{1}^{*}\pi_{1}^{-1}f + p_{2}^{*}
\pi_{3}^{-1}f = m^{*}\pi_{2}^{-1}f + d\rho + \pi_{Y_{123}}^{*}\b, 
\end{equation} 
on $\pi_{1}^{-1}Y\times_{X^{[3]}}
\pi_{3}^{-1}Y$ (Compare with equation~\ref{eq:8.2.1.10}). 
Here $\b$ is a one form on $X^{[3]}$ such that 
$\nabla_{A}(a) = \d(\b)\otimes a$.  This formula pullsback to 
$Y_{ijk}= Y_{jk}\times_{M}Y_{ij}$ to give the following 
equation in $\Omega^{2}(Y_{ijk})$: 
\begin{equation} 
\label{eq:9.2.5} 
p_{1}^{*}f_{jk} + p_{2}^{*}f_{ij} = 
m^{*}f_{ik} + d\rho_{ijk} + \pi_{Y_{ijk}}^{*}d\b_{ijk}, 
\end{equation} 
where $\rho_{ijk}$ is the pullback 
of $\rho$ to $Y_{ijk}$  
via the inclusion map 
$Y_{jk}\times_{M}Y_{ij}\to \pi_{1}^{-1}
Y\times_{X^{[3]}}\pi_{3}^{-1}Y$ and 
$\b_{ijk} = (s_{i},s_{j},s_{k})^{*}\b$.  
Pulling equation~\ref{eq:9.2.3} 
back to $(Y_{ijk})^{[2]}$ 
gives 
\begin{equation} 
\label{eq:9.2.6} 
\nabla_{jk} + \nabla_{ij} = 
\tilde{m}^{-1}\circ (m^{[2]})^{-1}
\nabla_{ik}\circ \tilde{m} + \d(\rho_{ijk}). 
\end{equation} 
Recall also from Section~\ref{sec:9.1} 
that we defined $\hat{L}_{ijk} = 
L_{jk}\otimes m^{-1}L_{ik}^{*}\otimes 
L_{ij}$ on $Y_{ijk}$.  
$\hat{L}_{ijk}$ has the tensor product 
connection 
$$
\hat{\nabla}_{ijk} = \tilde{\nabla}_{jk} 
- m^{-1}\tilde{\nabla}_{ik} + \tilde{\nabla}_{ij}. 
$$
As explained in Section~\ref{sec:9.1}, 
$\hat{L}_{ijk}$ descends to a bundle 
$L_{ijk}$ on $U_{ijk}$ and we may choose 
a non-vanishing section $\s_{ijk}$ of 
$L_{ijk}\to U_{ijk}$.  From equation~\ref{eq:9.2.6} 
we get 
$$
\d(\tilde{\nabla}_{jk} + \tilde{\nabla}_{ij}) 
= \tilde{m}^{-1}\circ \d(m^{-1}\tilde{\nabla}_{ik})
\circ \tilde{m} + \d(\rho_{ijk}), 
$$
which implies that the connection 
$\tilde{\nabla}_{jk} - m^{-1}\tilde{\nabla}_{ik} 
+ \tilde{\nabla}_{ij} - \rho_{ijk}$ on 
$\hat{L}_{ijk}$ descends to a connection 
on $L_{ijk}$ --- call it $\nabla_{ijk}$.  Let 
$\nabla_{ijk}(\s_{ijk}) = \a_{ijk}\otimes \s_{ijk}$.  
Then we have 
$$
\pi_{Y_{ijk}}^{*}d\a_{ijk} = \tilde{F}_{jk} - m^{*}\tilde{F}_{ik} 
+ \tilde{F}_{ij} - d\rho_{ijk}.  
$$
Therefore from equations~\ref{eq:9.2.1} 
and~\ref{eq:9.2.5} we have $\pi_{Y_{ijk}}^{*}d\a_{ijk} = 
\pi_{Y_{ijk}}^{*}(d\b_{ijk} - \d(\c_{ij})$ and so 
\begin{equation} 
\label{eq:9.2.7} 
dA_{ijk} = \c_{jk} - \c_{ik} + \c_{ij},  
\end{equation} 
where we have put $A_{ijk} = \b_{ijk} - 
\a_{ijk}$.  One can check that the isomorphism 
$$
L_{jkl}\otimes L_{ikl}^{*}\otimes 
L_{ijl}\otimes L_{ijk}^{*} \simeq A_{ijkl} 
$$
maps the tensor product connection $\nabla_{jkl} 
- \nabla_{ikl} + \nabla_{ijl} - \nabla_{ijk}$ 
on $L_{jkl}\otimes L_{ikl}^{*}\otimes L_{ijl}
\otimes L_{ijk}^{*}$ to the pullback connection 
$(s_{i},s_{j},s_{k},s_{l})^{-1}\nabla_{A} = 
\nabla_{ijkl}$ on $A_{ijkl}$.  Since $\nabla_{A}(a) = \d(\b)$, 
$\nabla_{ijkl}(a_{ijkl}) = \d(\b_{ijk})$ and so we have 
$$
\a_{jkl} - \a_{ikl} + \a_{ijl} - \a_{ijk} + 
g_{ijkl}^{-1}dg_{ijkl} = \d(\b_{ijk}), 
$$
or $\d(A_{ijk}) = g_{ijkl}^{-1}dg_{ijkl}$.  
Therefore we have the set of equations 
\begin{eqnarray*} 
A_{jkl} - A_{ikl} + A_{ijl} - A_{ijk} & = & 
g_{ijkl}^{-1}dg_{ijkl}                            \\ 
dA_{ijk} & = & \c_{jk}-\c_{ik}+\c_{ij}           \\ 
d\c_{ij} & = & K_{j} - K_{i}.   
\end{eqnarray*} 
It follows that the quadruple $(g_{ijkl},A_{ijk},\c_{ij},
K_{i})$ defines a class in the Deligne 
cohomology group $H^{3}(M:\underline{\C}^{\times}_{M}\to 
\underline{\Omega}^{1}_{M}\to \underline{\Omega}^{2}_{M}
\to \underline{\Omega}^{3}_{M})$.  Hence we have 
the following proposition. 

\begin{proposition} 
Every bundle 2-gerbe $(Q,Y,X,M)$ equipped with 
a bundle 2-gerbe connection $(\nabla,f)$ and 
2-curving $K$ gives rise to a class $D(Q,\nabla,f,K)$ 
in the Deligne cohomology group $H^{3}(M;\underline{\C}
^{\times}_{M}\to \underline{\Omega}^{1}_{M}\to 
\underline{\Omega}^{2}_{M}\to \underline{\Omega}^{3}_{M})$.  
\end{proposition} 

In particular it follows from this by standard double 
complex arguments (see eg. \cite{BotTu}) that the 
class defined by $g_{ijkl}$ is the same as the class 
defined by $\frac{1}{2\pi \sqrt{-1}}\Theta$.      
       
\begin{note} 
It is also possible to construct singular analogues 
of the 2-curving and four class of a bundle 2-gerbe 
in a similar way to that used in Chapter~\ref{chapter:4}.  
\end{note} 

\section{The first Pontryagin class} 
\label{sec:9.3} 

Recall the bundle 2-gerbe 
$(\tilde{Q},\tilde{P},P,M)$ constructed in 
Proposition~\ref{prop:11.3.2} associated to a 
principal $G$ bundle where $G$ is a compact, 
simple, simply connected Lie group.  We have the following 
proposition.  

\begin{proposition}[\cite{BryMcL1},\cite{BryMcL3}]   
The class in $H^{4}(M;\Z)$ associated 
to the bundle 2-gerbe $\tilde{Q}$ is the transgression of 
$[\nu]$, that is the first Pontryagin 
class $p_{1}$ of $P$.  
\end{proposition} 

\begin{proof} 
We calculate the \v{C}ech four class of the 
bundle 2-gerbe $\tilde{Q}$ and 
get it exactly equal to the \v{C}ech cocycle 
defined by Brylinski and McLaughlin in 
\cite{BryMcL1} and \cite{BryMcL3}.  We can then 
apply Theorem 6.2 of \cite{BryMcL1} to conclude 
that this \v{C}ech four class is $p_{1}$.  We set 
out to calculate the \v{C}ech cocycle $g_{ijkl}$.  
First choose an open cover $\{U_{i}\}_{i\in I}$ 
of $M$ relative to which $\pi:P\to M$ has 
local sections $s_{i}$.  Since $\tilde{P}\to 
P^{[2]}$ is a fibration, we can calculate the 
\v{C}ech cocycle representing the four class of 
$\tilde{Q}$ using the second construction 
of a \v{C}ech cocycle as detailed in Section~\ref{sec:9.1}.   
Recall that we first choose 
sections  $\sigma_{ij}:U_{ij}\to \tilde{P}_{ij} 
= (s_{i},s_{j})^{-1}\tilde{P}$.  This is 
equivalent to choosing maps $\c_{ij}:U_{ij}
\times I\to G$ such that $\c_{ij}(m,0) = 
1$ and $\c_{ij}(m,1) = g_{ij}(m)$.  Next 
we have to choose sections $\rho_{ijk}:U_{ijk}
\to (\sigma_{ik},\sigma_{jk}\circ \sigma_{ij})
^{-1}\tilde{Q}_{ik}$.  This amounts to 
choosing maps $\c_{ijk}:U_{ijk}\times I
\times I\to G$ such that $\c_{ijk}(m,0,t) = 
\c_{ik}(m,t)$, $\c_{ijk}(m,1,t) = 
(\c_{ij}\circ g_{ij}\c_{jk})(m,t)$, 
$\c_{ijk}(m,s,0) = 1$ and $\c_{ijk}(m,s,1) 
= g_{ij}(m)g_{jk}(m)$ (This uses the fact that 
$\pi_{1}(G) = 0$ and the cocycle condition 
$g_{ij}g_{jk} = g_{ik}$ satisfied by the 
transition functions $g_{ij}$).  
Define the section $t_{ijkl}$ of the 
bundle $(\sigma_{il},\sigma_{kl}\circ 
(\sigma_{jk}\circ \sigma_{ij}))^{-1}\tilde{Q}
_{il}$ by $t_{ijkl} = m_{\tilde{Q}}(
\hat{m}(e(\sigma_{kl})\otimes \rho_{ijk})
\otimes \rho_{ikl})$.  If we denote a 
homotopy with end points fixed $\phi$ from a path $\a$ to a 
path $\b$ by a 2-arrow $\phi:\a\Rightarrow 
\b$ then we see that $\hat{m}(e(\sigma_{kl})
\otimes \rho_{ijk})$ amounts to 
composing horizontally the 2-arrows 
in the following diagram 
$$
\diagram 
1 \rtwocell_{\c_{ij}\circ 
g_{ij}\c_{jk}}^{\c_{ik}} & 
g_{ik} \rtwocell^{g_{ik}\c_{kl}}
_{g_{ik}\c_{kl}} & g_{il}. \\ 
\enddiagram 
$$
The bundle gerbe product 
$m_{\tilde{Q}}(\hat{m}(e(\sigma_{kl})
\otimes \rho_{ijk})\otimes \rho_{ikl})$ 
is obtained by composing the 2-arrows   
$\c_{ikl}:\c_{il}\Rightarrow 
\c_{ik}\circ g_{ik}\c_{kl}$ 
and $\c_{ijk}\circ g_{ik}\c_{kl}   
:\c_{ik}\circ 
g_{ik}\c_{kl}\Rightarrow 
(\c_{ij}\circ g_{ij}\c_{jk})
\circ g_{ik}\c_{kl}$ vertically.  
In a similar manner, one constructs the 
section $s_{ijkl} = m_{\tilde{Q}}(\hat{m}(
\rho_{jkl}\otimes e(\sigma_{ij}))\otimes 
\rho_{ijl})$ of the bundle $(\sigma_{il},
\sigma_{kl}\circ (\sigma_{jk}\circ \sigma_{ij}))
^{-1}\tilde{Q}_{il}$.  We make $t_{ijkl}$ 
into a section of $(\sigma_{il},(\sigma_{kl}
\circ \sigma_{jk})\circ \sigma_{ij})^{-1}
\tilde{Q}_{il}$ by using the associator 
section: $m_{\tilde{Q}}(a(\sigma_{kl},
\sigma_{jk},\sigma_{ij})\otimes t_{ijkl})$.    
Finally we define the cocycle $g_{ijkl}$ 
by $s_{ijkl} = m_{\tilde{Q}}(a(\sigma_{kl},
\sigma_{jk},\sigma_{ij})\otimes t_{ijkl})
\cdot g_{ijkl}$.  We can get an explicit formula 
for $g_{ijkl}$ as follows: we choose a homotopy 
with endpoints fixed $H_{ijkl}:U_{ijkl}\times 
I\times I\times I\to G$ such that $H_{ijkl}(m,0,s,t) = 
(\c_{ij}\circ g_{ij}\c_{jkl})\c_{ijl}(s,t)$, $H_{ijkl}
(m,1,s,t) =  \bar{a}(\c_{kl},\c_{jk},\c_{ij})
(\c_{ijk}\circ g_{ik}\c_{kl})\c_{ikl}(s,t)$, 
$H_{ijkl}(m,r,0,t) = \c_{il}(m,t)$, $H_{ijkl}
(m,r,1,t) = (\c_{ij}\circ g_{jk}\c_{jk})
\circ g_{ik}\c_{kl}$, $H_{ijkl}(m,r,s,0) = 1$  
and $H_{ijkl}(m,r,s,1) = g_{il}$ and we set 
$g_{ijkl} = \exp(\int_{I^{3}}H_{ijkl}^{*}\nu)$.  
This is just the integral of $\nu$ over the 
tetrahedron described in \cite{BryMcL1} and 
\cite{BryMcL3}.  Thus our cocycle must agree 
with that defined by Brylinski and 
McLaughlin.    
\end{proof}

\setcounter{chapter}{11}
\chapter{Trivial bundle 2-gerbes} 
\label{chapter:12} 

\section{Trivial bundle 2-gerbes} 
\label{sec:12.1} 
Let $(Q,Y,X,M)$ be a bundle 2-gerbe.  
Choose an open cover $\mathcal{U}= \{U_{i}\}_{i \in I}$ 
of $M$.  Form the sets 
$X_{i} = \pi^{-1}(U_{i})$.  Then 
the collection $\{X_{i}\}_{i \in I}$ 
forms an open cover of $X$.  
Define maps $\hat{s}_{i}:X_{i} \to X^{[2]}$ 
by sending $x \in X_{i}$ to 
$(x,s_{i}(m)) \in X^{[2]}$ for 
$s_{i}:U_{i}\to X$ a local 
section of $\pi$ and 
$m = \pi(x)$.  

Pullback the bundle gerbe 
$(Q,Y,X^{[2]})$ to $X_{i}$ 
with the map $\hat{s}_{i}$ 
and define $Q_{i} = \hat{s}_{i}^{-1}Q$ 
and $Y_{i} = \hat{s}_{i}^{-1}Y$.  
Thus $Y_{i}$ has fibre 
$Y_{i_{x}}$ at $x \in X_{i}$ 
equal to $Y_{(x,s_{i}(m))}$ where 
again $m = \pi(x)$. 

We will assume for the 
time being that we can choose sections 
$\sigma_{ij}:U_{ij} \to Y_{ij}$ 
where $Y_{ij} \to U_{ij}$ is 
defined above.  
Define maps $\phi_{ij}:Y_{i}\to Y_{j}$ 
covering the identity map on 
$X_{ij}$ by sending 
$y \in Y_{i}$ to 
$m(\sigma_{ij},y) \in Y_{j}$.  
The $\phi_{ij}$ extend to 
bundle gerbe morphisms 
$\bar{\phi}_{ij}:Q_{i}\to Q_{j}$, 
$\bar{\phi}_{ij} = (\hat{\phi}_{ij}
,\phi_{ij},\text{id})$ 
with $\hat{\phi}_{ij}(u_{i}) = 
\hat{m}(e(\sigma_{ij})\otimes u_{i})$ 
for $u_{i} \in Q_{(y_{1},y_{2})}$ 
where $(y_{1},y_{2}) \in Y_{i}^{[2]}$.  

Let $Z = \coprod_{i \in I} Y_{i}$ and 
let $\pi_{Z}:Z\to X$ denote 
the obvious projection.   
The fibre product $Z^{[2]} =Z\times_{X} Z$ 
is equal to the disjoint 
union $\coprod_{i,j}Y_{i}\times_{X}Y_{j}$.  
Define maps $\chi_{ij}:Y_{i}\times_{X}Y_{j}\to Y_{j}^{[2]}$ 
by $\chi_{ij}(y_{i},y_{j}) = (\phi_{ij}(y_{i}),y_{j})$, 
for $y_{i}\in Y_{i}$ and 
$y_{j}\in Y_{j}$.  Let $L_{ij}$ 
denote the pullback line bundle 
$\chi_{ij}^{-1}Q_{j}$ on 
$Y_{i}\times_{X}Y_{j}$.  Define 
a line bundle $L\to Z^{[2]}$ by 
setting $L = \coprod_{i,j}L_{ij}$ 
with projection map $L\to Z^{[2]}$ 
induced by the projections 
$L_{ij}\to Y_{i}\times_{X}Y_{j}$.   

Our aim is to show that the triple 
$(L,Z,X)$ forms a bundle gerbe.  To do 
this we need to define an 
associative product 
$m_{L}:\pi_{1}^{-1}L\otimes \pi_{3}^{-1}L\to \pi_{2}^{-1}L$ 
on $Z^{[3]}$.  Since 
$Z^{[3]} = \coprod_{i,j,k}Y_{i}\times_{X}Y_{j}\times_{X}Y_{k}$ 
this amounts to finding 
a bundle map $L_{jk}\otimes L_{ij} \to L_{ik}$ 
which satisfies the associativity 
condition over 
$Y_{i}\times_{X}Y_{j}\times_{X}Y_{k}\times_{X}Y_{l}$.  

Let $u_{jk}\in L_{jk_{(y_{j},y_{k})}}$, 
$u_{ij}\in L_{ij_{(y_{i},y_{j})}}$ 
with $y_{i}\in Y_{i_{x}}$, 
$y_{j}\in Y_{j_{x}}$ and 
$y_{k}\in Y_{k_{x}}$.  So 
$u_{jk}\in Q_{k_(\phi_{jk}(y_{j}),y_{k})}$, 
$u_{ij}\in Q_{j_{(\phi_{ij}(y_{i}),y_{j})}}$.  
Then 
$$ 
\hat{m}(e(\s_{jk})\otimes u_{ij})
\in Q_{k_{(\phi_{jk}\circ \phi_{ij}(y_{i}),\phi_{jk}(y_{j}))}}. 
$$
Therefore 
$$
m_{Q}(u_{jk}\otimes \hat{m}(e(\s_{jk})\otimes u_{ij})) 
\in Q_{(\phi_{jk}\circ \phi_{ij}(y_{i}),y_{k})}.
$$
Denote as before by $Q_{ijk}$ the 
line bundle $(\s_{ik},m(\s_{jk},\s_{ij}))^{-1}Q$ 
on $U_{ijk}$ and choose as 
before a non-vanishing 
section $\s_{ijk}$ of $Q_{ijk}$.  
Define a section $\hat{\xi}_{ijk}$ of the 
line bundle $(\phi_{ik},\phi_{jk}\circ \phi_{ij})^{-1}Q_{k}$ 
on $Y_{i}$ by 
$$
\hat{\xi}_{ijk}(y_{i}) = 
m_{Q}(\hat{a}(\s_{jk},\s_{ij},y_{i})\otimes 
\hat{m}(\s_{ijk}\otimes e(y_{i}))),
$$
where $y_{i} \in Y_{i}$ and 
where $\hat{a}$ denotes the section 
of the pullback line bundle 
$(m_{1},m_{2})^{-1}Q_{14}$ which is mapped 
to the pullback of the associator section $\pi_{1234}^{-1}a$  
under the isomorphism of line bundles 
$(m_{1},m_{2})^{-1}Q_{14}\cong \pi_{1234}^{-1}A$.  
Here $\pi_{1234}:Y_{1234} \to X^{[4]}$ denotes the 
projection.   

Finally we define for $u_{jk}$, 
$u_{ij}$ as above, 
\begin{equation} 
\label{eq:12.4} 
m_{L}(u_{jk}\otimes u_{ij}) = 
m_{Q}(m_{Q}(u_{jk}\otimes \hat{m}(e(\s_{jk})\otimes u_{ij}))
\otimes \hat{\xi}_{ijk}(y_{i})).  
\end{equation} 
We have the following proposition.  
\begin{proposition}  
\label{prop:12.1.1} 
Suppose $(Q,Y,X,M)$ is a bundle 2-gerbe 
with trivial four class and suppose 
$L$ and $Z$ are defined as above.  Then 
$m_{L}$ is associative and the 
triple $(L,Z,X)$ forms a bundle 
gerbe.
\end{proposition}  

This uses crucially the fact that 
the cocycle $g_{ijkl}$ is trivial, 
without this the proof will fail.  

\begin{proof}  We need to show that 
$m_{L}$ is associative, that is 
$m_{L}\circ (1\otimes m_{L}) = m_{L}\circ (m_{L}\otimes 1)$.  
We have 
\begin{eqnarray*} 
&  &  m_{L}(u_{kl}\otimes m_{L}(u_{jk}\otimes u_{ij}))     \\              
&= & (u_{kl}(e(\s_{kl})\circ ((u_{jk}(e(\s_{jk})\circ u_{ij})) 
    \hat{\xi}_{ijk}(y_{i})))) \hat{\xi}_{ikl}(y_{i}),  
\end{eqnarray*}
where we have denoted composition 
by bundle gerbe product by   
juxtaposition and 
composition via the bundle  
2-gerbe product with $\circ$.  
We also have 
\begin{eqnarray*}
&   & m_{L}(m_{L}(u_{kl}\otimes u_{jk})\otimes u_{ij})        \\
& = & (((u_{kl}(e(\s_{kl})\circ u_{jk}))\hat{\xi}_{jkl}(y_{j}))  
    (e(\s_{jl})\circ u_{ij})) \hat{\xi}_{ijl}(y_{i}).  
\end{eqnarray*}
We have to show these two expressions are the 
same.  For this we will need the 
following lemma.  
\begin{lemma} 
\label{lemma:12.1.2} 
Under the hypotheses of the 
above proposition, we have 
$$
(e(\s_{kl})\circ \hat{\xi}_{ijk}(y_{i}))\hat{\xi}_{ikl}(y_{i}) 
= \hat{\xi}_{jkl}(\s_{ij}\circ y_{i})\hat{\xi}_{ijl}(y_{i}).  
$$
\end{lemma} 
To prove this one needs to use 
the associativity isomorphism 
to regroup brackets and use the 
coherency condition on $\hat{a}$.  

The first of the above two expressions is equal 
to 
\begin{eqnarray*} 
& = & (u_{kl}((e(\s_{kl})e(\s_{kl}))\circ ((u_{jk}(e(\s_{jk})\circ 
 u_{ij}))\hat{\xi}_{ijk}(y_{i}))))\hat{\xi}_{ikl}(y_{i})    \\
& = & (u_{kl}((e(\s_{kl})\circ u_{jk})(e(\s_{kl})\circ 
 ((e(\s_{jk})\circ u_{ij})\hat{\xi}_{ijk}(y_{i})))))
\hat{\xi}_{ikl}(y_{i})                                    \\
& = & (u_{kl}(e(\s_{kl})\circ u_{jk}))((e(\s_{kl})\circ 
((e(\s_{jk})\circ u_{ij})\hat{\xi}_{ijk}(y_{i}))) 
\hat{\xi}_{ikl}(y_{i}))                                  \\ 
& = & (u_{kl}(e(\s_{kl})\circ u_{jk}))((e(\s_{kl})
e(\s_{kl}))\circ ((e(\s_{jk})\circ u_{ij}) 
\hat{\xi}_{ijk}(y_{i})))\hat{\xi}_{ikl}(y_{i})            \\ 
& = & (u_{kl}(e(\s_{kl})\circ u_{jk}))((e(\s_{kl})\circ 
(e(\s_{jk})\circ u_{ij}))(e(\s_{kl})\circ 
\hat{\xi}_{ijk}(y_{i})))\hat{\xi}_{ikl}(y_{i})             \\ 
& = & (u_{kl}(e(\s_{kl})\circ u_{jk}))((e(\s_{kl})\circ 
(e(\s_{jk})\circ u_{ij}))(\hat{\xi}_{jkl}(\s_{ij}\circ 
y_{i})\hat{\xi}_{ijl}(y_{i}))                             \\ 
& = & (u_{kl}(e(\s_{kl})\circ u_{jk}))\hat{a}(\s_{kl}, 
\s_{jk},y_{j})((e(\s_{kl})\circ e(\s_{jk}))\circ 
u_{ij})                                                   \\ 
&   & \hat{a}(\s_{kl},\s_{jk},\s_{ij}\circ y_{i})^{-1} 
\hat{\xi}_{jkl}(\s_{ij}\circ y_{i})\hat{\xi}_{ijl}(y_{i})  \\ 
& = & (u_{kl}(e(\s_{kl})\circ u_{jk}))\hat{a}(\s_{kl},
\s_{jk},y_{j})(e(\s_{kl}\circ \s_{jk})\circ u_{ij}) 
(\s_{jkl}\circ e(\s_{ij}\circ y_{i}))\hat{\xi}_{ijl}(y_{i})  \\
& = & (u_{kl}(e(\s_{kl})\circ u_{jk}))\hat{a}(\s_{kl},
\s_{jk},y_{j})(\s_{jkl}\circ u_{ij})\hat{\xi}_{ijl}(y_{i}) \\ 
& = & (u_{kl}(e(\s_{kl})\circ u_{jk}))\hat{a}(\s_{kl},
\s_{jk},y_{j})((\s_{jkl}e(\s_{jl}))\circ 
(e(y_{j})u_{ij}))\hat{\xi}_{ijl}(y_{i})                     \\ 
& = & (u_{kl}(e(\s_{kl})\circ u_{jk}))\hat{a}(\s_{kl},
\s_{jk},y_{j})(\s_{jkl}\circ e(y_{j}))(e(\s_{jl})
\circ u_{ij})\hat{\xi}_{ijl}(y_{i})                          \\ 
& = & (u_{kl}(e(\s_{kl})\circ u_{jk}))\hat{\xi}_{jkl}(y_{i})
(e(\s_{jl})\circ u_{ij})\hat{\xi}_{ijl}(y_{i}).            
\end{eqnarray*}
This proves that $m_{L}$ is associative 
and hence the triple $(L,Z,X)$ is 
a bundle gerbe.  
\end{proof} 

\begin{definition}
\label{def:12.1.3} 
Let $(Q,Y,X,M)$ be a bundle  
2-gerbe.  We say that  
$Q$ is \emph{trivial} if there 
exists a bundle gerbe $(L,Z,X)$ on 
$X$ together with a bundle gerbe 
morphism $\bar{\eta}:\pi_{1}^{-1}L\otimes Q \to \pi_{2}^{-1}L$ 
over $X^{[2]}$. 
We also require that there is a transformation 
$\theta$ between the bundle gerbe  
morphisms in the diagram below.  
$$
\xymatrix{ 
\pi_{1}^{-1}\pi_{1}^{-1}L\otimes 
\pi_{1}^{-1}Q\otimes \pi_{3}^{-1}Q 
\ar[d]_{\pi_{1}^{-1}\bar{\eta}\otimes 1} 
\ar[rr]^{1\otimes \bar{m}} & & \pi_{1}^{-1}
\pi_{1}^{-1}L\otimes \pi_{2}^{-1}Q 
\ar @2{-}[d]                                \\ 
\pi_{1}^{-1}\pi_{2}^{-1}L\otimes \pi_{3}^{-1}
Q \ar @2{-}[d] 
 & & \pi_{2}^{-1}\pi_{1}^{-1}L
\otimes \pi_{2}^{-1}Q \ar[d]^{\pi_{2}^{-1}
\bar{\eta}}  \ar @<-1ex> @2{->}[ll]^{\theta}                \\ 
\pi_{3}^{-1}\pi_{1}^{-1}L\otimes \pi_{3}
^{-1}Q \ar[dr]_{\pi_{3}^{-1}\bar{\eta}} & & 
\pi_{2}^{-1}\pi_{2}^{-1}L \ar @2{-}[dl]     \\ 
& \pi_{3}^{-1}\pi_{2}^{-1}L                  } 
$$
So $\theta$ is a section of the 
$\cstar$-bundle 
$B=D_{\bar{\eta}_{1},\bar{\eta}_{2}}$ over $X^{[3]}$.  
Here $\bar{\eta}_{1} = \pi_{2}^{-1}\bar{\eta}\circ 
(\text{id}\otimes \bar{m})$   
and $\bar{\eta}_{2} = \pi_{3}^{-1}\bar{\eta} 
\circ (\pi_{1}^{-1}\bar{\eta} \otimes \text{id})$
and $D_{\bar{\eta}_{1},\bar{\eta}_{2}}$ denotes the descended 
bundle on $X^{[3]}$ whose lift to $Z_{4}\times_
{X^{[3]}}Y_{123}$ is isomorphic to $(\eta_{1},
\eta_{2})^{-1}L_{1}$ --- see Section~\ref{sec:3.4}.  
There is a canonical isomorphism  
$\pi_{1}^{-1}B\otimes \pi_{2}^{-1}B^{*}\otimes 
\pi_{3}^{-1}B\otimes \pi_{4}^{-1}B^{*}\simeq A$ 
over $X^{[4]}$, where 
$A$ is the associator bundle of $Q$.  
We require that $\deltatheta = a$ 
under this isomorphism.  We say that 
the bundle gerbe $L$ \emph{trivialises} 
$Q$.  
\end{definition}  

\begin{note} 
\label{note:12.1.4} 
The canonical isomorphism 
$\pi_{1}^{-1}B\otimes \pi_{2}^{-1}B^{*}\otimes 
\pi_{3}^{-1}B\otimes \pi_{4}^{-1}B^{*}\cong A$ 
is a result of the following.  We can define five 
bundle gerbe morphisms $\bar{f}_{i}:L_{4}\otimes 
Q_{1234}\to L_{1}$ for $i=1,\ldots ,5$ by 
\begin{eqnarray*} 
\bar{f}_{1} & = & \bar{\eta}\circ (\text{id}\otimes 
\bar{m})\circ (\text{id}\otimes \bar{m}\otimes 
\text{id}),                                              \\ 
\bar{f}_{2} & = & \bar{\eta}\circ (\text{id}\otimes 
\bar{m})\circ (\text{id}\otimes \text{id}\otimes 
\bar{m}),                                             \\ 
\bar{f}_{3} & = & \bar{\eta}\circ (\bar{\eta}
\otimes \text{id})\circ (\text{id}\otimes text{id}\otimes 
\bar{m}),                                              \\ 
\bar{f}_{4} & = & \bar{\eta}\circ (\bar{\eta}\otimes 
\text{id})\circ (\bar{\eta}\otimes \text{id}\otimes 
\text{id}),                                            \\ 
\bar{f}_{5} & = & \bar{\eta}\circ (\bar{\eta}\otimes 
\text{id})\circ (\text{id}\otimes \bar{m}\otimes 
\text{id}). 
\end{eqnarray*} 
Note that $\bar{f}_{3}$ can also be written as 
$$
\bar{f}_{3} = \bar{\eta}\circ (\text{id}\otimes \bar{m})
\circ (\bar{\eta}\otimes \text{id}\otimes \text{id}).  
$$
From this we can see that we have the 
following isomorphisms: 
\begin{eqnarray*} 
D_{\bar{f}_{1},\bar{f}_{2}} & \simeq & A,             \\  
D_{\bar{f}_{2},\bar{f}_{3}} & \simeq & \pi_{2}^{-1}B,   \\ 
D_{\bar{f}_{3},\bar{f}_{4}} & \simeq & \pi_{4}^{-1}B, \\ 
D_{\bar{f}_{4},\bar{f}_{5}} & \simeq & \pi_{1}^{-1}B^{*},  \\ 
D_{\bar{f}_{5},\bar{f}_{1}} & \simeq & \pi_{3}^{-1}B^{*},    
\end{eqnarray*} 
In each case, as in the notation used in Section
~\ref{sec:3.4}, $D_{\bar{f}_{i},\bar{f}_{j}}$ denotes 
the line bundle on $X^{[4]}$ whose 
lift to $Z_{4}\times_{X^{[4]}}Y_{1234}$ is the 
line bundle $(f_{i},f_{j})^{-1}L_{1}$.  The 
bundle gerbe product on $L_{1}$ provides 
a trivialisation of the line bundle 
$$
(f_{5},f_{1})^{-1}L_{1}\otimes (f_{4},f_{5})^{-1}
L_{1}\otimes (f_{3},f_{4})^{-1}L_{1}\otimes 
(f_{2},f_{3})^{-1}L_{1}\otimes (f_{1},f_{2})
^{-1}L_{1} 
$$
on $Z_{4}\times_{X^{[4]}}Y_{1234}$.   
This provides the canonical isomorphism 
$\pi_{1}^{-1}B\otimes \pi_{2}^{-1}B^{*}\otimes 
\pi_{3}^{-1}B\otimes \pi_{4}^{-1}B^{*}\simeq A$.  
The coherency condition on the section 
$\theta$ of $B$ is as follows.  Each bundle 
$(f_{i},f_{j})^{-1}L_{1}$ has a trivialising 
section $\hat{f}_{ij}$ induced from $\theta$ and 
given as follows: 
\begin{eqnarray*} 
\hat{f}_{12}(z_{4},y_{34},y_{23},y_{12}) & = & 
\hat{\eta}(e(z_{4})\otimes \hat{a}(y_{34},y_{23},y_{12})),   \\ 
\hat{f}_{23}(z_{4},y_{34},y_{23},y_{12}) & = & 
\hat{\theta}(z_{4},y_{34},m(y_{23},y_{12})), \\ 
\hat{f}_{34}(z_{4},y_{34},y_{23},y_{12}) & = & 
\hat{\theta}(\eta(z_{4},y_{34}),y_{23},y_{12}),     \\ 
\hat{f}_{45}(z_{4},y_{34},y_{23},y_{12}) & = & 
\hat{\eta}(\hat{\theta}(z_{4},y_{34},y_{23})\otimes e(y_{12}))^{-1}     \\ 
\hat{f}_{51}(z_{4},y_{34},y_{23},y_{12}) & = & 
\hat{\theta}(z_{4},m(y_{34},y_{23}),y_{12})^{-1}. 
\end{eqnarray*} 
where $\hat{\theta}$ denotes the lift of the 
section $\theta$ of $B = D_{\bar{\eta}_{1},\bar{\eta}_{2}}$ 
under the isomorphism $\hat{D}\simeq (\eta_{1},
\eta_{2})^{-1}L_{1}$.  The coherency condition 
$\pi_{1}^{-1}\theta\otimes \pi_{2}^{-1}\theta^{*}
\otimes \pi_{3}^{-1}\theta\otimes \pi_{4}^{-1}
\theta^{*} = a$ under the isomorphism 
$\pi_{1}^{-1}B\otimes \pi_{2}^{-1}B^{*}\otimes 
\pi_{3}^{-1}B\otimes \pi_{4}^{-1}B^{*}\cong A$ 
translates into 
$$
\hat{f}_{51} \hat{f}_{45} \hat{f}_{34}
 \hat{f}_{23} \hat{f}_{12} = 
e(\eta(z_{4},m(m(y_{34},y_{23}),y_{12}))), 
$$
where $e$ is the identity section of the 
bundle gerbe $(f_{1},f_{1})^{-1}L_{1}$.  
\end{note}  

We have the following lemma. 
\begin{lemma}
\label{12.1.5} 
Suppose $(Q,Y,X,M)$ is a trivial bundle 
2-gerbe with trivialising bundle 
gerbe $(L,Z,X)$.  Then the bundle  
2-gerbe $Q$ has zero four class in 
$H^{4}(M;\Z)$. 
\end{lemma} 

\begin{proof} Let 
$\U = \{U_{i}\}_{i \in I}$ be a 
`nice' open cover of $M$.  So each 
non-empty finite intersection 
$U_{i_{1}}\cap \cdots \cap U_{i_{p}}$ 
is contractible and there exist 
local sections $s_{i}:U_{i}\to X$ 
of $\pi_{X}:X\to M$.  Define maps 
$(s_{i},s_{j}):U_{ij}\to X^{[2]}$ by 
$m \mapsto (s_{i}(m),s_{j}(m))$.  
Pullback the bundle gerbe 
$Q = (Q,Y,X^{[2]})$ with this map 
to obtain a new bundle gerbe 
$Q_{ij} = (Q_{ij},Y_{ij},U_{ij})$.  
So we get a manifold $Y_{ij}$ 
and an induced projection map 
$\pi_{ij}:Y_{ij}\to U_{ij}$.  If 
$m \in U_{ij}$ then 
$\pi_{ij}^{-1}(m) = Y_{ij_{m}} = Y_{(s_{i}(m),s_{j}(m))}$.  
Since $U_{ij}$ is contractible 
the bundle gerbe $Q_{ij}$ is 
trivial.  So let $J_{ij}\to Y_{ij}$ 
be a $\cstar$ bundle trivialising $Q_{ij}$.  
Let $(L,Z,X)$ be a bundle gerbe 
trivialising the bundle 2-gerbe $Q$.  
Pullback the bundle gerbe $L$ with the 
sections $s_{i}$ to get a new bundle 
gerbe $L_{i}  = (L_{i},Z_{i},U_{i})$ 
on the contractible open set 
$U_{i}$.  Let $K_{i}\to Z_{i}$ be 
a $\cstar$ bundle trivialising the 
bundle gerbe $L_{i}$.  The 
bundle gerbe morphism 
$\bar{\eta}:\pi_{1}^{-1}L\otimes Q\to \pi_{2}^{-1}L$ 
induces a bundle gerbe morphism 
$\bar{\eta}_{ij}:L_{j}\otimes Q_{ij}\to L_{i}$ 
over $U_{ij}$.  Consider the $\cstar$ 
bundle $K_{j}\otimes J_{ij}\otimes 
\eta_{ij}^{-1}K_{i}^{*}$ on $Z_{j}\times_{M}Y_{ij}$.   
Since the bundles $K_{j}$, $J_{ij}$ and 
$K_{i}$ trivialise the bundle gerbes 
$L_{j}$, $Q_{ij}$ and $L_{i}$ respectively, 
and because $\bar{\eta}_{ij}$ is a bundle 
gerbe morphism and hence commutes with the 
bundle gerbe products on $L_{j}\otimes Q_{ij}$ 
and $L_{i}$, it follows that there 
is a section $\hat{\phi}_{ij}$ of 
$\d(K_{j}\otimes J_{ij}\otimes \eta_{ij}^{-1}K_{i}^{*})$ 
over $(Z_{j}\times_{M}Y_{ij})^{[2]}$ 
satisfying the coherency condition 
$\d(\hat{\phi}) = \underline{1}$.  Hence the 
$\cstar$ bundle $K_{j}\otimes J_{ij}\otimes 
\eta_{ij}^{-1}K_{i}^{*}$ descends to a bundle 
$C_{ij}$ on $U_{ij}$.  
From the construction of the \v{C}ech cocycle 
$g_{ijkl}$ associated to the bundle 2-gerbe 
$Q$ given in Section~\ref{sec:9.1} we have 
that the $\cstar$ bundle $\hat{L}_{ijk} = 
J_{jk}\otimes J_{ij}\otimes m^{-1}J_{ik}^{*}$ 
descends to a $\cstar$ bundle $L_{ijk}$ on 
$U_{ijk}$ and that furthermore there is an 
isomorphism 
$$
L_{jkl}\otimes L_{ikl}^{*}\otimes L_{ijl}
\otimes L_{ijk}^{*} \simeq A_{ijkl} 
$$
over $U_{ijkl}$, where $A_{ijkl} = (s_{i},s_{j},s_{k},s_{l})^{-1}A$.  
Also the \v{C}ech 3-cocycle $g_{ijkl}$ was 
defined by choosing sections $\s_{ijk}$ of 
$L_{ijk}$ and then comparing the section 
$a_{ijkl} = (s_{i},s_{j},s_{k},s_{l})^{-1}a$ 
with the image of the section $\s_{jkl}\otimes 
\s_{ikl}^{*}\otimes \s_{ijl}\otimes \s_{ijk}^{*}$ 
under the above isomorphism.  

It is not difficult to show that there is an 
isomorphism 
$$
C_{jk}\otimes C_{ij}\otimes B_{ijk} \simeq 
L_{ijk} \otimes C_{ik} 
$$
of $\cstar$ bundles over $U_{ijk}$, where 
$B_{ijk}$ denotes the pullback of the bundle 
$B$ on $X^{[3]}$ via the map $U_{ijk}\to X^{[3]}$ 
defined by sending a point $m$ of $U_{ijk}$ to 
$(s_{i}(m),s_{j}(m),s_{k}(m))$.  From here we see 
that if we choose sections $\rho_{ij}$ of $C_{ij}$ 
then we can define a section $\s_{ijk}$ of $L_{ijk}$ by 
$\s_{ijk} = \rho_{ij}\otimes \rho_{jk}\otimes \rho_{ik}^{*}
\otimes \theta_{ijk}$, where $\theta_{ijk} = 
(s_{i},s_{j},s_{k})^{-1}\theta$.  Comparing 
$\s_{jkl}\otimes \s_{ikl}^{*}\otimes \s_{ijl}\otimes 
\s_{ijk}^{*}$ with $a_{ijkl}$ and using the coherency 
condition $\d(\theta) = a$, we see that the cocycle 
$g_{ijkl}$ must be trivial.     
\end{proof} 

We shall now prove that the converse 
is also true, justifying the terminology.  
Let $(Q,Y,X,M)$ be a bundle 2-gerbe with 
zero four class in $H^{4}(M;\Z)$.  
We shall prove that the bundle gerbe 
$(L,Z,M)$ we constructed in 
Proposition~\ref{prop:12.1.1} trivialises the 
bundle 2-gerbe $(Q,Y,X,M)$.  
First of all we need to define the bundle 
gerbe morphism $\bar{\eta}:\pi_{1}^{-1}L\otimes Q\to \pi_{2}^{-1}L$.  
Recall $Z = \coprod_{i \in I}Y_{i}$.  
Therefore if $(x_{1},x_{2})\in X^{[2]}$, 
then $\pi_{1}^{-1}Z\times_{X^{[2]}}Y$ 
has fibre at $(x_{1},x_{2})$ 
equal to the disjoint union of the 
spaces $Y_{(x_{2},s_{i}(m))}\times Y_{(x_{1},x_{2})}$,   
where $m = \pi(x_{1}) = \pi(x_{2})$.  
Therefore the map 
$m:\pi_{1}^{-1}Y\times_{X^{[3]}}\pi_{3}^{-1}Y\to \pi_{2}^{-1}Y$ 
induces a map 
$\eta:\pi_{1}^{-1}Z\times_{X^{[2]}}Y\to \pi_{2}^{-1}Z$ 
which sends $(y_{i},y)$ to 
$m(y_{i},y)$ for $y_{i}\in \pi_{1}^{-1}Y_{i}$ and 
$y\in Y$.  We define a bundle map 
$\hat{\eta}:\pi_{1}^{-1}L\otimes Q\to \pi_{2}^{-1}L$ 
covering the map 
$\eta^{[2]}:(\pi_{1}^{-1}Z\times_{X^{[2]}}Y)^{[2]}\to \pi_{2}^{-1}Z^{[2]}$ 
as follows.  Let $u\in Q_{(y_{1},y_{2})}$ 
and $u_{ij}\in Q_{(m(\s_{ij},y_{i}),y_{j})}$ 
where $(y_{1},y_{2})\in Y^{[2]}$ 
and $y_{i}\in \pi_{1}^{-1}Y_{i}$, 
$y_{j}\in \pi_{1}^{-1}Y_{j}$.  
So $(y_{i},y_{j})\in \pi_{1}^{-1}Z^{[2]}$ 
and $(y_{i},y_{j},u_{ij})\in \pi_{1}^{-1}L_{(y_{i},y_{j})}$. 
Then $m((\hat{m}(u_{ij}\otimes u)\otimes 
\hat{a}(\s_{ij},y_{i},y_{1})^{-1})$ 
is in $Q_{(m(\s_{ij},m(y_{i},y_{1})),m(y_{j},y_{2}))}$.  
Therefore we define 
\begin{equation} 
\label{eq:12.5} 
\hat{\eta}(u_{ij}\otimes u) = 
m((\hat{m}(u_{ij}\otimes u)\otimes 
\hat{a}(\sigma_{ij}^{\a},y_{i},y_{1})^{-1}).
\end{equation} 
The problem that confronts us now is 
to show that $\bar{\eta}$ so 
defined with $\bar{\eta}= (\hat{\eta},\eta,\text{id})$  
is in fact a bundle 
gerbe morphism, that is it commutes 
with the respective bundle gerbe 
products.    

So let $((y_{i},y_{1}),(y_{j},y_{2}),(y_{k},y_{3}))$ 
be a point of 
$(\pi_{1}^{-1}Z\times_{X^{[2]}}Y)^{[3]}$ 
and let $u_{jk}\otimes u_{23}$ be 
an element of $(\pi_{1}^{-1}L\otimes Q)_{((y_{j},y_{2}),(y_{k},y_{3}))}$ 
and $u_{ij}\otimes u_{12}$ be an element of 
$(\pi_{1}^{-1}L\otimes Q)_{((y_{i},y_{1}),(y_{j},y_{2}))}$.  
Thus $u_{23} \in Q_{(y_{2},y_{3})}$ 
and $u_{jk}$ is an element of 
$Q_{(m(\s_{jk},y_{j}),y_{k})}$ and 
similarly for $u_{12}$ and $u_{ij}$.  
We show that 
$$
\hat{\eta}(u_{jk}\otimes u_{23})\cdot 
\hat{\eta}(u_{ij}\otimes u_{12}) = 
\hat{\eta}((u_{jk}\cdot u_{ij})\otimes (u_{23} u_{12})). 
$$
We start with the left hand side 
and for convenience of notation, 
we denote by $\circ$ both $m$ 
and $\hat{m}$.  
The left hand side is equal to 
$$
[(u_{jk}\circ u_{23}) \hat{a}(\s_{jk},y_{j},y_{2})^{-1}] 
\cdot [(u_{ij}\circ u_{12}) \hat{a}(\s_{ij},y_{i},y_{1})^{-1}],
$$
where the dot between the 
square brackets denotes multiplication in 
the bundle gerbe $L$.  Therefore the 
left hand side is equal to 
$$
[(u_{jk}\circ u_{23}) \hat{a}
(\s_{jk},y_{j},y_{2})^{-1}] [e(\s_{jk})
\circ ((u_{ij}\circ u_{12}) \hat{a}
(\s_{ij},y_{i},y_{1}))^{-1}] 
\hat{\xi}_{ijk}(y_{i}\circ y_{1}).
$$
Notice that since $e(\s_{jk})$ is 
just the identity section, we can 
write $e(\s_{jk}) = e(\s_{jk}) e(\s_{jk})$ 
and, since $\circ$ is a bundle gerbe 
morphism, we get the above expression 
is equal to 
$$
[(u_{jk}\circ u_{23}) \hat{a}
(\s_{jk},y_{j},y_{2})^{-1}] [e(\s_{jk})
\circ (u_{ij}\circ u_{12})] 
[e(\s_{jk})\circ \hat{a}(\s_{ij},y_{i},y_{1})^{-1}] 
 \hat{\xi}_{ijk}(y_{i}\circ y_{1}).
$$
The right hand side is more 
complicated:
\begin{eqnarray*}
\text{RHS} & = & 
\hat{\eta}([u_{jk} (e(\s_{jk})\circ u_{ij}) 
 \hat{\xi}_{ijk}(y_{i})]\otimes 
[u_{23} u_{12}])                                  \\
 & = & [(u_{jk}\circ (e(\s_{jk})\circ u_{ij}) 
 \hat{\xi}_{ijk}(y_{i}))\circ (u_{23} u_{12})]
 \hat{a}(\s_{ik},y_{i},y_{1})^{-1}.                     
\end{eqnarray*}
Next we use the fact that 
$\circ $ commutes with bundle 
gerbe products to see that the 
expression in square brackets is 
equal to 
$$
(u_{jk}\circ u_{23}) (((e(\s_{jk})
\circ u_{ij}) \hat{\xi}_{ijk}(y_{i})) 
\circ u_{12}). 
$$
Now we re-express 
$((e(\s_{jk})\circ u_{ij}) \hat{\xi}_{ijk}(y_{i}))\circ u_{12}$ 
by observing that it is equal to 
$$
((e(\s_{jk})\circ u_{ij}) \hat{\xi}_{ijk}(y_{i}))
\circ (u_{12} e(y_{1}))
$$
and then using again the fact that 
$\circ$ is a bundle gerbe morphism 
to get it equal to 
$$
((e(\s_{jk})\circ u_{ij})\circ u_{12}) 
 (\hat{\xi}_{ijk}(y_{i})\circ e(y_{1})).
$$
Next we use the associator 
section $\hat{a}$ to write 
$(e(\s_{jk})\circ u_{ij})\circ u_{12}$ 
as 
$$
\hat{a}(\s_{jk},y_{j},y_{2})^{-1} 
(e(\s_{jk})\circ (u_{ij}\circ u_{12})) 
\hat{a}(\s_{jk},\s_{ij}\circ y_{i},y_{1}).
$$
At this point, if we substitute 
all this into our original expression 
for the right hand side, 
then we see that if we can show that 
\begin{equation} 
\label{eq:12.6}
\hat{a}(\s_{jk},\s_{ij}\circ y_{i},y_{1})
 (\hat{\xi}_{ijk}(y_{i})\circ e(y_{1}))
 \hat{a}(\s_{ik},y_{i},y_{1})^{-1}
\end{equation} 
is equal to 
$(e(\s_{jk})\circ \hat{a}(\s_{ij},y_{i},y_{1})^{-1}) 
\hat{\xi}_{ijk}(y_{i}\circ y_{1})$  
then we are done.   
$\hat{\xi}_{ijk}(y_{i})\circ e(y_{1})$ 
is equal to 
\begin{eqnarray*} 
 & & (\hat{a}(\s_{jk},\s_{ij},y_{i}) 
(\s_{ijk}\circ e(y_{i})))\circ e(y_{1})        \\
 & = & (\hat{a}(\s_{jk},\s_{ij},y_{i})\circ 
e(y_{1})) ((\s_{ijk}\circ e(y_{i}))\circ 
e(y_{1})).                                          
\end{eqnarray*}
Using the associator section $\hat{a}$, 
we see that 
$(\s_{ijk}\circ e(y_{i}))\circ e(y_{1})$ 
is equal to 
$$
\hat{a}(\s_{jk}\circ \s_{ij},y_{i},y_{1})^{-1} 
(\s_{ijk}\circ (e(y_{i})\circ e(y_{1})) 
\hat{a}(\s_{ik},y_{i},y_{1}).
$$
Therefore 
\begin{multline*} 
\hat{\xi}_{ijk}(y_{i})\circ e(y_{1}) = (\hat{a}(
\s_{jk},\s_{ij},y_{i})\circ e(y_{1}))(\hat{a}(
\s_{jk}\circ \s_{ij},y_{i},y_{1})^{-1}          \\ 
(\s_{ijk}\circ 
(e(y_{i})\circ e(y_{1}))\hat{a}(\s_{ik},y_{i},y_{1})).  
\end{multline*} 
If we now substitute this into 
equation~\ref{eq:12.6} then we 
get
$$
\hat{a}(\s_{jk},\s_{ij}\circ y_{i},y_{1})
(\hat{a}(\s_{jk},\s_{ij},y_{i})\circ 
e(y_{1})) \hat{a}(\s_{jk}\circ \s_{ij},
y_{i},y_{1})^{-1} (\s_{ijk}\circ e(y_{i}\circ y_{1})),
$$
since $e(y_{i})\circ e(y_{1})$ 
equals $e(y_{i}\circ y_{1})$.  
Using the coherency condition on 
$\hat{a}$, we get that this is equal to 
$$
(e(\s_{jk})\circ \hat{a}(\s_{ij},y_{i},y_{1})^{-1})
 \hat{a}(\s_{jk},\s_{ij},y_{i}\circ y_{1})
 (\s_{ijk}\circ e(y_{i}\circ y_{1})),
$$
as required.  Thus $\bar{\eta}$ is 
a bundle gerbe morphism.  To complete 
the demonstration that $L$ 
trivialises the bundle 2-gerbe $Q$, 
we need to define the transformation 
section $\theta$ and show that it is 
coherent with $a$.  

First we need to define the section 
$\hat{\theta}$ of the $\cstar$ bundle 
$(\eta_{1},\eta_{2})^{-1}L_{1}$ on 
$Z_{3}\times_{X^{[3]}}Y_{23}\times_{X^{[3]}}Y_{12}$.  
If $(y_{i},y_{23},y_{12}) \in \zed$, then 
$\hat{J} = (\eta_{1},\eta_{2})^{-1}L_{1}$ 
has fibre equal to $Q_{(y_{i}\circ (y_{23}\circ y_{12}), 
\s_{ii}\circ ((y_{i}\circ y_{23})\circ y_{12}))}$ at 
$(y_{i},y_{23},y_{12})$.  Let $1(y_{i})$ 
denote the identity section of $L$ 
evaluated at $(y_{i},y_{i})$.  Then we can 
define a section $\hat{\theta}$ of $\hat{J}$ by 
\begin{equation} 
\label{eq:12.7} 
\hat{\theta}(y_{i},y_{23},y_{12}) = 
(e(\s_{ii})\circ \hat{a}(y_{i},y_{23},y_{12})^{-1})
\hat{a}(\s_{ii},y_{i},y_{23}\circ y_{12})(1(y_{i})^{-1}
\circ e(y_{23}\circ y_{12})).  
\end{equation} 
Here $1(y_{i})^{-1}$ means the inverse of the identity section 
$1(y_{i})$ of $L$ at $y_{i}$ calculated in 
the bundle gerbe $Q$ (recall $L_{(y_{i},y_{i})} = 
Q_{(\s_{ii}\circ y_{i},y_{i})}$ so $1(y_{i})^{-1}\in 
Q_{(y_{i},\s_{ii}\circ y_{i})}$).   
We need to show that $\hat{\theta}$ descends 
to a section $\theta$ of $J\to X^{[3]}$, that is we 
need to show that 
$\phi(\pi_{1}^{-1}\hat{\theta}) = \pi_{2}^{-1}\hat{\theta}$,
where $\phi:\pi_{1}^{-1}\hat{J} \to \pi_{2}^{-1}\hat{J}$ 
is the descent isomorphism.   
Recall that $\phi$ is constructed by 
choosing $u_{ij}\in Q_{(\s_{ij}\circ y_{i},y_{j})} = L_{(y_{i},y_{j})}$, 
$u_{23}\in Q_{(y_{23},y_{23}^{\prime})}$ 
and $u_{12}\in Q_{(y_{12},y_{12}^{\prime})}$ 
where $(y_{23},y_{23}^{\prime}) \in Y_{23}^{[2]}$, 
$(y_{12},y_{12}^{\prime})\in Y_{12}^{[2]}$ 
and $(y_{i},y_{j})\in Z_{3}^{[2]}$ and setting 
$$
\phi(v) = \hat{\eta}_{2}(u_{ij}\otimes u_{23}
\otimes u_{12})^{-1}\cdot v\cdot \hat{\eta}_{1}
(u_{ij}\otimes u_{23}\otimes u_{12})
$$
for $z \in J_{(y_{j},y_{23}^{\prime},y_{12}^{\prime})}$ 
and where the $\cdot$ denotes multiplication 
in the bundle gerbe $L$.   
It is very tedious to prove that $\theta$  
commutes with the descent isomorphism $\phi$ 
and so we will omit the calculation.  

The final step is to show that the 
descended section is coherent with 
$a$.  Recall from Note~\ref{note:12.1.4} 
the definitions of the bundle gerbe 
morphisms $\bar{f}_{i}:L_{4}\otimes 
Q_{1234}\to L_{1}$ for $i=1,\ldots,5$ 
and the sections $\hat{f}_{ij}$ of 
$(f_{i},f_{j})^{-1}L_{1}$ over 
$Z_{4}\times_{X^{[4]}}Y_{1234}$.  
Recall in particular that we had 
\begin{eqnarray*} 
\hat{f}_{12}(z_{4},y_{34},y_{23},y_{12}) & = & 
\hat{\eta}(1(z_{4})\otimes \hat{a}(y_{34},y_{23},y_{12})),     \\ 
\hat{f}_{23}(z_{4},y_{34},y_{23},y_{12}) & = & 
\hat{\theta}(z_{4},y_{34},m(y_{23},y_{12})),     \\ 
\hat{f}_{34}(z_{4},y_{34},y_{23},y_{12}) & = & 
\hat{\theta}(\eta(z_{4},y_{34}),y_{23},y_{12}),   \\   
\hat{f}_{45}(z_{4},y_{34},y_{23},y_{12}) & = & 
\hat{\eta}(\hat{\theta}(z_{4},y_{34},y_{23})\otimes e(y_{12}))^{-1},   \\ 
\hat{f}_{51}(y_{i},y_{34},y_{23},y_{12}) & = & 
\hat{\theta}(z_{4},m(y_{34},y_{23}),y_{12})^{-1}, 
\end{eqnarray*} 
where $(y_{i},y_{34},y_{23},y_{12})\in 
Z_{4}\times_{X^{[4]}}Y_{1234}$.  We need to show 
that 
$$
\hat{f}_{51} \hat{f}_{45} \hat{f}_{34}
 \hat{f}_{23} \hat{f}_{12} = 1(y_{i}\circ ((y_{34}\circ 
y_{23})\circ y_{12}))  
$$
in the bundle gerbe $(f_{1},f_{1})^{-1}L_{1}$.  
This will imply the coherency condition 
$\delta(\theta) = a$.  Again  
I will omit the proof as it is long and cumbersome.   
Combining the above results we have the 
following proposition.  

\begin{proposition}
\label{prop:12.1.6} 
Let $(Q,Y,X,M)$ be a bundle 2-gerbe.  
Then $Q$ has zero four class in 
$H^{4}(M;\Z)$ if and only if 
$Q$ is a trivial bundle 2-gerbe.
\end{proposition}

\begin{note} 
\label{note:12.1.7} 
Throughout we have assumed that it is possible 
to find a section $\s_{ij}$ of $Y_{ij} = 
(s_{i},s_{j})^{-1}Y\to U_{ij}$ over $U_{ij}$.  
In general though $\pi_{Y}:Y\to X^{[2]}$ 
only admits local sections and therefore it 
is only possible to find an open cover 
$\{U_{ij}^{\a}\}_{\a\in \Sigma_{ij}}$ of $U_{ij}$ 
and local sections $\s_{ij}^{\a}:U_{ij}^{\a}\to Y_{ij}$ 
of $Y_{ij}\to U_{ij}$.  It is possible to 
get around this restriction as we 
will now show.  Assume that $\pi_{Y}:Y\to 
X^{[2]}$ only admits local sections.  Then 
there is an open cover $\{U_{ij}^{\a}\}_{\a\in 
\Sigma_{ij}}$ of $U_{ij}$ such that there exist 
sections $\s_{ij}^{\a}:U_{ij}^{\a}\to Y_{ij}$.  
Let $X_{ij}^{\a} = X_{U_{ij}^{\a}}$.  We can 
define maps $\phi_{ij}^{\a}:Y_{i}\vert_{X_{ij}^{\a}}
\to Y_{j}\vert_{X_{ij}^{\a}}$ by sending 
$y\in Y_{i}\vert_{X_{ij}^{\a}}$ to 
$m(\s_{ij}^{\a},y)$.  As before the $\phi_{ij}^{\a}$ 
extend to define bundle gerbe morphisms 
$\bar{\phi}_{ij}^{\a}:Q_{i}\vert_{X_{ij}^{\a}}\to 
Q_{j}\vert_{X_{ij}^{\a}}$ with 
$\bar{\phi}_{ij}^{\a} = (\hat{\phi}_{ij}^{\a},
\phi_{ij}^{\a},\text{id})$ if we put 
$\hat{\phi}_{ij}^{\a}(u_{i}) = \hat{m}(
e(\s_{ij}^{\a})\otimes u_{i})$.    
Now we can use the $\phi_{ij}^{\a}$ to 
define maps $\chi_{ij}^{\a}:(Y_{i}\times_{X}
Y_{j})\vert_{X_{ij}^{\a}}\to Y_{j}^{[2]}\vert
_{X_{ij}^{\a}}$ by $\chi_{ij}^{\a}(y_{i},y_{j}) = 
(\phi_{ij}^{\a}(y_{i}),y_{j})$.  Define line 
bundles $L_{ij}^{\a}$ on $(Y_{i}\times_{X}Y_{j})
\vert_{X_{ij}^{\a}}$ by $L_{ij}^{\a} = 
(\chi_{ij}^{\a})^{-1}Q_{j}$.  

We now need to digress for a moment.  Let 
$U_{ij}^{\a\b} = U_{ij}^{\a}\cap U_{ij}^{\b}$.  
Let $\rho_{ij}^{\a\b}$ be a section of the 
pullback bundle $(\sigma_{ij}^{\a},\sigma_{ij}^{\b})
^{-1}Q_{ij}$ over $U_{ij}^{\a\b}$.  Observe that 
since $U_{ij}$ is contractible, the bundle gerbe 
$(Q_{ij},Y_{ij},U_{ij})$ is trivial and so it 
follows that we may choose the $\rho_{ij}^{\a\b}$ 
so that 
\begin{equation} 
\label{eq:12.8} 
m_{Q}(\rho_{ij}^{\b\c}\otimes \rho_{ij}^{\a\b}) 
= \rho_{ij}^{\a\c} 
\end{equation} 
on $U_{ij}^{\a\b\c} = U_{ij}^{\a}
\cap U_{ij}^{\b}\cap U_{ij}^{\c}$.  This is 
simply the statement that the Dixmier-Douady 
class of $Q_{ij}$ is zero.  Next choose sections 
$\sigma_{ijk}^{\a\b\c}$ of the pullback bundle 
$Q_{ijk}^{\a\b\c} = (\sigma_{jk}^{\b}\circ 
\sigma_{ij}^{\a},\sigma_{ik}^{\c})^{-1}Q_{ik}$ 
over $U_{ijk}^{\a\b\c}$.  Here $U_{ijk}^{\a\b\c} =
U_{ij}^{\a}\cap U_{jk}^{\b}\cap U_{ik}^{\c}$.  
If $\a^{'}$, $\b^{'}$ 
and $\c^{'}$ are more indices of $\Sigma_{ij}$, 
$\Sigma_{jk}$ and $\Sigma_{ik}$ respectively, we 
can define a map $h_{ijk}^{\a\b\c\a^{'}\b^{'}\c^{'}}:
U_{ijk}^{\a\b\c\a^{'}\b^{'}\c^{'}}\to \cstar$ 
by 
$$
m_{Q}(\rho_{ik}^{\c\c^{'}}\otimes \sigma_{ijk}
^{\a\b\c}) = m_{Q}(\sigma_{ijk}^{\a^{'}\b^{'}
\c^{'}}\otimes (\rho_{jk}^{\b\b^{'}}\circ 
\rho_{ij}^{\a\a^{'}})\cdot h_{ijk}^{\a\b\c\a^{'}
\b^{'}\c^{'}},  
$$
where $U_{ijk}^{\a\b\c\a^{'}\b^{'}\c^{'}} 
= U_{ijk}^{\a\b\c}\cap U_{ijk}^{\a^{'}\b^{'}
\c^{'}}$.  
If $\a^{''}$, $\b^{''}$ and $\c^{''}$ are 
more indices of $\Sigma_{ij}$, $\Sigma_{jk}$ 
and $\Sigma_{ik}$, then using 
Equation~\ref{eq:12.8} 
for the various 
indices of $I$ and the $\Sigma$, one can show 
that 
$$
h_{ijk}^{\a\b\c\a^{'}\b^{'}\c^{'}}
h_{ijk}^{\a^{'}\b^{'}\c^{'}\a^{''}
\b^{''}\c^{''}} = h_{ijk}^{\a\b\c
\a^{''}\b^{''}\c^{''}}.  
$$
Since $U_{ijk}$ is contractible and 
$\{U_{ijk}^{\a\b\c}\}_{\a\in \Sigma_{ij}, 
\b\in \Sigma_{jk}, \c\in \Sigma_{ik}}$ is a cover 
of $U_{ijk}$, it follows 
that there exist $c_{ijk}^{\a\b\c}:U_{ijk}
^{\a\b\c}\to \cstar$ such that 
$h_{ijk}^{\a\b\c\a^{'}\b^{'}\c^{'}} = 
c_{ijk}^{\a^{'}\b^{'}\c^{'}}c_{ijk}^{\a\b\c\ -1}$.  
It follows that we can define the sections 
$\sigma_{ijk}^{\a\b\c}$ of $Q_{ijk}^{\a\b\c}$ 
so that 
\begin{equation} 
\label{eq:12.9} 
m_{Q}(\rho_{ik}^{\c\c^{'}}\otimes 
\sigma_{ijk}^{\a\b\c}) = m_{Q}(
\sigma_{ijk}^{\a^{'}\b^{'}\c^{'}}
\otimes (\rho_{jk}^{\b\b^{'}}\circ 
\rho_{ij}^{\a\a^{'}})). 
\end{equation} 
Now we can define a section 
$s_{ijkl}^{\a\b\c\eta\epsilon}$ of the 
pullback bundle $(\sigma_{il}^{\epsilon},
(\sigma_{kl}^{\d}\circ \sigma_{jk}^{\b})
\circ \sigma_{ij}^{\a})^{-1}Q_{il}$ 
in the same fashion as in 
equation~\ref{eq:9.1}, namely we set 
\begin{equation} 
\label{eq:12.10} 
s_{ijkl}^{\a\b\c\eta\epsilon} = 
(\sigma_{jkl}^{\b\c\eta}\circ e(\sigma_{ij}
^{\a}))\sigma_{ijl}^{\a\eta\epsilon}. 
\end{equation} 
We also define a section $t_{ijkl}
^{\a\b\c\d\epsilon}$ of the pullback 
bundle $(\sigma_{il}^{\epsilon},
\sigma_{kl}^{\d}\circ (\sigma_{jk}^{\b}
\circ \sigma_{ij}^{\a}))^{-1}Q_{il}$ 
in the analogous way to equation~\ref{eq:9.2}, 
\begin{equation} 
\label{eq:12.11} 
t_{ijkl}^{\a\b\c\d\epsilon} = 
((e(\sigma_{kl}^{\d})\circ \sigma_{ijk}
^{\a\b\c})\sigma_{ikl}^{\c\d\epsilon}. 
\end{equation} 
Finally, after equation~\ref{eq:9.3}, 
we define a map $g_{ijkl}^{\a\b\c\d\eta
\epsilon}:U_{ijkl}^{\a\b\c\d\eta\epsilon}
\to \cstar$ by 
\begin{equation} 
\label{eq:12.12} 
t_{ijkl}^{\a\b\c\d\epsilon} = 
a(\sigma_{kl}^{\d},\sigma_{jk}^{\b},
\sigma_{ij}^{\a})s_{ijkl}^{\a\b\d\eta
\epsilon}\cdot g_{ijkl}^{\a\b\c\d\eta
\epsilon}.   
\end{equation} 
Using equation~\ref{eq:12.9}, one 
can show that if $\a^{'}$, $\b^{'}$, 
$\c^{'}$, $\d^{'}$, $\eta^{'}$ and 
$\epsilon^{'}$ are some more indices 
of $\Sigma$, then 
\begin{equation} 
\label{eq:12.13} 
g_{ijkl}^{\a\b\c\d\eta\epsilon} = 
g_{ijkl}^{\a^{'}\b^{'}\c^{'}\d^{'}
\eta^{'}\epsilon^{'}}. 
\end{equation} 
Therefore the various maps 
$g_{ijkl}^{\a\b\c\d\eta\epsilon}$ are 
the restrictions to $U_{ijkl}^{\a\b\c
\d\eta\epsilon}$ of a map $g_{ijkl}:
U_{ijkl}\to\cstar$.  It is not hard 
to see that $g_{ijkl}$ satisfies the 
cocycle condition 
$$
g_{jklm}g_{iklm}^{-1}
g_{ijlm}g_{ijkm}^{-1}g_{ijkl} = 1
$$ 
and 
that $g_{ijkl}$ is the \v{C}ech 3-cocycle 
constructed in Section~\ref{sec:9.1}.  
     
Returning to the problem at hand, 
we will use the $\rho_{ij}^{\a\b}$ to define 
maps $\psi_{ij}^{\a\b}:L_{ij}^{\a}\to L_{ij}^{\b}$ 
covering the identity on $(Y_{i}\times_{X}Y_{j})\vert_
{X_{ij}^{\a}\cap X_{ij}^{\b}}$  
as follows.  Recall that $L^{\a}_{ij_{(y_{i},y_{j})}} 
= Q_{j_{(\phi_{ij}^{\a}(y_{i}),y_{j})}}$.  Therefore if $u_{ij}^{\a}\in 
L^{\a}_{ij_{(y_{i},y_{j})}}$ 
then $u_{ij}^{\a}\in Q_{j_{(\phi_{ij}^{\a}(y_{i}),y_{j})}}$.  
Therefore  
$u_{ij}^{\a}\hat{m}(\rho_{ij}^{\a\b}\otimes e(y_{i}))$  
lives in $Q_{j_{(\phi_{ij}^{\b}(y_{i}),y_{j})}}$.  So we 
let $\psi_{ij}^{\a\b}(u_{ij}^{\a}) = 
u_{ij}^{\a}\hat{m}(\rho_{ij}^{\a\b}\otimes e(y_{i}))$.  
Equation~\ref{eq:12.8}  
shows that 
$\psi_{ij}^{\b\c}\circ \psi_{ij}^{\a\b} = 
\psi_{ij}^{\a\c}$.  Hence we can glue the $L_{ij}^{\a}$ 
together using the standard clutching construction 
to form a line bundle $L_{ij}$ on the whole of 
$Y_{i}\times_{X}Y_{j}$.  

One can define a bundle gerbe product 
$m_{ijk}^{\a\b\c}:L_{jk}^{\b}\otimes 
L_{ij}^{\a}\to L_{ik}^{\c}$  using the analogue 
of equation~\ref{eq:12.4} --- providing, of 
course, that the cocycle $g_{ijkl}$ is trivial : 
\begin{equation} 
\label{eq:12.14} 
m_{ijk}^{\a\b\c}(u_{jk}^{\b}\otimes u_{ij}^{\a}) 
= m_{Q}(m_{Q}(u_{jk}^{\b}\otimes \hat{m}(
e(\sigma_{ij}^{\a})\otimes u_{ij}^{\a}))\otimes 
\hat{\xi}_{ijk}^{\a\b\c}(y_{i})), 
\end{equation} 
where $\hat{\xi}_{ijk}^{\a\b\c}$ is the 
section of the bundle $(\phi_{ik}^{\c},
\phi_{jk}^{\b}\circ \phi_{ij}^{\a})^{-1}Q_{k}$ 
on $Y_{i}\vert_{X_{ijk}^{\a\b\c}}$ given 
by 
$$
\hat{\xi}_{ijk}^{\a\b\c}(y_{i}) = 
m_{Q}(\hat{a}(\sigma_{jk}^{\b},\sigma_{ij}^{\a},
y_{i})\otimes \hat{m}(\sigma_{ijk}^{\a\b\c}
\otimes e(y_{i}))) 
$$
where $y_{i}\in Y_{i}\vert_{X_{ijk}^{\a\b\c}}$.  
We want to show that the various bundle 
gerbe products $m_{ijk}^{\a\b\c}$ extend to define   
a bundle gerbe product $m_{L}:L_{jk}\otimes 
L_{ij}\to L_{ik}$.  For this, we need to 
know that the following diagram commutes: 
$$ 
\xymatrix{ 
L_{jk}^{\b}\otimes L_{ij}^{\a} \ar[r]^-{m_{ijk}^{\a\b\c}} 
\ar[d]_-{\psi_{jk}^{\b\b^{'}}\otimes \psi_{ij}^{\a\a^{'}}} 
& L_{ik}^{\c} \ar[d]^-{\psi_{ik}^{\c\c^{'}}}                \\ 
L_{jk}^{\b^{'}}\otimes L_{ij}^{\a^{'}} 
\ar[r]_-{m_{ijk}^{\a^{'}\b^{'}\c^{'}}} & L_{ik}^{\c^{'}}.    } 
$$
The commutativity of this diagram follows from 
equation~\ref{eq:12.9}.  We also need to define the 
bundle gerbe morphism $\bar{\eta}:\pi_{1}^{-1}L
\otimes Q\to \pi_{2}^{-1}L$.  Again we define 
$\bar{\eta}$ locally and then show that it 
is compatible with the glueing isomorphisms 
$\psi_{ij}^{\a\b}$.  We define $\bar{\eta}^{\a}_{ij}:
\pi_{1}^{-1}L^{\a}_{ij}\otimes Q\to \pi_{2}^{-1}
L_{ij}^{\a}$ by putting $\hat{\eta}_{ij}^{\a}(
u_{ij}^{\a}\otimes u) = (u_{ij}^{\a}\circ u)
\hat{a}(\s_{ij}^{\a},y_{i},y_{1})^{-1}$ where 
$u_{ij}^{\a}\in L_{ij}^{\a}$ and $u\in Q$.  One can 
then show that the diagram below commutes: 
$$
\xymatrix{ 
\pi_{1}^{-1}L_{ij}^{\a}\otimes Q 
\ar[r]^-{\hat{\eta}^{\a}_{ij}} 
\ar[d]_-{\pi_{1}^{-1}\psi_{ij}^{\a\b}\otimes 
1_{Q}} & \pi_{2}^{-1}L_{ij}^{\a} 
\ar[d]^-{\pi_{2}^{-1}\psi_{ij}^{\a\b}}      \\ 
\pi_{1}^{-1}L_{ij}^{\b}\otimes Q \ar[r]_-{
\hat{\eta}^{\b}_{ij}} & \pi_{2}^{-1}L_{ij}^{\b}. } 
$$
Hence the various $\hat{\eta}^{\a}_{ij}$ glue 
together to form a bundle gerbe morphism 
$\eta_{ij}:\pi_{1}^{-1}L_{ij}\otimes Q\to 
\pi_{2}^{-1}L_{ij}$.  Similarly, one can define 
the transformation $\theta$ locally and 
show that it glues together to form a globally defined 
transformation which satisfies the required coherency condition.  
\end{note} 

\section{Trivial stable bundle 2-gerbes}
\label{sec:12.2} 

We first define the notion of a 
\emph{trivial} stable bundle 2-gerbe.  
So let $(T,Z,X)$ be a bundle gerbe 
on $X$ where $\pi:X\to M$ is a 
surjection admitting local sections.  
Form the bundle gerbe $(\d(T),\d(Z),X^{[2]})$ 
where $\d(T) = \pi_{1}^{-1}T\otimes \pi_{2}^{-1}
T^{*}$ and where $\d(Z) = \pi_{1}^{-1}Z\times_{
X^{[2]}}\pi_{2}^{-1}Z$.   
We want to show that the quadruple 
$(\d(T),\d(Z),X,M)$ 
defines a stable bundle 2-gerbe.  

We first need to define a stable morphism 
$$ 
m:\pi_{1}^{-1}(\d(T))\otimes \pi_{3}^{-1}(\d(T))
\to \pi_{2}^{-1}(\d(T))
$$
Such a stable morphism will correspond to a 
trivialisation $(L_{m},\phi_{m})$ of the bundle gerbe 
$$
(\pi_{1}^{-1}(\d(T))
\otimes \pi_{3}^{-1}(\d(T)) 
)^{*}\otimes \pi_{2}^{-1}(
\d(T)).  
$$
Let us use the notation of Chapter 7 to denote 
the bundle gerbe $(\pi_{1}^{-1}T,\pi_{1}^{-1}Z,
X^{[2]})$ as $(T_{2},Z_{2},X^{[2]})$ and the 
bundle gerbe $(\pi_{1}^{-1}\pi_{2}^{-1}T, 
\pi_{1}^{-1}\pi_{2}^{-1}Z,X^{[3]})$ as $(T_{3},
Z_{3},X^{[3]})$ and so on.   
Hence the bundle gerbe 
above can be rewritten in the new notation as 
$$
(T_{3}\otimes T_{2}^{*}\otimes T_{2}\otimes T_{1}^{*})
^{*}\otimes (T_{3}\otimes T_{1}^{*}).  
$$
We seek a trivialisation $(L_{m},\phi_{m})$ 
of this bundle gerbe.  So $L_{m}$ will 
be a $\cstar$ bundle on 
$$
Z_{3}\times_{\pi}Z_{2}\times_{\pi}Z_{2}
\times_{\pi}Z_{1}\times_{\pi}Z_{3}\times_{\pi}
Z_{1}.  
$$
The obvious candidate for $L_{m}$ is 
$T_{3}\otimes T_{2}\otimes T_{1}^{*}$.  It is easy 
to check that this trivialises the bundle gerbe.  
One can also check that there is a transformation 
$a:m\circ (m\otimes 1)\Rightarrow m\circ (1\otimes 
m)$ of stable morphisms of bundle gerbes over 
$X^{[4]}$ which satisfies the coherency condition 
over $X^{[5]}$.  

We now define what it means for a stable bundle 
2-gerbe to be trivial.  
\begin{definition} 
\label{def:trivialstableb2g} 
We say that a stable bundle 2-gerbe $(Q,Y,X,M)$ 
is \emph{trivial} if there is a bundle gerbe 
$(T,Z,X)$ on $X$ and a stable morphism $f:Q\to 
\d(T)$ over $X^{[2]}$.  We also require that there 
is a transformation of stable morphisms of bundle 
gerbe over $X^{[3]}$ as pictured in the following diagram.  
$$
\xymatrix{ 
\pi_{1}^{-1}Q\otimes \pi_{3}^{-1}Q 
\ar[d]_-{m} \ar[rr]^-{\pi_{1}^{-1}f\otimes 
\pi_{3}^{-1}f} & & \pi_{1}^{-1}\d(T)
\otimes \pi_{3}^{-1}\d(T) \ar[d] 
\ar @2{->}[dl]^-{\theta}                    \\ 
\pi_{2}^{-1}Q \ar[rr]^-{\pi_{2}^{-1}f} & & 
\pi_{2}^{-1}\d(T).                            } 
$$
Finally we require that $\theta$ satisfies 
the compatibility condition $\d(\theta) = a$ 
with $a$. 
\end{definition}    

\begin{note} 
\label{note:stableb2gmorphism} 
This definition suggests that we could 
define a stable bundle 2-gerbe morphism 
from a stable bundle 2-gerbe $(P,Y,
X,M)$ to a stable bundle 2-gerbe $(Q,Z,X,M)$ 
to be a stable morphism $f:P\to Q$ from 
the bundle gerbe $(P,Y,X^{[2]})$ to the 
bundle gerbe $(Q,Z,X^{[2]})$ together with a 
transformation $\theta_{f}$ as pictured in the 
following diagram 
$$
\xymatrix{ 
\pi_{1}^{-1}P\otimes \pi_{3}^{-1}P \ar[r]^-{m} 
\ar[d]_-{\pi_{1}^{-1}f\otimes \pi_{3}^{-1}f} & 
\pi_{2}^{-1}P \ar[d]^-{\pi_{2}^{-1}f} 
\ar @2{->}[dl]^-{\theta_{f}}                    \\ 
\pi_{1}^{-1}Q\otimes \pi_{3}^{-1}Q \ar[r]_-{m} & 
\pi_{2}^{-1}Q.                                   } 
$$
$\theta$ is required to satisfy the compatibility 
condition $\d(\theta) = a_{P}^{*}\otimes a_{Q}$.  
One could then go further and show that any 
stable bundle 2-gerbe morphism was invertible 
up to a coherent `stable bundle 2-gerbe transformation'.  
\end{note} 
  
The following Lemma is not difficult to 
prove.  

\begin{lemma} 
\label{lemma:trivialstableb2gclass} 
The class in $\check{H}^{3}(M;\underline{\C}^{\times}_{M})$ 
represented by the cocycle $g_{ijkl}$ associated 
to a trivial stable bundle 2-gerbe 
$(Q,Y,X,M)$ is zero.  
\end{lemma}

\setcounter{chapter}{12}
\chapter{The relationship of bundle 2-gerbes with 
2-gerbes}
\label{chapter:13} 
\section{Simplicial gerbes and 2-gerbes}
\label{section:13.1} 
We first recall the definition of  
a simplicial gerbe (see \cite{BryMcL1} 
and \cite{BryMcL2}).  A simplicial  
gerbe is defined as follows.  

\begin{definition}[\cite{BryMcL1}, \cite{BryMcL2}] 
\label{def:13.1.1}  
Let $X=\{X_{p}\}$ be a simplicial 
manifold.  A \emph{simplicial gerbe} on 
$X$ consists of the 
following data.
\begin{enumerate} 
\item  A gerbe $\mathcal{Q}$ on 
$X_{1}$.  We will only be 
interested in the case where 
$\mathcal{Q}$ has band equal to $\underline{\C}^{\times}_{X_{1}}$.  

\item A morphism of gerbes 
$$
m:d_{0}^{*}\mathcal{Q}\otimes 
d_{2}^{*}\mathcal{Q}
\to d_{1}^{*}\mathcal{Q} 
$$
over $X_{2}$ where the $d_{i}$ 
denote the face operators of 
$X_{\cdot}$. 

\item An invertible natural transformation $\theta$  
between morphisms of gerbes 
over $X_{3}$ as 
pictured in the following 
diagram 

$$ 
\xymatrix{ 
d_{1}^{*}d_{0}^{*}\mathcal{Q}\otimes 
d_{3}^{*}(d_{0}^{*}\mathcal{Q}\otimes 
d_{2}^{*}\mathcal{Q}) \ar[r]^-{\simeq}  
\ar[d]_{1\otimes d_{3}^{*}m} & 
d_{0}^{*}(d_{0}^{*}\mathcal{Q}\otimes 
d_{2}^{*}\mathcal{Q})\otimes d_{2}^{*}
d_{2}^{*}\mathcal{Q} \ar[d]^{d_{0}^{*}m
\otimes 1}                                   \\ 
d_{1}^{*}d_{0}^{*}\mathcal{Q}\otimes 
d_{3}^{*}d_{1}^{*}\mathcal{Q} \ar[d]_-{\simeq}  
 & 
d_{0}^{*}d_{1}^{*}\mathcal{Q}\otimes 
d_{2}^{*}d_{2}^{*}\mathcal{Q} \ar[d]^-{\simeq}   
\ar @<4.2ex> @2{->}[l]^-{\theta}     \\ 
d_{1}^{*}(d_{0}^{*}\mathcal{Q}\otimes 
d_{2}^{*}\mathcal{Q}) \ar[d]_{d_{1}^{*}m} 
& d_{2}^{*}(d_{0}^{*}\mathcal{Q}\otimes 
d_{2}^{*}\mathcal{Q}) \ar[d]^{d_{2}^{*}m}   \\ 
d_{1}^{*}d_{1}^{*}\mathcal{Q} \ar[r]^-{\simeq} & 
d_{2}^{*}d_{1}^{*}\mathcal{Q}.              } 
$$
Here the equivalences of gerbes denoted by 
$\stackrel{\simeq}{\to}$ are the equivalences of 
gerbes $\phi_{f,g}$ of Lemma~\ref{lemma:6.3.2}. 

\item $\theta$ is required to satisfy a   
coherency condition on 
$X_{4}$.  This is $\d(\theta) = 1$, 
where we use the notation of 
Note~\ref{note:6.6.6}.  
\end{enumerate} 
\end{definition} 

Note that this differs from the definition 
of simplicial gerbe given in \cite{BryMcL1} 
and \cite{BryMcL2} in that $m:d_{0}^{*}\mathcal{Q}
\otimes d_{2}^{*}\mathcal{Q}\to d_{1}^{*}
\mathcal{Q}$ is only a morphism of gerbes, 
rather than an equivalence of gerbes.  

We also recall the definition of a 
2-gerbe (see \cite{Bre}, \cite{BryMcL1} 
and \cite{BryMcL2}).  

\begin{definition}[\cite{Bre},\cite{BryMcL1},\cite{BryMcL2}]  
\label{def:13.1.2} 
A \emph{2-gerbe} on a manifold $M$ consists 
of the following data: 
\begin{enumerate} 
\item  For each open subset $U \subset M$ 
there is an assignment of a bicategory   
$\B_{U}$ to $U$, $U \mapsto \B_{U}$.  
Given another open subset 
$V \subset U \subset M$ there is a 
restriction bifunctor 
$\rho_{U,V}:\B_{U} \to \B_{V}$ 
such that if we have a third open 
subset $W$ with 
$W \subset V \subset U \subset M$ then 
$\rho_{V,W}\circ \rho_{U,V} = \rho_{U,W}$. 
Such an assignment is called a 
pre-sheaf of bicategories on 
$M$. 

\item  Various glueing laws for 
objects, 1-arrows and 2-arrows are 
required to hold.  The glueing 
conditions for objects are as follows.  
Let $U \subset M$ and let 
$\{U_{i}\}_{i \in I}$ be a covering 
of $U$ by open subsets $U_{i}$ of 
$M$.  Suppose we are given objects 
$A_{i}$ in each $\B_{U_{i}}$ together 
with 1-arrows 
$\a_{ij}:\rho_{U_{i},U_{ij}}(A_{i}) \to \rho_{U_{j},U_{ij}}$ 
between the restrictions of the 
objects $A_{i}, A_{j}$ to $\B_{U_{ij}}$.  
So the $\a_{ij}$ belong to $\B_{U_{ij}}$.  
Suppose also that we are given 
2-arrows $\phi_{ijk}$ of $\B_{U_{ijk}}$ 
such that 
$$ 
\phi_{ijk} : \rho_{U_{jk},U_{ijk}}(\a_{jk}) 
\circ \rho_{U_{ij},U_{ijk}}(\a_{ij}) \Rightarrow  
\rho_{U_{ik},U_{ijk}}(\a_{ik}), 
$$
where $\circ$ denotes the 
application of the composition 
functor in $\B_{U_{ijk}}$.  We try 
and indicate this in the following 
diagram which lives in $\B_{U_{ijk}}$ 
(we abuse notation and omit the 
restriction bifunctors $\rho$).  
$$ 
\diagram 
A_{i} \rrtwocell^{\a_{jk}\circ \a_{ij}}_
{\a_{ik}}{\ \ \ \phi_{ijk}}  &  & A_{k}.        \\  
\enddiagram 
$$
We can now construct two 2-arrows 
$\a_{kl}\circ \a_{jk}\circ \a_{ij} \Rightarrow \a_{il}$ 
in $\B_{U_{ijkl}}$ as pictured in the 
following diagram 
$$
\xymatrix{ 
& & \a_{kl}\circ (\a_{jk}\circ \a_{ij}) 
\ar @2{->}[drr]^-{1_{\a_{kl}}\circ \phi_{ijk}} & & \\  
(\a_{kl}\circ \a_{jk})\circ \a_{ij} 
\ar @2{->}[urr]^-{a} \ar @2{->}[dr]_-{\phi_{jkl}\circ 
1_{\a_{ij}}} & & & & \a_{kl}\circ \a_{ik} 
\ar @2{->}[dl]^-{\phi_{ikl}}                       \\   
& \a_{jl}\circ \a_{ij} \ar @2{->}[rr]_-{\phi_{ijl}} 
& & \a_{il} &                                       } 
$$
Finally we can state the 
glueing axiom for objects by saying that if these 
two 2-arrows coincide   
then there is an object $A$ belonging 
to $\B_{U}$ and a 1-arrow 
$\b$ in $\text{Hom}(\rho_{U,U_{i}}(A),A_{i})$ which 
is compatible with the glueing 1-arrows 
$\a_{ij}$ up to a coherent 2-arrow in the 
sense of Proposition~\ref{prop:6.6.1}.   

The glueing conditions for 1-arrows 
and 2-arrows are as follows.  
For an open set $U\subset M$ and 
objects $A_{1}$ and $A_{2}$ of 
$\B_{U}$ we require that 
$\text{Hom}(A_{1},A_{2})$ is a 
sheaf of categories.  
Thus, given an open cover 
$\{V_{i}\}_{i\in I}$ of $U$ 
by open sets $V_{i}\subset M$ together  
with 1-arrows $\a_{i}:A_{1}|_{V_{i}}\to A_{2}|_{V_{i}}$ 
in $\text{Hom}(A_{1},A_{2})|_{V_{i}}$ 
and 2-arrows $\phi_{ij}:\a_{i}\Rightarrow \a_{j}$ 
which satisfy the cocycle condition 
$\phi_{jk}\phi_{ij} = \phi_{ik}$, 
then there is a 1-arrow $\a$ in 
$\text{Hom}(A_{1},A_{2})$ which restricts 
to each $\a_{i}$.  The glueing 
conditions for 2-arrows are handled 
similarly.  

Such an assignment $U\mapsto \B_{U}$ 
of a pre-sheaf of bicategories to 
$M$ which satisfies the effective descent 
or glueing conditions on objects, 1-arrows 
and 2-arrows is called a 2-stack or a 
sheaf of bicategories.  

\item For each point $m \in M$ there 
is an open neighbourhood $U$ of 
$m$ such that $\B_{U}$ is nonempty; 
that is, the set of objects of $\B_{U}$ 
is nonempty. 

\item Given an open subset $U \subset M$,  
any two objects $A_{1}$ and $A_{2}$ 
of $\B_{U}$ can be locally joined 
by a 1-arrow.  That is there is 
a covering $\{U_{i}\}_{i \in I}$ of 
$U$ by open sets of $M$ such that 
there exist 1-arrows 
$\a_{i} :\rho_{{U},U_{i}}(A_{1})\to \rho_{U,U_{i}}(A_{2})$ 
in $\B_{U_{i}}$.  

\item We also require that 1-arrows 
are weakly invertible. 
This means that given an open 
set $U\subset M$ and objects $A_{1}$ 
and $A_{2}$ of $\B_{U}$ together 
with a 1-arrow $\a:A_{1}\to A_{2}$, 
there exists a  1-arrow 
$\b:A_{2}\to A_{1}$   
of $\B_{U}$ together with 
2-arrows $\phi:\b\circ \a\Rightarrow Id_{A_{1}}$ 
in $\B_{U}$. 

\item Finally we require that all 
two arrows are invertible.     
\end{enumerate} 
\end{definition} 

We will also need to have the 
notion of a 2-gerbe on $M$ \emph{bound} 
by $\underline{\C}^{\times}_{M}$.  This means that any 
two 1-arrows of $\mathcal{B}_{U}$ 
can be locally joined by a 2-arrow 
and for every 1-arrow $\a$ there 
is a given isomorphism $\underline{\C}^{\times}_{M} 
\simeq \underline{Aut}(\a)$ 
which is required to be compatible 
with restrictions to smaller open 
sets and is compatible with 2-arrows. 
There is also a notion of 
equivalence of 2-gerbes --- we 
will not spell this out but refer 
instead to \cite{Bre}.  
One can associate to a 2-gerbe 
$\B$ on $M$ bound by $\underline{\C}^{\times}_{M}$ 
a \v{C}ech three cocycle 
$g_{ijkl}$ in $\check{H}^{3}(M;\underline{\C}^{\times}_{M})$.  
We will not do this here and refer to 
\cite{Bre} and \cite{BryMcL1} for a 
discussion of this.  There is an 
isomorphism between equivalence 
classes of 2-gerbes and 
$\check{H}^{3}(M;\underline{\C}^{\times}_{M})$ via this 
cocycle.  As is the case for gerbes, 
a 2-gerbe with zero class in 
$\check{H}^{3}(M;\underline{\C}^{\times}_{M})$ is trivial 
precisely when it has a global 
object.  For more details we refer 
to \cite{Bre}. 

Let $\U = \{U_{i}\}_{i\in I}$ be an 
open cover of $M$ such that each 
non-empty intersection $U_{i_{0}}\cap 
\cdots \cap U_{i_{p}}$ is contractible.  
Let $\mathcal{B}$ be a 2-gerbe on $M$ 
bound by $\underline{\C}^{\times}_{M}$.  
Let $\mathcal{B}_{i}$ denote the 
restriction of $\mathcal{B}$ to $U_{i}$ 
for each $i\in I$.  Since $U_{i}$ is 
assumed contractible, the 2-gerbe $\mathcal{B}_{i}$ 
has a global object $A_{i}$ say.  Denote 
by $A_{i}|_{U_{ij}}$ and $A_{j}|_{U_{ij}}$ 
the restrictions of $A_{i}$ and $A_{j}$ 
to $U_{ij}$ respectively.  We have the 
presheaf of categories $\mathcal{B}_{ij}$ 
on $U_{ij}$ where $\mathcal{B}_{ij}$ 
has fibre at an open set $V\subset U_{ij}$ 
$\mathcal{B}_{ij}(V)$ equal to $\text{Hom}
(A_{i}|_{V},A_{j}|_{V})$ (we ignore any 
complexities that may result from the 
objects $A_{i}|_{U_{ij}}|_{V}$ and 
$A_{i}|_{V}$ not being equal but isomorphic 
by a coherent 1-arrow --- otherwise we get  
that $\mathcal{B}_{ij}$ is a fibred 
category (\cite{Bre}) or a pre-sheaf 
of categories in the sense of \cite{Bry}).  
By definition of the glueing 
laws in $\mathcal{B}$, $\mathcal{B}_{ij}$ 
is a sheaf of categories.  In fact, since 
$\mathcal{B}$ is bound by $\underline{\C}
^{\times}_{M}$, $\mathcal{B}_{ij}$ is a gerbe 
on $U_{ij}$ bound by $\underline{\C}^{\times}
_{M}$.  The composition functor 
$m:\text{Hom}(A_{j}|_{U_{ijk}},A_{k}|_{U_{ijk}})
\times \text{Hom}(A_{i}|_{U_{ijk}},A_{j}|_
{U_{ijk}})\to \text{Hom}(A_{i}|_{U_{ijk}},
A_{k}|_{U_{ijk}})$ gives rise to a morphism 
$\mathcal{B}_{jk}\times \mathcal{B}_{ij}\to 
\mathcal{B}_{ik}$ of gerbes `bound' by 
the morphism of sheaves $\underline{\C}^
{\times}_{M}\times \underline{\C}^{\times}_{M}
\to \underline{\C}^{\times}_{M}$.  One can 
show that this morphism in turn defines a 
morphism of gerbes $\mathcal{B}_{jk}\otimes 
\mathcal{B}_{ij}\to \mathcal{B}_{ik}$ over 
$U_{ijk}$.  The associator natural 
transformation in the bicategory $\mathcal{B}_
{U_{ijkl}}$ then provides a natural transformation 
$\phi$ between morphisms of gerbes as pictured in the 
following diagram (we have omitted any 
restriction functors and equivalences of 
gerbes that may result by passing between 
them for convenience): 
$$
\xymatrix{ 
\mathcal{B}_{kl}\otimes \mathcal{B}_{jk}
\otimes \mathcal{B}_{ij} \ar[r]^-{1\otimes 
m} \ar[d]_-{m\otimes 1} & \mathcal{B}_{jl}
\otimes \mathcal{B}_{ij} \ar[d]^-{m} 
\ar @2{->}[dl]^-{\phi}                     \\ 
\mathcal{B}_{kl}\otimes \mathcal{B}_{ik} 
\ar[r]_-{m} & \mathcal{B}_{il}.            } 
$$
It follows that we have defined a 
simplicial gerbe on the simplicial manifold 
$X = \{X_{p}\}$ with $X_{0} = \coprod U_{i}$, 
$X_{1} = \coprod U_{ij}$ and so on.  It is not 
hard to construct a \v{C}ech 3-cocycle $g_{ijkl}$ 
representing a class in $\check{H}^{3}(M;\underline{\C}^{\times}_{M}) 
= H^{4}(M;\Z)$ associated to this simplicial 
gerbe (compare with Theorem 5.7 of \cite{BryMcL1}).     
It is also easy to see that the class in 
$H^{4}(M;\Z)$ associated to this simplicial gerbe is 
equal to the class of the 2-gerbe in $H^{4}(M;\Z)$.  

This simplicial gerbe on $X = \{X_{p}\}$ with 
$X_{p} = \coprod U_{i_{0}}\cap \cdots \cap 
U_{i_{p}}$ is an example of a simplicial 
gerbe on the special class of simplicial manifolds 
formed by taking iterative fibre products of 
a surjection $\pi:X\to M$.  I had hoped to relate 
bundle 2-gerbes to 2-gerbes by first showing that 
a bundle 2-gerbe gave rise to a simplicial gerbe 
on the simplicial manifold $X = \{X_{p}\}$ 
with $X_{p} = X^{[p+1]}$ and then showing that 
any simplicial gerbe on $X$ gave rise to a 2-gerbe 
on $M$.  Unfortunately, as noted in Chapter~\ref{chapter:6} 
(Note~\ref{note:6.5.2}), checking directly that the construction which 
associates a gerbe $\mathcal{G}(P)$ to a given 
bundle gerbe $P$ commutes with equivalences 
of gerbes induced by successive tensor products 
and transitivity of pullbacks gets very messy 
and I was unable to produce a convincing proof 
that this was the case.  Note that any simplicial 
gerbe $\Q$ on the simplicial manifold $X=\{X_{p}\}$ formed by taking 
iterated fibre products of a  
surjection $\pi:X\to M$ admitting local sections 
gives rise to a `local description' of a 2-gerbe 
as described above --- ie gerbes $\Q_{ij}$ on 
$U_{ij}$ and morphisms of gerbes $m:\Q_{jk}\otimes 
\Q_{ij}\to \Q_{ik}$ and so on --- by pulling back the      
simplicial gerbe on $X$ to the simplicial gerbe on 
simplicial manifold associated to the nerve of the 
open cover $\U = \{U_{i}\}$ using the simplicial map 
$\coprod U_{i_{1}}\cap \cdots \cap U_{i_{p}}\to X^{[p]}$, 
$(i_{1},\ldots,i_{p},m)\mapsto (s_{i_{1}}(m),\ldots 
,s_{i_{p}}(m))$.  Hence every simplicial gerbe on 
$X = \{X_{p}\}$ has a class in $H^{4}(M;\Z)$ 
associated to it.   

Let $\pi:X\to M$ be a surjection 
admitting local sections.  Form the 
simplicial manifold $X=\{X_{p}\}$ 
with $X_{p} = X^{[p+1]}$.  
Let $(\Q,m,\theta)$ 
be a simplicial gerbe on $X$.  
We define a sheaf of 2-categories $\B$ on $M$ as 
follows.   

For an open set $U\subset M$, the 
objects of the 2-category $\B_{U}$ 
will be gerbes $\T_{U}$ 
on $X_{U} = \pi^{-1}(U)$ which `trivialise' 
$\Q$ on $X^{[2]}_{U}$.  Thus an 
object will consist of a gerbe 
$\T_{U}$ on $X_{U}$ together 
with a morphism of gerbes 
$\eta:\pi_{1}^{*}\mathcal{T}_{U}\otimes \Q 
\to \pi_{2}^{*}\mathcal{T}$ plus an 
invertible natural transformation 
$\chi$ as pictured in the following 
diagram: 
$$
\xymatrix{ 
\pi_{1}^{*}\pi_{1}^{*}\T_{U}\otimes 
\pi_{1}^{*}\Q\otimes \pi_{3}^{*}\Q 
\ar[rr]^-{1\otimes m} \ar[d]_-{\pi_{1}^{*}\eta
\otimes 1} & & \pi_{1}^{*}\pi_{1}^{*}\T_{U} 
\ar[d]^-{\simeq}                               \\ 
\pi_{1}^{*}\pi_{2}^{*}\T_{U}\otimes 
\pi_{3}^{*}\Q \ar[d]_-{\simeq} \ar @<-4.2ex>  
@2{->}[rr]^-{\chi} & & \pi_{2}^{*}\pi_{1}^{*}\T_{U}
\otimes \pi_{2}^{*}\Q \ar[d]^-{\pi_{2}^{*}
\eta}                                              \\
\pi_{3}^{*}\pi_{1}^{*}\T_{U}\otimes \pi_{3}^{*}\Q 
\ar[dr]_-{\pi_{3}^{*}\eta} & & \pi_{2}^{*}\pi_{2}^{*}
\T_{U} \ar[dl]_-{\simeq}                            \\ 
& \pi_{3}^{*}\pi_{2}^{*}\T_{U}.                     } 
$$
$\chi$ is required to satisfy the 
compatibility condition $\d(\chi) = \theta$ 
with the natural transformation $\theta$ --- 
see Note~\ref{note:6.6.6} for the notation 
`$\d$'.   

Given objects $(\T_1,\eta_1,\chi_1)$ and 
$(\T_2,\eta_2,\chi_2)$ a 1-arrow  
$(\phi,\tau_{\phi}):(\mathcal{T}_{1},
\eta_{1},\chi_{1})\to (\mathcal{T}_{2},
\eta_{2},\chi_{2})$ consists of an 
equivalence of gerbes $\phi:\mathcal{T}_{1}
\to \mathcal{T}_{2}$ over $X_{U}$ and an 
invertible natural transformation 
of gerbes $\tau_{\phi}$ as pictured in the 
following diagram: 
$$
\xymatrix{ 
\pi_{1}^{*}\mathcal{T}_{1}\otimes \mathcal{Q} 
\ar[r]^{\eta_{1}}\ar[d]_{\pi_{1}^{*}\phi\otimes 
\text{id}} & \pi_{2}^{*}\mathcal{T}_{1} 
\ar[d]^{\pi_{2}^{*}\phi}                         \\ 
\pi_{1}^{*}\mathcal{T}_{2}\otimes \mathcal{Q} 
\ar[r]_{\eta_{2}} \ar @2{->}[ur]^{\tau_{\phi}} & 
\pi_{2}^{*}\mathcal{T}_{2}.                       } 
$$
We require that $\tau_{\phi}$ satisfies 
$\d(\tau_{\phi}) = \chi_{2}\otimes \chi_{1}^{-1}$.    
A 2-arrow $\theta:(\phi_{1},\tau_{\phi_{1}})
\Rightarrow (\phi_{2},\tau_{\phi_{2}}):(
\mathcal{T}_{1},\eta_{1},\chi_{1})\to 
(\mathcal{T}_{2},\eta_{2},\chi_{2})$ is an 
invertible natural transformation 
$\theta$ as pictured in the following 
diagram: 
$$
\xymatrix{ 
\mathcal{T}_{1} \ar[r]^{\phi_{1}} 
\ar @2{-}[d] & \mathcal{T}_{2} \ar @2{-}[d]    \\ 
\mathcal{T}_{1} \ar @2{->}[ur]^{\theta} 
\ar [r]_{\phi_{2}} & \mathcal{T}_{2}.          } 
$$
$\theta$ is required to be compatible with 
$\tau_{\phi_{1}}$ and $\tau_{\phi_{2}}$ in that 
$\d(\theta) = \tau_{\phi_{2}}\tau_{\phi_{1}}^{-1}$.  
Given another open subset $V$ with 
$V\subset U\subset M$, there is a natural 
restriction 2-functor $\rho_{U,V}:
\B_{U}\to \B_{V}$.  It follows that we 
have defined a presheaf of 2-categories.  
In actual fact the assignment $U\to \mathcal{B}_{U}$ 
is a sheaf of 2-categories as we will now 
see.  Firstly we check that the glueing axioms 
for objects are satisfied.  Let $U\subset M$ 
be an open set and let $\{U_{i}\}_{i\in I}$ be 
a covering of $U$ by open subsets $U_{i}$ of 
$M$.  Let $(\mathcal{T}_{i},\eta_{i},\chi_{i})$ 
be objects of $\mathcal{B}_{U_{i}}$ for $i\in I$. 
Suppose we have 1-arrows $(\phi_{ij},\tau_{\phi_{ij}}):
(\mathcal{T}_{i},\eta_{i},\chi_{i})
\to (\mathcal{T}_{j},\eta_{j},\chi_{j})$, 
where for convenience of notation we have omitted  
the restriction 2-functors $\rho_{U_{i},U_{ij}}$.   
Suppose we also have 
2-arrows $\theta_{ijk}:(\phi_{jk}\circ \phi_{ij},
\tau_{\phi_{jk}}\circ \tau_{\phi_{ij}})\Rightarrow 
(\phi_{ik},\tau_{\phi_{ik}})$ which satisfy the 
glueing condition 
$$
\theta_{ijl}(\theta_{jkl}\circ \text{id}_{\phi_{ij}}) 
= \theta_{ikl}(\text{id}_{\phi_{kl}}\circ \theta_{ijk}) 
$$
of (2) in Definition~\ref{def:13.1.2} above and where 
$\text{id}_{\phi_{ij}}$ is the identity 2-arrow 
at $(\phi_{ij},\tau_{\phi_{ij}})$, that is 
$\text{id}_{\phi_{ij}}$ is the identity natural 
transformation from $\phi_{ij}$ to itself.           
$\mathcal{T}_{i}$, $\phi_{ij}$ and $\theta_{ijk}$ 
provide 2-descent data for the 2-stack of 
gerbes on $X_{U}$ --- see \cite{Bre}[pages 32--33] 
and also Proposition~\ref{prop:6.6.1}.   
Consequently, there is a gerbe  
$\mathcal{T}$ on $X_{U}$ which 
restricts to $\mathcal{T}_{i}$ on $X_{U_{i}}$.  
In other words there are equivalences of gerbes 
$\psi_{i}:\mathcal{T}|_{X_{i}}\stackrel{\simeq}
{\to}\mathcal{T}_{i}$ on $X_{i}$ which are compatible 
with the glueing 1-arrows $\phi_{ij}$ up 
to a coherent natural transformation $\xi_{ij}$ 
(the coherency condition that the $\xi_{ij}$ satisfy 
is $\xi_{ik}(\theta_{ijk}\circ 1_{\psi_{i}}) = 
\xi_{jk}(1_{\phi_{jk}}\circ \xi_{ij})$).    
We have to show also that the gerbe morphisms 
$\eta_{i}$ glue together to give a gerbe 
morphism $\eta:\pi_{1}^{*}\mathcal{T}\otimes 
\mathcal{Q}\to \pi_{2}^{*}\mathcal{T}$ and 
similarly that the $\chi_{i}$ glue together 
to give an invertible natural transformation 
$\chi$.  The gerbe morphisms $\eta_{i}:\pi_{1}^{*}
\mathcal{T}_{i}\otimes \mathcal{Q}|_{X_{i}^{[2]}}
\to \pi_{2}^{*}\mathcal{T}_{i}$ give rise 
to morphisms of gerbes $\tilde{\eta}_{i}:(\pi_{1}^{*}
\mathcal{T}\otimes \mathcal{Q})|_{X_{i}^{[2]}}\to 
(\pi_{2}^{*}\mathcal{T})|_{X_{i}^{[2]}}$ over 
$X_{i}^{[2]}$.  $\tilde{\eta}_{i}$ is defined by 
the following diagram. 
$$
\xymatrix{ 
(\pi_{1}^{*}\mathcal{T}\otimes \mathcal{Q})|_{X_{i}^{[2]}} 
\ar[d]_-{\simeq} \ar[r]^-{\tilde{\eta}_{i}} & 
(\pi_{2}^{*}\mathcal{T})|_{X_{i}^{[2]}}  \\  
\pi_{1}^{*}\mathcal{T}|_{X_{i}^{[2]}}  
\otimes \mathcal{Q}|_{X_{i}^{[2]}} \ar[d]_-{\simeq}  &  \\ 
\pi_{1}^{*}(\mathcal{T}_{X_{i}})\otimes \mathcal{Q}|_{X_{i}^{[2]}} 
\ar[d]_-{\pi_{1}^{*}\psi\otimes \text{id}_{\mathcal{Q}|_{X_{i}^{[2]}}}} & 
\pi_{2}^{*}(\mathcal{T}|X_{i}) \ar[uu]                  \\ 
\pi_{1}^{*}\mathcal{T}_{i}\otimes 
\mathcal{Q}|_{X_{i}^{[2]}} \ar[r]^-{\eta_{i}} & 
\pi_{2}^{*}\mathcal{T}_{i} \ar[u]^-{\pi_{2}^{*}\psi_{i}^{-1}}  } 
$$
We can define natural transformations 
$\tilde{\tau}_{ij}:\tilde{\eta}_{i}\Rightarrow 
\tilde{\eta}_{j}$ by $\tilde{\tau}_{ij} = 
\pi_{2}^{*}\xi_{ij}^{-1}\tau_{\phi_{ij}}\pi_{1}^{*}\xi_{ij}$.  
Because of the coherency conditions on 
$\tau_{\phi_{ij}}$, $\xi_{ij}$ and $\theta_{ijk}$ 
one can show that $\tilde{\tau}_{jk}\tilde{\tau}_{ij} = 
\tilde{\tau}_{ik}$.  Therefore by Proposition~\ref{prop:6.4.1} 
there is a morphism of gerbes $\eta:\pi_{1}^{*}\mathcal{T}
\otimes \mathcal{Q}\to \pi_{2}^{*}\mathcal{T}$ over 
$X_{U}^{[2]}$ and a natural transformation 
$\tau_{i}:\eta|_{X_{i}^{[2]}}\Rightarrow \eta_{i}$ 
between the morphisms of gerbes over $X_{i}^{[2]}$.  
Similarly one can show that the natural transformations 
$\chi_{i}$ glue together to give a natural 
transformation $\chi$ which is compatible with $\phi$ 
in the required sense.  
 
We also have to check that the glueing 
conditions for 1-arrows and 2-arrows are 
satisfied.  The glueing laws for 2-arrows 
are clearly satisfied, to show that the same 
is true for 1-arrows requires a little more work.  
So suppose that we have objects $(\mathcal{T}_{1},
\eta_{1},\chi_{1})$ and $(\mathcal{T}_{2},\eta_{2},
\chi_{2})$ of $\mathcal{B}_{U}$ and suppose that 
there are 1-arrows $(\phi_{i},\tau_{\phi_{i}})$ from 
the restriction of $(\mathcal{T}_{1},\eta_{1},\chi_{1})$ 
to $X_{U_{i}}$ to the restriction of $(\mathcal{T}_{2},
\eta_{2},\chi_{2})$ to $X_{U_{i}}$ for an open cover 
$\{U_{i}\}$ of $U$ and glueing 2-arrows $\theta_{ij}:
(\phi_{i},\tau_{\phi_{i}})\Rightarrow (\phi_{j},
\tau_{\phi_{j}})$.  Thus we have equivalences of gerbes 
$\phi_{i}:\mathcal{T}_{1}|_{X_{i}}\to \mathcal{T}_{2}|_{X_{i}}$ 
and natural transformations $\theta_{ij}:\phi_{i}|_{X_{ij}}
\Rightarrow \phi_{j}|_{X_{ij}}$ which satisfy the glueing 
condition $\theta_{jk}\theta_{ij} = \theta_{ik}$.  
Therefore by Proposition~\ref{prop:6.4.1} there 
is an equivalence of gerbes $\phi:\mathcal{T}_{1}\to 
\mathcal{T}_{2}$ over $X_{U}$ and natural transformations  
$\a_{i}:\phi|_{X_{i}}\Rightarrow \phi_{i}$ which are compatible 
with the $\theta_{ij}$.  We can use the natural transformations 
$\a_{i}$ and $\tau_{\phi_{i}}$ to construct natural 
transformations $\tilde{\tau}_{i}$ between the morphisms 
of gerbes bounding the  
following diagram    
$$
\xymatrix{ 
(\pi_{1}^{*}\mathcal{T}_{1}\otimes \mathcal{Q})|_{X_{i}^{[2]}} 
\ar[r]^-{\eta_{1}|_{X_{i}^{[2]}}} \ar[d]_{(\pi_{1}^{*}
\phi\otimes 1_{\mathcal{Q}})|_{X_{i}^{[2]}}} & 
(\pi_{2}^{*}\mathcal{T}_{1})|_{X_{i}^{[2]}} 
\ar[d]^-{(\pi_{2}^{*}\phi)|_{X_{i}^{[2]}}} 
\ar @2{->}[dl]^-{\tilde{\tau}_{i}}                  \\ 
(\pi_{1}^{*}\mathcal{T}_{2}\otimes \mathcal{Q})|_{X_{i}^{[2]}} 
\ar[r]_-{\eta_{2}|_{X_{i}^{[2]}}} & (\pi_{2}^{*}\mathcal{T}_{2})
|_{X_{i}^{[2]}}.                                         } 
$$ 
Because of the coherency conditions satisfied by 
$\a_{i}$ and $\tau_{\phi_{i}}$ with $\theta_{ij}$ one 
can show that $\tilde{\tau}_{i}|_{X_{ij}^{[2]}} = 
\tilde{\tau}_{j}|_{X_{ij}^{[2]}}$ and hence there 
is a natural transformation $\tau_{\phi}:\pi_{2}^{*}
\phi\circ \eta_{1}\Rightarrow \eta_{2}\circ (\pi_{1}^{*}
\phi\otimes 1_{\mathcal{Q}})$ with $\tau_{\phi}|_{X_{i}^{[2]}}
= \tilde{\tau}_{i}$.  Again one can show that the 
natural transformation $\tau_{\phi}$ satisfies the 
required coherency condition with $\chi_{1}$ and $\chi_{2}$.     
Thus we have shown that the assignment $U\mapsto \mathcal{B}_{U}$ 
is a sheaf of 2-categories.  
We have the following conjecture. 
\begin{conjecture} 
\label{prop:13.2.2} 
The sheaf of 2-categories $\mathcal{B}$ 
described above is actually a 2-gerbe.   
Furthermore the four 
class of the 2-gerbe $\B$ equals 
the four class of the simplicial 
gerbe $\Q$. 
\end{conjecture} 

Unfortunately I have been unable to give a concrete 
proof of this.  However one can prove a version 
of this conjecture with gerbes replaced by 
bundle gerbes and gerbe morphisms replaced by 
stable morphisms (bundle gerbe morphisms will 
not do here as they do not glue together 
properly).  

\begin{proposition} 
Let $(Q,Y,X,M)$ be a stable bundle 2-gerbe.  
Then there is a 2-gerbe $\mathcal{Q}$ on $M$ 
bound by $\underline{\C}^{\times}_{M}$ 
such that the class in $H^{4}(M;\Z)$ of 
$\mathcal{Q}$ is equal to the class of $Q$ 
in $H^{4}(M;\Z)$. 
\end{proposition} 

Let $U\subset M$ be an open subset of $M$.  
We define a bicategory $\mathcal{Q}_{U}$ 
as follows.  $\mathcal{Q}_{U}$ has as objects 
the bundle gerbes $(T,Z,X_{U})$ on 
$X_{U}$ such that there is a stable morphism 
$f:\pi_{1}^{-1}T\otimes Q\to \pi_{2}^{-1}
T$ of bundle gerbes over $X_{U}^{[2]}$ and 
a transformation $\chi$ between the stable 
morphisms of bundle gerbes over $X_{U}^{[3]}$ as 
pictured in the following diagram.  
$$
\xymatrix{ 
\pi_{1}^{-1}\pi_{1}^{-1}T\otimes 
\pi_{1}^{-1}Q\otimes \pi_{3}^{-1}Q 
\ar[d]_{\pi_{1}^{-1}f\otimes 1} 
\ar[rr]^-{1\otimes m} & & 
\pi_{1}^{-1}\pi_{1}^{-1}T\otimes 
\pi_{2}^{-1}Q \ar @2{-}[d]              \\ 
\pi_{1}^{-1}\pi_{2}^{-1}T\otimes 
\pi_{3}^{-1}Q \ar @2{->}[rr]^-{\chi} 
\ar @2{-}[d] & & \pi_{2}^{-1}\pi_{1}
^{-1}T\otimes \pi_{2}^{-1}Q 
\ar[d]^-{\pi_{2}^{-1}f}                 \\ 
\pi_{3}^{-1}\pi_{1}^{-1}T\otimes 
\pi_{3}^{-1}Q \ar[dr]_-{\pi_{3}^{-1}f} 
& & \pi_{2}^{-1}\pi_{2}^{-1}T 
\ar @2{-}[dl]                           \\ 
& \pi_{3}^{-1}\pi_{2}^{-1}T.        } 
$$ 
Therefore $\chi$ is a section of the $\cstar$ 
bundle $B$ on $X_{U}^{[3]}$ such that 
$$
L_{\pi_{2}^{-1}f\circ (1_{\pi_{1}^{-1}\pi_{1}^{-1}
T}\otimes m)} = L_{\pi_{3}^{-1}f\circ (\pi_{1}^{-1}
f\otimes 1_{\pi_{3}^{-1}Q})}\otimes \pi^{-1}B, 
$$
where 
$$
\pi:\pi_{1}^{-1}\pi_{1}^{-1}Z\times_{X_{U}^{[3]}}
\pi_{1}^{-1}Y\times_{X_{U}^{[3]}}\pi_{3}^{-1}Y
\times_{X_{U}^{[3]}}\pi_{3}^{-1}\pi_{2}^{-1}Z\to 
X_{U}^{[3]}
$$
is the projection.  One can show, see Chapter~\ref{chapter:12}, 
that there is an isomorphism $\d(B)\simeq A$ of 
$\cstar$ bundles over $X_{U}^{[4]}$.  We finally 
require that the induced section $\d(\chi)$ of 
$\d(B)$ maps to $a$ under this isomorphism.   
Given two such objects $(T,f_{T},\chi_{T})$ 
and $(S,f_{S},\chi_{S})$ a 1-arrow $(T,f_{T},
\chi_{T})\to (S,f_{S},\chi_{S})$ consists 
of a stable morphism $g:T\to S$ and a 
transformation $\a_{g}$ between the stable 
morphisms of bundle gerbes over $X_{U}^{[2]}$ 
as pictured in the following diagram: 
$$
\xymatrix{ 
\pi_{1}^{-1}T\otimes Q \ar[r]^-{f_{T}} 
\ar[d]_-{\pi_{1}^{-1}g\otimes 1} & 
\pi_{2}^{-1}T \ar[d]^-{\pi_{2}^{-1}g} 
\ar @2{->}[dl]_-{\a_{g}}               \\ 
\pi_{1}^{-1}S\otimes Q \ar[r]_-{f_{S}} 
& \pi_{2}^{-1}S.                       } 
$$
We also require that $\a_{g}$ satisfies the 
coherency condition $\d(\a_{g}) = \chi_{T}
\otimes \chi_{S}^{*}$ over $X_{U}^{[3]}$.  
Finally a 2-arrow 
$\theta:(g,\a_{g})\Rightarrow (h,\a_{h}):
(T,f_{T},\chi_{T})\to (S,f_{S},\chi_{S})$ is 
a transformation $\theta:g\Rightarrow h$ 
between stable morphisms of bundle gerbes 
over $X_{U}$ which satisfies the coherency 
condition $\d(\theta) = \a_{g}\otimes \a_{h}^{*}$ 
over $X_{U}^{[2]}$.  By Proposition~\ref{prop:stablebigrpd},  
$\mathcal{Q}_{U}$ is a bicategory.  
The proof that $U\mapsto \mathcal{B}_{U}$ is a 
sheaf of 2-categories given above carries through 
to the case of bundle gerbes and stable morphisms 
(although now $\mathcal{Q}$ is a sheaf of bicategories).  
Therefore we have only have to show that the 
sheaf of bicategories $\mathcal{Q}$ is (i) locally 
non-empty, (ii) any two objects of $\mathcal{Q}_{U}$ 
may locally be joined by a 1-arrow, (iii) 1-arrows are 
weakly invertible and (iv) 2-arrows are invertible.  

\begin{proof} 
We first prove that the sheaf of bicategories 
$\mathcal{Q}$ is locally non-empty.  Choose 
an open cover $\U = \{U_{i}\}_{i\in I}$ of 
$M$ such that there exist local sections 
$s_{i}:U_{i}\to X$ of $\pi:X\to M$.  Form 
the maps $\hat{s}_{i}:X_{i}\to X^{[2]}$ 
sending $x\in X_{i}$ to $(x,s_{i}(\pi(x)))$.  
Let $Q_{i}=\hat{s}_{i}^{-1}Q$.  So $Q_{i} = 
(Q_{i},Y_{i},X_{i})$ is a bundle gerbe on $X_{i}$.  
Now pullback the stable morphism 
$m:\pi_{1}^{-1}Q\otimes \pi_{3}^{-1}Q\to 
\pi_{2}^{-1}Q$ using the map $X_{i}^{[2]}\to 
X^{[3]}$, $(x_{1},x_{2})\mapsto (x_{1},x_{2},
s_{i}(m))$ where $m=\pi(x_{1})=\pi(x_{2})$.  
Let $f_{i}$ denote this pullback stable 
morphism $f_{i}:\pi_{1}^{-1}Q_{i}\otimes Q
\to \pi_{2}^{-1}Q_{i}$.  To construct a 
transformation  $\chi_{i}:\pi_{2}^{-1}f_{i}
\circ (1_{\pi_{1}^{-1}\pi_{1}^{-1}Q_{i}}
\otimes m)\Rightarrow \pi_{3}^{-1}f_{i}\circ 
(\pi_{1}^{-1}f_{i}\otimes 1_{\pi_{3}^{-1}Q})$ 
we pullback the transformation $a^{-1}$ using 
the map $X_{i}^{[3]}\to X^{[4]}$, $(x_{1},
x_{2},x_{3})\mapsto (x_{1},x_{2},x_{3},s_{i}
(m))$ where again $m=\pi(x_{1})=\pi(x_{2})=
\pi(x_{3})$.  The coherency condition $\d(a)=
1$ satisfied by $a$ ensures that $\d(\chi_{i})
=a$.  Therefore for each $i\in I$ we have 
constructed objects $(Q_{i},f_{i},\chi_{i})$ 
of $\mathcal{Q}_{U_{i}}$.  Therefore $\mathcal{Q}$ 
is locally non-empty.  

We now show that $\mathcal{Q}$ is locally 
connected.  Let $(T,f,\chi)$ and $(S,g,\rho)$ 
be objects of $\mathcal{Q}_{U}$.  So $T$ 
and $S$ are bundle gerbes $(T,Z,X)$ and 
$(S,W,X)$.  We need 
to show that at least locally on $U$ these 
two objects are connected by a 1-arrow.  
We have stable morphisms $f:\pi_{1}^{-1}T\otimes 
Q\to \pi_{2}^{-1}T$ and $g:\pi_{1}^{-1}S\otimes 
Q\to \pi_{2}^{-1}S$ corresponding to trivialisations 
$(L_{f},\phi_{f})$ and $(L_{g},\phi_{g})$ 
of the bundle gerbes $(\pi_{1}^{-1}T\otimes Q)
^{*}\otimes \pi_{2}^{-1}T$ and $(\pi_{1}^{-1}S\otimes 
Q)^{*}\otimes \pi_{2}^{-1}S$ respectively.  
Form the bundle gerbe $(\pi_{1}^{-1}(S^{*}\otimes 
T))^{*}\otimes \pi_{2}^{-1}(S^{*}\otimes T)
\otimes (Q^{*}\otimes Q)$.  This is trivialised by 
the $\cstar$ bundle $L_{f}\otimes L_{g}^{*}$ 
which lives over $V\times_{\pi}Y^{[2]}$, where 
$$
V = \pi_{1}^{-1}(Z\times_{\pi}W)\times_{\pi} 
\pi_{2}^{-1}(Z\times_{\pi}W).   
$$
Consider the $\cstar$ bundle $\hat{L} = L_{f}\otimes 
L_{g}^{*}\otimes Q^{*}$ which lives over $V\times_{\pi}Y^{[2]}$.   
One can show that there is a section $\hat{\phi}$ 
of $\d_{Y^{[2]}}(\hat{L})$ satisfying 
the coherency condition $\d_{Y^{[2]}}(\hat{\phi}) = 
\underline{1}$.  $\hat{\phi}$ corresponds to the 
isomorphism given fibrewise by composition 
\begin{eqnarray*} 
&    & L_{f}(z,y_{1},z^{'})\otimes 
L_{g}^{*}(w,y_{2},w^{'})\otimes Q^{*}(y_{1},y_{2}) \\ 
& \simeq & L_{f}(z,y_{1},z^{'})\otimes 
L_{g}^{*}(w,y_{2},w^{'})\otimes Q^{*}(y_{1},y_{3})
\otimes Q^{*}(y_{3},y_{4})\otimes Q^{*}(y_{4},y_{2}) \\ 
& \simeq & L_{f}(z,y_{3},z^{'})\otimes 
L_{g}^{*}(w,y_{2},w^{'})\otimes Q(y_{2},y_{4})
\otimes Q^{*}(y_{3},y_{4})                           \\ 
& \simeq & L_{f}(z,y_{3},z^{'})\otimes 
L_{g}^{*}(w,y_{4},w^{'})\otimes Q^{*}(y_{3},y_{4}).  
\end{eqnarray*} 
One can show this isomorphism satisfies a coherency 
condition given a third point $(z,z^{'},w,w^{'},y_{5},y_{6})$ 
of $V\times_{\pi}Y^{[2]}$ lying in the same fibre as 
the points $(z,z^{'},w,w^{'},y_{1},y_{2})$ and 
$(z,z^{'},w,w^{'},y_{3},y_{4})$.  Therefore by 
Lemma~\ref{lemma:gendescent} $\hat{L}$ descends 
to a $\cstar$ bundle $L$ on $V$ and it is 
not difficult to show that $L$ trivialises 
the bundle gerbe $\pi_{1}^{-1}(S^{*}\otimes T)^{*}
\otimes \pi_{2}^{-1}(S^{*}\otimes T)$ over 
$V^{[2]}$.  Therefore the trivialisation $L$ 
corresponds to a stable morphism $h:\pi_{1}^{-1}(
S^{*}\otimes T)\to \pi_{2}^{-1}(S^{*}\otimes T)$.  
One can also show that there is a transformation     
$\theta:\pi_{1}^{-1}h\circ \pi_{3}^{-1}h\Rightarrow 
\pi_{2}^{-1}h$ satisfying the coherency condition 
$\d(\theta) = 1$.  We will omit the proof of this as 
the calculation gets rather messy.  The calculation 
exploits the fact that the transformations $\chi$ 
and $\rho$ satisfy $\d(\chi) = a$ and $\d(\rho) = a$ and 
so $\d(\chi\otimes \rho^{*}) = 1$ in a sense one 
has to make precise.   
Therefore we can apply Proposition~\ref{prop:6.6.3}
to conclude there is a bundle gerbe $(P,X_{P},U)$ on 
$U$ and a stable morphism $\eta:\pi^{-1}P
\to S^{*}\otimes T$ over $X_{U}$ corresponding 
to a trivialisation $(L_{\eta},\phi_{\eta})$.  There is also a 
transformation $\xi:\phi\circ \pi_{1}^{-1}\eta
\Rightarrow \pi_{2}^{-1}\eta$ compatible 
with $\psi$.  By shrinking $U$ if necessary, 
we may assume that the bundle gerbe $P$ is 
trivial --- say $P \simeq \d(J)$.  Therefore 
$S^{*}\otimes T$ is trivial as well.  In fact, 
if one considers the $\cstar$ bundle 
$\pi^{-1}J^{*}\otimes L_{\eta}$ on 
$\pi^{-1}X_{P}\times_{\pi}W\times_{\pi}Z$ then one 
can show that there is a coherent section 
$\hat{\phi}$ of $\d_{\pi^{-1}X_{P}}(\pi^{-1}J^{*}\otimes L_{\eta})$ 
and hence $\pi^{-1}J^{*}\otimes L_{\eta}$ descends 
to a $\cstar$ bundle $L_{k}$ on $W\times_{\pi}Z$.  
One can also show that there is a 
trivialisation $\phi_{k}:\d(L_{k})\to S^{*}\otimes T$.  
So we have a stable morphism $k:S\to T$.  As well, 
one can show that the transformations $\xi$ 
give rise to transformations $\a_{k}:\pi_{2}^{-1}k\circ 
g\Rightarrow f\circ (\pi_{1}^{-1}k\otimes 1)$ and 
which satisfy the compatibility condition with 
$\chi$ and $\rho$.  Therefore $(k,\a_{k}):(S,g,\rho)\to 
(T,f,\chi)$ is a 1-arrow.  So $\mathcal{Q}$ is 
locally connected.    

Since the bicategory of bundle gerbes and 
stable morphisms is a bigroupoid (Proposition~\ref{prop:stablebigrpd}), 
all 1-arrows of $\mathcal{Q}$ are coherently invertible and 
all 2-arrows of $\mathcal{Q}$ are invertible.  

We now have to show that the 2-gerbe 
$\mathcal{Q}$ is bound by the sheaf 
of abelian groups $\underline{\C}^{\times}_{M}$.  
For this, we first need to show that any 
two 1-arrows can be locally connected by 
a 2-arrow.  So let $(S,g,\rho)$ and 
$(T,f,\chi)$ be objects of $\mathcal{Q}_{U}$ 
joined by 1-arrows $(h,\a_{h}), (k,\a_{k}):
(S,g,\rho)\to (T,f,\chi)$.  So $h$ and 
$k$ correspond to trivialisations $(L_{h},\phi_{h})$ 
and $(L_{k},\phi_{k})$ respectively of $S^{*}\otimes T$.  
Therefore there is a $\cstar$ bundle $\hat{D}$ 
on $X_{U}$ such that $L_{h} \simeq L_{k}
\otimes \pi^{-1}\hat{D}$.  It is easy to check 
that $\a_{h}$ and $\a_{k}$ provide a section 
$\a_{k}\otimes \a_{h}^{*}$ of the bundle 
$\d(\hat{D})$ over $X_{U}^{[2]}$.  The coherency 
conditions satisfied by $\a_{h}$ and $\a_{k}$ 
ensure that $\d(\a_{k}\otimes \a_{h}^{*}) = 1$ and hence by 
Lemma~\ref{lemma:3.4.5} the bundle $\hat{D}$ 
descends to a $\cstar$ bundle $D$ on $U$.  
By shrinking $U$ if necessary we may assume that 
$D$ is trivial.  Hence there is a section of 
$D$ which lifts to a section $\theta$ of $\hat{D}$ 
which is compatible with $\a_{k}\otimes \a_{h}^{*}$.  
In other words, $\theta$ is a 2-arrow $(h,\a_{h})
\Rightarrow (k,\a_{k})$.  In a similar way one can 
show that for a given 1-arrow $(h,\a_{h})$ the 
sheaf of 2-arrows $(h,\a_{h})\Rightarrow 
(h,\a_{h})$ is isomorphic to $\underline{\C}^{\times}_{M}$ and 
that this isomorphism is natural with respect to 
composition of 2-arrows.  Therefore 
$\mathcal{Q}$ is a 2-gerbe on $M$ bound 
by $\underline{\C}^{\times}_{M}$.  

If one calculates the \v{C}ech cocycle $g_{ijkl}$ 
representing the class of $\mathcal{Q}$ in 
$\check{H}^{3}(M;\underline{\C}^{\times}_{M})$ then it is not hard to 
show that $g_{ijkl}$ represents the same class 
as the four class of the stable bundle 2-gerbe $Q$.            
\end{proof} 
 
As a consequence of this Proposition we 
have the following corollary: 

\begin{corollary} 
A stable bundle 2-gerbe $(Q,Y,X,M)$ has 
zero class in $H^{4}(M;\Z)$ precisely 
when it is trivial.  
\end{corollary} 

\begin{proof} 
We have already seen (Lemma~\ref{lemma:trivialstableb2gclass}) that a trivial 
stable bundle 2-gerbe has zero class 
in $H^{4}(M;\Z)$.  Conversely, given 
a stable bundle 2-gerbe with zero 
class in $H^{4}(M:\Z)$, the associated 
2-gerbe also has zero class in $H^{4}(M;\Z)$ 
and hence has a global object.  This global 
object trivialises the stable bundle 
2-gerbe $Q$.   
\end{proof} 

Let $(Q,Y,X,M)$ be a stable bundle 2-gerbe.  
We can form a new stable bundle 2-gerbe $(Q^{*},Y,X,M)$ from 
$Q$ which we could think of as the `dual' 
of $Q$ by letting the bundle gerbe $(Q^{*},Y,X^{[2]})$ 
be the dual of the bundle gerbe $(Q,Y,X^{[2]})$.  
So the $\cstar$ bundle $Q^{*}\to Y^{[2]}$ is 
$Q$ with the action of $\cstar$ changed to its 
inverse.  It is clear that the stable morphism 
$m:\pi_{1}^{-1}Q\otimes \pi_{3}^{-1}Q\to 
\pi_{2}^{-1}Q$ induces a stable morphism 
$m^{*}:\pi_{1}^{-1}Q^{*}\otimes \pi_{3}^{-1}Q^{*}\to 
\pi_{2}^{-1}Q^{*}$ and it is not difficult to show 
$m^{*}$ satisfies the requirements of 
Definition~\ref{def:7.2.1}.  It is also not 
difficult to show that the four class of the 
stable bundle 2-gerbe $Q^{*}$ is the negative of the 
four class of the stable bundle 2-gerbe $Q$.  
Similarly, given stable bundle 2-gerbes $(Q_{1},
Y_{1},X_{1},M)$ and $(Q_{2},Y_{2},X_{2},M)$ there 
is a notion of their product.  This is the 
stable bundle 2-gerbe $(Q_{1}\otimes Q_{2},Y_{12},X_{12},M)$ 
where $X_{12} = X_{1}\times_{M}X_{2}$, $Y_{12} = 
(p_{1}^{[2]})^{-1}Y_{1}\times_{X_{12}^{[2]}}(p_{2}
^{[2]})^{-1}Y_{2}$ (here $p_{1}:X_{1}\times_{M}X_{2}
\to X_{1}$ and $p_{2}:X_{1}\times_{M}X_{2}\to X_{2}$ 
denote the projections on the first and second factors 
respectively).  $Q_{1}\otimes Q_{2}$ is the $\cstar$ 
bundle on $Y_{12}^{[2]}$ given by $(q_{1}^{[2]})^{-1}Q_{1}
\otimes (q_{2}^{[2]})^{-1}Q_{2}$ where now 
$q_{1}:Y_{12}^{[2]}\to Y_{1}^{[2]}$ and 
$q_{2}:Y_{12}^{[2]}\to Y_{2}^{[2]}$ are the 
projections on the first and second factors 
respectively.  Again it is not difficult to show 
that $Q_{1}\otimes Q_{2}$ is a stable bundle 2-gerbe 
and that the four class of $Q_{1}\otimes Q_{2}$ 
is the sum of the four classes of the stable 
bundle 2-gerbes $Q_{1}$ and $Q_{2}$.  

Let us now apply the preceding discussion to the 
following situation.  Suppose we are given stable 
bundle 2-gerbes $(Q_{1},Y_{1},X_{1},M)$ and 
$(Q_{2},Y_{2},X_{2},M)$ with four classes 
$D_{4}(Q_{1})$ and $D_{4}(Q_{2})$ respectively.  
Suppose moreover that $D_{4}(Q_{1}) = 
D_{4}(Q_{2})$.  Form the stable bundle 2-gerbe 
$Q_{1}\otimes Q_{2}^{*}$.  Then by the above 
discussion $D_{4}(Q_{1}\otimes Q_{2}^{*}) = 0$.  
Therefore by the corollary above 
$Q_{1}\otimes Q_{2}^{*}$ is trivial.  Hence there 
is a bundle gerbe $(T,Z,X_{12})$ on $X_{12}$ and 
a stable morphism $f:Q_{1}\otimes Q_{2}^{*}\to 
\d(T)$ satisfying the requirements of 
Definition~\ref{def:trivialstableb2g}.  Therefore 
$f$ induces a stable bundle 2-gerbe morphism 
(see Note~\ref{note:stableb2gmorphism}) 
$Q_{1}\otimes Q_{2}^{*}\otimes Q_{2}\to \d(T)\otimes Q_{2}$.  
The stable bundle 2-gerbe $Q_{2}^{*}\otimes Q_{2}$ 
can then be written as $\d(T^{'})$ and so we have 
a stable bundle 2-gerbe morphism $Q_{1}\otimes \d(T^{'})\to 
Q_{2}\otimes \d(T)$.  It should then be the case 
that this notion of `stable 2-isomorphism' for 
stable bundle 2-gerbes introduces an equivalence 
relation on the set of all stable bundle 2-gerbes 
on $M$.  Since every class in $H^{4}(M;\Z)$ appears 
as the four class of some stable bundle 2-gerbe (for 
example the stable bundle 2-gerbe arising from the 
strict bundle 2-gerbe associated to the $BB\cstar$ 
bundle with the four class in question --- see 
Section~\ref{sec:11.1}) this leads us to the 
following result which we will state as a conjecture 
since we have not given a definite proof.  

\begin{conjecture} 
There is a bijection between the set of all 
stable 2-isomorphism classes of stable bundle 2-gerbes 
on $M$ and $H^{4}(M;\Z)$.  
\end{conjecture}

\bibliography{booklist} 

\begin{thebibliography}{10}

\bibitem{Ben}
J.~B{\'e}nabou.
\newblock {\em Introduction to Bicategories}.
\newblock Number~47 in Lecture Notes in Mathematics. Springer-Verlag, Berlin,
  1967.

\bibitem{Bor}
A.~Borel.
\newblock {\em Topics in the {H}omology {T}heory of {F}ibre {B}undles}.
\newblock Number~36 in Lecture Notes in Mathematics. Springer-Verlag, Berlin,
  1967.

\bibitem{BotShuSta}
R.~Bott, H.~Shulman, and J.~Stasheff.
\newblock On the de {R}ham {T}heory of {C}ertain {C}lassifying {S}paces.
\newblock {\em Advances in Mathematics}, 20:43--56, 1976.

\bibitem{BotTu}
R.~Bott and L.~W. Tu.
\newblock {\em Differential Forms in Algebraic Topology}.
\newblock Springer-Verlag, New York, 1982.

\bibitem{Bre}
L.~Breen.
\newblock {\em On the Classification of Two Gerbes and Two Stacks}.
\newblock Number 225 in Ast\'{e}risque. Soci\'{e}t\'{e} de Math\'{e}matique de
  France, 1994.

\bibitem{Bry1}
J.-L. Brylinski.
\newblock Private communication.

\bibitem{Bry}
J.-L. Brylinski.
\newblock {\em Loop Spaces, Characteristic Classes and Geometric Quantization},
  volume 107 of {\em Progress in Mathematics}.
\newblock Birkhauser, Berlin, 1992.

\bibitem{BryMcL1}
J.-L. Brylinski and D.~Mc{L}aughlin.
\newblock Geometry of {D}egree {F}our {C}haracteristic {C}lasses and {L}ine
  {B}undles on {L}oop {S}paces {I}.
\newblock {\em Duke Mathematical Journal}, 75(3):603--638, 1994.

\bibitem{BryMcL2}
J.-L. Brylinski and D.~Mc{L}aughlin.
\newblock Geometry of {D}egree {F}our {C}haracteristic {C}lasses and of {L}ine
  {B}undles on {L}oop {S}paces {I}{I}.
\newblock {\em Duke Mathematical Journal}, 83(1):105--139, 1996.

\bibitem{BryMcL3}
Jean~Luc Brylinski and D.~Mc{L}aughlin.
\newblock A {G}eometric {C}onstruction of the {F}irst {P}ontryagin {C}lass.
\newblock In Louis~H. Kauffman and Randy~A. Baadhio, editors, {\em Series on
  Knots and Everything}, pages 209--220, Singapore, 1992. World Scientific.

\bibitem{CaCroMur}
A.~L. Carey, D.~Crowley, and M.~K. Murray.
\newblock {P}rincipal {B}undles and the {D}ixmier-{D}ouady {C}lass.
\newblock {\em Communications in Mathematical Physics}, 193(1):171--196, 1997.

\bibitem{CaMicMur}
A.~L. Carey, J.~Mickelsson, and M.~K. Murray.
\newblock {I}ndex {T}heory, {G}erbes and {H}amiltonian {Q}uantization.
\newblock {\em Communications in Mathematical Physics}, 183(3):702--722, 1997.

\bibitem{CaMur}
A.~L. Carey and M.~K. Murray.
\newblock {F}addeev's {A}namoly and {B}undle {G}erbes.
\newblock {\em Letters in Mathematical Physics}, 37:29--36, 1996.

\bibitem{CaMuWa}
A.~L. Carey, M.~K. Murray, and B.~L. Wang.
\newblock {H}igher {B}undle {G}erbes and {C}ohomology {C}lasses in {G}auge
  {T}heories.
\newblock {\em Journal of Geometry and Physics}, 21:183--197, 1997.

\bibitem{Car}
E.~Cartan.
\newblock La topologie des espaces repr\'{e}sentatifs des groupes de {L}ie.
\newblock In {\em Oeuvres Compl\`{e}tes, Part 1, Vol. 2}, pages 1307--1330.
  Gauthier-Villars, Paris, 1952.

\bibitem{Che}
K-T. Chen.
\newblock {I}terated {I}ntegrals of {D}ifferential {F}orms and {L}oop {S}pace
  {H}omology.
\newblock {\em Ann. of Math.}, 97:217--246, 1973.

\bibitem{Dup}
J.~L. Dupont.
\newblock {\em Curvature and Characteristic Classes}.
\newblock Number 640 in Lecture Notes in Mathematics. Springer-Verlag, 1978.

\bibitem{Gaj}
P.~Gajer.
\newblock {G}eometry of {D}eligne {C}ohomology.
\newblock {\em Inventiones Math.}, 127:155--207, 1997.

\bibitem{Gir}
J.~Giraud.
\newblock {\em Cohomologie non-ab\'{e}lienne}.
\newblock Number Band 179 in Die Grundlehren der Mathematischen Wissenschaften.
  Springer-Verlag, New York-Berlin, 1971.

\bibitem{GorPowStr}
R.~Gordon, A.~J. Power, and R.~Street.
\newblock {\em Coherence in Tri-categories}.
\newblock Number 558 in Memoirs of the A.M.S. Wiley Interscience, 1995.

\bibitem{GriHar}
P.~Griffiths and J.~Harris.
\newblock {\em Principles of Algebraic Geometry}.
\newblock Wiley (Interscience), New York, 1978.

\bibitem{Hit}
N.~Hitchin.
\newblock {L}ectures on {S}pecial {L}agrangian {S}ubmanifolds.
\newblock In {\em Proceedings of School on Differential Geometry}, Trieste,
  1999. International Centre for Theoretical Physics.

\bibitem{Hus}
D.~Husemoller.
\newblock {\em Fibre Bundles}.
\newblock Mc{G}raw Hill Book Company, New York, 1966.

\bibitem{Kal}
J.~Kalkkinen.
\newblock Gerbes and {M}assive {T}ype {I}{I} {C}onfigurations.
\newblock \texttt{hep-th/9905018}.

\bibitem{KelStr}
G.~M. Kelly and R.~Street.
\newblock {\em Review of the Elements of 2-categories}.
\newblock Number 420 in Lecture Notes in Mathematics. Springer Verlag, 1974.

\bibitem{KoNo}
S.~Kobayishi and K.~Nomizu.
\newblock {\em Foundations of Differential Geometry}, volume~1.
\newblock Wiley (Interscience), New York, 1963.

\bibitem{Mac}
S.~Mac{L}ane.
\newblock {\em Categories for the Working Mathematician}.
\newblock Springer-Verlag, New York, 1971.

\bibitem{MaMo}
S.~Mac{L}ane and I.~Moerdijk.
\newblock {\em Sheaves in {G}eometry and {L}ogic}.
\newblock Universitext. Springer-Verlag, New York, 1992.

\bibitem{May1}
J.~P. May.
\newblock {\em Simplicial Objects in Algebraic Topology}.
\newblock Van Nostrand, 1967.

\bibitem{May2}
J.~P. May.
\newblock {\em Geometry of Iterated Loop Spaces}.
\newblock Number 271 in Lecture Notes in Mahtematics. Springer-Verlag, 1970.

\bibitem{MosPer}
M.~Mostow and J.~Perchick.
\newblock {N}otes on {G}el'fand-{F}uks {C}ohomology and {C}haracteristic
  {C}lasses ({L}ectures by {R}. {B}ott).
\newblock In {\em Eleventh Holiday Symposium}. New Mexico State University,
  December 1973.

\bibitem{Mos}
M.~A. Mostow.
\newblock {T}he {D}ifferentiable {S}pace {S}tructures of {M}ilnor {C}lassifying
  {S}paces, {S}implicial {C}omplexes, and {G}eometric {R}ealisation.
\newblock {\em Journal of Differential Geometry}, 14(2):255--293, 1979.

\bibitem{Mur}
M.~K. Murray.
\newblock Bundle {G}erbes.
\newblock {\em Journal of the London Mathematical Society}, 54(2):403--416,
  1996.

\bibitem{MurSte}
M.~K. Murray and D.~Stevenson.
\newblock Bundle {G}erbes: {S}table {I}somorphism and {L}ocal {T}heory.
\newblock Preprint.

\bibitem{Seg2}
G.~Segal.
\newblock Classifying {S}paces and {S}pectral {S}equences.
\newblock {\em Institut des Hautes \'{E}tudes Scientifiques Publications
  Math\'{e}matiques}, 34:105--112, 1968.

\bibitem{Seg3}
G.~Segal.
\newblock Cohomology of {T}opological {G}roups.
\newblock In {\em Symposia Mathematica}, volume~IV, pages 377--387, 1970.

\bibitem{Seg1}
G.~Segal.
\newblock Categories and {C}ohomology {T}heories.
\newblock {\em Topology}, 13:293--312, 1974.

\bibitem{Spa}
E.H. Spanier.
\newblock {\em Algebraic Topology}.
\newblock Mc{G}raw-Hill Book Company, New York, 1966.

\bibitem{Str}
R.~Street.
\newblock Categorical {S}tructures.
\newblock In M.~Hazewinkel, editor, {\em Handbook of Algebra}, pages 529--577.
  Elsevier Science, Amsterdam, 1996.

\end{thebibliography}

\end{document}